\documentclass[a4paper,11pt]{amsart}
\usepackage{amsmath,amsthm,amssymb}
\usepackage{graphicx}
\usepackage{subcaption}
\usepackage{float}
\usepackage{cancel}
\usepackage{slashed}
\usepackage{comment}
\newcommand{\mychoice}[3]{#3
}
\mychoice{
\newcommand{\plabel}[1]{ \label{#1}}
\newcommand{\gbibitem}[1]{ \bibitem{#1}}
\newcommand{\snewpage}{}

\specialcomment{commentx}{}{}\excludecomment{commentx}
\specialcomment{commenty}{}{}\excludecomment{commenty}
\newcommand{\rechoicecomm}[1]{}
}{
\usepackage{xcolor}
\newcommand{\plabel}[1]{ \label{#1}\rlap{\smash{${}^{^{[#1]}}$}}}
\newcommand{\gbibitem}[1]{ \bibitem{#1}\rlap{\smash{${}^{^{[#1]}}$}}}
\newcommand{\snewpage}{\newpage}

\newenvironment{commentx}{\color{magenta} }{\color{black} }
\newenvironment{commenty}{\color{blue} }{\color{black} }
\newcommand{\rechoicecomm}[1]{#1}
}{
\usepackage{xcolor}
\newcommand{\plabel}[1]{ \label{#1}}
\newcommand{\gbibitem}[1]{ \bibitem{#1}}
\newcommand{\snewpage}{}

\specialcomment{commentx}{}{}\excludecomment{commentx}
\newenvironment{commenty}{\scriptsize }{\normalsize }
\newcommand{\rechoicecomm}[1]{}
}
\DeclareMathOperator{\Tan}{\slashed{\mathrm T}an}
\DeclareMathOperator{\ad}{ad}
\DeclareMathOperator{\Ad}{Ad}
\DeclareMathOperator{\AC}{AC}
\DeclareMathOperator{\ATT}{AT}
\DeclareMathOperator{\LR}{LR}
\DeclareMathOperator{\sgn}{sgn}
\DeclareMathOperator{\Id}{Id}
\DeclareMathOperator{\id}{id}

\DeclareMathOperator{\des}{des}

\DeclareMathOperator{\artanh}{artanh}
\DeclareMathOperator{\sech}{sech}

\DeclareMathOperator{\arcosh}{arcosh}

\DeclareMathOperator{\mul}{mul}

\DeclareMathOperator{\asc}{asc}\DeclareMathOperator{\conv}{conv}
\DeclareMathOperator{\MP}{\mathcal M}
\DeclareMathOperator{\intD}{\mathring D}
\newcommand{\real}{\mathrm{real}}
\newcommand{\anal}{\mathrm{an}}
\DeclareMathOperator{\Dbar}{D}
\DeclareMathOperator{\Rea}{Re}
\DeclareMathOperator{\Ima}{Im}
\DeclareMathOperator{\Rexp}{exp_{R}}
\DeclareMathOperator{\Lexp}{exp_{L}}
\DeclareMathOperator{\BCH}{BCH}
\DeclareMathOperator{\spec}{sp}

\theoremstyle{definition}
\newtheorem{point}{}[section]
\newtheorem{defin}[point]{Definition}
\newtheorem{deco}[point]{}
\newtheorem{remark}[point]{Remark}
\newtheorem{example}[point]{Example}
\newtheorem{discussion}[point]{Discussion}
\theoremstyle{plain}

\newtheorem{lemma}[point]{Lemma}
\newtheorem{cor}[point]{Corollary}
\newtheorem{theorem}[point]{Theorem}
\theoremstyle{remark}

\newcommand{\bo}{\boldsymbol}
\newcommand{\leaveout}[1]{}

\newcommand{\qedexer}{\renewcommand{\qedsymbol}{$\diamondsuit$}\qed\renewcommand{\qedsymbol}{$\Box$}}
\newcommand{\qedremark}{\renewcommand{\qedsymbol}{$\triangle$}\qed\renewcommand{\qedsymbol}{$\Box$}}
\newcommand{\qedno}{\renewcommand{\qedsymbol}{}}

\newcommand{\eqed}{\pushQED{\qed}\qedhere\popQED}
\newcommand{\eqedexer}{\renewcommand{\qedsymbol}{$\diamondsuit$}\pushQED{\qed}\qedhere\popQED\renewcommand{\qedsymbol}{$\Box$}}
\newcommand{\eqedremark}{\renewcommand{\qedsymbol}{$\triangle$}\pushQED{\qed}\qedhere\popQED\renewcommand{\qedsymbol}{$\Box$}}

\newcommand{\marginextend}[1]{ \addtolength{\oddsidemargin}{-#1}  \addtolength{\evensidemargin}{-#1}
  \addtolength{\textwidth}{#1}\addtolength{\textwidth}{#1}}
\newcommand{\updownextend}[1]{ \addtolength{\topmargin}{-#1}  \addtolength{\textheight}{#1}
\addtolength{\textheight}{#1}}
\marginextend{1cm}
\updownextend{0cm}
\allowdisplaybreaks[4]
\title{Convergence estimates for the Magnus expansion I. Banach algebras}
\author{Gyula Lakos}
\address{Alfréd Rényi Institute of Mathematics, Reáltanoda utca 13-15, Budapest, H--1053, Hungary}
\email{lakos@renyi.hu}
\begin{document}
\begin{abstract}
We  review and provide simplified proofs related to the Magnus expansion, and improve convergence estimates.
Observations and improvements concerning the Baker--Campbell--Hausdorff expansion are also made.

In this Part I, we consider the general Banach algebraic setting.
We show that the (cumulative) convergence radius of the Magnus expansion is $2$;
and of the Baker--Campbell--Hausdorff series is $\mathrm C_2=2.89847930\ldots$.
More generally, the resolvent method is developed in the analytic setting.
\end{abstract}
\maketitle
\section*{Introduction}
\textbf{General introduction.}
The
Baker--Campbell--Hausdorff expansion and its continuous generalization, the Magnus expansion,
attract attention from time to time.
This is probably partly due to their aesthetically pleasing nature.
There are many rediscoveries in this area;
%For example,
just broadly related to the Magnus expansion,
see Magnus \cite{M} and Chen \cite{Ch};  Mielnik, Pleba\'nski \cite{MP} and Strichartz \cite{St} and
Vinokurov \cite{V}; Goldberg \cite{G} and Helmstetter \cite{H}.

In this series of articles, we review and provide simplified proofs related to the Magnus expansion, and improve convergence estimates.
Observations and improvements concerning the Baker--Campbell--Hausdorff expansion are also made.
See Blanes, Casas, Oteo, Ros \cite{BCOR} for a general review of the Magnus expansion for aspects we do not discuss.
A recent reference with respect to vector fields is Beauchard,   Le Borgne,  Marbach \cite{BLBM}.

We refer to Bourbaki \cite{B}, Reutenauer \cite{R}, and Bonfiglioli, Fulci \cite{BF}
for general background in Lie theory, including the Baker--Campbell--Hausdorff expansion,
the Poincar\'e--Birkhoff--Witt theorem, Friedrichs' criterion for commutator polynomials,
the Dynkin--Specht--Wever lemma,
and the Dynkin formula, i.~e.~the ``standard universal Lie algebraic tools'';
for these we do not make particular references.
(They can mostly be circumvented for our purposes, but they are unavoidable  for putting the material into context.)
Especially \cite{BF} can be useful:
It provides a detailed discussion of the above mentioned topics and a guide to further literature,
but the  discussion itself ends around the same place where our discussion begins.
Dunford, Schwartz \cite{DS} still provides a reasonably good reference for functional analysis.
We cite Flajolet, Sedgewick \cite{FS} for  generating function
techniques, including a discussion of the role of Pringsheim's theorem.
We refer to Coddington, Levinson \cite{CL}, Hsieh, Sibuya \cite{HS}, Teschl \cite{T} for the theory
of ordinary differential equations, including
existence and uniqueness theorems, and some standard solution techniques.
%(Of course, the above mentioned works contain much more than the necessary background material.)
Some other background references will be introduced topically.

Our viewpoints are synthetic; and
in our discussion we will allow relatively more general integrands than usual;
and but we always do try to relate to the original papers.

\snewpage
\textbf{Introduction to the Banach algebraic setting (Part I).}
In this part, we aim to deal with general Banach algebras, although some liberties are taken to extend the range of algebras.
We sometimes use the locally convex algebras $\mathrm F^{1,\mathrm{loc}}_{\mathbb K}$ for computational purposes.
Moreover, we may also consider the ``algebraic'' topology of formal power series,
which we use for algebraic purposes without much comment.

Regarding earlier results for convergence in general Banach algebras,
 the most naive power series based estimate for the BCH series shows convergence
 for cumulative norm (i.~e.~the norm variation of the noncommutative mass to be expanded)
 less than $\log 2$, cf.~Dynkin \cite{Dy}, which bound extends to the Magnus case without difficulty, cf.~remarks
 in Pechukas, Light \cite{PL}, or Karasev, Mosolova \cite{KM}.
Beyond such naive power series estimates, the relevant papers are Strichartz \cite{St}, Thompson \cite{Th},
 Day, So, Thompson \cite{DST}, Vinokurov \cite{V}, Moan, Oteo \cite{MO}.
In particular, although their results are formulated in slightly weaker settings,
 Thompson \cite{Th} (1989) and Moan, Oteo \cite{MO} (2001) essentially state
 that the BCH expansion and the Magnus expansion, respectively, converge if the cumulative norm is less than $2$.
Divergence was studied by Wei \cite{W} (who presents counterexamples but does not optimize
for cumulative norm) and Michel \cite{Mi} (who presents sophisticated analysis
but for the BCH case of Banach--Lie algebras), but the most relevant information is given by  Vinokurov \cite{V} (1991),
 that there are simple counterexamples for the convergence of the BCH expansion with cumulative norm greater than $\pi$;
 and by Sch\"affer \cite{Scha} (1964) (indirectly), and by Moan \cite{Ma} (2002), Moan, Niesen \cite{MN} (2008)
 showing that there is a counterexample for the convergence of the Magnus expansion with cumulative norm $\pi$.
These latter counterexamples are, however, in the $C^*$-algebraic (i.~e.~Hilbert space operator) setting
where stronger convergence results hold;
 and they turn out to be still not representative with respect to the general Banach algebraic setting.
The above mentioned papers altogether put the guaranteed convergence radius of the Magnus expansion (and the BCH expansion)
 between $2$ and $\pi$ in the Banach algebraic setting.
\snewpage

\textbf{Outline of content (Part I).}
In this article, we show that the Magnus expansion converges if the cumulative norm is $<2$;
 moreover,  it converges even if the  cumulative norm is $2$ and the integrand is of Lebesgue--Bochner type;
 and there are counterexample measures with cumulative norm $2$, and there are counterexamples of
 multivariable BCH type with cumulative norm $>2$.
In the above mentioned convergent cases we show that the exponential formula holds even in the stronger logarithmic sense.
In establishing these results, the key features are the consideration of the Eulerian generating function and the resolvent method.
(These were introduced by Mielnik, Pleba\'nski \cite{MP} into the study of the Magnus expansion on the formal level.)
We also show that the Baker--Campbell--Hausdorff series converges if the cumulative norm
 is less than $\mathrm C_2=2.89847930\ldots$, and there are counterexamples with cumulative norm $\mathrm C_2$.
More specifically:
Section \ref{sec:MagnusBanach} provides an introduction to the Magnus expansion in the setting of Banach algebras.
Apart from the discussion of divergence, the emphasis is on recovering some known results, organized efficiently.
The BCH expansion and certain algebraic aspects are also discussed.
Section \ref{sec:MagnusResolvent} uses the resolvent technique to refine the results concerning the convergence of the Magnus expansion;
 we prove the logarithmic version of the Magnus formula, and we discuss the critical case when the cumulative norm is $2$.
In Section \ref{sec:DiscRes}, the resolvent method is reformulated in a more discrete setting, and we
 deal with the convergence of the Baker--Campbell--Hausdorff expansion.
Section \ref{sec:adj} treats the Magnus expansion in the adjoint representation,
encountering a situation where the naive application of the resolvent method is less effective than the algebraic
combinatorial
approach.
In Appendix \ref{sec:mea} we discuss measure-theoretic terminology.
In Appendix \ref{sec:ellone}, we describe the spaces $\mathrm F^{1,\mathrm{loc}}_{\mathbb K}(\Omega,\mathfrak F,\omega)$,
 which are essentially formal power series with continuous variables.
The algebras introduced here are technically quite important for us (even if they are quite natural);
 their discussion is put here not to interrupt the flow of the discussion of the Magnus expansion.
In Appendix  \ref{sec:sc}, we recall information on the spectral calculus of the adjoint representation.
In Appendix  \ref{sec:schuranal}, we consider certain analytic aspects of Schur's formulae.
\begin{commenty}
In what follows, in smaller print, we include some supplementary information.
\end{commenty}
\snewpage

\textbf{General notation and terminology.}
The principal point of the Magnus expansion is that (in appropriate circumstances)
 the (sum of the) time-ordered exponential expansion of an ``ordered non-commutative mass'' $\phi$
 is the ordinary exponential of the (sum of the) Magnus expansion of this ``ordered non-commutative mass'' $\phi$.
Traditionally such a $\phi$ is presented as
\[\phi=f\mathbf 1_{[a,b)}\]
or
 \[\text{`$t\mapsto$'\,}\phi(t)=\text{`$t\mapsto$'\,}f(t)\,\mathrm dt|_{[a,b)} ;\]
 where $f$ is a continuous or at least Lebesgue--Bochner integrable Banach algebra (or Banach--Lie algebra) valued function,
 and $\mathbf 1$ denotes the Lebesgue measure.
In what follows, a bit more generally, we take any continuous Banach algebra
 (or Banach--Lie algebra) valued ordered measure of finite variation for $\phi$.
(The notation of measures, as it is usual in the mathematical literature, may sometimes be a sloppy
but it is always done for the ease of notation.)

All Banach algebras will be (real or complex) unital algebras.
The unit element of the algebra $\mathfrak A$ will be denoted by $1_{\mathfrak A}$ or,
a bit haphazardly, just by $1$; but in case of $n\times n$ matrices, by $\Id_n$.
However, the identity operator on the linear space $\mathfrak A$ will be denoted by $\Id_{\mathfrak A}$
(and $\id_{\mathfrak A}$ denotes the identity function just on any set $\mathfrak A$).
For the purposes of spectral calculus, it is customary to complexify the real Banach algebras.
If $\mathfrak A$ is a real Banach algebra with norm $|\cdot|_{\mathfrak A}$, then for the complexification
 $\mathbb C\mathfrak A$, our choice of norm is
 \[|x|_{\mathbb C\mathfrak A}=\inf\left\{\sum_{i\in I} |t_i|\cdot|x_i|_{\mathfrak A}\, :\,
 x=\sum_{i\in I} t_ix_i,\, x_i\in\mathfrak A, t_i\in\mathbb C \right\}.\]
This is the maximal (i.~e.~most pessimistic) isometric complexification of $(\mathfrak A,|\cdot|_{\mathfrak A})$.
Other complexifications are possible, and, in fact, useful in certain places.
(For example, the complexification of operators on real Hilbert spaces, that is the $C^*$-algebraic complexification.)
If $A$ is an element of a Banach algebra, then $\spec(A)$ will denote its spectrum, and $\mathrm r(A)$ will denote its spectral radius.
The value $1/\mathrm r(A)$ can be called as the convergence radius of $A$, as it is equal to the
 radius of the maximal open disc of center $0$ on the complex plane such that the power series $z\mapsto \sum_{n=0}^\infty A^nz^n$
 extends to an analytic function on that disk.
On the complex plane, $\Dbar(z_0,r)$ denotes the (possibly degenerate) closed disk with center $z_0$ and radius $r$.
$\intD(z_0,r)$ denotes the corresponding open disk.
Sometimes $\mathbb C$ is identified with $\mathbb R^2$, which yields the notation  $\Dbar((a,b),r)$.

In the arithmetic context, we may use `$A\times B$' when it can be understood either as $AB$ or $BA$.
Also, `$\frac AB$' may be used  when it can be understood either as $AB^{-1}$ or $B^{-1}A$.
Even more generally, `$\frac AB$' may be used for $X$, if it is the unique solution of the system of equations $A=BX=XB$.
\snewpage

When we consider generating functions, we may understand them either in formal, real function theoretic, or analytic sense.
For example, the series $\sum_{n=0}^\infty\frac{(2n)!}{4^nn!(n+1)!} x^n$ can be
 understood as a formal power series; or as a real function which is
 equal to $\frac2{1+\sqrt{1-x}}$ for $0\leq x\leq1$ and equal to $+\infty$ for $x>1$;
 or as the analytic function $\frac2{1+\sqrt{1-x}}$ defined at least on $\intD(0,1)$, but rather on $\mathbb C\setminus[1,\infty) $.
We try to be clear in what sense we understand the generating function
 at hand, but, of course, much of the power of generating functions comes from change between various viewpoints.
Using power series, one often encounters the ``problem of specification''.
Most typically:
When a formal series is given not in series form, but, say, as a result of division, e. g. $\frac{\mathrm e^x-1}{x}$, then
 it can be well-defined as a formal power series, but substituting a zero divisor to the place of $x$
 naively will lead to an expression which does no make sense as such.
For that reason we prefer to use concrete function names, e. g. $\alpha(x)$ in the previous case.
We do not do this consistently, but we do this when the a divisional form might cause confusion.
However, we do not eschew the divisional arithmetic of formal power series
 (or the use of analytical methods) either, as they are very useful to establish identities.
If $A(X)=\sum_{n=0}^\infty a_nX^n$, $B(X)=\sum_{n=0}^\infty b_nX^n$ are formal power series, then we use the notation
 \[A(X)\stackrel{\forall X}\leq B(X)\]
 to indicate that $a_n\leq b_n$ for all $n$.
If the coeffients in $A(X)$ are from a Banach algebra $\mathfrak A$, then we use the notation
$|A(X)|_{\mathfrak A}^{\forall X}=\sum_{n=0}^\infty |a_n|_{\mathfrak A}X^n$.
These notations also extend to the case of several, possibly non-commuting variables.

If $A(X_\omega\,:\,\omega\in \Omega)$ is a formal power series in several  non-commuting variables, then
 we might consider its version $A(x_\omega\,:\,\omega\in \Omega)$ with commuting variables.
If $A(X_\omega\,:\,\omega\in \Omega)$ has only nonnegative coefficients, then we can
 consider $A_\real(x_\omega\,:\,\omega\in \Omega)$ with $x_\omega\in[0,\infty]$, taking possibly infinite values.
For us, by default, power series are formal; but we can substitute elements from a
 more general algebra $\mathfrak A$, in which case existence is understood as at least order-independent
 convergence in the topology of $\mathfrak A$, but preferably absolute convergence (in locally convex spaces);
 except in the (nonnegative) `real' case, when we allow $+\infty$ as the result or even as the argument.
 $\Sigma_S$ denotes the permutations of the set $S$, $\Sigma_n$ denotes the permutations of $\{1,\ldots,n\}$.
\snewpage

For a Banach algebra element $A$, if the spectrum of $A$ is contained in $\mathbb C\setminus (-\infty,0] $,
 then we define the logarithm as
\begin{equation}
\qquad
\log A=\int_{\lambda=0}^1 \frac{A-1}{\lambda +(1-\lambda)A }\,d\lambda
\qquad\left(=\int^{0}_{s=-\infty}\frac{A-1}{(1-s)(A-s)}\, \mathrm ds\right).
\plabel{eq:logdef0}
\end{equation}
If the spectrum of $A$ is contained in $\intD(1,1)$, then we define
\begin{equation}
\qquad
\log^{(\mathrm{pow})} A=\sum_{n=1}^\infty \frac{(-1)^{n+1}}{n}(A-1)^n.
\plabel{eq:logpow0}
\end{equation}
It is well-known that in the latter case $\log A= \log^{(\mathrm{pow})} A$, thus
 $\log^{(\mathrm{pow})}$ can be considered as a restricted use of $\log$.
In what follows, it goes without saying that  $\log^{(\mathrm{pow})}$ can be used instead of $\log$ if $A$ gets close to $1$.
These definitions can also work in the formal case, when we consider formal perturbations of $1$.
(In fact, $\log^{(\mathrm{pow})}$ can be used for formal perturbations of $c1$ with $c\in \intD(1,1)$, and
 $\log$ can be used even for formal perturbations of $c1$ with $c\in \mathbb C\setminus (-\infty,0] $,
 but the latter extension can also be achieved by prescription, and it is better not to use anyway.)
Formula \eqref{eq:logdef0} is very natural from analytical viewpoint; it makes sense even on the level of elementary calculus:
 we integrate the derivative of $\lambda\mapsto \log(\lambda +(1-\lambda)A)$ (times $-1$, due to the direction);
 moreover it is spectral, which is emphasized by the change of variables $\lambda=\frac{s}{s-1}$, $s=\frac{\lambda}{\lambda-1}$
 (which are orientation-reversing on the indicated domains).
Note that $\lambda\mapsto \log(\lambda +(1-\lambda)A)$ can be defined and studied by any reasonable holomorphic
 functional calculus validating the ``$\lambda$-integral'' in \eqref{eq:logdef0}; there is no need
 to start from the sectorial ``$s$-integral''.
The form \eqref{eq:logpow0} is better for  combinatorial purposes.
It is known that $\log$ and $\exp$ induces a bijection of the spectral domains
 $\{A\,:\,\spec A\subset \mathbb C\setminus(-\infty,0]\}
 \overset{\log}{\underset{\exp}\rightleftarrows}\{A\,:\,\spec A\subset \{z\in\mathbb C\,:\,|\Ima z|<\pi\}\}.$

We emphasize that for us $\log$ or the ``canonical'' logarithm is the ``spectral logarithm'' of \eqref{eq:logdef0}.
If $\exp C=A$ holds, then we may say that $C$ is an ``accidental logarithm'' of $A$
 (but even this terminology is dangerous).
The switch  $\lambda=\frac{s}{s-1}$, $s=\frac{\lambda}{\lambda-1}$ can be termed
in general as passing from the ``additively presented'' spectrum (resolvent) via $s$
to  ``the multiplicatively presented'' spectrum (resolvent) via $\lambda$.
This has no mathematical depth in itself but it is useful for us.

In what follows, it goes without much saying that how analytical formulas have formal consequences by application to truncated polynomial algebras.
In general, once a Banach algebraic result is obtained, we do not spend much time on how to obtain a formal version.
In the other direction, one must be more careful.

\snewpage
\begin{commentx}
\tableofcontents
\end{commentx}
\snewpage
\section{The combinatorial approach}
\plabel{sec:MagnusBanach}
\subsection{The Magnus expansion}
~\\

Let $\mathfrak A$ be a Banach algebra (real or complex). For $X_1,\ldots,X_k\in\mathfrak A$, we define the Magnus commutator
 (or Dynkin commutator) by
\begin{equation}\mu_k(X_1,\ldots,X_k)=
\frac{\partial^k}{\partial t_1\cdot\ldots\cdot\partial t_k}\log(\exp(t_1X_1)\cdot\ldots\cdot\exp(t_kX_k))\Bigl|_{t_1=\ldots=t_k=0};
\plabel{eq:mudef}
\end{equation}
 where $\exp$ is defined as usual, and $\log$ can be defined either by power series around $1$,
 or with spectral calculus, the complex plane cut along the negative real axis.
Algebraically, in terms of formal power series,
\begin{equation}\mu_k(X_1,\ldots,X_k)=\log(\exp(X_1)\cdot\ldots\cdot\exp(X_k))_{\text{the variables }
X_1\ldots X_k\text{ has multiplicity }1},\plabel{eq:sinmu}\end{equation}
i.~e.~the part of $\log(\exp(X_1)\cdot\ldots\cdot\exp(X_k))$ where every variable $X_i$ has multiplicity $1$.
This implies that, in terms of formal power series,
\begin{equation}\log(\exp(X_1)\cdot\ldots\cdot\exp(X_k))=\sum_{\substack{\{i_1,\ldots, i_l\}\subset \{1,\ldots,k\}  \\ i_1<\ldots <i_l }}
\mu_l(X_{i_1},\ldots,X_{i_l}) +H(X_1,\ldots,X_k),\plabel{eq:bare}
\end{equation}
where $H(X_1,\ldots,X_k)$ collects the terms where some variable have multiplicity more than one.
Taking the exponential, and detecting the terms where every variable has multiplicity $1$, this yields
\begin{equation}
X_1\cdot\ldots\cdot X_k=\sum_{\substack{I_1\dot\cup\ldots \dot\cup I_s=\{1,\ldots,k\}\\I_j=\{i_{j,1},\ldots,i_{j,l_j}\}\neq
\emptyset \\i_{j,1}<\ldots<i_{j,l_j}}}\frac1{s!}\cdot
\mu_{l_1}(X_{i_{1,1}},\ldots ,X_{i_{1,l_1}})\cdot \ldots \cdot\mu_{l_s}(X_{i_{s,1}},\ldots ,X_{i_{s,l_s}}),
\plabel{eq:ppod}
\end{equation}
where we sum over ordered partitions.

A continuous recast of the identities \eqref{eq:ppod} is  given by
\begin{theorem}[Magnus \cite{M} (1954), Chen \cite{Ch} (1957)]\plabel{th:MagnusPre}
Let $\phi$ be a continuous $\mathfrak A$-valued measure of finite variation on the interval $I$.
If $T$ is a formal commutative variable, then
\begin{equation}
1+\sum_{k=1}^\infty T^k\cdot\int_{t_1\leq\ldots\leq t_k\in I}\phi(t_1)\cdot\ldots\cdot\phi(t_k)
=\exp\sum_{k=1}^\infty T^k\cdot\int_{t_1\leq\ldots\leq t_k\in I}\mu_k(\phi(t_1),\ldots,\phi(t_k)).
\plabel{eq:magPre}
\end{equation}
\begin{proof} Take the exponential on the RHS, and contract the terms of order $k$ according to \eqref{eq:ppod}.
(This is understood so that, in each degree, $X_j$ can be corresponds to the $j$th chronological contribution under the integral sign.)
\end{proof}
\end{theorem}
In a truly formal setting, if $\phi$ is ``sufficiently formal'', then $T=1$ can be substituted.
(In particular, this applies if $\mathfrak A$ is graded, complete in appropriate sense, and $\phi$ is graded positively.)
In the Banach algebraic setting, convergence is a more delicate matter.
\begin{theorem}[Magnus \cite{M} (1954), Chen \cite{Ch} (1957)]\plabel{th:Magnus}
Let $\phi$ be a continuous $\mathfrak A$-valued measure of finite variation on the interval $I$.
If
\begin{equation}
\sum_{k=1}^\infty\left|\int_{t_1\leq\ldots\leq t_k\in I}\mu_k(\phi(t_1),\ldots,\phi(t_k))\right|<+\infty,
\plabel{eq:sersum}
\end{equation}
then
\begin{equation}
1+\sum_{k=1}^\infty\int_{t_1\leq\ldots\leq t_k\in I}\phi(t_1)\cdot\ldots\cdot\phi(t_k)
=\exp\sum_{k=1}^\infty\int_{t_1\leq\ldots\leq t_k\in I}\mu_k(\phi(t_1),\ldots,\phi(t_k)).
\plabel{eq:mag}
\end{equation}
\begin{proof} Substitute $T=1$ into \eqref{eq:magPre}, which is possible due to absolute convergence.
\end{proof}
\end{theorem}
\renewcommand{\qedsymbol}{$\triangle$}
\begin{proof}[Note]
More precisely, Magnus \cite{M} has the information of this amount, however
he is more interested in the commutator recursion for $\mu_k$ than starting with \eqref{eq:mudef} directly.
While being roughly equivalent, Magnus \cite{M} has more of the ODE viewpoint, Chen \cite{Ch} is more algebraical.
Chen \cite{Ch} also exhibits the commutator recursion but less practically.
On the other hand, the idea of locally induced composition, without integrals,
(applied to the BCH expansion) is already present in Dynkin \cite{Dyy} (1949).
\end{proof}
\renewcommand{\qedsymbol}{$\Box$}
Let us digress on some possible logarithmic versions of the theorems above.
\begin{theorem}\plabel{th:MagnusPreLogg}
Let $\phi$ be a continuous $\mathfrak A$-valued measure of finite variation on the interval $I$.
If $T$ is a formal commutative variable, then
\begin{equation}
\log\left(1+\sum_{k=1}^\infty T^k\cdot\int_{t_1\leq\ldots\leq t_k\in I}\phi(t_1)\cdot\ldots\cdot\phi(t_k)\right)
=\sum_{k=1}^\infty T^k\cdot\int_{t_1\leq\ldots\leq t_k\in I}\mu_k(\phi(t_1),\ldots,\phi(t_k)).
\plabel{eq:formmag}
\end{equation}
\begin{proof}
As we deal with formal perturbations of $1$, we can take the logarithm without trouble.
\end{proof}
\end{theorem}
\begin{theorem}\plabel{th:MagnusLogg}
(Logarithmic Magnus formula by analytic continuation.)
Let $\phi$ be a continuous $\mathfrak A$-valued measure of finite variation on the interval $I$.
If \eqref{eq:sersum} holds, and
\begin{equation}
\log\left(1+t^k\sum_{k=1}^\infty\int_{t_1\leq\ldots\leq t_k\in I}\phi(t_1)\cdot\ldots\cdot\phi(t_k)\right)
 \text{ exists for any $t\in[0,1]$, }
\plabel{eq:logcont}
\end{equation}
then
\begin{equation}
\log\left(1+\sum_{k=1}^\infty\int_{t_1\leq\ldots\leq t_k\in I}\phi(t_1)\cdot\ldots\cdot\phi(t_k)\right)
=\sum_{k=1}^\infty\int_{t_1\leq\ldots\leq t_k\in I}\mu_k(\phi(t_1),\ldots,\phi(t_k))
\plabel{eq:maglog}
\end{equation}
holds.
\begin{proof}
For $t\in[0,1]$, let $E(t)=1+t^k\sum_{k=1}^\infty\int_{t_1\leq\ldots\leq t_k\in I}\phi(t_1)\cdot\ldots\cdot\phi(t_k)$ and
 $L(t)=\sum_{k=1}^\infty t^k\int_{t_1\leq\ldots\leq t_k\in I}\mu_k(\phi(t_1),\ldots,\phi(t_k))$.
We know that $E(t)=\exp L(t)$.
For $t\sim0$, we have $E(t)\sim 1$ and $L(t)\sim0$, which also holds in spectral sense.
In particular, for $t\sim0$, we have $\spec E(t)\subset \mathbb C\setminus(-\infty,0]$ and $\spec L(t)\subset\{z\in\mathbb C:|\Ima z|<\pi\}$.
As on these spectral domains $\log$ and $\exp$ are well defined and invert each other, we have $L(t)=\log E(t)$ for $t\sim 0$.
As both $\log E(t)$ and $L(t)$ are analytic in $t\in(0,1)$ and continuous in $t\in[0,1]$, we obtain $\log E(1)=L(1)$ by
 the unicity of analytical continuation.
\end{proof}
\end{theorem}
\begin{remark}\plabel{eq:remlogver}
Without condition \eqref{eq:logcont}, both the existence of $\log$ and the equality in \eqref{eq:maglog} may fail, cf.
 the cases $\phi=\pi\mathrm i \mathbf 1_{[0,1)}$ and $\phi=2\pi\mathrm i \mathbf 1_{[0,1)}$ for $\mathfrak A=\mathbb C$.
\qedremark
\end{remark}

Using a more compact notation, the right Magnus expansion of $\phi$ is
\begin{equation}
\mu_{\mathrm R}(\phi)=\sum_{k=1}^\infty\mu_{k,\mathrm R}(\phi),
\plabel{eq:MagnusR}
\end{equation}
where
\[\mu_{k,\mathrm R}(\phi)=\int_{t_1\leq\ldots\leq t_k\in I}\mu_k(\phi(t_1),\ldots,\phi(t_k));\]
and Theorem \ref{th:Magnus} says that it exponentiates to the right time-ordered exponential
\[\exp_{\mathrm R}(\phi)=1+\sum_{k=1}^\infty\int_{t_1\leq\ldots\leq t_k\in I}\phi(t_1)\cdot\ldots\cdot\phi(t_k)\]
of $\phi$, as long as the Magnus expansion \eqref{eq:MagnusR} is absolute convergent as a series.

One can, of course, formulate  similar statements in terms of the left Magnus expansion $\mu_{\mathrm L}(\phi)$ and
left time-ordered exponential $\exp_{\mathrm L}(\phi)$ by reversing the ordering on $I$.
The two formalisms have the same capabilities (we just have to reverse the measures $\phi\leftrightarrow\phi^\dag$),
but the `$\mathrm R$' formalism is better suited for purely algebraic manipulations,
and the `$\mathrm L$' formalism is better suited to the standard setup of differential equations.

In what follows, instead of 'continuous measure on an interval' we may simply write `ordered measure'.
(Cf. Appendix \ref{sec:mea}.)

\begin{remark}\plabel{rem:simplelog}
A nontrivial estimate, but which involves only a simple trick,
provides the absolute convergence and  the validity of the logarithmic version of the Magnus formula for
\begin{equation}
\int|\phi|<F=1.505241496\ldots,
\plabel{eq:Fbound}
\end{equation}
where $F$ is the positive solution of the equation $\mathrm e^F=3+F$.
($F=-3-W_{-1}(-\mathrm e^{-3})$
in the Lambert function terminology,
which we will occasionally use;
see
Corless,  Gonnet,   Hare,   Jeffrey,    Knuth \cite{CGHJK} or Mező \cite{MeI}
for a description.)

Indeed, when we are to take the logarithm of the time ordered exponential $A=\exp_{\mathrm R}(\phi)$,
we have to consider the meaningfulness of the resolvent expressions  $\frac{A-1}{\lambda+(1-\lambda)A}$
where $\lambda\in[0,1]$
The trick is to (try to) rewrite it as
\[\frac{A-1}{\lambda+(1-\lambda)A}
=\frac{(A-1)(\lambda+(1-\lambda)A^{-1})}{(\lambda+(1-\lambda)A)(\lambda+(1-\lambda)A^{-1})}=
\frac{(A-1)(\lambda+(1-\lambda)A^{-1})}{1+\lambda(1-\lambda)\left(A+A^{-1}-2\right)}.\]

Here $A^{-1}$ exists, and, in fact, $A^{-1}=\exp_{\mathrm R}(-\phi^\dag)$, where the negative reverse measure was taken.
Also note, that when $A+A^{-1}=\exp_{\mathrm R}(\phi)+\exp_{\mathrm R}(-\phi^\dag)$ is taken,
then $\exp_{\mathrm R,1}(\phi)$ and $\exp_{\mathrm R,1}(-\phi^\dag)$ cancel each other.
Then it is easy to see that $ |\lambda(1-\lambda)\left(A+A^{-1}-2\right)|\leq\frac14\cdot 2(\exp(\smallint\phi)-1-\smallint\phi)$.
This latter expression is less than $1$ if the initial norm inequality is satisfied, therefore the resolvent expression makes sense, thus $\log \exp_{\mathrm R}(\phi)$ exists.

Regarding the formal estimate
\begin{align*}
\mu_{\mathrm R}(T\cdot\phi)
&=\int_{\lambda=0}^1 \frac{\exp_{\mathrm R}(T\cdot \phi)-1}{\lambda+(1-\lambda)\exp_{\mathrm R}(T\cdot \phi)}\,\mathrm d\lambda
\\
&=\int_{\lambda=0}^{\frac12} \frac{\exp_{\mathrm R}(T\cdot \phi)-1}{\lambda+(1-\lambda)\exp_{\mathrm R}(T\cdot \phi)}
+\frac{\exp_{\mathrm R}(T\cdot \phi)-1}{(1-\lambda)+ \lambda \exp_{\mathrm R}(T\cdot \phi)}\,\mathrm d\lambda
\\
&=\int_{\lambda=0}^{\frac12} \frac{\exp_{\mathrm R}(T\cdot \phi)-\exp_{\mathrm R}(T\cdot \phi)^{-1} }{
1+\lambda(1-\lambda)\left(\exp_{\mathrm R}(T\cdot \phi)+\exp_{\mathrm R}(T\cdot \phi)-2\right)}
\,\mathrm d\lambda.
\end{align*}
Therefore
 \begin{equation}
\sum_{k=1}^\infty T^k \cdot|\mu_{k,\mathrm R}(\phi)|\stackrel{\forall T}\leq \int_{\lambda=0}^{\frac12}
\frac{2\left(\exp (T\cdot \smallint|\phi|)-1\right)  }{
1-\lambda(1-\lambda)2\left(\exp T\cdot \smallint|\phi|)-1-\smallint|\phi|\right)
}
\,\mathrm d\lambda
\plabel{eq:arc1}
\end{equation}
\[\Biggl(=\AC\Bigl(2+T\cdot(\smallint|\phi|)-\exp(T\cdot\smallint|\phi|)\Bigr)\cdot \left(\exp (T\cdot \smallint|\phi|)-1\right)
\qquad\text{[cf. \eqref{eq:ACdef} and \eqref{eq:ATdef}]}\]
\[=\ATT\left(\frac{\frac12\left(\exp (T\cdot \smallint|\phi|)-1-T\cdot\smallint|\phi|\right)
}{1-\frac12\left(\exp (T\cdot \smallint|\phi|)-1-T\cdot\smallint|\phi|\right) } \right)\cdot
\frac{\left(\exp T(\cdot \smallint|\phi|)-1\right)}{1-\frac12\left(\exp (T\cdot \smallint|\phi|)-1-T\cdot\smallint|\phi|\right) }\quad\Biggr)\]
\[
\stackrel{\forall T}\leq
\frac{\left(\exp T(\cdot \smallint|\phi|)-1\right)}{1-\frac12\left(\exp (T\cdot \smallint|\phi|)-1-T\cdot\smallint|\phi|\right) }  .
\]
Under \eqref{eq:Fbound}, $T=1$ can be plugged in, which shows the absolute convergence of the Magnus expansion.

Ultimately, under \eqref{eq:Fbound}, both the absolute convergence of the Magnus expansion and the existence of the
 logarithms has be demonstrated.
Taking Theorem \ref{th:MagnusLogg} into consideration, this also implies the logarithmic Magnus formula
 under the indicated norm restriction.
\qedremark
\end{remark}

\begin{remark}\plabel{rem:philos}
Here is some philosophy.
The Magnus expansion was obviously developed as the natural solution to the
 logarithm problem of the time-ordered exponentials.
Indeed, if it converges (absolutely) then it provides an accidental logarithm to the time-ordered exponential.
However, its success is at most partial, because the existence of logarithms is more typical than the
 convergence of the Magnus expansion.
 More precisely:
For example, considering large matrix algebras we see the following
``hierarchy'':

(i) The existence of the Magnus expansion.
This is the most universal one.
If $\phi$ is a time-ordered measure of bounded variation, then the associated Magnus expansion exists.

(ii) The existence of accidental logarithms to the time-ordered exponential.
This is ``typical'' in the case of complex matrices, although it is scarcer, or rather ``regional''
in the case of real matrices.

(iii) The existence of the (canonical) logarithm to the time-ordered exponential.
This is obviously more special than (ii) but the same observations apply regarding typicality.

(iv) The (absolute) convergence of the Magnus expansion.
This behaviour is certainly more special than (ii).
It is not more special than (iii) but less typical.
Indeed, (iv) holds ``regionally'' but untypically.
Its nature is not just about the existence of $\log$ but rather about the convergence radius
of the function $t\mapsto \sum_{k=1}^\infty t^k\mu_{k,\mathrm R}(\phi)$,
 implying the existence of a bunch of accidental logarithms for the $\Rexp(t\cdot\phi)$,  in a complex analytical way.
This is the least typical behaviour among the cases here.

One sometimes hears frustration over the lack of universal convergence of the Magnus series but which I find somewhat unreasonable.
$\frac1{1-x}=\sum_{k=1}^{\infty}x^k$, mother of all power series is not universally convergent,
and universal convergence is not typical in complex analysis.
It is like seeing `$\exp$' but not seeing `$\log$'.

For these reasons, in present applications the Magnus expansion provides
a theoretical framework rather for small-time and perturbative developments, cf.
Iserles,   Munthe-Kaas,  Nørsett,  Zanna   \cite{IMKNZ}, Blanes, Casas, Oteo, Ros \cite{BCOR}.
One, however, should not ignore questions of macroscopic convergence.
Partly, because these are connected to microscopical ones, partly because
in their investigations there is a plethora of analytical methods
which are most instructive.  \qedremark
\end{remark}

Although we have already touched upon the logarithmic Magnus formula and the spectral method by
Theorem \ref{th:MagnusLogg}, it will be addressed fully only by the resolvent method later.
In the rest of the section, we concentrate on the (absolute) convergence of the Magnus series.

We can give an estimate for the convergence of the Magnus expansion as follows:
Let $\Theta_k$ be $\frac1{k!}$ times the sum of the absolute value of the coefficients
in the monomial expansion of $\mu_k(X_1,\ldots,X_k)$.
(Cf. \eqref{eq:mpform} later.)
We define the absolute Magnus characteristic as
\[\Theta(x)=\sum_{k=1}^\infty\Theta_kx^k,\]
i.~e.~as the  exponential generating function associated
to the sum of the absolute value of the coefficients in the
monomial expansion of the Magnus commutators.
Then
\begin{equation}
\sum_{k=1}^\infty T^k\cdot\left|\mu_{k,\mathrm R}(\phi)\right|\stackrel{\forall T}
\leq
\sum_{k=1}^\infty T^k\cdot\int_{t_1\leq\ldots\leq t_k\in I}
\left|\mu_k(\phi(t_1),\ldots,\phi(t_k))\right|
\stackrel{\forall T}\leq
\sum_{k=1}^\infty T^k\cdot \Theta_k\cdot\left(\textstyle{\int|\phi|} \right)^k,
\plabel{eq:Tmeasum}
\end{equation}
and, consequently,
\begin{equation}
\sum_{k=1}^\infty\left|\mu_{k,\mathrm R}(\phi)\right|
\leq\sum_{k=1}^\infty\int_{t_1\leq\ldots\leq t_k\in I}\left|\mu_k(\phi(t_1),\ldots,\phi(t_k))\right|
\leq\Theta_\real\left(\textstyle{\int|\phi|}\right).
\plabel{eq:measum}
\end{equation}

Note that we have equalities in the case $\phi=c\cdot\mathrm Z^1_{[0,1]}$, $c\in[0,\infty)$,
the ordered totally noncommutative continuous mass of norm $c$.
See the Appendix \ref{sec:ellone} for a detailed explanation of this measure.
(Actually, the equalities also hold for $c\in\mathbb C$, where $ c\cdot\mathrm Z^1_{[0,1]}$
is equivalent to the ordered totally noncommutative continuous mass of norm $|c|$.)
Even if we do not care about the $\Theta_k$ much, this particular testing example shows that in the realm of general Banach algebras
universal estimates involving the $\left|\mu_{k,\mathrm R}(\phi)\right|$
automatically yield the estimates for the $\int_{t_1\leq\ldots\leq t_k\in I}\left|\mu_k(\phi(t_1),\ldots,\phi(t_k))\right|$.
(This is not necessarily true for more special classes of algebras.)

\begin{remark}
\plabel{rem:multgen}
One can formulate Theorem \ref{th:Magnus}
and the subsequent discussion for not necessarily continuous measures, too.
Then, we have to insert the multiplicity terms
\[\int_{t_1\leq\ldots\leq t_k\in I}\ldots\quad \rightsquigarrow\quad
\int_{t_1\leq\ldots\leq t_k\in I}\frac{1}{\mul(t_1,\ldots,t_n)!}\ldots\quad
\]
everywhere. This is in accordance to taking the continuous blowup $\phi^*$ of $\phi$.
Thus, it is left to the reader to decide whether he wants to incorporate non-continuous measures under the
name of `ordered measure', or even to relax the intervals to ordered sets.
\qedremark
 \end{remark}
$\Theta(x)$ can be worked out as follows.
Firstly, one has a direct expression for $\mu_k$.
In order to have it, let us introduce some terminology.
If $\sigma=(\sigma(1),\ldots,\sigma(k) )$ is a finite sequences of real numbers,
let $\asc(\sigma)$ denote the number of its ascents,
i.~e.~the number of pairs such that $\sigma(i)<\sigma(i+1)$;
and let $\des(\sigma)$ denote the number of its descents,
i.~e.~the number of pairs such that $\sigma(i)>\sigma(i+1)$.
This naturally applies in the special case when
$\sigma$ is a permutation from the symmetric group $\Sigma_k$.
Then $\asc(\sigma)+\des(\sigma)=k-1$.
\begin{theorem}[Mielnik, Pleba\'nski \cite{MP} (1970), Dynkin \cite{Dyy} (1949)]\plabel{th:mp}
\begin{equation}\mu_k(X_1,\ldots,X_k)=
\sum_{\sigma\in\Sigma_k}(-1)^{\des(\sigma)}
\frac{\asc(\sigma)!\des(\sigma)!}{k!}X_{\sigma(1)}\cdot\ldots\cdot X_{\sigma(k)} .
\plabel{eq:mpform}\end{equation}
\begin{proof}(Dynkin \cite{Dyy},  Strichartz \cite{St} (1987).)
Considering \eqref{eq:sinmu}, and the power series of $\log$, one can see that $\mu_k(X_1,\ldots,X_k)$ is a sum of terms
\[\frac{(-1)^{j-1}}{j}
X_{\sigma(1)}\cdot\ldots\cdot X_{\sigma(l_1)}\mid  X_{\sigma(l_1+1)}\cdot\ldots\ldots
\cdot X_{\sigma(l_{j-1})}\mid X_{\sigma(l_{j-1}+1)}\cdot\ldots\cdot X_{\sigma(k)},\]
where separators show the ascendingly indexed components which enter into power series of $\log$.
Considering a given permutation $\sigma$, the placement of the separators is either necessary (in case $\sigma(i)>\sigma(i+1)$),
or optional (in case  $\sigma(i)<\sigma(i+1)$).
Summing over the $2^{\asc(\sigma)}$ many optional possibilities, the coefficient of $X_{\sigma(1)}\cdot\ldots\cdot X_{\sigma(k)}$
is
\[\sum_{p=0}^{\asc(\sigma)}\frac{(-1)^{\des(\sigma)+p}}{\des(\sigma)+1+p}
\begin{pmatrix}\asc(\sigma)\\p\end{pmatrix}.\]
This simplifies  according to the combinatorial identity
$\sum_{p=0}^{r}\frac{(-1)^{p}}{d+1+p}\binom rp=\frac{r!d!}{(d+1+r)!}$.
%(But see the approach of Mielnik and Pleba\'nski later.)
\end{proof}
\end{theorem}
\renewcommand{\qedsymbol}{$\triangle$}
\begin{proof}[Note]
The $\mu_k$ were defined and given explicitly by Dynkin \cite{Dyy}, entirely in the way we presented them;
but, at that time, he applies them only in the setting of the Baker--Campbell--Hausdorff expansion
(cf. Achilles, Bonfiglioli \cite{AB}).
Nevertheless, the `Magnus commutator' could also be called as `Dynkin commutator' or `Dynkin--Magnus commutator'.
The RHS of \eqref{eq:mpform} is redeveloped and scrutinized by Solomon \cite{S} (1968), however not directly
in the context of the Magnus expansion, but as a canonical projection.
(That idea will be explained in Subsection \ref{ss:sol}.)
The explicit formula \eqref{eq:mpform} was introduced into the study of Magnus expansion only by Mielnik, Pleba\'nski \cite{MP},
whose results were announced in Bialynicki-Birula, Mielnik, Pleba\'nski \cite{BMP} (1969) somewhat earlier.
They use the resolvent method, which we discuss in Subsection \ref{ss:IntroResolvent}.
A rediscoverer is Strichartz \cite{St}, who presents the combinatorial argument again.
In this series of papers we simply use the term `Magnus commutator'
(as the relevant term in the Magnus expansion).
\end{proof}
\renewcommand{\qedsymbol}{$\Box$}

The following is very classical, cf. Comtet \cite{Co}, Graham, Knuth, Patashnik \cite{GKP}, or Petersen \cite{P}.
However, we indicate its proof, because it is a prototype argument.
\begin{theorem}[Euler, 1755] \plabel{th:Euler}
Let $A(n,m)$ denote the number of permutations $\sigma\in\Sigma_n$ such that $\des(\sigma)=m$.
(These are the Eulerian numbers.)
 Consider the exponential generating function
\begin{equation}G(u,v)\equiv\sum_{0\leq m< n}^\infty \frac{A(n,m)}{n!} u^{n-1-m}v^m.\plabel{eq:mpgen}\end{equation}
Then: As formal power series in $u,v$,
\[%\begin{align*}
G(u,v)
%&
=
\dfrac{\frac{\tanh\frac{u-v}2}{\frac{u-v}2}}{1-\frac{u+v}2\frac{\tanh\frac{u-v}2}{\frac{u-v}2}}=
\frac{1-\frac{1}{3}\left(\frac{u-v}2\right)^2+\frac{2}{15}\left(\frac{u-v}2\right)^4-\frac{17}{315}\left(\frac{u-v}2\right)^6+\ldots}
{1-\left(\frac{u+v}2\right)+\frac{1}{3}\left(\frac{u+v}2\right)\left(\frac{u-v}2\right)^2-\frac{2}{15}\left(\frac{u+v}2\right)
\left(\frac{u-v}2\right)^4+\ldots}
\]
\begin{commentx}
\[%\\&
=\dfrac{\frac{\sinh\frac{u-v}2}{\frac{u-v}2}}{\cosh\frac{u-v}2-\frac{u+v}2\frac{\sinh\frac{u-v}2}{\frac{u-v}2}}=
\frac{
1-\frac{1}{3!}\left(\frac{u-v}2\right)^2+\frac{1}{5!}\left(\frac{u-v}2\right)^4+\ldots
}
{
1-\left(\frac{u+v}2\right)+\frac{1}{2!}\left(\frac{u-v}2\right)^2-\frac{1}{3!}\left(\frac{u+v}2\right)
\left(\frac{u-v}2\right)^2+\ldots
}
\]
\end{commentx}
\[%\\&
=\frac{\mathrm e^u-\mathrm e^v}{u\mathrm e^v-v\mathrm e^u}=\frac{\cancel{(u-v)}\left(1+\frac{u+v}{2!}
+\frac{u^2+uv+v^2}{3!}+\frac{u^3+u^2v+uv^2+v^3}{4!}+\ldots\right)}{
\cancel{(u-v)}\left(1-\frac{uv}{2!}-\frac{uv(u+v)}{3!}-\frac{uv(u^2+uv+v^2)}{4!}-\ldots\right)}.
\]%\end{align*}
This generating function can be interpreted as an analytic function at $(0,0)$.

\begin{proof}[Note]
The division(al simplification) is valid, because the formal power series form a ring without zero-divisors.
\qedno
\end{proof}
\begin{proof}[Sketch of proof]
Consider the extended generating function
\[G(u,v;x)=\sum_{0\leq m< n}^\infty \frac{A(n,m)}{n!} u^{n-1-m}v^m x^n.\]
Based on the image of "1" in the permutations,  one can develop a recursion for the Eulerian numbers,
which, in terms of the generating function, formally reads as the differential equation (in $x$)
\[ G'(u,v;x)=\left(1+u G(u,v;x)(1+vG(u,v;x))\right),\]
\[G(u,v;0)=0.\]

Now, $G(u,v;x)=xG(ux,vx)$ solves this differential equation not only formally
but also analytically for small $(ux,vx)$.
Substituting $ux\mapsto u, vx\mapsto v, x\mapsto 1$ formally gives the indicated $G(u,v)$
as formal power series, but also analytically for small $(u,v)$.
The analytical behaviour at $(0,0)$ is also transparent from the first, ``tangential'', explicit form.
(Originally, Euler computes $1+tG(t,tz)=\frac{1-z}{\mathrm e^{t(z-1)}-z} $, but uses the same idea.)
\end{proof}
\end{theorem}

\begin{remark}\plabel{rem:Euler}
Using  (the meromorphic extension of) the notation \eqref{eq:defTan}, we have the more compact and less singular presentation
\[G(u,v)=\frac{\Tan\left(-\left(\frac{u-v}2\right)^2\right)}{1-\frac{u+v}2\Tan\left(-\left(\frac{u-v}2\right)^2\right)}.\]
Such notation is more advantageous analytically or in demonstrably avoiding some singularities.
Yet, such compact presentations might obscure the arithmetical combinatorial nature of the given objects.
Thus they will be used less than possible.
\qedremark
\end{remark}

The following lemma is a natural companion statement concerning $G(u,v)$.
\begin{lemma}
\plabel{lem:coEuler}
(a) The generating function
is analytically equal to, i.~e.~it has the same Taylor expansion at $(0,0)$ as, the function
\[G_{\anal}(u,v)=\begin{cases}
\frac{\mathrm e^u-\mathrm e^v}{u\mathrm e^v-v\mathrm e^u}&u\neq v,\quad \frac{\mathrm e^u}u\neq\frac{\mathrm e^v}v\\
\infty&u\neq v,\quad \frac{\mathrm e^u}u=\frac{\mathrm e^v}v\\
\frac1{1-u}=\frac1{1-v}&u= v\neq1,\\
\infty&u=v=1,
\end{cases}\]
which is meromorphic (meant: analytic with respect to the range Riemann sphere).

This function is analytic around $(0,0)$, with poles in the real quadrant $u,v>0$  at
\[u+v=\begin{cases}\frac{\frac uv+1}{\frac uv-1}\log \frac uv&\frac uv\neq1,\\2&\frac uv=1,\end{cases}\]
using the `rational polar' coordinates $u+v,\frac uv$.

(b)  Regarding absolute convergence, for $u,v\geq 0$, the generating function gives
the real function $G_{\real}(u,v)$, which is finite and equal to $G_\anal(u,v)$ if
 \[u=0 \quad \text{or}\quad v=0\quad\text{or}\quad  u+v<
 \begin{cases}\frac{\frac uv+1}{\frac uv-1}\log \frac uv&\frac uv\neq1,\\2&\frac uv=1,\end{cases}\]
 and $ G_\real(u,v)=+\infty$ otherwise.

(c) In particular it can be said that $G(u,v)$ is analytic and absolute convergent for $|u|+|v|<2$
but these fail to be true for $u=v=1$.
\begin{proof}
(a) The ``exponential form'' $\frac{\mathrm e^u-\mathrm e^v}{u\mathrm e^v-v\mathrm e^u}$ is a quotient of analytic functions and cannot yield $\frac00$ unless $u=v$.
For (near) $u=v$, however, $\frac{\tanh\frac{u-v}2}{\frac{u-v}2}$ makes sense as an analytical
function with value (near) $1$, thus there the ``tangential form'' is a quotient of analytic functions and cannot yield $\frac00$.
Hence it is safe to call the function meromorphic.
The nature of the poles follows from elementary considerations.
(b) The statement about absolute convergence follows from restricting
to radial directions and considering the general properties of power series with nonnegative coefficients.
(c) follows from (b).
\end{proof}
\end{lemma}
Mielnik and Pleba\'nski, who have the generating function $G(u,v)$,
are not interested in metric estimates, thus they miss to give the following
\begin{theorem} \plabel{th:mpest} The absolute Magnus characteristic is given by
\begin{equation}
\Theta(x)=\int_{\lambda=0}^1 x G(\lambda x,(1-\lambda)x)\,\mathrm d\lambda;
\plabel{eq:mpest}
\end{equation}
and
\[\Theta_\real(x)=\int_{\lambda=0}^1 x G_\real(\lambda x,(1-\lambda)x)\,\mathrm d\lambda
=\int_{y=0}^x G_\real(y,x-y)\,\mathrm dy.\]
In particular, $\Theta_\real(x)<+\infty$ if $0\leq x<2$; and $\Theta_\real(x)=+\infty$ if $2\leq x$.
\begin{proof} According to \eqref{eq:mpform},
\[\Theta(x)=\sum_{0\leq m< n}^\infty \frac{A(n,m)}{n!} \cdot \frac{(n-1-m)!m!}{n!} x^n.\]
Then the first equality in the statement follows from the beta function identity
\[ \frac{(n-1-m)!m!}{n!}=\int_{\lambda=0}^1 \lambda^{n-1-m} (1-\lambda)^m\,\mathrm d\lambda.\]
In the real case, summation is valid, as the terms are nonnegative.
Now, $\Theta_\real(x)<+\infty$ if $0\leq x<2$,
because the integrand is finite and continuous (it is easy to see that for the poles $u+v\geq2$).
On the other, hand,
\[\Theta_\real(2)=\int_{\lambda=0}^12G_\real(2\lambda,2(1-\lambda))\,d\lambda= \int_{\nu=-1}^{1}G_\real(1+\nu,1-\nu)\,d\nu.\]
The integrand is continuous away from $\nu=0$, but it has a double pole there,
\[G_\real(1+\nu,1-\nu)=\frac{\tanh  \nu}{\nu-\tanh  \nu}\sim3\nu^{-2}+\text{holomorphic in }\nu.\]
This implies $\Theta_\real(2)=+\infty$.
\end{proof}
\end{theorem}
As corollaries, we have
\begin{theorem}[Moan, Oteo \cite{MO}, 2001]
\plabel{cor:moo}
In the Banach algebra setting, the Magnus expansion converges absolutely
%(even as a time-ordered integral, cf. \eqref{eq:measum})
if $\int|\phi|<2$.
\qed
\end{theorem}
\begin{commentx}
\renewcommand{\qedsymbol}{$\triangle$}
\begin{proof}[Note]
In effect, Moan and Oteo generalize  the estimate of Thompson \cite{Th}, who is using the results of Goldberg \cite{G} in the BCH case,
to the Magnus case but using the results of Mielnik, Pleba\'nski \cite{MP} .
They formulate convergence on the unit interval, for $L^\infty$ norm $<2$.
\end{proof}
\renewcommand{\qedsymbol}{$\Box$}
\end{commentx}
\begin{theorem}\plabel{cor:mooexact}
In the Banach algebra setting, the Magnus expansion diverges
if $\phi=2\cdot\mathrm Z^1_{[0,1)}$, in which case $\int|\phi|=2$.
\begin{proof}
In this case the Magnus expansion converges in $\mathrm F^{1,\mathrm{loc}}_{\mathbb R}([0,1))$ but to an element whose
conventional norm is $+\infty$. This implies divergence in $\mathrm F^1_{\mathbb R}([0,1))$.
\end{proof}
\renewcommand{\qedsymbol}{$\triangle$}
\begin{proof}[Note]
This counterexample is a measure and a not nice function (times the Lebesgue measure).
This is not accidental; it will be addressed subsequently.
\end{proof}
\renewcommand{\qedsymbol}{$\Box$}
\end{theorem}

The $n$-summand in \eqref{eq:mpgen},
\begin{equation}
G_n(u,v)\equiv\sum_{m=0}^{n-1}\frac{A(n,m)}{n!} u^{n-1-m}v^m,
\plabel{eq:Eulerpol}
\end{equation}
is, by definition, an Eulerian polynomial (in bivariate form) divided by $n!$.
It gives the individual coefficient
\[\Theta_k=\int_{\lambda=0}^1 G_k(\lambda,1-\lambda)\,\mathrm d\lambda.\]

For later reference, we give here some terms:
\begin{equation}
\Theta(x)=x+\frac{1}{2}x^2+\frac{2}{9}x^3+\frac{7}{72}x^4+\frac{13}{300}x^5
+\frac{71}{3600}x^6+\frac{67}{7350}x^7
+O(x^8).\plabel{eq:thetaexp}
\end{equation}
We also give a crude but explicit estimate for $\Theta_k$ as follows.

For $\lambda\in[0,1]$, let
\begin{equation}\Theta^{(\lambda)}(x)=xG(\lambda x,(1-\lambda)x).
\plabel{eq:ThetaLdefpre}
\end{equation}

Note that $\Theta^{(\lambda)}(x)$ solves the differential equation (IVP)
\begin{equation}
\Theta^{(\lambda)\prime}(x)=(1+\lambda \Theta^{(\lambda)}(x) )(1+(1-\lambda) \Theta^{(\lambda)}(x) ),
\plabel{eq:Thetaeqpre}
\end{equation}
\[\Theta^{(\lambda)}(0)=0.\]
(Of course, $\Theta^{(\lambda)}(x)$  could be expressed completely expicitly, cf. later.)

In this notation, as for a formal power series in $x$,
\[\Theta(x)=\int_{\lambda=0}^1\Theta^{(\lambda)}(x)\,\mathrm d\lambda.\]
\begin{lemma}\plabel{lem:ThetaMest} For $\lambda\in[0,1]$,
\begin{equation}\Theta^{(\lambda)}(x)\stackrel{\forall x}\leq\Theta^{(1/2)}(x)\equiv\frac{x}{1-\frac12x}.\plabel{eq:ThetaLest}\end{equation}
Consequently,
\begin{equation}\Theta(x)\stackrel{\forall x}\leq\Theta^{(1/2)}(x)\equiv\frac{x}{1-\frac12x};\plabel{eq:ThetaMest}\end{equation}
which can be written in other terms as
\begin{equation}\Theta_k\leq2^{1-k}.\plabel{eq:ThetaMestdis}\end{equation}
\begin{proof}
Then $\Theta^{(\lambda)}(x)$ solves the differential equation (IVP)
\[\Theta^{(\lambda)\prime}(x)=1+ \Theta^{(\lambda)}(x)+\lambda(1-\lambda) \Theta^{(\lambda)}(x)^2 ,\]
\[\Theta^{(\lambda)}(0)=0.\]
Thinking about this as a recursion for the nonnegative coefficients of the formal variable $x$,
it is easy to prove that the greatest growth is achieved when $\lambda(1-\lambda)\in[0,\frac14]$ is maximal,
which is at $\lambda=\frac12$.
This proves \eqref{eq:ThetaLest}.
Integrated, it yields \eqref{eq:ThetaMest}.
In terms of the individual coefficients of the power series,
we obtain \eqref{eq:ThetaMestdis}.
\end{proof}
\end{lemma}
The individual estimate \eqref{eq:ThetaMestdis} appears in Moan, Oteo \cite{MO} first; they use it to prove Theorem \ref{cor:moo}.
This closes the discussion of known results regarding the convergence of the Magnus expansion in the Banach algebraic case.
In Section \ref{sec:MagnusResolvent}, using the resolvent method, we reobtain these results in sharper form.
\\~

\snewpage
\subsection{The Baker--Campbell--Hausdorff expansion}
~\\

After discussing the convergence of the Magnus expansion, let us consider the relationship to the Baker--Campbell--Hausdorff expansion.
A special case of the Magnus expansion is when
$\phi=X\mathbf 1_{[0,1)}\boldsymbol.Y\mathbf 1_{[0,1)}$, i.~e.~when we take the constant function $X$ on one
unit interval, take the constant function $Y$ on another one, and concatenate them.
(The first term, $X\mathbf 1_{[0,1)}$ comes in lower, the second term $Y\mathbf 1_{[0,1)}$ comes in
higher in the ordering of the new interval.)
Then the Magnus integral of order $n$ immediately specifies to
\begin{equation}\BCH_n(X,Y)=\sum_{j=0}^n\frac1{j!(n-j)!}\mu_n(\underbrace{X,\ldots,X}_{j\,\text{ terms}},
\underbrace{Y,\ldots,Y}_{n-j\text{ terms}}).\plabel{eq:ABCH}\end{equation}
(Another equivalent viewpoint is that we consider measures supported only at two points, one with
mass $X$ and one with mass $Y$.)
Thus, from the Magnus formula we obtain
\[\exp(X)\exp(Y)=\exp\left(\sum_{n=1}^\infty \BCH_n(X,Y) \right),\]
which is valid, as long as the sum in the RHS converges absolutely.
This situation also adapts to formal variables $X,Y$;
just take certain nilpotent elements in truncated associative free algebras to see that
\[\log(\exp(X)\exp(Y))=\BCH(X,Y)\equiv\sum_{n=1}^\infty \BCH_n(X,Y)\]
in formal sense. $\BCH_n(X,Y)$ collects exactly the terms of degree $n$.
(We must note that formula \eqref{eq:ABCH} was already obtained by Dynkin \cite{Dyy}.)
This yields some explicit formulas for the terms of the BCH expansion,
see Theorem \ref{th:Goldberg}.

One can  discuss the BCH expansion $\BCH(X_1,\ldots,X_k)=\log((\exp X_1)\ldots (\exp X_k))$
with several variables in similar manner; the formulas are analogous.

Even without the explicit knowledge of the terms
$\BCH(X_1,\ldots,X_k)$, we can derive some qualitative statements about the
convergence of the Magnus and BCH expansions.

\textbf{(O)}
Let $\Gamma(X,Y)$ be the same noncommutative formal series as $\BCH(X,Y)$ but all coefficients turned into nonnegative.
We can apply this series with commutative variables, or with nonnegative real variables.
Note that
\[\Gamma_\real(x_1,\ldots,x_k)<+\infty\]
($x_i\in[0,+\infty$)) implies the $\BCH(X_1,\ldots,X_k)$ will surely converge as long $|X_i|\leq x_i$.
On the other hand, if
\[\Gamma_\real(x_1,\ldots,x_k)=+\infty,\]
then there is counterexample, where $|X_i|=x_i$, and  $\BCH(X_1,\ldots,X_k)$ does not converge;
in fact there is no element whose exponential is the product $(\exp X_1)\ldots(\exp X_k)$.
Indeed, $X_i=x_iY_i$  can  be taken in the setting of $\mathrm F^1_{\mathbb R}[ Y_i\,:\,1\leq i \leq k]$.

\textbf{(I)}
The first observation, related to convergence of the Magnus expansion,
is that the (multivariable) BCH expansions are special cases of the Magnus expansion,
thus the estimates of the Magnus expansion also apply to them.
In particular, the formal expansion of $\log(\exp(X)\exp(Y))$ converges absolutely in the Banach algebra setting if $|X|+|Y|<2$.

\textbf{(II)}
The second observation, related to divergence of the Magnus expansion, is that while the counterexample $\phi=2\cdot\mathrm Z^1_{[0,1)}$
is only a general measure (and not a function times a Lebesgue measure), it is possible to divide it to by small
intervals into parts $2\cdot\mathrm Z^1_{[c_i,c_{i+1})}$, and then replace
those parts by $\mu_{\mathrm R}(2\cdot\mathrm Z^1_{[c_i,c_{i+1})})\cdot\mathbf 1_{[c_i,c_{i+1})}$.
In that way, the Magnus expansion of the resulted measure $\phi_C$ is still convergent in
$\mathrm F^{1,\mathrm{loc}}_{\mathbb R}([0,1))$ to the very same element as of $\phi$,
thus it is divergent $\mathrm F^1_{\mathbb R}([0,1))$.
Then $\phi_C$ is of (multivariable) BCH type, which can be actually be smoothed out completely; yet the total variation of
$\phi_C$ is $\sum_i\Theta_{\real}(2(c_{i+1}-c_i) )$, which can be arbitrarily close to $2$, as
$\Theta_\real(x)\sim x$ for small $x$.
Thus we have a nice set of counterexamples for the Magnus expansion although not with cumulative norm $2$ but with $2+\varepsilon$.

In particular, a numerical consequence is
\begin{lemma}
\plabel{rem:GammaDec}
 For $x_i\in[0,+\infty)$,
\[\Gamma_\real(x_1,\ldots,x_k)\leq \Theta_\real(x_1+\ldots+x_k)\leq \Gamma_\real(\Theta_\real(x_1),\ldots,\Theta_\real(x_k)).\]
\begin{proof}
The first inequality follows from applying the Magnus expansion estimate
to $Y_1\mathbf 1_{[0,x_1)}\bo.\ldots\bo.Y_k\mathbf 1_{[x_1+\ldots+x_{k-1},x_1+\ldots+x_k)}$
in $\mathrm F^{1,\mathrm{loc}}_{\mathbb R}[ Y_i\,:\,1\leq i \leq k]$.
The second one follows from the trivial $\Gamma$-estimate for the BCH expansion with respect to
$X_1=\mu_{\mathrm R}(\mathrm Z^1_{[0,x_1)}),\ldots,$ $X_k=\mu_{\mathrm R}(\mathrm Z^1_{[x_1+\ldots+x_{k-1},x_1+\ldots+x_k)})$
in $\mathrm F^{1,\mathrm{loc}}_{\mathbb R}([0,x_1+\ldots+x_k))$.
\end{proof}
\end{lemma}
%\begin{commentx}
As a special case, we find that
\[\Gamma_\real\Biggl(\underbrace{\Theta_\real\left(\frac2n\right),\ldots,\Theta_\real\left(\frac2n\right)}_{n\text{ times}}\Biggr)=+\infty.\]

If  $n\nearrow\infty$, then the cumulative norm $n\Theta\left(\frac2n\right)\searrow2$.
By \eqref{eq:ThetaMest}, $\Theta\left(\frac2n\right)$ can be replaced  by $\frac{2}{n-1}$,
leading to nicely quantified counterexamples for the Magnus expansion.
(Sharper results will be met in Section \ref{sec:DiscRes}.)

Taking a closer look, from the explicit form of Magnus expansion one can obtain
\begin{theorem}[Goldberg \cite{G}, 1956] \plabel{th:Goldberg}
The coefficient of the monomial
\[M=(X \vee Y)^{k_1}\cdot\ldots\cdot X^{k_i}Y^{k_{i+1}}\cdot\ldots\cdot(X \vee Y)^{k_p}\]
($X$ and $Y$ alternating) in $\log(\exp(X)\exp(Y))$ is
\begin{equation}
c_M=\int_{t=0}^1 t^{\asc(M)} (t-1)^{\des(M)}G_{k_1}(t,t-1)\cdot\ldots\cdot G_{k_p}(t,t-1)\,\mathrm dt,
\plabel{eq:goldberg}
\end{equation}
where
$\asc(M)$ is number of consecutive $XY$ pairs in $M$,
$\des(M)$ is number of consecutive $YX$ pairs in $M$,
and $G_n(u,v)$ is as in \eqref{eq:Eulerpol}.
\begin{proof} (Cf.~Helmstetter, \cite{H}.)
Let $\deg_X(M)$ be sum of the exponents $k_i$ belonging to $X$, and let $\deg_Y(M)$ be sum of the exponents $k_i$ belonging to $Y$.
Consider \eqref{eq:mpform} in the case when the first $\deg_X(M)$ many variables
are substituted by $X$, and the remaining  $\deg_Y(M)$ many variables are substituted by $Y$.
Examine those permutations which lead to $M$, and compute the generating polynomial of their
ascents ($u$) and descents ($v$).
$M$ itself introduces ordered partitions of the variables. There are
$\frac{\deg_X(M)!\deg_Y(M)! }{k_1!\cdot\ldots\cdot k_p!}$ many possible partitions.
Inside each partition set coming from $(X\vee Y)^{k_i}$, the generating polynomial is $k_i!G_{k_i}(u,v)$;
and there are $u^{\asc(M)}$ and $v^{\des(M)}$ coming from the boundaries between the partitions.
Thus the generating polynomial is
\[\deg_X(M)!\deg_Y(M)!\, u^{\asc(M)}v^{\des (M)}G_{k_1}(u,v)\cdot\ldots\cdot G_{k_p}(u,v).\]
We obtain the coefficient of $M$ in $\mu_k(X,\ldots,X,Y,\ldots,Y)$
by replacing the terms $u^av^b$ with $(-1)^b\frac{a!b!}{(a+1+b)!}$, respectively.
According to the beta function identity, this corresponds exactly to integration
of the $u\mapsto t$ and $v\mapsto t-1$ substituted expression as in the statement.
Then, according to \eqref{eq:ABCH}, we have to divide by $\deg_X(M)!\deg_Y(M)! $ in order to get the coefficients in the BCH expansion.
\end{proof}
\end{theorem}
Originally, in \cite{G}, although truly well-informed about the relationship to the Eulerian polynomials,
Goldberg primarily develops the polynomials $G_n(\lambda,\lambda-1)$ via the following recursion, which we now supply as a practical lemma.
\begin{lemma}
\plabel{lem:goldgen}
The polynomials $G_n(\lambda,\lambda-1)$ are generated by the recursion
\[G_1(\lambda,\lambda-1)=1;\]
\[G_n(\lambda,\lambda-1)=\frac1n\,\frac{\mathrm d}{\mathrm d\lambda}\left(\lambda(\lambda-1)G_{n-1}(\lambda,\lambda-1)\right)
\quad\text{ for }\quad n>1.\]

In particular,
\[\int_{\lambda=0}^1G_n(\lambda,\lambda-1)\,\mathrm d\lambda=\delta_{n,1}.\]
\begin{proof}
Recall that
\[n!G_n(\lambda,\lambda-1)
=\sum_{m=0}^{n-1}A(n,m)\lambda^{n-1-m}(\lambda-1)^m
=\sum_{\sigma\in\Sigma_n}\lambda^{\asc(\sigma)}(\lambda-1)^{\des(\sigma)}.\]
Let us call this now as the generating polynomial.
Let us identify the permutation $\sigma\in\Sigma_n$ by the sequence $\sigma(1),\ldots,\sigma(n)$.
Let us consider how to obtain a permutation $\tilde \sigma\in\Sigma_{n+1}$ by inserting a cipher `$n+1$'.
We can put `$n+1$' to the front, that introduces a descent, the generating polynomial is multiplied $\lambda-1$.
We can put `$n+1$' to the end, that introduces an ascent, the generating polynomial is multiplied $\lambda$.
We insert `$n+1$' to the middle, that replaces a factor $\lambda$ or $\lambda-1$ by $1$  (that is the effect of derivation)
but after that the polynomial is multiplied $\lambda(\lambda-1)$.
Ultimately, we obtain
\begin{multline*}
(n+1)!G_{n+1}(\lambda,\lambda-1)=\\=\lambda n!G_n(\lambda,\lambda-1) +(\lambda-1)n!G_n(\lambda,\lambda-1)
+\lambda(\lambda-1)\frac{\mathrm d}{\mathrm d\lambda}\left(n!G_n(\lambda,\lambda-1)\right),
\end{multline*}
which is equivalent to the recursion step above.
The statement about the integral follows from the Newton--Leibniz formula.
\end{proof}
\end{lemma}

Considering the explicit nature of \eqref{eq:goldberg}, and the fact that Goldberg \cite{G} also computes the generating functions
of the coefficients (with fixed $p$), in theory, one could obtain much sharper estimates for $\Gamma_\real(x,y)$.
However, in practice, this is not entirely straightforward.
For example,   Thompson \cite{Th} (in some sense, the predecessor of Moan, Oteo \cite{MO})  obtains convergence for $|X|,|Y|<1$ only.
Arguments in Section \ref{sec:MagnusResolvent} will already imply that the bound $2$ is not sharp in the BCH case,
but, in Section \ref{sec:DiscRes} we obtain more precise information about the convergence domain.

\begin{remark}
\plabel{rem:comput}
In the literature, one can find clever ways
to compute the coefficients of the BCH expansion; see Reinsch \cite{Re} or Van-Brunt, Visser \cite{VBV}.

At this point, however, we would like to argue that the computation of the BCH expansion up to order $n$
is (in relative sense) not extremely difficult or resource-consuming  even for computing through $\log(\exp(x)\exp(y))$ naively.
The following observations  will be imprecise, as
under `space' we mean essentially the number of coefficients organized,
and under `time' we mean the number rational arithmetical operations without penalizing the actual complexity
of the rational numbers; but they may be of practical value.

In the author's experience BCH expansion is best to be treated as a noncommutative polynomial (as we do in this paper).
The write-time (and write-space) of the BCH expansion up to order $n$ can be considered as $O(2^n)$.
We call an algorithm relatively polynomial (for this write-time of the result) if it is of time  $O(p(n)2^n)$
where $p$ is a polynomial.
It is easy to see that computing $\log(\exp(x)\exp(y))$ is of this class.
(In general, both $\log$ and $\exp$ and multiplication are this class if we consider non-commutative
polynomials in $2$ variables, up to order $n$.
Furthermore, these all can be organized in space  $O(2^n)$.)
For non-commutative series, these relatively polynomial algorithms are quite acceptable.
(The recovery time of the individual coefficients is, of course, tiny compared to this.)
However, making optimalizations on the BCH series, or specifying them from the Magnus expansion by \eqref{eq:ABCH}
is indeed much more resource-consuming.
\qedremark
\end{remark}
~

\snewpage
\subsection{Some algebraic properties of the Magnus commutators}
\plabel{ss:sol}
~\\

Here we discuss some algebraic properties of the Magnus commutators.
The content of this subsection will not be used for our main convergence estimates;
and it will be used only Subsection \ref{ss:adjaction} via
Theorem \ref{th:magcomm}, which will be redeveloped in Part III later.

\begin{lemma}\plabel{lem:magchar}
The polynomials $\mu_n(X_1,\ldots,X_n)$ ($n\geq1$) satisfy the properties

(P0) $\mu_n$ is linear in its variables.

(P1) $\mu_1(X)=X$; and
\begin{equation}\mu_n(X,\ldots,X)=0 \qquad\text{for $n>1$, }\plabel{eq:muhigh}\end{equation}
or, equivalently,
\begin{equation}\sum_{\sigma\in\Sigma_n}\mu_n(X_{\sigma(1)},\ldots,X_{\sigma(n)})=0
\qquad\text{for $n>1$. }\plabel{eq:muhighver}\end{equation}

(P2) For $1\leq k<n$,
\[\mu_n(\ldots_1,X,Y,\ldots_2)-\mu_n(\ldots_1,Y,X,\ldots_2)=\mu_{n-1}(\ldots_1,[X,Y],\ldots_2). \]
\begin{proof}[Note]
Formula \eqref{eq:muhighver} follows from \eqref{eq:muhigh} by polarization,
i.~e.~by replacing $X$ with $t_1X_1+\ldots+t_nX_n$ (the $t_i$ are commuting ``scalars'')
and taking the coefficient of $t_1\cdot\ldots\cdot t_n$.
Conversely, \eqref{eq:muhigh} follows from \eqref{eq:muhighver} by replacing all $X_i$ with $X$ and dividing by $n!$.
\qedno
\end{proof}
%\snewpage
\begin{proof}~[This proof uses the definition of $\mu_n$ through the exponentials.]
Recall that
\begin{equation}
\mu_n(X_1,\ldots,X_n)=\text{the coefficient of $t_1\cdot\ldots\cdot t_n$ in } \log(\exp(t_1X_1)\cdot\ldots\cdot\exp(t_nX_n))
\plabel{eq:mron}
\end{equation}

(P1) Specifying all $X_i$ to $X$, we find that
\[\mu_n(X ,\ldots,X )=\text{the coefficient of $t_1\cdot\ldots\cdot t_n$ in }  (t_1+\ldots+t_n)X.  \]
This immediately implies the statement.

(P2) Due to the nature of the problem, we can work modulo $t_1^2,\ldots,t_n^2$.
In particular, in \eqref{eq:mron}, we can replace
\[\exp(t_kX_k)\exp(t_{k+1}X_{k+1}) \rightsquigarrow
\exp(t_{k+1}X_{k+1})\exp(t_kX_k)\exp(t_kt_{k+1}[X_k,X_{k+1}]).\]
Indeed, it is sufficient to check this up to order $2$ as in higher orders the term will be multiples
of either $t_k^2$ or $t_{k+1}^2$.
Let us consider now
\begin{multline}
\text{the coefficient of $t_1\cdot\ldots\cdot t_n$ in }\log\bigl(\exp(t_1X_1)\ldots\exp(t_{k-1}X_{k-1}) \\
\cdot
\exp(t_{k+1}X_{k+1})\exp(t_kX_k)\exp(t_kt_{k+1}[X_k,X_{k+1}])\cdot
\exp(t_{k+2}X_{k+2})\ldots\exp(t_{n}X_{n})
\bigr).
\plabel{eq:mronver}
 \end{multline}
As $\exp$ and $\log$ can be interpreted as formal power series the logarithmic expression
is a formal series of $t_1 X_1,\ldots, t_nX_n$ and $t_kt_{k+1}[X_k,X_{k+1}]$ (that is a sum of monomials of these terms).
If the very latter term $t_kt_{k+1}[X_k,X_{k+1}]$ makes a contribution then, it will render the terms $t_k X_k$ and $ t_{k+1}X_{k+1}$
inactive, that is  $\exp(t_k X_k)$ and $\exp( t_{k+1}X_{k+1})$ can be replaced by $1$.
If the  term $t_kt_{k+1}[X_k,X_{k+1}]$ makes no contribution then $\exp(t_kt_{k+1}[X_k,X_{k+1}]  )$ can be replaced by $1$.
Thus, \eqref{eq:mronver} is equal to
\begin{multline}
\text{the coefficient of $t_1\cdot\ldots\cdot t_n$ in }\log\bigl(\exp(t_1X_1)\ldots\exp(t_{k-1}X_{k-1}) \\
\cdot
\exp(t_{k+1}X_{k+1})\exp(t_kX_k)\cdot
\exp(t_{k+2}X_{k+2})\ldots\exp(t_{n}X_{n})
\bigr)
\plabel{eq:mronver2}
 \end{multline}
\begin{multline*}
+\text{the coefficient of $t_1\cdot\ldots\cdot t_n$ in }\log\bigl(\exp(t_1X_1)\ldots\exp(t_{k-1}X_{k-1}) \\
\cdot
\exp(t_kt_{k+1}[X_k,X_{k+1}])\cdot
\exp(t_{k+2}X_{k+2})\ldots\exp(t_{n}X_{n})
\bigr).
 \end{multline*}
This is then equal to
\begin{multline}
\mu_n(X_1,\ldots,X_{k-1},X_{k+1},X_{k },X_{k+2},\dots,X_n)+\\+
\mu_{n-1}(X_1,\ldots,X_{k-1},[X_{k},X_{k+1 }],X_{k+2},\dots,X_n).
\plabel{eq:mronf}
 \end{multline}
Ultimately, we find that \eqref{eq:mron} is equal to \eqref{eq:mronf}, which is same as the statement of (ii).
\end{proof}
%\snewpage
\begin{proof}[Alternative proof] [This proof uses the explicit formula \eqref{eq:mpform}.]
Using the beta function identity, we can write
\[\mu_n(X_1,\ldots,X_n)=\int_{\lambda=0}^1\sum_{\sigma\in\Sigma_n}
\lambda^{\asc(\sigma)}(\lambda-1)^{\des(\sigma)}X_{\sigma(1)}\cdot \ldots\cdot X_{\sigma(n)}\,\mathrm d\lambda.  \]

(P1)
\begin{commentx}
We can argue that (formally)
\begin{multline*}
\sum_{n=1}^\infty\frac{\mu_n(X ,\ldots,X )}{n!}=
\int_{\lambda=0}^1\left(\sum_{n=1}^\infty G_{n}(\lambda ,\lambda-1)X^n\right)\,\mathrm d\lambda=
\int_{\lambda=0}^1   G(\lambda ,\lambda-1;X)  \,\mathrm d\lambda=\\=
\int_{\lambda=0}^1 \frac{(\exp X)-1}{\lambda+(1-\lambda)(\exp X)}  \,\mathrm d\lambda=
\bigl[-\log\left(\lambda+(1-\lambda)\exp X\right)\bigr]_{\lambda=0}^1=X.
\end{multline*}

Or, alternatively:
\end{commentx}
Applying Lemma \ref{lem:goldgen},
\[
\frac{\mu_n(X ,\ldots,X )}{n!}=\left(\int_{\lambda=0}^1G_{n-1}(\lambda ,\lambda-1)\,\mathrm d\lambda\right)X^n=\delta_{n,1}X^n.
\]

(P2) Let $\iota$ be the transposition sending `$k$' to `$k+1$' and vice versa, but leaving the other ciphers the same.
Then
\begin{align*}
\mu_n(X_1,\ldots,X_{k-1},&X_{k+1},X_{k },X_{k+2},\dots,X_n)=\\
=&\int_{\lambda=0}^1\sum_{\sigma\in\Sigma_n}
\lambda^{\asc(\sigma)}(\lambda-1)^{\des(\sigma)}X_{\iota\circ\sigma(1)}\cdot \ldots\cdot X_{\iota\circ\sigma(n)}\,\mathrm d\lambda\\
=&\int_{\lambda=0}^1\sum_{\sigma\in\Sigma_n}
\lambda^{\asc(\iota\circ\sigma)}(\lambda-1)^{\des(\iota\circ\sigma)}X_{\sigma(1)}\cdot \ldots\cdot X_{\sigma(n)}\,\mathrm d\lambda.
\end{align*}
Now, the ascent and descent numbers of $\sigma$ and $\iota\circ\sigma$ are the same except when $k$ and $k+1$ are taken in neighboring
positions (i.~e.~$|\sigma^{-1}(k)-\sigma^{-1}(k+1)|=1$.
In this case let $\sigma'$ be the same as the sequence $\sigma(1),\ldots,\sigma(n)$ but the consecutive terms
$k,k+1$ or $k+1,k$ are contacted to $k+\frac12$.
Also, in such cases,
\begin{multline*}
\lambda^{\asc(\sigma)}(\lambda-1)^{\des(\sigma)}-\lambda^{\asc(\iota\circ\sigma)}(\lambda-1)^{\des(\iota\circ\sigma)}=
\\=\lambda^{\asc(\sigma')}(\lambda-1)^{\des(\sigma')}\cdot\begin{cases}
+1\qquad\text{if }\sigma^{-1}(k)<\sigma^{-1}(k+1)\\
-1\qquad\text{if }\sigma^{-1}(k)>\sigma^{-1}(k+1).
\end{cases}
\end{multline*}
Thus we find
\begin{multline}
\mu(X_1,\ldots,X_n)-\mu_n(X_1,\ldots,X_{k-1},X_{k+1},X_{k },X_{k+2},\dots,X_n)=\\
 =\int_{\lambda=0}^1\sum_{\sigma'}
\lambda^{\asc(\sigma')}(\lambda-1)^{\des(\sigma')}X_{\sigma'(1)}\cdot \ldots\cdot X_{\sigma'(n-1)}\,\mathrm d\lambda
\plabel{eq:wrom}
\end{multline}
where $X_{k+\frac12}$ is an abbreviation for $X_kX_{k+1}-X_{k+1}X_k=[X_k,X_{k+1}]$.
But then, the RHS of \eqref{eq:wrom} is just $\mu_{n-1}(X_1,\ldots,X_{k-1},[X_{k },X_{k+1}],X_{k+2},\dots,X_n)$.

(The previous argument foreshadows the resolvent method;
a purely combinatorial argument is furnished by Solomon \cite{S}.)
\end{proof}
\end{lemma}
The previous lemma can be interpreted as the Magnus commutators realizing a canonical projection.
We can explain this as follows:

For $A_1,\ldots,A_k$, let the symmetrically distributed product these $n$ many terms be
\begin{equation}
A_1\cdot_\Sigma\ldots \cdot_\Sigma A_k=\frac1{k!}\sum_{\sigma\in\Sigma_k}A_{\sigma(1)}\cdot\ldots\cdot A_{\sigma(k)}.
\plabel{eq:Asymm}
\end{equation}
Let us now consider any homogeneous non-commutative polynomial of $P$ order of order $n$ ($n\geq1$).
Then we can bring $P$ into form
\begin{equation}
P=H_n+\ldots+H_1
\plabel{eq:Psymm}
\end{equation}
such that $H_i$ is a linear combination of symmetrically distributed products of $i$ many commutator monomials of the variables of $P$.
This will be termed as a `symmetrically distributed decomposition'.

A possible way to achieve this is as follows:
$P$ is a linear combination of product of $n$ many variables.
Then interchanging \textit{neighbouring} terms in the products, we can symmetrize $P$ at the sufferance of creating some commutators
\[\ldots_1XY\ldots_2\rightsquigarrow\ldots_1YX\ldots_2+\ldots_1[X,Y]\ldots_2.\]
Thus we have $P=H_n+P'$, where is $P'$ is a linear combination of products of $n-1$ many terms,
some are variables, some are commutators of variables.
Doing this inductively, we arrive to \eqref{eq:Psymm}.
We call this as a `(non-deterministic) elementary symmetrization process'.

\begin{lemma}[Solomon \cite{S} (1968), special case] \plabel{lem:magsol}
Assume that $P$ is a homogeneous non-commutative polynomial of order $n$;
\[P=\sum_{\lambda:\{1,\ldots,n\}\rightarrow \{1,\ldots,k\}}c_\lambda X_{\lambda_1}\cdot\ldots\cdot X_{\lambda_n}.\]
Assume that \eqref{eq:Psymm} is a result of a (non-deterministic) elementary symmetrizing decomposition.
Then,
\[ \sum_{\lambda:\{1,\ldots,n\}\rightarrow \{1,\ldots,k\}}c_\lambda \mu_n(X_{\lambda_1},\ldots, X_{\lambda_n})=H_1.\]
\begin{proof}
Let us consider the elementary symmetrization process but the products of the monomials
is replaced by Magnus commutators.
Then by the previous lemma: As the Magnus commutators are linear in the variables  (P1), and the creation  of commutators
is compatible to  taking the Magnus commutators (P2), we obtain
\[\mu_n(P)=\mu_n(H_n)+\ldots +\mu_1(H_1)\]
(with some abuse of notation).
However, by (P1), the higher Magnus commutators vanish, yielding $\mu_n(P)= \mu_1(H_1)\equiv H_1$.
\end{proof}
\end{lemma}
Despite its simplicity, this statement, and the method of its proof have already non-trivial consequences:
\begin{theorem}
[Dynkin \cite{Dyy}, Magnus \cite{M}, Solomon \cite{S}]
\plabel{th:magcomm}
$\mu_n(X_1,\ldots,X_n)$ is a commutator polynomial.
\begin{proof}
We can apply Lemma \ref{lem:magsol} for $P=X_1\cdot \ldots\cdot X_n$.
As $H_1$ is a commutator polynomial, so is $\mu_n(X_1,\ldots,X_n)$.
\end{proof}
\end{theorem}
\renewcommand{\qedsymbol}{$\triangle$}
\begin{proof} [Note]
This statement occurs in all of the indicated sources, but one should recognize that it applies to different constructs.
(See some comments later.)
By the time we arrive here, it is already clear  those constructs are the same.
\end{proof}
\renewcommand{\qedsymbol}{$\Box$}

\begin{remark}\plabel{rem:Eichler}
In the argument above the symmetrization process was used, in which the multivariable
nature of the Magnus commutators played a fundamental role.
Therefore, such a symmetrization argument cannot be applied as such in the specialized setup of  \eqref{eq:ABCH}.
This is, however, not completely true, see Eichler \cite{E}.
\qedremark
\end{remark}

\begin{theorem}[Solomon \cite{S}]\plabel{th:magchar}
The properties (P0)--(P2) characterize the Magnus commutators algebraically.
\begin{proof}
We can apply Lemma \ref{lem:magsol} for $P=X_1\cdot \ldots\cdot X_n$.
As its proof used only properties (P0)--(P2), the result is determined by them.
\end{proof}
\end{theorem}
\renewcommand{\qedsymbol}{$\triangle$}
\begin{proof} [Note]
We have arrived to Lemma \ref{lem:magchar} in two different ways.
In light of Theorem \ref{th:magchar}, that yields a universal Lie algebraic proof of Theorem \ref{th:mp}.
\end{proof}
\renewcommand{\qedsymbol}{$\Box$}
Another consequence which we could mention is that Lemma \ref{lem:magsol} shows that the part $H_1$ in \eqref{eq:Psymm}
does not depend on the non-deterministic elementary symmetrization process.
It is, however, very easy to see by other methods that there is nothing non-deterministic about the decomposition \eqref{eq:Psymm},
elementary symmetrization process or not:

Assume that $P$ is a non-commutative polynomial of some variables $X_i$.
We define the co-shuffle $F(P)$ of $P$ as follows: we get
$X_i$ replaced by $X_i\otimes 1+1\otimes X_i$
($1$ should be replaced by $1\otimes1$), and then, the resulted `coproduct' $\Delta(P)=\sum_I P_1^I\otimes P_2^I$
is replaced by $F(P)=\sum_I P_1^I  P_2^I$.

Let us recall the easy part of Friedrichs' criterion:
If $P$ is a commutator polynomial, then $\Delta(P)=P\otimes1+1\otimes P$.
Indeed, this is sufficient to prove for commutator monomials, which is very easy by induction.
Then we can prove
\begin{theorem}\plabel{th:Frie}
Let us  consider any homogeneous non-commutative polynomial of $P$ order of order $n$ ($n\geq1$).
Then for any symmetrically distributed decomposition \eqref{eq:Psymm},
\[F(H_i)=2^iH_i.\]
I.~e.~the symmetrically distributed decomposition is just the natural eigenspace decomposition of $F$
(with eigenvalues $2^i$), completely determined by linear algebraic reasons.

(If $n=0$, then $F$ is the identity operation with eigenvalue $1$.)
\begin{proof}
If we take a symmetrically distributed product \eqref{eq:Asymm} of commutator monomials,
the $\Delta$ may distribute  it in $2^k$ many ways.
As order of the components was random  anyway, the multiplicative reassembling just introduces a multiplier   $2^k$.
\end{proof}
\end{theorem}
\snewpage
Immediate consequence are:
\begin{cor}
\plabel{cor:Frie}
(Friedrichs' criterion) $P$ is a commutator polynomial if and only if $\Delta(P)=P\otimes1+1\otimes P$.
\begin{proof} We have already discussed the easy direction $(\Rightarrow)$.
Conversely, if $\Delta(P)=P\otimes1+1\otimes P$, then $F(P)=2^1P$.
Thus, in each homogeneity degree, $P$ turns out to be a commutator polynomial ($H_1$ component).
\end{proof}
\end{cor}
\begin{cor}
\plabel{cor:magshuff}
Let us  consider any homogeneous non-commutative polynomial of $P$ order of order $n$ ($n\geq2$).
Then (with some abuse of notation),
\[\mu_n(P)= \frac{( F-2^n\Id)\ldots( F-2^2\Id)}{(2-2^n)\ldots(2-2^2)}(P).\]
\begin{proof} This is immediate from the eigenvalue decomposition.
\end{proof}
\end{cor}
\begin{remark}\plabel{rem:magshuff}
At first sight there is little use for computing Magnus commutators via co-shuffles.
Yet, one can imagine situations where it can be useful:
If there are only two variables $X,Y$ and one wants to compute
Magnus commutators of shape $\mu_n(X\vee Y,\ldots,X\vee Y)$, then making a co-shuffle is not terribly expensive.
By successive substitutions, it can be realized in space $O(n2^n)$ with relatively polynomial time $O(p(n)2^n)$.
This is certainly better than summing up $n!$ many permutations.
Also, we can apply co-shuffles when we take the canonical projection of an already large non-commutative polynomial.
\qedremark
\end{remark}
\begin{theorem}[Solomon \cite{S} (1968), special case] \plabel{th:magsol}
Assume that $P$ is a homogeneous non-commutative polynomial of order $n$;
\[P=\sum_{\lambda:\{1,\ldots,n\}\rightarrow \{1,\ldots,k\}}c_\lambda X_{\lambda_1}\cdot\ldots\cdot X_{\lambda_n}.\]
Assume that \eqref{eq:Psymm} is a / the symmetrically distributed decomposition.
Then,
\[ \sum_{\lambda:\{1,\ldots,n\}\rightarrow \{1,\ldots,k\}}c_\lambda \mu_n(X_{\lambda_1},\ldots, X_{\lambda_n})=H_1.\]
\begin{proof}
By its unicity, it does not matter how we have obtained $H_1$.
\end{proof}
\end{theorem}

\begin{remark}\plabel{rem:algrem}
(i) As the symmetrically distributed decompositions turn out to be unique, one can see that they can be obtained from
\eqref{eq:ppod} by linear extension, with number $s$ being contributing to $H_s$.

(ii) Theorem \ref{th:Frie} and its consequences extend to universal enveloping algebras
over fields of characteristic $0$ without trouble.
\qedremark
\end{remark}
\snewpage
Let us continue with some remarks on Theorem \ref{th:magcomm}.
The simplest way to prove Theorem \ref{th:magcomm} is, of course,
\begin{point}\plabel{po:poalt1}
\begin{proof}[Alternative proof for Theorem \ref{th:magcomm}]
We can write  formally
\[\log(\exp(t_1X_1)\cdot\ldots\cdot\exp(t_kX_k))=\BCH(\ldots\BCH(\BCH(t_1X_1,t_2X_2),t_3X_3)\ldots,t_nX_n).\]
$\BCH$ is well-known to be a commutator series (cf. Discussion \ref{disc:tradBCH}),
thus the coefficient of $t_1\cdot\ldots\cdot t_n$ will be a commutator polynomial.
\end{proof}
\end{point}
The coproduct for $\mu_n(X_1,\ldots,X_n)$ can also be obtained easily:
As the exponential of a sum of commuting elements is the product of exponentials;
and, conversely, the logarithm of a product of commuting perturbations of the unity element
is the sum of the logarithms; one quickly finds
\begin{multline*}
\log(\exp(t_1(X_1\otimes1+1\otimes X_1))\cdot\ldots\cdot\exp(t_k(X_k\otimes1+1\otimes X_k) ))=
\\
=\log(\exp(t_1X_1)\cdot\ldots\cdot\exp(t_kX_k))\otimes 1
+1\otimes\log(\exp(t_1   X_1)\cdot\ldots\cdot\exp(t_kX_k)).
\end{multline*}
Detecting the coefficients of $t_1\cdot\ldots\cdot t_k$ immediately yields
 \begin{multline}
\mu_n(X_1\otimes1+1\otimes X_1,\ldots, X_n\otimes1+1\otimes X_n)=
\\=
 \mu_n(X_1 ,\ldots, X_n )\otimes1+1\otimes\mu_n(X_1 ,\ldots, X_n ).
\plabel{eq:Magco}
\end{multline}

\begin{point}
\plabel{po:poalt2}
\begin{proof}[Another proof for Theorem \ref{th:magcomm}]
Identity \eqref{eq:Magco}, by  Friedrichs' criterion, implies that $\mu_n(X_1 ,\ldots, X_n )$ is a commutator polynomial.
\end{proof}
\end{point}
\renewcommand{\qedsymbol}{$\triangle$}
\begin{proof}[Note]
\plabel{rem:alg}
Dynkin \cite{Dyy} is certainly aware of the argument of Proof \ref{po:poalt1},
 but he wants to use prove the BCH result and not use it;
 thus, his proof of Theorem \ref{th:magcomm} is entirely combinatorial, cf.~Achilles, Bonfiglioli \cite{AB}.
Proof \ref{po:poalt1} is also implicit at Magnus \cite{M} but, due to the analytical circumstances,
 he prefers to obtain that $\mu_k$ is a commutator polynomial from his recursion formula, which is based on Schur's formula.
This recursion formula of Magnus, and anything more on  explicit commutator expressions for the $\mu_k$
 will, however, be discussed only later, in Part III.
Chen \cite{Ch} is similar to Magnus \cite{M}.
Interestingly, Magnus \cite{M} contains a proof of Friedrichs' criterion, but
 he uses only to obtain a quick proof for the Baker--Campbell--Hausdorff result.
Solomon \cite{S}, starting from the direction of the standard universal Lie algebraic tools, is in pursuit
 of some naturally defined Lie elements in the first place; only an explicit form is needed, which he finds.
Strichartz \cite{St} works out the Magnus expansion with explicit formula in two different ways,
 meanwhile both Proof \ref{po:poalt1} and Proof \ref{po:poalt2} appear.
Helmstetter \cite{H} extends Solomon \cite{S}, uses the coalgebra structure with respect
 to $\mu_k$ as defined by Dynkin (cf.~Proof \ref{po:poalt2}), rediscovers Goldberg's theorem; he obtains several formulas.
As we see, several of these authors arrive  to $\mu_k$ from different directions and contexts
(and typically being only partially aware of the related works).
The related algebra (up to ordinary Lie algebras)  is discussed the most comprehensively in Reutenauer \cite{R}.
In terms of the associated combinatorics, see Mielnik, Pleba\'nski \cite{MP}, Garsia \cite{Gar}, Reutenauer \cite{R}.
\end{proof}
\renewcommand{\qedsymbol}{$\Box$}

\snewpage
\begin{commenty}
\subsection{On the power-series based holomorphic calculus}\plabel{ss:CombHol}
~\\

Here we consider what one can do without the resolvent method but using power series calculus.
Algebraic consequences are also examined.
\begin{theorem}[Cf.~Dynkin \cite{Dy} (BCH case);  Pechukas, Light \cite{PL},  Karasev, Mosolova \cite{KM}]\plabel{th:dynest}
Let $\phi$ be a continuous $\mathfrak A$-valued measure of finite variation on the interval $I$.
If
\begin{equation}
\int|\phi|<\log2,
\plabel{eq:sersub}
\end{equation}
or even just
\begin{equation}
\sum_{k=1}^\infty\left|\int_{t_1\leq\ldots\leq t_k\in I}\phi(t_1)\cdot\ldots\cdot\phi(t_k)\right|<1,
\plabel{eq:sersubb}
\end{equation}
holds, then
\begin{equation*}
\log^{(\mathrm{pow})}
\left(1+\sum_{k=1}^\infty\int_{t_1\leq\ldots\leq t_k\in I}\phi(t_1)\cdot\ldots\cdot\phi(t_k)\right)
=\sum_{k=1}^\infty\int_{t_1\leq\ldots\leq t_k\in I}\mu_k(\phi(t_1),\ldots,\phi(t_k)).
%\plabel{eq:magsub}
\end{equation*}
\begin{proof}
It is easy to see that \eqref{eq:sersub} implies \eqref{eq:sersubb}:
Indeed, in this case
\[\sum_{k=1}^\infty\left|\int_{t_1\leq\ldots\leq t_k\in I}\phi(t_1)\cdot\ldots\cdot\phi(t_k)\right|<
\sum_{k=1}^\infty\frac1{k!}(\log 2)^k=\exp(\log 2)-1=1.\]
In turn, \eqref{eq:sersubb} implies that we can substitute $T=1$ into \eqref{eq:formmag} understood with $\log^{(\mathrm{pow})}$.
\end{proof}
\end{theorem}

The following digression will also demonstrate that how changing the kind of logarithm influences matters.
One can define the inverse-symmetric power series logarithm as
\[\log^{\mathrm{ipow}} A=\log^{\mathrm{pow}}\left(\frac{A+1}2\right)- \log^{\mathrm{pow}}\left(\frac{A^{-1}+1}2\right).\]
(It is considered to be meaningful if both summands are.)
Considering its domain, $\log^{\mathrm{ipow}}$ is incomparable $\log^{(\mathrm{pow})}$
 but it is still a special case of our default spectral $\log$.
\snewpage
\begin{theorem}\plabel{th:log3}
Let $\phi$ be a continuous $\mathfrak A$-valued measure of finite variation on the interval $I$.
If
\begin{equation}
\int|\phi|<\log3,
\plabel{eq:sersub3}
\end{equation}
or even just
\begin{equation}
\sum_{k=1}^\infty\left|\int_{t_1\leq\ldots\leq t_k\in I}\phi(t_1)\cdot\ldots\cdot\phi(t_k)\right|<2
\,\,\text{and}\,\,
\sum_{k=1}^\infty\left|\int_{t_1\leq\ldots\leq t_k\in I}\phi(t_k)\cdot\ldots\cdot\phi(t_1)\right|<2,
\plabel{eq:sersubb3}
\end{equation}
holds, then
\begin{equation*}
\log^{(\mathrm{ipow})}
\left(1+\sum_{k=1}^\infty\int_{t_1\leq\ldots\leq t_k\in I}\phi(t_1)\cdot\ldots\cdot\phi(t_k)\right)
=\sum_{k=1}^\infty\int_{t_1\leq\ldots\leq t_k\in I}\mu_k(\phi(t_1),\ldots,\phi(t_k)).
\end{equation*}
\begin{proof}
It is easy to see that \eqref{eq:sersub3} implies \eqref{eq:sersubb3}:
Indeed, in this case
\[\sum_{k=1}^\infty\left|\int_{t_1\leq\ldots\leq t_k\in I}\phi(t_1)\cdot\ldots\cdot\phi(t_k)\right|<
\sum_{k=1}^\infty\frac1{k!}(\log 3)^k=\exp(\log3)-1=2.\]
In this case $\log^{(\mathrm{pow})}\frac{\Rexp(\phi)+1}2$ makes sense.
The argument is similar with respect to $\log^{(\mathrm{pow})}\frac{\Rexp(\phi)^{-1}+1}2\equiv
\log^{(\mathrm{pow})}\frac{\Rexp(-\phi^\dag) +1}2$, etc.
\end{proof}
\end{theorem}
\begin{commentx}
\begin{proof}[Note]
 $\log^{(\mathrm{pow})}$ ceases to work already for
$\phi=(\log 2)X\mathbf 1_{[0,1)}$ in the setting formal power series with $\ell^1$ norm, where $\int|\phi|=\log 2$.
 $\log^{(\mathrm{ipow})}$ ceases to work already for
$\phi=(\log 3)X\mathbf 1_{[0,1)}$ in the setting formal power series with $\ell^1$ norm, where $\int|\phi|=\log 3$.
(Meanwhile, the extension with the ordinary $\log$ will work for $\int|\phi|< 2$.)
\renewcommand{\qedsymbol}{$\triangle$}
\end{proof}
\renewcommand{\qedsymbol}{$\Box$}
\end{commentx}
\begin{remark}
\plabel{rem:logimp}
The argument of Theorem \ref{th:dynest} yields
\begin{equation}
\sum_{n=1}^\infty T^n\cdot|\mu_{\mathrm R,n}(\phi)| \leq
 \log\frac1{ 2-  \exp\left(T\cdot\textstyle{\int|\phi| }\right) } ;\plabel{eq:Tlog2}
\end{equation}
and the  argument of Theorem \ref{th:log3} yields
\begin{equation}
\sum_{n=1}^\infty T^n\cdot|\mu_{\mathrm R,n}(\phi)| \leq
 2\log\frac2{ 3-  \exp\left(T\cdot\textstyle{\int|\phi| }\right) } .\plabel{eq:Tlog3}
\end{equation}
As one expects, \eqref{eq:Tlog2} is weaker than \eqref{eq:Tlog3}, which is weaker than
\eqref{eq:arc1}, which is weaker than
\eqref{eq:Tmeasum}/\eqref{eq:mpest}/\eqref{eq:ThetaMest}.
\qedremark
\end{remark}

Assume that $F(w)$ is a holomorphic function near $w\sim1$.
Let $R_F$ be the convergence radius of $F$ around $w=1$
(but there is no danger in choosing a smaller positive number); such that $F(w)=\sum_{n=0}^\infty F_i\cdot(w-1)^n$.
If the spectrum of $A$ is contained in $\intD(1,R_F)$, then we set
\[F^{(\mathrm{pow})}(A)=\sum_{n=0}^\infty F_i\cdot(A-1)^n.\]
Let $\mu_k^{\langle F\rangle}(X_1,\ldots,X_k)$ be the
universally 1-homogeneous part of $F( \exp(X_1)\cdot\ldots\cdot\exp(X_k))$ (formally).
Then
\begin{theorem}\plabel{th:holdes}
Let $\phi$ be a continuous $\mathfrak A$-valued measure of finite variation on the interval $I$.
Then, for
\[\int|\phi|<\log(1+R_F),\]
or even just for
\[
\sum_{k=1}^\infty\left|\int_{t_1\leq\ldots\leq t_k\in I}\phi(t_1)\cdot\ldots\cdot\phi(t_k)\right|< R_F,
\]
it holds that
\[
F^{(\mathrm{pow})}
\left(1+\sum_{k=1}^\infty\int_{t_1\leq\ldots\leq t_k\in I}\phi(t_1)\cdot\ldots\cdot\phi(t_k)\right)
=\sum_{k=1}^\infty\int_{t_1\leq\ldots\leq t_k\in I}\mu_k^{\langle F\rangle}(\phi(t_1),\ldots,\phi(t_k)),
\]
where the sum on the RHS is absolutely convergent.
\begin{proof}
Let us take the formal variable $T$. Take
\[F^{(\mathrm{pow})}\left(1+\sum_{k=1}^\infty T^k\int_{t_1\leq\ldots\leq t_k\in I}\phi(t_1)\cdot\ldots\cdot\phi(t_k)\right);\]
expand it and collect it in powers of $T$.
Due to the combinatorial nature of the process
(again, we can imagine that $X_j$ is a placeholder for the $j$th chronological contribution), what we obtain must be
\[\sum_{k=0}^\infty T^k\int_{t_1\leq\ldots\leq t_k\in I}\mu_k^{\langle F\rangle}(\phi(t_1),\ldots,\phi(t_k)).\]
Then, due to the convergence relations, $T=1$ can be substituted.
\end{proof}
\end{theorem}
This argument can also be made on the formal level (practically, with the insertion of a `$T$'),
starting from a formal power series `$F(1+x)$'.
The combinatorics of this situation was studied by Mielnik, Pleba\'nski \cite{MP}.

In what follows, we may use the notation
$\mu_{k,\mathrm R}^{\langle F\rangle}(\phi)\equiv\int_{t_1\leq\ldots\leq t_k\in I}\mu_k^{\langle F\rangle}(\phi(t_1),\ldots,\phi(t_k))$,
etc.
Let us apply the theorem above in some special situations.

\snewpage
For a Banach algebra element $A$ whose spectrum lies in $\mathbb C\setminus(-\infty,0]$,
 and for $\alpha\in\mathbb C\setminus\mathbb Z$, we define the `fractional' power by
\[A^\alpha=\exp(\alpha\log A).\]
With the given spectral condition this is also valid $\alpha\in\mathbb Z$, but then $A^\alpha$ is defined more generally.
(The fractional powers can, in general, also be defined a bit more generally, but we apply them only for
invertible operators, thus the definition above is sufficient.)
If the spectrum of $A$ is contained in $\intD(1,1)$, then set
\[A^{\alpha\,{(\mathrm{pow})}}=\sum_{n=0}^\infty\binom{\alpha}{n}(A-1)^n.\]
The domain $A^{\alpha\,{(\mathrm{pow})}}$ is the same as the domain of $\log^{(\mathrm{pow})} A$;
$A^{\alpha\,{(\mathrm{pow})}}=\exp(\alpha\log^{(\mathrm{pow})} A)$ holds in every major sense.
For $k\in\mathbb N$, let $\mu^{\langle\langle\alpha\rangle\rangle}(X_1,\ldots,X_k)$ be the universally $1$-homogeneous part of
\[\left((\exp X_1)\cdot\ldots\cdot(\exp X_k)\right)^\alpha \equiv
\exp\left(\alpha\log\left((\exp X_1)\cdot\ldots\cdot(\exp X_k)\right)\right)\]
 in the non-commutative variables $X_1,\ldots,X_k$.
Then we obtain
\begin{theorem}\plabel{th:mppow}
Let $\phi$ be a continuous $\mathfrak A$-valued measure of finite variation on the interval $I$.
If \eqref{eq:sersub}, or even just \eqref{eq:sersubb} holds, then
\begin{equation*}
\left(1+\sum_{k=1}^\infty\int_{t_1\leq\ldots\leq t_k\in I}\phi(t_1)\cdot\ldots\cdot\phi(t_k)\right)^{\alpha\,{(\mathrm{pow})}}
=\sum_{k=0}^\infty\int_{t_1\leq\ldots\leq t_k\in I}\mu_k^{\langle\langle\alpha\rangle\rangle}(\phi(t_1),\ldots,\phi(t_k)),
%\plabel{eq:magsubpow}
\end{equation*}
where the sum on the RHS is absolutely convergent.
\qed
\end{theorem}

\begin{theorem}[Mielnik, Pleba\'nski \cite{MP}]\plabel{th:powex}
\begin{equation}
\mu_k^{\langle\langle\alpha\rangle\rangle}(X_1,\ldots,X_k)
=\sum_{\sigma\in\Sigma_k}\binom{\asc(\sigma)+\alpha}{k} X_{\sigma(1)}\ldots X_{\sigma(k)}.
\plabel{eq:powex}
\end{equation}
\begin{proof} We model the proof after the proof of Theorem \ref{th:mp}.
Considering  the formal power series of $(1+x)^\alpha$, one can see that
$\mu_k^{\langle\langle\alpha\rangle\rangle}(X_1,\ldots,X_k)$ is a sum of terms
\[\binom{\alpha}{k}
X_{\sigma(1)}\cdot\ldots\cdot X_{\sigma(l_1)}\mid  X_{\sigma(l_1+1)}\cdot\ldots\ldots
\cdot X_{\sigma(l_{j-1})}\mid X_{\sigma(l_{j-1}+1)}\cdot\ldots\cdot X_{\sigma(k)},\]
where separators show the ascendingly indexed components which enter into power series of $(1+x)^\alpha$.
Considering a given permutation $\sigma$, the placement of the separators is either necessary (in case $\sigma(i)>\sigma(i+1)$),
or optional (in case  $\sigma(i)<\sigma(i+1)$).
Summing over the $2^{\asc(\sigma)}$ many optional possibilities, the coefficient of $X_{\sigma(1)}\cdot\ldots\cdot X_{\sigma(k)}$
is
\[\sum_{p=0}^{\asc(\sigma)}\begin{pmatrix}\asc(\sigma)\\p\end{pmatrix}\binom{\alpha}{p+\des\sigma+1}.\]
This simplifies  according to the combinatorial identity
$\sum_{p=0}^{r}\binom rp\binom\alpha{p+d+1}= \binom {r+\alpha}{r+d+1}$.
\end{proof}
\end{theorem}
\snewpage
Assume that $n\geq1$, and let
\[X_1\cdot\ldots\cdot X_n=\mu^{[n]}_n(X_1,\ldots, X_n)+\ldots+\mu^{[1]}_n(X_1,\ldots, X_n)\]
be a symmetrically distributed decomposition as in \eqref{eq:Psymm}.
Then $\mu^{[1]}_n(X_1,\ldots, X_n)\equiv\mu_n(X_1,\ldots, X_n)$ is the first canonical projection.
Taking
\[\exp(X_1)\cdot\ldots\cdot \exp(X_n)=\sum_{s=0}^\infty \frac1{s!}\left(\log\left( \exp(X_1)\cdot\ldots\cdot \exp(X_n)\right)\right)^s\]
and restricting to the  universally $1$-homogeneous part;
from the uniqueness of the symmetrically distributed decompositions,
 the term $\mu^{[s]}_n(X_1,\ldots, X_n)$ must be the universally $1$-homogeneous part of
 the ``$s$-symmetric term''
\[\frac1{s!}\left(\log\left( \exp(X_1)\cdot\ldots\cdot \exp(X_n)\right)\right)^s.\]

\begin{theorem}\plabel{th:logepo}
Let $\phi$ be a continuous $\mathfrak A$-valued measure of finite variation on the interval $I$.
If \eqref{eq:sersub}, or even just \eqref{eq:sersubb} holds, then
\begin{equation}
\frac1{s!}\log^{(\mathrm{pow})}
\left(1+\sum_{k=1}^\infty\int_{t_1\leq\ldots\leq t_k\in I}\phi(t_1)\cdot\ldots\cdot\phi(t_k)\right)^s
=\sum_{k=s}^\infty\int_{t_1\leq\ldots\leq t_k\in I}\mu_k^{[s]}(\phi(t_1),\ldots,\phi(t_k)),
\plabel{eq:magsuber}
\end{equation}
where the sum on the RHS is absolutely convergent.
\qed
\end{theorem}
Then, Theorem \ref{th:holdes} can be augmented by
\begin{theorem}[Mielnik, Pleba\'nski \cite{MP}]\plabel{th:logdes}
Assume that $f(z)=\sum_{i=0}^\infty f_iz^i$ is a holomorphic function near $z\sim0$ such that $F(w)=f(\log w)$
for $w\sim 1$.
Then
\[\mu_k^{\langle F\rangle}(X_1,\ldots,X_k)=\sum_{s=0}^k f_s \cdot s!\cdot \mu_k^{[s]}(X_1,\ldots,X_k).\]
\begin{proof}
This follows from considering the expansion in the formal setting,
and comparing Theorem \ref{th:holdes} and  Theorem \ref{th:logepo}.
\end{proof}
\end{theorem}
This, in particular, yields
\begin{theorem}[Mielnik, Pleba\'nski \cite{MP}]\plabel{th:mppowgen}
\begin{equation}
\mu_k^{\langle\langle\alpha\rangle\rangle}(X_1,\ldots,X_k)=\sum_{s=0}^k\alpha^s\mu_k^{[s]}(X_1,\ldots,X_k).
\plabel{eq:mppowgen}
\end{equation}
\begin{proof}
Apply Theorem \ref{th:logdes} with respect to $F(w)=w^\alpha$ and $f(z)=\exp(\alpha z)$.
\end{proof}
\end{theorem}

As a corollary, we obtain
\snewpage
\begin{theorem}[Mielnik, Pleba\'nski \cite{MP}, Garsia \cite{Gar}]\plabel{th:mpgarsia}
The ``generating function of the canonical projections'' is given as RHS\eqref{eq:mppowgen}$=$RHS\eqref{eq:powex}.
\begin{commentx}
\[\sum_{s=0}^kx^s\mu_k^{[s]}(X_1,\ldots,X_k)
=\sum_{\sigma\in\Sigma_k}\binom{x+\asc(\sigma)}{k} X_{\sigma(1)}\ldots X_{\sigma(k)}.\]
\end{commentx}
In particular
\[\mu_k^{[s]}(X_1,\ldots,X_k)=
\sum_{\sigma\in\Sigma_k}
\left.\frac1{s!}\frac{\mathrm d^s}{\mathrm dx^s}\binom{x+\asc(\sigma)}{k}\right|_{x=0}
X_{\sigma(1)}\ldots X_{\sigma(k)}.\]
\begin{proof}
Compare Theorem \ref{th:powex} and \ref{th:mppowgen}; and replace $\alpha$ by $x$.
\end{proof}
\end{theorem}

See Reutenauer \cite{R} for more about canonical projections (from other viewpoints).

Let us consider another example.
If the spectrum of $A$ is does not contain $\frac{\lambda}{\lambda-1}$, then we can take
\[\mathcal R^{(\lambda)}(A)=\frac{A-1}{\lambda+(1-\lambda)A}.\]
If the spectrum of $A$ is contained in $\intD\left(1,\frac1{|1-\lambda|}\right)$, then we set
\[\mathcal R^{(\lambda)\,(\mathrm{pow})}(A)=\sum_{n=1}^\infty (\lambda-1)^{n-1}A^n.\]
Let $\mu_k^{(\lambda)}(X_1,\ldots,X_k)$ be the universally $1$-homogeneous part of
\[\mathcal R^{(\lambda)}\left( (\exp X_1)\cdot\ldots\cdot(\exp X_k)\right) \]
 in the non-commutative variables $X_1,\ldots,X_k$ (in formal sense).
Then we obtain
\begin{theorem}\plabel{th:mppowr}
Let $\phi$ be a continuous $\mathfrak A$-valued measure of finite variation on the interval $I$.
If
\begin{equation}
\int|\phi|<\log\left( 1+\frac1{|1-\lambda|}\right),
\plabel{eq:sersubr}
\end{equation}
or even just
\begin{equation}
\sum_{k=1}^\infty\left|\int_{t_1\leq\ldots\leq t_k\in I}\phi(t_1)\cdot\ldots\cdot\phi(t_k)\right|< \frac1{|1-\lambda|},
\plabel{eq:sersubbr}
\end{equation}
holds, then
\begin{equation*}
\mathcal R^{(\lambda)\,(\mathrm{pow})}\left(1+\sum_{k=1}^\infty\int_{t_1\leq\ldots\leq t_k\in I}\phi(t_1)\cdot\ldots\cdot\phi(t_k)\right)
=\sum_{k=0}^\infty\int_{t_1\leq\ldots\leq t_k\in I}\mu_k^{(\lambda)}(\phi(t_1),\ldots,\phi(t_k)),
%\plabel{eq:magsubpow}
\end{equation*}
where the sum on the RHS is absolutely convergent.
\qed
\end{theorem}

See Theorem \ref{th:powexr} for the concrete form of $\mu_k^{(\lambda)}(X_1,\ldots,X_k)$.

\begin{commentx}
\begin{remark}\plabel{rem:ipow}
One can define $A^{\alpha \,(\mathrm{ipow}})$ and
$\mathcal R^{(\lambda)\,(\mathrm{ipow})}$ for $\lambda\in[0,1]$ under the same
spectral conditions as $\log^{\mathrm{ipow}}$.
Under the condition \eqref{eq:sersubb3}, the statements of
 Theorem \ref{th:mppow} and Theorem \ref{th:mppowr} for $\lambda\in[0,1]$  hold.
This is, however, not quite the right direction, as the resolvent should be treated the most fundamental quantity.
\qedremark
\end{remark}
\end{commentx}

\begin{commentx}
\begin{remark}\plabel{rem:ripow}
If $F(w)$ is holomorphic for $w\sim1$, then $F^{(\mathrm{pow})}(A)$ was obtained as follows:
Take $F^{(\mathrm{pow})}_T(A)=F(1+T\cdot(A-1))$, where $T$ is a formal infinitesimal variable,
expand it as a series of $T$, and take $T=1$ in the already expanded series.

For us, taking  $F^{(\mathrm{ipow})}(A)$ is an ad hoc procedure.
nevertheless, a more systematic version, but is quite acceptable in several situations is as follows:
Take $F^{(\mathrm{ipow})}_T(A)=F\left(\frac{1+T\cdot A}{1+T\cdot A^{-1}}\right)$, where $T$ is a formal infinitesimal variable,
expand it as a series of $T$, and take $T=1$ in the already expanded series.
\qedremark
\end{remark}
\end{commentx}

\begin{theorem}[Mielnik, Pleba\'nski \cite{MP}]\plabel{th:gemcan}
\[\mu_k^{(\lambda)}(X_1,\ldots,X_k)=\sum_{s=0}^k G_s(\lambda,\lambda-1) \cdot s!\cdot \mu_k^{[s]}(X_1,\ldots,X_k).\eqed\]
\begin{proof}
Apply Theorem \ref{th:logdes}
with respect to $F(w)=\frac{w-1}{\lambda+(1-\lambda) w}$ and $f(z)= \frac{\exp(z)-1}{\lambda+(1-\lambda) \exp (z)}$.
By Theorem \ref{th:Euler}, $\frac{(\exp x)-1}{\lambda+(1-\lambda)(\exp x)}=G(\lambda,\lambda-1;x)
\equiv\sum_{k=1}^\infty G_k(\lambda,\lambda-1)x^k $.
\end{proof}
\end{theorem}

If $p(\lambda)$ is a polynomial of $\lambda$, then we may define the operation $W_{/\lambda}$ by
\begin{equation}
W_{/\lambda}\,p(\lambda)=-\left(\frac1\lambda\int_{t=0}^\lambda p(t)\,\mathrm dt +\frac1{\lambda-1}\int_{\lambda}^1p(t)\mathrm dt\right).
\plabel{eq:Wdef}
\end{equation}
While this may look complicated at first sight, it can be algebraized very simply:
For $n\in\mathbb N$,
\begin{equation}
W_{/\lambda}\,\left(\lambda^n\right)=\frac1{n+1}\left(1+\ldots+\lambda^{n-1}\right).
\plabel{eq:Walt}
\end{equation}

\begin{lemma}\plabel{lem:cogenres} For $k\geq 1$,
\[W_{/\lambda}\,G_k(\lambda,\lambda-1)=\frac1kG_{k-1}(\lambda,\lambda-1),\]
understood such that $G_{0}(\lambda,\lambda-1)=0$.
\begin{proof}
Let $W_{_{/\lambda},\mathrm{pol}}$ be the linear operation which acts on the $G_k(\lambda,\lambda-1)$'s as expected.
By Lemma \ref{lem:goldgen}, $W_{_{/\lambda},\mathrm{pol}}$ can be realized as follows:
Correct $p(\lambda)$ to $p(\lambda)-\int_{t=0}^1p(t)\, dt$, so that its integral on $[0,1]$ should vanish.
Take a primitive function which vanishes at $\lambda=0$, and henceforth at $\lambda=1$.
Divide by $\lambda(\lambda-1)$.
The algebraic formulation of this process yields $W_{/\lambda}$ as in \eqref{eq:Wdef}.
(A combinatorial proof based on \eqref{eq:Walt} can also be given.)
\end{proof}
\end{lemma}
This gives another generating method for the higher canonical projections:
\begin{theorem}\plabel{th:cogenres}
 For $s\geq1$,
\[\mu_k^{[s]}(X_1,\ldots,X_k)=\int_{\lambda=0}^1 (W_{/\lambda})^{s-1}\mu_k^{(\lambda)}(X_1,\ldots,X_k)\,\mathrm d\lambda.\]
\begin{proof}
By Lemma \ref{lem:cogenres}, $(W_{/\lambda})^{s-1}\mu_k^{(\lambda)}(X_1,\ldots,X_k)=G_1(\lambda,\lambda-1)\mu_k^{[s]}(X_1,\ldots,X_k)+$
 a linear combination of $G_r(\lambda,\lambda-1)$'s with $r>1$.
Then, by Lemma \ref{lem:goldgen}, integration gives the desired result.
\end{proof}
\end{theorem}

\snewpage

In terms of the canonical logarithm, powers and integral powers of $\log$ do not depart much from the Magnus expansion:
\begin{theorem}\plabel{th:undep}
Let $\phi$ be a continuous $\mathfrak A$-valued measure of finite variation on the interval $I$.
Assume that
\begin{equation}
\log\Rexp(\phi)=\sum_{k=1}^\infty\mu_{k,\mathrm R}(\phi)
\plabel{eq:logconc}
\end{equation}
(so that the LHS exists and the RHS is absolutely convergent).
Then, for $s\in\mathbb N$,
\begin{equation}
\frac1{s!}\left(\log\Rexp(\phi)\right)^s=\sum_{k=s}^\infty\mu_{k,\mathrm R}^{[s]}(\phi);
\plabel{eq:logpowconc}
\end{equation}
and, for $\alpha\in\mathbb C$,
\begin{equation}
\left(\Rexp(\phi)\right)^\alpha=\sum_{k=0}^\infty\mu_{k,\mathrm R}^{\langle\langle\alpha\rangle\rangle}(\phi)
\plabel{eq:powconc}
\end{equation}
hold (so that the LHS's exist and the RHS's  are absolutely convergent).
\begin{proof}
Simply take
$\frac1{s!}\left(\sum_{k=1}^\infty\mu_{k,\mathrm R}(\phi)\right)^s$
and
$\sum_{s=0}^\infty\alpha^s\left(\sum_{k=1}^\infty\mu_{k,\mathrm R}(\phi)\right)^s;$
expand them to linear combinations of monomials of $\mu_{k_1,\mathrm R}(\phi)\ldots \mu_{k_r,\mathrm R}(\phi)$, and contract them
according to homogeneity (i.~e.~in $k_1+\ldots+k_r$).
Absolute convergence makes this process OK.
\end{proof}
\end{theorem}

However, $\mathcal R^{(\lambda)}(\Rexp(\phi))$ cannot not be treated as simple derivative of $\log(\Rexp(\phi))$; but
the other way around.
In fact, $\mathcal R^{(\lambda)}$ is a game changer from analytic viewpoint,
 because it is a function which is directly connected to the spectral properties of its argument.
That prompts the direct study of  $\mu_k^{(\lambda)}(X_1,\ldots,X_k)$, which is the subject of the next section.
Nevertheless, one can already state
\begin{theorem}\plabel{eq:powecan}
Let $\phi$ be a continuous $\mathfrak A$-valued measure of finite variation on the interval $I$.
If
$\log(\Rexp(t\cdot\phi))$  exists for all $ t\in\Dbar(0,1)$,
then \eqref{eq:logconc}, and \eqref{eq:logpowconc}, and \eqref{eq:powconc} hold.
Moreover, for $\lambda\in[0,1]$,
\begin{equation}
\mathcal R^{(\lambda)}\left(\Rexp(\phi)\right)=\sum_{k=0}^\infty\mu_{k,\mathrm R}^{(\lambda)}(\phi)
.
\plabel{eq:resconc}
\end{equation}
\begin{proof}
As the spectrum is closed, $t\in\Dbar(0,1)$ can be replaced by $t\in\Dbar(0,1+\varepsilon)$.
Then $\log(\Rexp(t\cdot\phi))$, $\frac1{s!}\log(\Rexp(t\cdot\phi))^s$,  $(\Rexp(t\cdot\phi))^\alpha$,
 $\mathcal R^{(\lambda)}(\Rexp(t\cdot\phi))$
 exist and  holomorphic in $t\in\Dbar(0,1+\varepsilon)$; thus the corresponding power series representations exist.
By the previous theorems, they are given as $\sum_{k=0}^\infty t^k\mu_{k,\mathrm R}(\phi) $, \ldots,
 $\sum_{k=0}^\infty t^k\mu_{k,\mathrm R}^{(\lambda)}(\phi) $ for $t\sim0$.
However, by the power series representation, those are also valid for $t=1$.
\end{proof}
\end{theorem}
The condition `$\log(\Rexp(t\cdot\phi))$  exists for all $ t\in\Dbar(0,1)$'
 will be termed later as being `$M$-controlled'.
\end{commenty}

\snewpage
\section{The (infinitesimal) resolvent approach}\plabel{sec:MagnusResolvent}
\subsection{Introduction to the resolvent method}\plabel{ss:IntroResolvent}
~\\

As we have seen, in the Banach algebraic setting,  the Magnus expansion can be treated quite directly.
The resolvent approach of Mielnik, Pleba\'nski \cite{MP} was not needed in its full power
(although the Eulerian generating function was used).
However, their resolvent method comes handy when we inquire about finer analytical details.

Recall, in Banach algebras, we can define the logarithm of $A\in\mathfrak A$ by
\begin{equation}
\qquad
\log A=\int_{\lambda=0}^1 \frac{A-1}{\lambda +(1-\lambda)A }\,d\lambda
\qquad\left(=\int^{0}_{s=-\infty}\frac{A-1}{(1-s)(A-s)}\, \mathrm ds\right),
\plabel{eq:logdef}
\end{equation}
which we consider well-defined if and only if $\spec(A)$ is disjoint from the closed negative real axis.
We may call  $A$ $\log$-able, if it has this spectral property.

Notice the (modified) resolvent expression
\begin{equation}\mathcal R^{(\lambda)}(A)\equiv\frac{A-1}{\lambda+(1-\lambda)A}.\plabel{eq:resdef}\end{equation}
Here we have used $\lambda$ to parametrize the resolvent.
Technically, $\nu=2\lambda-1$ or $\xi=1-\lambda$ are equally good choices, or even better;
but the use of $\lambda$ fits to the earlier discussions.
(More specifically, Mielnik, Pleba\'nski \cite{MP} chooses our $\nu$ as the parameter, and the
their normalization of the resolvent expression is also slightly different.
Our conventions fit more to Goldberg's formalism.)
\begin{point}\plabel{rem:Res1}
Note: The resolvent can be defined as any element $\mathcal R^{(\lambda)}(A)$ such that
\begin{equation}
(\lambda+(1-\lambda)A)\mathcal R^{(\lambda)}(A) =A-1\quad\text{and}\quad\mathcal R^{(\lambda)}(A)(\lambda+(1-\lambda)A) =A-1
\plabel{eq:reseq}
\end{equation}
hold (``equational'' definition).
Indeed, then it is easy to show that
\[(\lambda+(1-\lambda)A)(1-(1-\lambda)\mathcal R^{(\lambda)}(A)) =1\quad\text{and}
\quad(1-(1-\lambda)\mathcal R^{(\lambda)}(A))(\lambda+(1-\lambda)A) =1\]
 hold.
Hence, $(\lambda+(1-\lambda)A)^{-1}$ exists.
Then \eqref{eq:reseq} implies that the quantity $\mathcal R^{(\lambda)}(A)$
is also given as the product of  $A-1$ and $(\lambda+(1-\lambda)A)^{-1}$ (in arbitrary order; in particular, the terms commute).
Thus, the ``equational'' definition \eqref{eq:reseq} is equivalent to
\eqref{eq:resdef}.
\qedremark
\end{point}
\begin{remark}\plabel{rem:Res2}
The identities
\begin{equation}
(1-(1-\lambda)\mathcal R^{(\lambda)}(A))=\frac1{\lambda+(1-\lambda)A}
\plabel{eq:rres1}
\end{equation}
and
%\[
\begin{equation}
(1-\lambda\mathcal R^{(\lambda)}(A))=\frac A{\lambda+(1-\lambda)A}
\plabel{eq:rres2}
\end{equation}
%\]
are often useful.
\eqref{eq:rres1}, in particular, shows (again) that the existence of $\mathcal R^{(\lambda)}(A)$
is equivalent to invertibility of $\lambda+(1-\lambda)A$.
\qedremark
\end{remark}

Now $\mu_k^{(\lambda)}(X_1,\ldots,X_k)$ is the universally $1$-homogeneous part of
\[\mathcal R^{(\lambda)}\left( (\exp X_1)\cdot\ldots\cdot(\exp X_k)\right) \]
 in the non-commutative variables $X_1,\ldots,X_k$.

 \begin{theorem}[Mielnik, Pleba\'nski \cite{MP}]\plabel{th:powexr}
\begin{equation}
\mu_k^{(\lambda)}(X_1,\ldots,X_k)
=\sum_{\sigma\in\Sigma_k}\lambda^{\asc\sigma}(\lambda-1)^{\des\sigma} X_{\sigma(1)}\ldots X_{\sigma(k)}.
\plabel{eq:muresdef}
\end{equation}
\begin{proof}[Proof]
It is well-known that if $A$ is a formal perturbation of $1$, and $\boldsymbol\varepsilon$
is a formal perturbation of $0$, then
\[(A+\boldsymbol\varepsilon)^{-1}-A^{-1}= -A^{-1}\boldsymbol\varepsilon A^{-1}+O(\boldsymbol\varepsilon)^2,\]
where $O(\boldsymbol\varepsilon)^2$ simply means terms which can be expressed with higher multiplicative
multiplicity in $\boldsymbol\varepsilon$.
In similar manner, one can prove
\[
\mathcal R^{(\lambda)}(A+\boldsymbol\varepsilon)-\mathcal R^{(\lambda)}(A)=
(1+(\lambda-1) \mathcal R^{(\lambda)}(A))\boldsymbol\varepsilon(1+(\lambda-1) \mathcal R^{(\lambda)}(A))+O(\boldsymbol\varepsilon)^2.
\]
Furthermore, substituting, $\boldsymbol\varepsilon\mapsto A\boldsymbol\varepsilon$, one finds
\begin{equation}
\mathcal R^{(\lambda)}(A(1+\boldsymbol\varepsilon))-\mathcal R^{(\lambda)}(A)=
(1+\lambda \mathcal R^{(\lambda)}(A))\boldsymbol\varepsilon(1+(\lambda-1) \mathcal R^{(\lambda)}(A))+O(\boldsymbol\varepsilon)^2.
\plabel{eq:formresdiff}
\end{equation}
Applying this to $A=\exp(X_1)\cdot\ldots\cdot\exp(X_{k-1})$ and $\boldsymbol\varepsilon=\exp(X_k)-1$, we find
\begin{align}\mathcal R^{(\lambda)}\,&(\exp(X_1)\cdot\ldots\cdot\exp(X_k))=\notag\\
=\,&\mathcal R^{(\lambda)}(\exp(X_1)\cdot\ldots\cdot\exp(X_{k-1}))\notag\\
&+X_k\notag\\
&+\lambda\, \mathcal R^{(\lambda)}(\exp(X_1)\ldots\cdot\exp(X_{k-1}))X_k \notag\\
&+(\lambda-1)\, X_k\mathcal R^{(\lambda)}(\exp(X_1)\cdot\ldots\cdot\exp(X_{k-1}))\notag\\
&+\lambda(\lambda-1)\, \mathcal R^{(\lambda)}(\exp(X_1)\cdot\ldots
\cdot\exp(X_{k-1})) X_k\mathcal R^{(\lambda)}(\exp(X_1)\cdot\ldots\cdot\exp(X_{k-1}))\notag\\
&+H(X_1,\ldots,X_k),\notag
\end{align}
where $H(X_1,\ldots,X_k)$ contains some terms where some variables $X_i$ have multiplicity more than $1$.

Using this as an induction step, one can prove that,
in terms of formal variables,
\begin{multline}\mathcal R^{(\lambda)}(\exp(X_1)\cdot\ldots\cdot\exp(X_k))=\\=
\sum_{\substack{\mathbf i=(i_1,\ldots, i_l)\in \{1,\ldots,k\}^l
  \\ i_a\neq i_b,\,l\geq1 }}
\lambda^{\asc(\mathbf i)}(\lambda-1)^{\des(\mathbf i)} X_{i_1}\cdot \ldots\cdot X_{i_l} +H(X_1,\ldots,X_k),\plabel{eq:resbare}
\end{multline}
where $H(X_1,\ldots,X_k)$ collects the terms with multiplicities in the variables.
Taking this in degree $k$, yields the statement.
\end{proof}

\begin{proof}[Alternative proof] (We model this proof after the proof of Theorem \ref{th:mp}.)
Considering  the formal power series $\mathcal R^{(\lambda)}(1+x)$, one can see that
$\mu_k^{(\lambda)}(X_1,\ldots,X_k)$ is a sum of terms
\[(\lambda-1)^{j-1}
X_{\sigma(1)}\cdot\ldots\cdot X_{\sigma(l_1)}\mid  X_{\sigma(l_1+1)}\cdot\ldots\ldots
\cdot X_{\sigma(l_{j-1})}\mid X_{\sigma(l_{j-1}+1)}\cdot\ldots\cdot X_{\sigma(k)},\]
where separators show the ascendingly indexed components which enter into power series of $\mathcal R^{(\lambda)}(1+x)$.
Considering a given permutation $\sigma$, the placement of the separators is either necessary (in case $\sigma(i)>\sigma(i+1)$),
or optional (in case  $\sigma(i)<\sigma(i+1)$).
Summing over the $2^{\asc(\sigma)}$ many optional possibilities, the coefficient of $X_{\sigma(1)}\cdot\ldots\cdot X_{\sigma(k)}$
is
\[\sum_{p=0}^{\asc(\sigma)}\begin{pmatrix}\asc(\sigma)\\p\end{pmatrix}(\lambda-1)^{p+\des\sigma}.\]
This simplifies  according to the identity
$\sum_{p=0}^{r}\binom rp(\lambda-1)^{p+d}= \lambda^r(\lambda-1)^d$ (which is but the binomial identity).
\end{proof}
\end{theorem}

\begin{proof}[Alternative proof for Theorem \ref{th:mp}]
If we integrate $\mathcal R^{(\lambda)}\left( (\exp X_1)\cdot\ldots\cdot(\exp X_k)\right)$ in $\lambda\in[0,1]$,
then we obtain $\log\left( (\exp X_1)\cdot\ldots\cdot(\exp X_k)\right)$
Thus, if we integrate
$\mu_k^{(\lambda)}(X_1,\ldots,X_k)$ in $\lambda\in[0,1]$, then we obtain $\mu_k(X_1,\ldots,X_k)$.
The actual integration (up to signs) is just the beta function identity; yielding $\mu_k(X_1,\ldots,X_k)$ as indicated.
\end{proof}

Assume that $\phi$ is a continuous  $\mathfrak A$-valued measure of finite variation on the interval $I$.
Let us define the Mielnik--Pleba\'nski integrals
\begin{align}
\mu^{(\lambda)}_{k,\mathrm R}(\phi)
&=\int_{t_1\leq\ldots\leq t_k\in I}\mu_k^{(\lambda)}(\phi(t_1),\ldots,\phi(t_k))
\plabel{eq:MPinter}\\\notag
&=\int_{\mathbf t=(t_1,\ldots,t_k)\in I^k}\lambda^{\asc(\mathbf t)}(\lambda-1)^{\des(\mathbf t)}\phi(t_1)\ldots\phi(t_k),
\end{align}
analogously to the terms of the Magnus expansion.
Then
\begin{equation}
\mu_{k,\mathrm R}(\phi)=\int_{\lambda=0}^1\mu^{(\lambda)}_{k,\mathrm R}(\phi)\,\mathrm d\lambda.
\plabel{eq:resmag}
\end{equation}

\begin{lemma}\plabel{lem:reslem}
Let $\phi$ be a continuous $\mathfrak A$-valued measure of finite variation on the interval $I$.
Let $T$ be a formal variable.
Then
\begin{equation}
\underbrace{\left(\lambda 1+(1-\lambda)\sum_{m=0}^\infty T^{m}\exp_{m,\mathrm R}(\phi)
\right)}_{\equiv\lambda+(1-\lambda)\Rexp(T\cdot\phi)}
\times
\left(\sum_{k =1}^\infty T^{k }\mu_{k ,\mathrm R}^{(\lambda)}(\phi)\right)
=\underbrace{\sum_{n =1}^\infty T^{n}\exp_{n ,\mathrm R}(\phi)}_{\equiv \Rexp(T\cdot\phi)-1},
\plabel{eq:reslem}
\end{equation}
where $\times$ means product in arbitrary order.
\begin{proof}
The very fact that $\mu_k^{(\lambda)}(X_1,\ldots,X_k)$ is the universally $1$-homogeneous part of
$\mathcal R^{(\lambda)}\left( (\exp X_1)\cdot\ldots\cdot(\exp X_k)\right)$
sufficiently many identities to prove the contraction \eqref{eq:reslem}:

In degree $k$, that follows from the identity $\lambda1+(1-\lambda)\exp\left( (\exp X_1)\cdot\ldots\cdot(\exp X_k)\right)$
times $\mathcal R^{(\lambda)}\left( (\exp X_1)\cdot\ldots\cdot(\exp X_k)\right)$ equals $\exp\left( (\exp X_1)\cdot\ldots\cdot(\exp X_k)\right)-1$, but the universally $1$-homogeneous part taken in the variables $X_1,\ldots,X_n$,
and then $X_i$ is replaced by $\phi(t_i)$, and then integrated in $t_1\leq\ldots\leq t_k\in I$.
\end{proof}
\end{lemma}
\begin{theorem}\plabel{th:MPresF} (The resolvent formula of Mielnik, Pleba\'nski \cite{MP}, a formal version.)
Let $\phi$ be a continuous $\mathfrak A$-valued measure of finite variation on the interval $I$.
Let $T$ be a formal variable.
Then
\begin{equation}\mathcal R^{(\lambda)}(\Rexp(T\cdot\phi))=
\sum_{k=1}^\infty \mu^{(\lambda)}_{k,\mathrm R}(T\cdot\phi).\plabel{eq:MPresT}
\end{equation}
\begin{proof}
This is immediate from the previous lemma, as the equational definition to the resolvent expression
(cf. Remark \ref{rem:Res1}) is satisfied.
\end{proof}
\end{theorem}

\begin{theorem}\plabel{th:MPres} (The resolvent formula of Mielnik, Pleba\'nski \cite{MP}, analytic version.)

Let $\phi$ be a continuous $\mathfrak A$-valued measure of finite variation on the interval $I$.
If
\begin{equation}
\sum_{k=1}^\infty \left|\mu^{(\lambda)}_{k,\mathrm R}(\phi)\right| <+\infty ,\plabel{eq:preMPres}
\end{equation}
then
\begin{equation}\mathcal R^{(\lambda)}(\Rexp(\phi))=
\sum_{k=1}^\infty \mu^{(\lambda)}_{k,\mathrm R}(\phi).\plabel{eq:MPres}
\end{equation}
\begin{proof}
Due to absolute convergence $T=1$ can be substituted into \eqref{eq:reslem}.
Then the equational definition to the resolvent expression is satisfied.
\end{proof}
\end{theorem}
We may use the notation $\mu^{(\lambda)}_{ \mathrm R}(\phi)$ for \eqref{eq:MPres}
(if it exist as an absolute convergent sum).
\begin{point}\plabel{rem:resseries}
Let us apply the previous statement with $\phi=X\mathbf 1_{[0,1)}$ in the formal setting.
Then \eqref{eq:MPres} and \eqref{eq:muresdef} give
\[\mathcal R^{(\lambda)}(\exp X)=\sum_{n=1}^\infty G_n(\lambda,\lambda-1)X^n\]
(by the very definition \eqref{eq:Eulerpol}, with the factorial coming from the integration domain).

Alternatively, considering $XG(uX,vX)$ with $u=\lambda$, $v=\lambda-1$, Theorem \ref{th:Euler} yields
\[\sum_{n=1}^\infty G_{n}(\lambda,\lambda-1)X^{n}
=XG(\lambda X,(\lambda-1)X))=\frac{(\exp X)-1}{\lambda+(1-\lambda)(\exp X)}\equiv\mathcal R^{(\lambda)}(\exp X).\]

It is easy to see that the knowledge of $G(u,v)$ is equivalent to the one of $G(\lambda,\lambda-1;X)$;
thus, keeping the first argument and
reversing the direction in the latter argument, we can view it as an alternative proof for Theorem \ref{th:Euler};
we just have to substitute  $X=u-v$, $\lambda=\frac{u}{u-v}$
into $\frac1XG(\lambda,(\lambda-1)X)$, which happens to work due to the special shape of the generating series.
\end{point}

Let us take another look at the Eulerian generating series.
Let $\Theta^{(\lambda)}_n$ be $1/n!$ times the sum of the absolute value of the coefficients in $\mu^{(\lambda)}(X_1,\ldots,X_n)$;
and let $\Theta^{(\lambda)}(x)=\sum_{n=1}^{\infty}\Theta^{(\lambda)}_n(x)$  be the associated generating function.
Note that
\begin{equation}\Theta^{(\lambda)}(x)=xG(|\lambda|x,|1-\lambda|x)
=\sum_{0\leq m<n}\frac{A(n,m)}{n!} |\lambda|^{n-1-m} |1-\lambda|^mx^n,
\plabel{eq:ThetaLdef}
\end{equation}
cf. Theorem \ref{th:Euler}.
This a conservative extension of the notation \eqref{eq:ThetaLdefpre} from $\lambda\in[0,1]$,
but it has now a direct meaning concerning the resolvent expressions.

Also note that $\Theta^{(\lambda)}(x)$ solves the differential equation (IVP)
\begin{equation}
\Theta^{(\lambda)\prime}(x)=(1+|\lambda| \Theta^{(\lambda)}(x) )(1+|1-\lambda| \Theta^{(\lambda)}(x) ),
\plabel{eq:Thetaeq}
\end{equation}
\[\Theta^{(\lambda)}(0)=0;\]
which shows that the naive majorizing estimate obtainable from  \eqref{eq:formresdiff} is actually exact.

\begin{lemma}\plabel{lem:Thetaspecpre}
If $\lambda\in[0,1]$, then $\Theta^{(\lambda)}(x)=xG(\lambda x,(1-\lambda)x)$, and
\newcommand{\tanhti}{\frac{\tanh \frac {2\lambda-1}2 x}{\frac {2\lambda-1}2  }}
\[
\Theta^{(\lambda)}(x)=
\begin{cases}
\dfrac{\tanhti}{1-\frac12\tanhti}
=\dfrac{\mathrm e^{(1-\lambda)x}-\mathrm e^{\lambda x}}{(1-\lambda)\mathrm e^{\lambda x}
-\lambda\mathrm e^{(1-\lambda)x}}&\text{if}\quad \lambda\in[0,1]\setminus \{\frac12\},\\
\dfrac x{1-\frac 12x }&\text{if}\quad \lambda=\frac12.
\end{cases}
\]
In terms of Taylor series,
\begin{equation}
\Theta^{(\lambda)}(x)=x+\frac12x^2+\frac{1 +2\lambda(1-\lambda)}{6}x^3+\ldots\qquad\text{if}\quad \lambda\in[0,1].
\plabel{eq:ThetaTay}
\end{equation}
More generally, for $\lambda\in\mathbb C$,
\begin{equation}
\Theta^{(\lambda)}(x)=\frac{\Theta^{\left(\frac{|\lambda|}{|\lambda|+|1-\lambda|}\right)}((|\lambda|+|1-\lambda|)x)}{|\lambda|+|1-\lambda|}.
\plabel{eq:Thetaspec}
\eqed
\end{equation}
\end{lemma}
\begin{commentx}
In particular,
\[
\Theta^{(\lambda)}(x)=\dfrac{2\tanh\frac x2}{1-|2\lambda-1|\tanh\frac x2}
=\dfrac{\mathrm e^x-1}{|1-\lambda|\mathrm e^x+|\lambda|}\qquad\text{if}\quad \lambda\in\mathbb R\setminus[0,1].
\eqedremark
\]
\begin{remark}\plabel{rem:Thetaspecpre}
As a consequence of \eqref{eq:Thetaspec} and \eqref{eq:ThetaTay},
\[\Theta^{(\lambda)}(x)=x+\frac{|\lambda|+|1-\lambda|}2x^2+\frac{(|\lambda|+|1-\lambda|)^2+2|\lambda|\cdot|1-\lambda|}6x^3+\ldots\,.\]

Equation \eqref{eq:Thetaspec} and Lemma \ref{lem:ThetaMest} shows that
\[\Theta^{(\lambda)}(x)\stackrel{\forall x}\leq\frac{\Theta^{(1/2)}((|\lambda|+|1-\lambda|)x)}{|\lambda|
+|1-\lambda|}\equiv\frac{x}{\frac1{|\lambda|+|1-\lambda|}-\frac12x},\]
i. e.
\[\Theta^{(\lambda)}_{k}\leq\left(\frac{|\lambda|+|1-\lambda|}{2}\right)^{k-1}.\eqedremark\]
\end{remark}
\end{commentx}

%\end{commentx}

Then, analogously to the Magnus case, from \eqref{eq:ThetaLdef},
\[\sum_{k=1}^\infty T^k\cdot\left|\mu^{(\lambda)}_{k,\mathrm R}(\phi)\right|\stackrel{\forall T}\leq
\sum_{k=1}^\infty T^k\cdot\int_{t_1\leq\ldots
\leq t_k\in I}\left|\mu^{(\lambda)}_k(\phi(t_1),\ldots,\phi(t_k))\right|\stackrel{\forall T}\leq
\sum_{k=1}^\infty T^k\cdot\Theta^{(\lambda)}_k\cdot\left(\int|\phi| \right)^k,\]
and,
\begin{equation}
\sum_{k=1}^\infty\left|\mu^{(\lambda)}_{k,\mathrm R}(\phi)\right|
\leq\sum_{k=1}^\infty\int_{t_1\leq\ldots\leq t_k\in I}\left|\mu^{(\lambda)}_k(\phi(t_1),\ldots,\phi(t_k))\right|
\leq\Theta^{(\lambda)}_\real\left(\int|\phi|\right);
\plabel{eq:measummu}
\end{equation}
with equalities for $\phi=c\cdot \mathrm Z^1_{[0,1)}$.

\begin{theorem} \plabel{cor:MPresshort}
(Resolvent formula via norm control.)
If $\phi$ is a continuous $\mathfrak A$-valued measure, and $\Theta^{(\lambda)}\left(\int |\phi|\right)<+\infty$, then
 $\mathcal R^{(\lambda)}\left(\Rexp(\phi)\right)$ exists,
 $\sum_{k=1}^\infty \left|\mu_{k,\mathrm R}^{(\lambda)}(\phi)\right|<+\infty $, and
\begin{equation}
\mathcal R^{(\lambda)}\left(\Rexp(\phi)\right)= \sum_{k=1}^\infty \mu_{k,\mathrm R}^{(\lambda)}(\phi).
\plabel{eq:calado}
\end{equation}
\begin{proof} In this case,
$\sum_{k=1}^\infty\left|\mu^{(\lambda)}_{k,\mathrm R}(\phi)\right|\leq \Theta^{(\lambda)}_\real(x)<+\infty$,
and Theorem \ref{th:MPres} applies.
\end{proof}
\end{theorem}

%\snewpage
\begin{lemma}
\plabel{lem:convradplain}
For $\lambda\in[0,1]$, the convergence radius of $\Theta^{(\lambda)}(x)$ around $x=0$ is
\[{\mathrm C}_\infty^{(\lambda)}=\begin{cases}
2&\text{if }\lambda=\frac12
\\
\dfrac{2\artanh (1-2\lambda)}{1-2\lambda}=\dfrac{\log\dfrac{1-\lambda}{\lambda}}{1-2\lambda}&\text{if }\lambda\in(0,1)
\setminus\{\frac12\}\\
+\infty&\text{if }\lambda\in\{0,1\}.
\end{cases}\]
This is a strictly convex, nonnegative function in $\lambda\in(0,1)$,
symmetric for $\lambda\mapsto1-\lambda$;  its minimum is ${\mathrm C}_\infty^{(1/2)}=2$.
In particular, in $\lambda\in[0,1]$, it yields a $[2,+\infty]$-valued strictly convex continuous function.

For $\lambda\in[0,1]$, let
\[w^{(\lambda)}=1/{\mathrm C}_\infty^{(\lambda)}.\]
In $\lambda\in[0,1]$, it is a $[0,1/2]$-valued strictly convex continuous function, symmetric for $\lambda\mapsto1-\lambda$;
its maximum is $w^{(1/2)}=1/2$.
\begin{proof}
${\mathrm C}_\infty^{(\lambda)}$ is just the convergence radius of $G(\lambda x,(1-\lambda)x)$ in $x$, cf. Theorem \ref{th:Euler}.
The rest is elementary calculus.
\end{proof}
\end{lemma}
\begin{theorem}\plabel{th:gecco}
For $\lambda\in\mathbb C\setminus\left\{0,1\right\}$, the convergence radius of $\Theta^{(\lambda)}(x)$ around $x=0$ is
\[
{\mathrm C}_\infty^{(\lambda)}=
\begin{cases}
\dfrac{2\artanh  \frac{|1-\lambda|-|\lambda|}{|\lambda|+|1-\lambda|}}{|1-\lambda|-|\lambda|}
=\dfrac{\log{|1-\lambda|}-\log{|\lambda|}}{|1-\lambda|-|\lambda|} &\text{if }\Rea\lambda\neq\frac12,\\
\frac2{|\lambda|+|1-\lambda|}=\frac2{\sqrt{1+(2\Ima \lambda)^2}}&\text{if }\Rea\lambda=\frac12.
\end{cases}
\]
The measure $\mathrm Z^1_{[0,1)}$ does also play an extremal role here: For $r\geq0$,
\[\frac{\Rexp(r\cdot\mathrm Z^1_{[0,1)})-1}{\lambda+(1-\lambda)\Rexp(r\cdot\mathrm Z^1_{[0,1)})}\]
exists if and only if
\[{\mathrm C}_\infty^{(\lambda)}>r .\]
\begin{proof}
The first part follows from Lemma \ref{lem:Thetaspecpre}.
The resolvent expression always exists in $\mathrm F^{1,\mathrm{loc}}([0,1))$;
the point is when it falls into $\mathrm F^{1 }([0,1))$.
\end{proof}
\end{theorem}

\begin{theorem}\plabel{th:cogecco}
For $s\in\mathbb C$, let us set
\[
{\mathrm C}_\infty^{[s]}:={\mathrm C}_\infty^{(\frac{s}{s-1})}=
\begin{cases}
\frac{(\log|s|)|s-1|}{|s|-1} &\text{if }|s|\neq1,0\\
 |s-1| &\text{if }|s|=1,\\
 +\infty &\text{if }s=0.
\end{cases}
\]
Then, for $r\geq0$,
\[\frac{1}{s-\Rexp(r\cdot\mathrm Z^1_{[0,1)})}\]
exists if and only if
\[{\mathrm C}_\infty^{[s]}>r .\]
In other  terms, for $r\geq0$, regarding the spectrum,
\[\spec \Rexp(r\cdot\mathrm Z^1_{[0,1)})=\{s\in\mathbb C\,:\,{\mathrm C}_\infty^{[s]}\leq r \}.\]

Furthermore, if for any other ordered measure $\phi$, if $\smallint|\phi|\leq r$, then
\[\spec \Rexp(\phi)\subset\spec \Rexp(r\cdot\mathrm Z^1_{[0,1)}). \]
\begin{proof}
This is just the previous statement ``transformed''.
\end{proof}
\end{theorem}
By that, we have a full spectral inclusion theory for the Banach algebra valued time-ordered exponentials
in terms of the cumulative norm.

\begin{figure}[t]
\centering
\includegraphics[width=.639\textwidth]{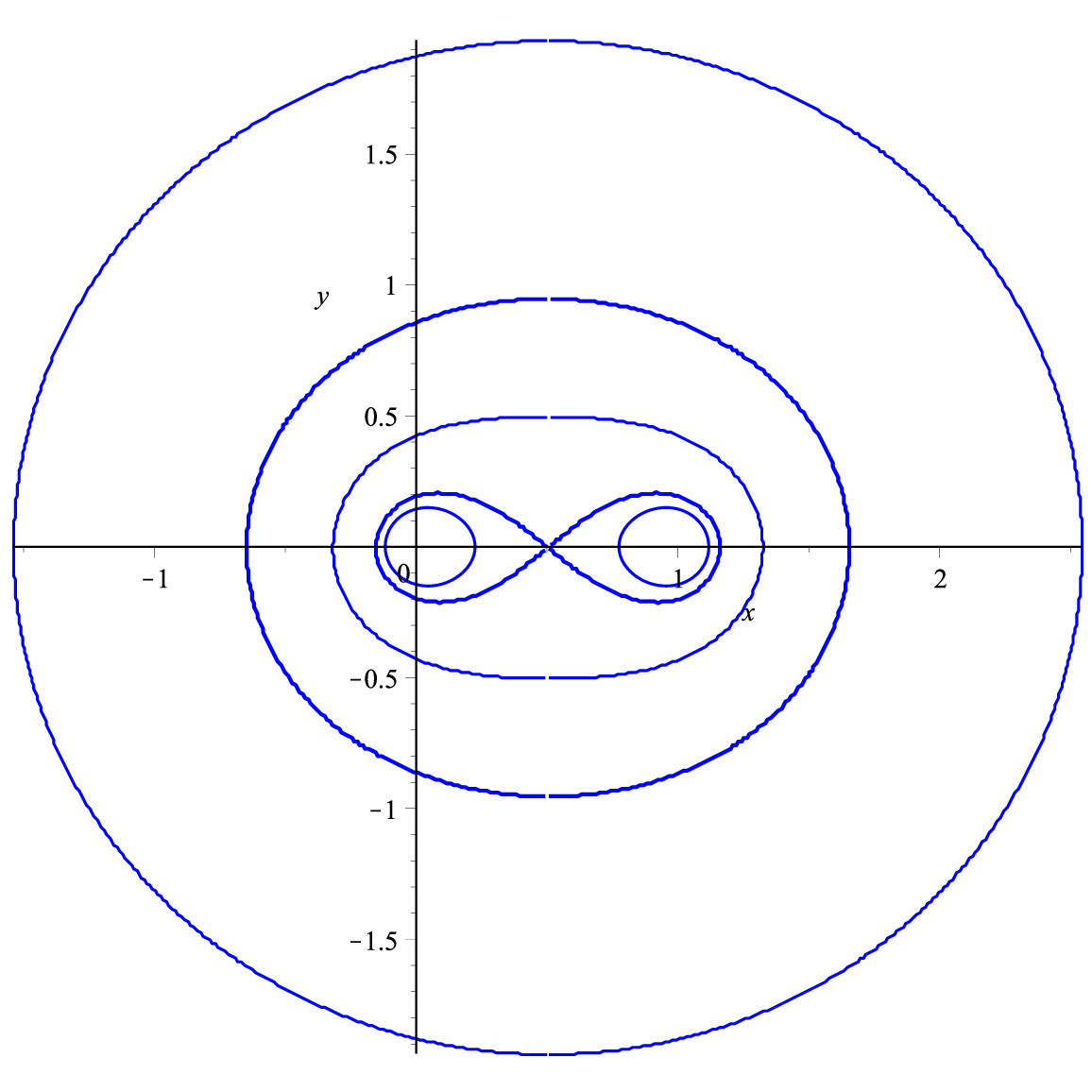}
\caption{The sets ${\mathrm C}_\infty^{(\lambda)}=r$ with
$\lambda=x+\mathrm iy$ and $r=\frac12$, $\boldsymbol{F_1}$, $\sqrt2$, $\boldsymbol2$, $F_2$.\plabel{fig:1}}
\vspace{3mm}
\centering
\includegraphics[width=.639\textwidth]{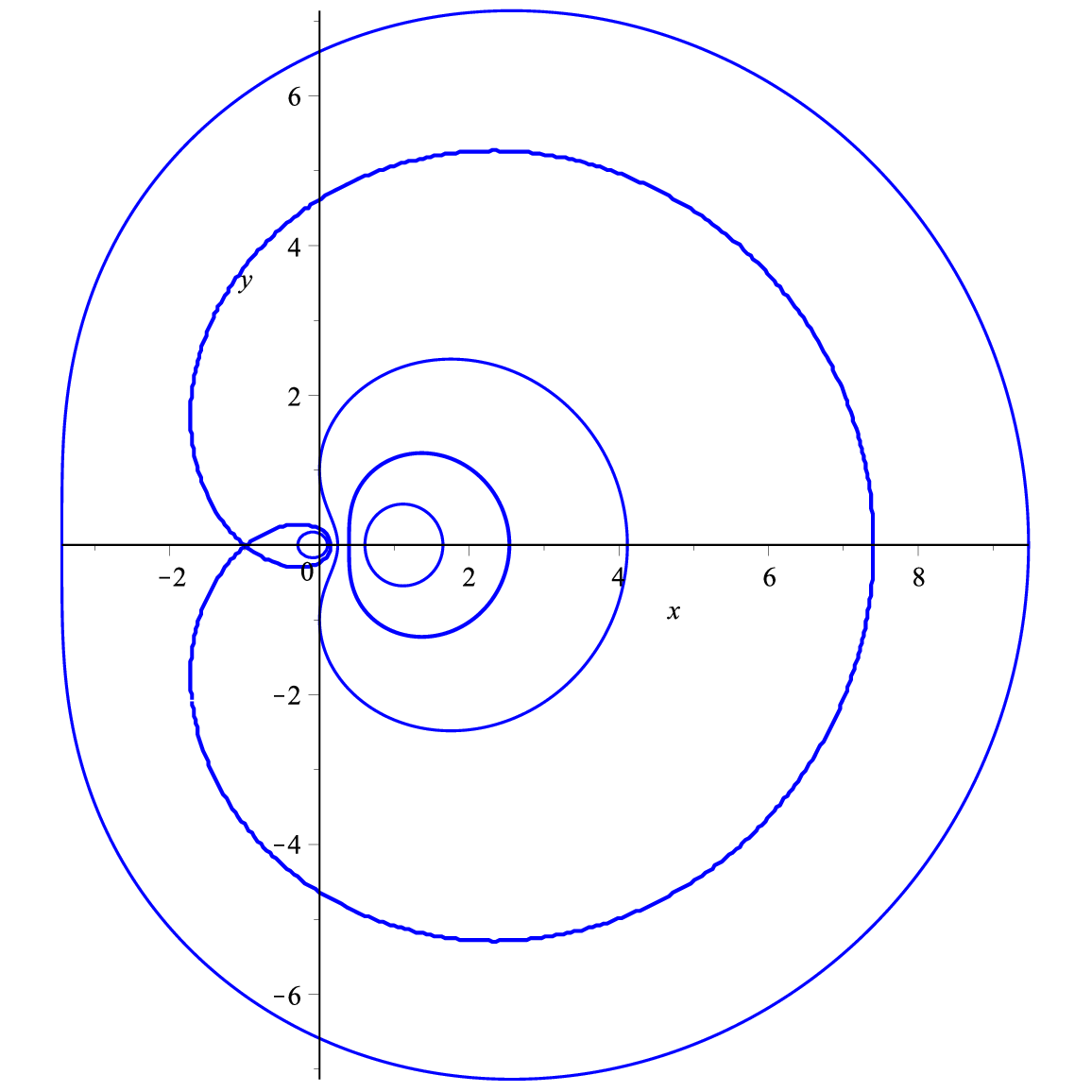}
\caption{The sets ${\mathrm C}_\infty^{[s]}=r$ with
$s=x+\mathrm iy$ and $r=\frac12$, $\boldsymbol{F_1}$, $\sqrt2$, $\boldsymbol2$, $F_2$.\plabel{fig:2}}
\end{figure}

In Figures \ref{fig:1} and \ref{fig:2}, we show the level sets ${\mathrm C}_\infty^{(\lambda)}=r$
 and ${\mathrm C}_\infty^{[s]}=r$ regarding some selected valued of $r$.
The selection of the $r$ is primarily according to the actual spectrum.
$r=2$ is when the spectrum first hits $(-\infty,0]$ (at $-1$);
$r=\sqrt2$ is when the spectrum first hits $\mathrm i\mathbb R$ (at $\pm\mathrm i$).
The values $r=F_1$ and $r=F_2$ are when when the critical convexity changes happen in the spectrum:
$F_1=-\log G_1 =0.930821\ldots$ is such that
$G_1=0.3942298\ldots$ is the solution of the equation $-\log(G_1)=G_1-2+\frac1{G_1}$ for $G_1\in(0,1)$;
$F_2=\frac{-G_2+1}{-G_2-1}\log (-G_2) =2.2484627\ldots$ is such that
$G_2=-3.44271909\ldots$ is the solution of the equation $\log(-G_2)=\frac{(G_2+1)(G_2-1)^2}{G_2(1-3G_2)} $ for $G_2\in(-\infty,-1)$.
(In appropriate circumstances the convexity of the level sets of $f(x,y)\leq c$
can be read off from the curvature function $\frac{-1}{\left(f_x)^2+(f_y)^2\right)^{3/2}}
\left|
\begin{smallmatrix}
f_{xx}&f_{xy}&f_x\\f_{xy}&f_{yy}&f_y\\f_x&f_y&0
\end{smallmatrix}
\right|$.)
The most interesting case is, of course, $r=2$.
${\mathrm C}_\infty^{(\lambda)}=2$ yields a ``logarithmic quasi lemniscate'';
${\mathrm C}_\infty^{[s]}=2$, the inverse image from $1$, is a ``logarithmic quasi prolate cardioid''.

Information about the (maximal) spectrum $\spec \Rexp(r\cdot\mathrm Z^1_{[0,1)})$ is certainly
 useful if one is to consider spectral calculus regarding time-ordered exponentials.
(Cf. Dunford, Schwartz \cite{DS}, Mart\'inez Carracedo, Sanz Alix \cite{MS}, or Haase \cite{Haa}
for the relevant spectral calculus.)
It is an obvious but suggestive property of this maximal spectrum that it is invariant for $s\mapsto\frac1s$.
Here we examine only the sectoriality of $\spec \Rexp(r\cdot\mathrm Z^1_{[0,1)})$ in details.
For a Banach algebra element $A$, its angle of sectoriality is the smallest $\omega\geq0$ such that
$\spec A\subset \{0\}\,\cup\, \{\exp z\,:\, |\Ima z|\leq \omega\}$.
The shape of this spectrum for $r>2$ is a (transformed) annulus around $0$.
\begin{theorem}\plabel{th:secto}
For $0\leq r \leq 2$, the angle of sectoriality of $\spec \Rexp(r\cdot\mathrm Z^1_{[0,1)})$ is
\[\omega(r)=2\arcsin\frac r2.\]

The critical rays $\mathrm e^{\pm \mathrm i\alpha(r)}\cdot[0,+\infty)$ meet
$\spec \Rexp(r\cdot\mathrm Z^1_{[0,1)})$ (only) in the unit circle.
\begin{proof}
Using the notation $s=x+\mathrm iy$, let us consider the radial (or, more precisely, dilatational) derivative
\[R(s):= \left.\frac{\mathrm d{\mathrm C}_\infty^{[ts]}}{\mathrm dt}\right|_{t=1}=x\frac{\partial}{\partial x} {\mathrm C}_\infty^{[x+\mathrm iy]}+y\frac{\partial}{\partial y} {\mathrm C}_\infty^{[x+\mathrm iy]}.\]
This is a smooth function but undefined for $s=0$ and $s=1$. More closely, we find:

(i) If $|s|\neq 1$ and $s\notin[0,+\infty)$, then
\[R(s)=\frac{\left(1+\tfrac{\Rea s}{|s|}\right)}{2|s-1|} (|s|-1)+\frac{\left(1-\tfrac{\Rea s}{|s|}\right)}{2|s-1|}
\frac{(|s|+1)(|s|^2-2|s|\log |s|-1)}{(|s|-1)^2}.\]
The fist summand is positive for $|s|>1$ and negative for $|s|<1$,
the second summand  nonnegative for $|s|>1$ and nonpositive for $|s|<1$.
Ultimately, $R(s)$ is positive for $|s|>1$ and negative for $|s|<1$.

(ii)  If $|s|= 1$ and $s\neq1$, then smoothness and point (i) implies that
\[R(s)=0.\]

(iii) If $s>1$, then
\[R(s)=1.\]

(iv) If $0<s<1$, then
\[R(s)=-1.\]

(v) If $s=1$, then $R(1)$ is undefined. (In the limit, it is the interval $[-1,1]$.)

(vi) If $s=0$, then $R(0)$ is undefined. (In the limit, it is $-1$.)

As $s=0$ cannot occur in the spectrum, this implies that the (radial) sectorial supporting lines can meet
the spectrum only in the unit circle $|s|=1$.
Then, for $s=\cos\omega+\mathrm i\sin\omega$, it yields  ${\mathrm C}_\infty^{[s]}=|s-1|=2\left|\sin\frac\omega2\right|$,
from which the statement about the angle of sectoriality follows.
\end{proof}
\end{theorem}
In fact, it is easy to say more: For $u=x+\mathrm iy\in\mathbb C$, let us set
\[
{\mathrm C}_\infty^{\{u\}}\equiv{\mathrm C}_\infty^{\{ x+\mathrm iy\}}
:=\frac{x}{\mathrm e^x-1}\sqrt{((\mathrm e^x\cos y)-1)^2+(\mathrm e^x\sin y)^2}
=\frac{x\sqrt{\mathrm e^{2x}-2\mathrm e^x\cos y+1}}{\mathrm e^x-1}
,
\]
where $\frac{x}{\mathrm e^x-1}$ is resolved as $1$ for $x=0$, etc.
This is continuous, and it is smooth except for $x=0$ and $y\in\pi\mathbb Z$.
For $r\geq0$, let
\[\mathfrak M_r=\{s\in\mathbb C\,:\,|\Ima u|\leq\pi\text{ and } {\mathrm C}_\infty^{\{u\}}\leq r\}.\]

Then it is easy to see that
\[\spec \Rexp(r\cdot\mathrm Z^1_{[0,1)})=\{s\in\mathbb C\,:\,{\mathrm C}_\infty^{[s]}\leq r \}=\exp \mathfrak M_r.\]
(This requires any effort only for $s\in(-\infty, 0]$, i. e. for $u=\pm\pi\mathrm i+\mathbb R$.)
As for special values,
\begin{equation}
{\mathrm C}_\infty^{\{x\}}=|x|,
\qquad
{\mathrm C}_\infty^{\{\mathrm iy\}}=\left|2\sin\frac y2\right|,
\qquad
{\mathrm C}_\infty^{\{x\pm\mathrm i\pi\}}=\frac{\mathrm e^x+1}{\mathrm e^x-1}x.
\plabel{eq:mspecial}
\end{equation}
In Figure \ref{fig:3}, we show the level sets ${\mathrm C}_\infty^{\{u\}}=r$ regarding some selected values of $r$.
\begin{figure}[ht]
\centering
\includegraphics[width=.639\textwidth]{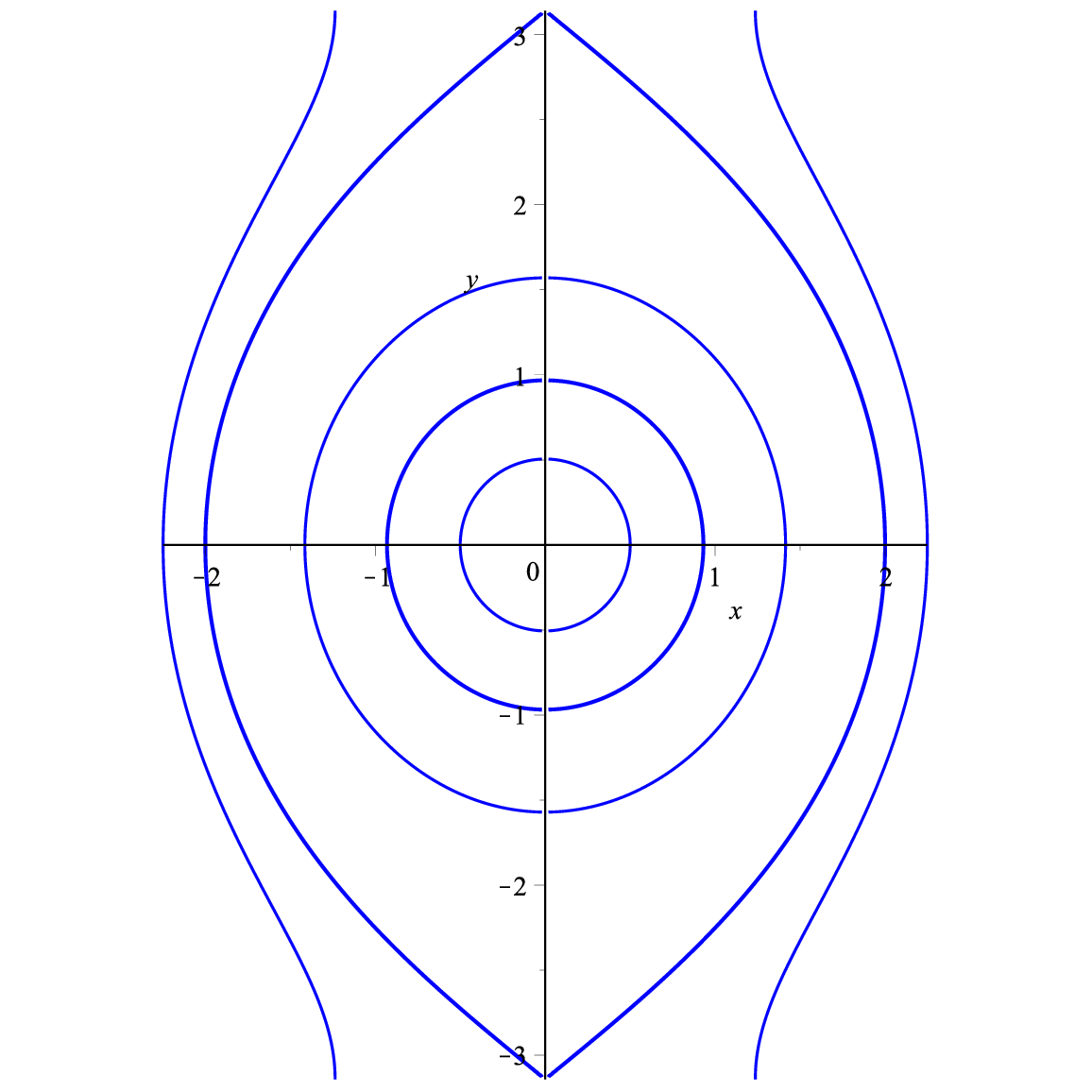}
\caption{The sets ${\mathrm C}_\infty^{\{u\}}=r$ with
$u=x+\mathrm iy$ and $r=\frac12$, $\boldsymbol{F_1}$, $\sqrt2$, $\boldsymbol2$, $F_2$.\plabel{fig:3}}
\end{figure}

\begin{lemma}\plabel{lem:logto}
(a) $\mathfrak M_r$ is star-like with respect to its center $0$.

(b) Its farthest points from $0$ are taken on $y$ axis or on the boundary $y=\pm\pi$.

(c) For $r<2$, $\mathfrak M_r$ lies in $\{z\in\mathbb C\,:\,\}$
the farther points from $0$ are taken on the $y$-axis with $|y|<\pi$;
$\mathfrak M_r=\log \spec \Rexp(r\cdot\mathrm Z^1_{[0,1)})$.

(e) For $r\geq2$,
the farther points of $\mathfrak M_r$ from $0$ are taken on the boundary $|y|=\pi$.

(e) In general,  $\mathfrak M_r$ is the closure of $\log \spec \Rexp(r\cdot\mathrm Z^1_{[0,1)})$.
\begin{proof}
For the radial derivative of ${\mathrm C}_\infty^{\{x+\mathrm iy\}}$,
\begin{multline*}
 x\frac{\partial}{\partial x} {\mathrm C}_\infty^{\{x+\mathrm iy\}}
+y\frac{\partial}{\partial y} {\mathrm C}_\infty^{\{x+\mathrm iy \}}
=  {\frac {{x}^{2}{{\mathrm e}^{\frac32\,x}}}{\sqrt {{{\mathrm e}^{2\,x}}-2\,{{\mathrm e}
^{x}}\cos \left( y \right) +1} \left( {{\mathrm e}^{x}}-1 \right) ^{2}}}
\cdot\Biggl(
\\
\frac {\left( {{\mathrm e}^{x}}-1 \right){ {\mathrm e}^{-\frac x2}} }{x}
\cdot
\left( y\sin y \right)
+
\frac { \left( x{{\mathrm e}^{x}}-2\,{{\mathrm e}^{x}}+x+2 \right){{\mathrm e}^{-\frac x2}} }{x}
\cdot
\left( 1+\cos y \right)
+
\\
+ {\frac { \left( {{\mathrm e}^{3\,x}}-2\,{{\mathrm e}^{2\,x}}x+{{\mathrm e}^{2\,x}}
-2\,x{{\mathrm e}^{x}}-{{\mathrm e}^{x}}-1 \right) {{\mathrm e}^{-\frac32\,x}}}{x}}
\Biggr).
\end{multline*}
One can check that in the big bracket each summand is nonnegative for $|y|<\pi$;
 the first one is positive if $y\neq 0$,
 the second and third ones are positive if $x\neq 0$.
Ultimately, one can deduce that the rotational derivative is positive for $|y|<\pi$ unless  $x=0$ and $y=0$.
This already implies that $\mathfrak M_r$ is star-like, and also the closure property.

The rotational derivative of ${\mathrm C}_\infty^{\{x+\mathrm iy\}}$ is
\begin{multline*}
 -y\frac{\partial}{\partial x} {\mathrm C}_\infty^{\{x+\mathrm iy\}}
+x\frac{\partial}{\partial y} {\mathrm C}_\infty^{\{x+\mathrm iy \}}
=-xy\frac{x\sqrt{\mathrm e^{2x}-2\mathrm e^x\cos y+1}}{\mathrm e^x-1}\cdot\Biggl(
\\
{\frac {-{x}^{2}{{\mathrm e}^{x}}+{{\mathrm e}^{2\,x}}-2\,{{\mathrm e}^{x}}+1}{{x}
^{2} \left( {{\mathrm e}^{x}}-1 \right) ^{2}}}
+
\frac{\left(  2\,y-3\,\sin   y  + y\cos   y  \right)
}{3y}\cdot  {\frac { {{\mathrm e}^{x}}}{ \left( {{\mathrm e}^{2\,x}}-2\,{{\mathrm e}^{x}}
\cos   y   +1 \right)  }}+
\\
 {\frac {x{{\mathrm e}^{-x}}+x{{\mathrm e}^{x}}+3\,{{\mathrm e}^{-x}}-3\,{
{\mathrm e}^{x}}+4\,x}{3x}}\cdot\left(  1-\cos   y\right)
\cdot
{\frac {{{\mathrm e}^{2\,x}}}{ \left( {{\mathrm e}^{x}}-1 \right) ^{2} \left(
{{\mathrm e}^{2\,x}}-2\,{{\mathrm e}^{x}}\cos   y   +1 \right) }}\Biggr).
\end{multline*}
One can check that in the big bracket each summand is nonnegative, and actually positive if $x\neq0$ and $y\neq0$;
therefore the same applies to the big bracket.
Ultimately, one can deduce that the rotational derivative never vanishes unless  $x=0$ or $y=0$.
That means that the farthest points of the level sets from the origin are on the coordinate axes or on the boundary $y=\pm\pi$.
From the sign behaviour of the rotational derivative, however, we can deduce that this cannot happen on the $x$-axis for $x\neq0$.
The latter thing can also be checked from the special values \eqref{eq:mspecial}, which also tell about the location of the
points of extremal distance.
\end{proof}
\end{lemma}

\begin{theorem}\plabel{th:logto}
For $0\leq r \leq 2$, the closure of $\log\spec \Rexp(r\cdot\mathrm Z^1_{[0,1)})$ is contained in the disc
$\Dbar(0, \omega(r))$, where
\[\omega(r)=2\arcsin\frac r2.\]

The closure of $\log\spec \Rexp(r\cdot\mathrm Z^1_{[0,1)})$ has common point with the boundary of the indicated
disk only at $\pm\mathrm i \omega(r)$.
\begin{proof}
This follows from the previous lemma and the special values \eqref{eq:mspecial}.
\end{proof}
\end{theorem}
An immediate consequence of the spectral mapping theorem is that if $\log\Rexp(\phi)$ is defined, then
its spectrum must lie in $\mathfrak M_{\int|\phi|}$.

In general, the most notable case for the resolvent is
$ \mathcal R^{(1/2)}(A)=-2 \frac{1-A}{1+A}$, which is $-2$ times the real involutive Cayley transform.
(This, in various contexts, also plays the role of ``the poor man's logarithm''.)
The existence of $\mathcal R^{(1/2)}(A)$, or, equivalently, $(A+1)^{-1}$ is the ``turn back'' problem,
as it is, classically, involves the eigenvalue $-1$ for $A$ (in an appropriate representation space).

\snewpage
\subsection{The logarithmic Magnus formula}\plabel{ss:MagLog}
~\\
\begin{theorem} \plabel{cor:MPshort} (Logarithmic Magnus formula via norm control.)
If $\phi$ is a continuous $\mathfrak A$-valued measure, and $\int |\phi|<2$, then
then $\log\Rexp(\phi)$ exists,
$\sum_{n=1}^\infty |\mu_{n,\mathrm R}(\phi)|<+\infty$ holds,
$\sum_{n=1}^\infty \max_{\lambda\in[0,1]}|\mu_{n,\mathrm R}^{(\lambda)}(\phi)|<+\infty$ holds, and
\begin{equation}
\log \Rexp(\phi)
=\int_{\lambda=0}^1\sum_{k=1}^\infty \mu_{k,\mathrm R}^{(\lambda)}(\phi)\,\mathrm d\lambda
=\sum_{k=1}^\infty \mu_{k,\mathrm R}(\phi)
.\plabel{eq:MPloge}
\end{equation}
Moreover, in these circumstances,regarding the spectral radius,
\[\mathrm r\left(\log\Rexp(\phi)\right)\leq\omega(r)\equiv2\arcsin\frac r2.\]
\begin{proof}
By \eqref{eq:measummu} and \eqref{eq:ThetaLest}, for $n\geq1$,
$|\mu_{n,\mathrm R}(\phi)|,|\mu_{n,\mathrm R}^{(\lambda)}(\phi)| \leq \Theta_{\real}^{(1/2)}(\int|\psi|)=2^{1-n}\cdot(\int|\psi|)^n$
holds (with any $\lambda\in[0,1]$), which imply the boundedness properties.
These estimates are sufficient to justify the application of Theorem \ref{th:MPres}, and
thus to establish the existence of $\log\Rexp(\phi)$.
Furthermore, the estimes are also sufficient to justify  the Beppo Levi / Fubini theorems in
\begin{multline}
\sum_{k=1}^\infty  \mu_{k,\mathrm R} (\phi)
=
\sum_{k=1}^\infty \int_{\lambda=0}^1\mu_{k,\mathrm R}^{(\lambda)}(\phi)\,\mathrm d\lambda
=\\=
\int_{\lambda=0}^1\sum_{k=1}^\infty \mu_{k,\mathrm R}^{(\lambda)}(\phi)\,\mathrm d\lambda
=
\int_{\lambda=0}^1 \mathcal R^{(\lambda)}(\Rexp(\phi))\,\mathrm d\lambda
=\log \Rexp(\phi),
\plabel{eq:corMPshort}
\end{multline}
yielding \eqref{eq:MPloge}.
The statement about the spectral radius is a consequence of spectral mapping theorem and Theorem \ref{th:logto}.
\end{proof}
\end{theorem}

\snewpage
Let $\phi$ be a continuous $\mathfrak A$-valued measure of finite variation on the interval $I$.
We define
\begin{equation}
\Theta^{(\lambda)}(\phi)=\sum_{k=1}^\infty \left|\mu^{(\lambda)}_{k,\mathrm R}(\phi)\right|.
\end{equation}
By Theorem \ref{th:MPres}, $\Theta^{(\lambda)}(\phi)<+\infty$ implies the existence of $\mathcal R^{(\lambda)}(\Rexp(\phi))$.
We also know that $\Theta^{(\lambda)}(\phi)\leq\Theta^{(\lambda)}_{\real}\left(\int|\phi|\right)$.
If $x\geq0$ equality is realized by $\Theta^{(\lambda)}\left(x\cdot\mathrm Z^1_{[0,1]}\right)=\Theta^{(\lambda)}_{\real}\left(x\right)$.
\begin{lemma}\plabel{lem:ThetaCont}
Let $\lambda,\nu\in\mathbb C$.

(a) If $\Theta^{(\lambda)}(\phi)<+\infty$ and $|\nu-\lambda|\Theta^{(\lambda)}(\phi)<1$, then

\[ |\Theta^{(\nu)}(\phi)-\Theta^{(\lambda)}(\phi)|\leq
 \frac{|\nu-\lambda|\Theta^{(\lambda)}(\phi)^2}{1-|\nu-\lambda|\Theta^{(\lambda)}(\phi) }. \]

(b)  If $\Theta^{(\lambda)}(\phi)=+\infty$, then
\[\Theta^{(\nu)}(\phi)\geq \frac1{|\nu-\lambda|}.\]
\begin{proof}
(a) If $T$ is a formal commutative variable, then on the formal level,
\[\mathcal R^{(\nu)}(\Rexp(T\cdot \phi))=\sum_{n=0}^\infty (\nu-\lambda)^n\mathcal R^{(\lambda)}(\Rexp(T\cdot \phi))^{n+1};\]
or, the resolvent expression expanded,
\[\sum_{k=1}^\infty T^k\mu_{k,\mathrm R}^{(\nu)}(\phi)=
\sum_{n=0}^\infty (\nu-\lambda)^n\left(\sum_{k=1}^\infty T^k\mu_{k,\mathrm R}^{(\lambda)}(\phi)\right)^{n+1}.\]
Slightly rewritten,
\[\left(\sum_{k=1}^\infty T^k\mu_{k,\mathrm R}^{(\nu)}(\phi)\right)
-
\left(\sum_{k=1}^\infty T^k\mu_{k,\mathrm R}^{(\lambda)}(\phi)\right)
=
\sum_{n=1}^\infty (\nu-\lambda)^n\left(\sum_{k=1}^\infty T^k\mu_{k,\mathrm R}^{(\lambda)}(\phi)\right)^{n+1}.\]
Now, applying straightforward norm estimates, the statement follows.

(b) This follows from (a) interchanging the role of $\lambda$ and $\nu$.
\end{proof}
\end{lemma}
\begin{cor}\plabel{cor:ThetaCont}
The function $\lambda\in\mathbb C\mapsto \Theta^{(\lambda)}(\phi)\in[0,+\infty]$ is a continuous functions
(as a $[0,+\infty]$-valued function).
\qed
\end{cor}
We can also define
\[\Theta^{([0,1]\mathrm m)}(\phi)=\sup_{\lambda\in[0,1]}\Theta^{(\lambda)}(\phi);\]
\[\Theta^{([0,1]\mathrm i)}(\phi)=\int_{\lambda=0}^1\Theta^{(\lambda)}(\phi)\,\mathrm d\lambda
=\int_{\lambda=0}^1
\sum_{k=1}^\infty \left|\mu^{(\lambda)}_{k,\mathrm R}(\phi)\right|
\,\mathrm d\lambda;
\]
\[\Theta(\phi)=
\sum_{k=1}^\infty \left| \int_{\lambda=0}^1\mu^{(\lambda)}_{k,\mathrm R}(\phi)\,\mathrm d\lambda\right|
=
\sum_{k=1}^\infty \left| \mu_{k,\mathrm R}(\phi)\right|.\]
Then $\Theta(\phi)<+\infty$ is just the convergence of the Magnus expansion.
We know that
$\Theta(\phi)\leq \Theta^{([0,1]\mathrm i)}(\phi)\leq\Theta^{([0,1]\mathrm m)}(\phi)\leq \Theta_{\real}^{(1/2)}(\int|\phi|)$,
and $\Theta(\phi)\leq \Theta(\int|\phi|) \leq\Theta_{\real}^{(1/2)}(\int|\phi|)$.
\snewpage
\begin{theorem}\plabel{th:control}
Let $\phi$ be a continuous $\mathfrak A$-valued measure of finite variation on the interval $I$.
Then the following conditions are equivalent:

(i)
\[\underbrace{\mathcal R^{(\lambda)}(\Rexp(t\cdot\phi))\text{ exists for any }\lambda\in[0,1]
}_{\equiv \log\Rexp(t\cdot\phi)\text{ exists}}, t\in\Dbar(0,1);\]

(ii)
\[\underbrace{\mathcal R^{(\lambda)}(\Rexp(t\cdot\phi))\text{ exists for any }\lambda\in[0,1]
}_{\equiv \log\Rexp(t\cdot\phi)\text{ exists}}, t\in\Dbar(0,1+\varepsilon),
\text{ with some }\varepsilon>0 ;\]

(iii)
\[\mathcal R^{(\lambda)}(\Rexp(t\cdot\phi))\text{ exists for any }\lambda\in[0,1], t\in\Dbar(0,1+\varepsilon),
\text{ with some }\varepsilon>0 ,\]
\[\text{ given by the locally uniformly convergent power series }\sum_{k=1}^\infty t^k \mu_{k,\mathrm R}^{(\lambda)}(\phi). \]

(iv)
\[\sum_{k=1}^\infty t^k\max_{\lambda\in[0,1]}\left|\mu_{k,\mathrm R}^{(\lambda)}(\phi)\right|<+\infty;
\text{ for } t\in[0,1+\varepsilon),
\text{ with some }\varepsilon>0,\]

(v)
\[\sum_{k=1}^\infty \max_{\lambda\in[0,1]}\left|\mu_{k,\mathrm R}^{(\lambda)}(\phi)\right|<+\infty
\text{\quad(i.~e.~uniformly absolute convergent in $\lambda\in[0,1]$)} ;\]

(vi) $\Theta^{([0,1]\mathrm m)}(\phi)<+\infty$, i. e.
\[\Theta^{(\lambda)}(\phi)\equiv\sum_{k=1}^\infty \left|\mu_{k,\mathrm R}^{(\lambda)}(\phi)\right|\text{ is uniformly bounded in }\lambda\in[0,1];\]

(vii) $\Theta^{([0,1]\mathrm i)}(\phi)<+\infty$, i. e.
\begin{equation}
\int_{\lambda=0}^1\Theta^{(\lambda)}(\phi) \,\mathrm d\lambda
\equiv
\int_{\lambda=0}^1\sum_{k=1}^\infty \left|\mu_{k,\mathrm R}^{(\lambda)}(\phi)\right|\,\mathrm d\lambda<+\infty,
\plabel{eq:MPpre}
\end{equation}

(viii)
\[\Theta^{(\lambda)}(\phi)\equiv\sum_{k=1}^\infty \left|\mu_{k,\mathrm R}^{(\lambda)}(\phi)\right|\text{ is pointwise bounded for }\lambda\in[0,1];\]

\begin{proof}

As $t\mapsto\exp(t\cdot\phi)$ is continuous, (i)$\Rightarrow$(ii) follows from a simple compactness argument.
Assume now (ii). Then $(t,\lambda)\mapsto\mathcal R^{(\lambda)}(\exp(t\cdot\phi))$ is analytic,
 and it has power series expansion in $t\in\Dbar(0,1+\varepsilon)$, uniformly in $\lambda$.
By Theorem \ref{cor:MPresshort}, for any $\lambda\in[0,1]$, it is true that for small $t$ equality \eqref{eq:calado} holds.
This identifies the power series of $(t,\lambda)\mapsto\mathcal R^{(\lambda)}(\exp(t\cdot\phi))$
 as $\sum_{k=1}^\infty t^k\mu_{k,\mathrm R}^{(\lambda)}(\phi)$.
This proves (ii)$\Rightarrow$(iii).
A standard application of the Cauchy formula gives (iii)$\Rightarrow$(iv).
The implications (iv)$\Rightarrow$(v)$\Rightarrow$(vi)$\Rightarrow$(vii) are obvious.
Assume that (vii) holds.
If $\Theta^{(\lambda)}(\phi)=+\infty$ for some $\lambda\in[0,1]$, then the function
$\nu\in[0,1]\mapsto\Theta^{(\nu)}(\phi)$ has ``at least a pole singularity'' by Lemma \ref{lem:ThetaCont}(b), making
 $\int_{\lambda=0}^1\Theta^{(\lambda)}(\phi)\,\mathrm d\lambda=+\infty$.
This contradiction implies  (vii)$\Rightarrow$(viii).
Assume that (viii) holds.
By Theorem \ref{th:MPres}, $\mathcal R^{(\lambda)}(\Rexp(\phi))$ exists for every $\lambda\in[0,1]$.
If $\phi$ is replaced by $t\cdot\phi$ ($t\in\Dbar(0,1)$), then $\mu_{k,\mathrm R}(\phi)$
 gets replaced by  $t^k\cdot\mu_{k,\mathrm R}(\phi)$, where $|t^k\cdot\mu_{k,\mathrm R}(\phi)|\leq \mu_{k,\mathrm R}(\phi)$.
Applying the previous argument for that we obtain (vi).
This demonstrates the implication (viii)$\Rightarrow$(i).
These implications together prove the statement.
\end{proof}
\end{theorem}
If the measure $\phi$ satisfies (any of) these equivalent conditions above, then we say that $\phi$ is $M$-controlled.
Then we have the following generalization of Theorem \ref{cor:MPshort}:
\begin{theorem}\plabel{th:logMagnus}(Logarithmic Magnus formula via analytic control.)
Let $\phi$ be a continuous $\mathfrak A$-valued measure of finite variation on the interval $I$.

If $\phi$ is $M$-controlled then $\log\Rexp(\phi)$ exists,
$\sum_{n=1}^\infty |\mu_{n,\mathrm R}(\phi)|<+\infty$ holds,
moreover,
$\sum_{n=1}^\infty \max_{\lambda\in[0,1]}|\mu_{n,\mathrm R}^{(\lambda)}(\phi)|<+\infty$ holds, and
\begin{equation}
\log \Rexp(\phi)
= \int_{\lambda=0}^1\sum_{k=1}^\infty \mu_{k,\mathrm R}^{(\lambda)}(\phi)\,\mathrm d\lambda
=\sum_{k=1}^\infty \mu_{k,\mathrm R}(\phi)
.\plabel{eq:MPlog}
\end{equation}
\begin{proof}
In the setting of Theorem \ref{th:control}(iii)/(iv), both the existence of $\log$,
and the bounds, and the sequence as in \eqref{eq:corMPshort} are straightforward.
\end{proof}
\end{theorem}

\begin{remark}\plabel{rem:Cayley}
If $\int |\phi|<2$, then
 $\Theta^{(\lambda)}(\phi)\leq \Theta^{(\lambda)}_\real(\int|\phi|)<+\infty$ holds  for $\lambda\in[0,1]$.
For  $\int |\phi|=2$, this also holds except at $\lambda=\frac12$.
Thus, in the critical case $\int |\phi|=2$, the possible divergence of the Magnus expansion
 comes from the case of $\mathcal R^{(1/2)}(\Rexp(\phi))$,
 which is $-2$ times the real involutive Cayley transform of $\Rexp(\phi)$.
If this Cayley transform exists, then $\log \Rexp(\phi)$ exists.
If $\Theta^{(1/2)}(\phi)<+\infty$, then, by Theorem \ref{th:logMagnus},
 the Magnus expansion is convergent and the logarithmic  Magnus formula holds.
\qedremark
\end{remark}

Now, the situation of Theorem \ref{th:logMagnus} is more general than the situation of Theorem \ref{cor:MPshort},
 yet it can be applied relatively flexibly, due to Theorem \ref{th:control}.
The situation of Theorem \ref{th:MagnusLogg} is even more general, but less practical.

\snewpage
\begin{remark}\plabel{rem:res}
Similarly to the Magnus expansion, one can formulate
Theorem \ref{th:MPres} and the subsequent discussion for not necessarily continuous measures, too.
Then, in \eqref{eq:MPinter}, we have to insert the multiplicity terms
\[\int_{t_1\leq\ldots \leq t_n\in I}\ldots\quad \rightsquigarrow\quad
\int_{t_1\leq\ldots \leq t_n\in I}\frac1{\mul(t_1,\ldots,t_n)!}\ldots\]
or, depending on viewpoint,
\begin{multline}\notag
\int_{\mathbf t=(t_1,\ldots,t_n)\in I^n}\lambda^{\asc(\mathbf t)}(\lambda-1)^{\des(\mathbf t)}
\ldots\quad \rightsquigarrow\\\rightsquigarrow\quad
\int_{\mathbf t=(t_1,\ldots,t_n)\in I^n}\lambda^{\asc(\mathbf t)}(\lambda-1)^{\des(\mathbf t)}
G_{\mul(t_1,\ldots,t_n)}(\lambda,\lambda-1)
\ldots
\end{multline}
throughout, where
\[G_{(m_1,\ldots,m_s)}(u,v)=G_{m_1}(u,v)\cdot\ldots\cdot G_{m_s}(u,v)\]
by definition.
Again, this is in accordance to taking the continuous blowup $\phi^*$ of $\phi$.

Then the Goldberg presentation reads as
the expression for the Magnus expansion via \eqref{eq:MPinter}/\eqref{eq:resmag}
in the extreme case when the measure is supported in two points.
\qedremark
\end{remark}

\begin{remark}
\plabel{rem:MPrem}
Mielnik, Pleba\'nski \cite{MP} contains several innovations.
On the combinatorial level, they introduce a symbolic ``chronological'' differential calculus.
Using that, they approach $\mu_k$ through taking linear combinations of partially selected ensembles of permutation monomials.
In those they are followed by Saenz, Suarez \cite{SS} but not by many.
(In particular, the combinatorial arguments we used are also entirely conventional,
 and we obtain the permutation-wise controlled formula \eqref{eq:mpform} directly.)
Then, they recognize the significance of the resolvent expression in the Magnus expansion.
(They are not much followed in that either, but Moan, Oteo \cite{MO} utilizes the related Eulerian generating function.)
Moreover,  they  develop  the corresponding formal holomorphic calculus; and they (re)obtain several elegant combinatorial identities.
\qedremark
\end{remark}

For $A\in\mathfrak A$, we define its Magnus exponent as
\[\MP_{\mathfrak A}(A):=\inf\left\{ \int |\phi|\,:\,\Rexp(\phi)=A \right\}. \]

Here it was left unexplained what kind ordered measures $\phi$ we consider at all.
Now, the point is that as long as we consider ordered measures of finite
variation (and this definition aims those) those measures can be modified
to multivariable BCH type, as it was explained in point \textbf{(II)}, at arbitrarily
small expense. Those can be smoothed out, etc. Ultimately, the Magnus exponent
is the same for  reasonable measure classes.
\begin{theorem}
\plabel{th:MPtheta}
If $\MP_{\mathfrak A}(A)<2$, then $A$ is $\log$-able, and
\[|\log A|\leq \Theta_\real(\MP_{\mathfrak A}(A)).\]
\begin{proof}
 Whenever $A=\Rexp(\phi)$, $\int |\phi|<2$,
by Theorems \ref{th:logMagnus} and \ref{th:mpest}, we know that $|\log A|=|\mu_{\mathrm R}(\phi)|\leq
\Theta_\real(\int |\phi|)$ holds.
\end{proof}
\end{theorem}

\snewpage

\subsection{Chronological decompositions the Magnus integrand}\plabel{ss:deco}
~

\begin{deco}\plabel{po:deco}
As a formal noncommutative power series, let us define
\begin{align}
\mathrm M^{(\lambda)}(X,Y)=&X(1-\lambda(\lambda-1) YX)^{-1}+Y(1-\lambda(\lambda-1) XY)^{-1}\plabel{eq:mult2}\\
&+\lambda XY(1-\lambda(\lambda-1) XY)^{-1}+(\lambda-1) YX(1-\lambda(\lambda-1) YX)^{-1}\notag\\
=&\sum_{k=0}^\infty\biggl(
\lambda^k(\lambda-1)^kY(XY)^{k}
+\lambda^k(\lambda-1)^kX(YX)^{k}\notag\\
&+\lambda^{k+1}(\lambda-1)^{k}(XY)^{k+1}
+\lambda^{k}(\lambda-1)^{k+1}(YX)^{k+1} \biggr),\notag
 \end{align}
and then, inductively,
\[\mathrm M^{(\lambda)}(X_1,\ldots,X_n)=\mathrm M^{(\lambda)}(\mathrm M^{(\lambda)}(X_1,\ldots,X_{n-1}),X_{n} ).\]

Then it is easy to show that
\[
\mathrm M^{(\lambda)}(X_1,\ldots,X_n)=
\sum_{\substack{i_1,\ldots,i_s\in\{1,\ldots,n\}\\ s\geq1,\,i_j\neq i_{j+1}}}
\lambda^{\asc(i_1,\ldots,i_s)}(\lambda-1)^{\des(i_1,\ldots,i_s)}
X_{i_1}\cdot\ldots\cdot X_{i_s};
\]
and  $\mathrm M^{(\lambda)}$ is an associative operation formally.

Let $\mathrm P^{(\lambda)}(X_1,\ldots,X_n)$ denote $\mathrm M^{(\lambda)}(X_1,\ldots,X_n)$
but here every coefficient of the monomials $X_{j_1}\cdot\ldots\cdot X_{j_s}$ is replaced by its absolute value.
It is notable that in the expansion of \eqref{eq:mult2} there are no monomials where a variable occurs twice consecutively,
and there are no multiplicities
among the monomials coming from the four terms of \eqref{eq:mult2}. This allow us to conclude that
\begin{multline}
\mathrm P^{(\lambda)}(X,Y)=Y(1-|\lambda|\cdot|\lambda-1| XY)^{-1}+(1-|\lambda|\cdot|\lambda-1| XY)^{-1}X\\
+|\lambda| XY(1-|\lambda|\cdot|\lambda-1| XY)^{-1}+|\lambda-1| Y(1-|\lambda|\cdot|\lambda-1| XY)^{-1}X;\notag
\end{multline}
and
\[\mathrm P^{(\lambda)}(X_1,\ldots,X_n)=\mathrm P^{(\lambda)}(\mathrm P^{(\lambda)}(X_1,\ldots,X_{n-1}),X_{n} ) ;\]
and  $\mathrm P^{(\lambda)}$ is still an associative operation formally.
Putting the formal commutative variables $x,y$ into the place of $X,Y$ we find
\[\mathrm P^{(\lambda)}(x,y)=\frac{x+y+(|\lambda|+|1-\lambda|)xy}{1-|\lambda|\cdot|1-\lambda|xy},\]
If $\lambda\in[0,1]$, then
\[\mathrm P^{(\lambda)}(x,y)=\frac{x+y+xy}{1-\lambda(1-\lambda)xy}.
\eqedexer\]
\end{deco}

\begin{lemma}\plabel{lem:Mcomputer}
 (a) $\mathrm P^{(\lambda)}(x_1,\ldots,x_n)$ applied with commutative variables is symmetric.

(b)\[\mathrm P^{(1/2)}(\underbrace{x,\ldots,x}_{\text{$n$ terms}})=\frac{2nx}{2-(n-1)x}.\]
\begin{proof}
(a) follows from the identities
$\mathrm P^{(\lambda)}(\mathrm P^{(\lambda)}(x_1,x_2),x_3)=\mathrm P^{(\lambda)}(x_1,\mathrm P^{(\lambda)}(x_2,x_3))$
and
$\mathrm P^{(\lambda)}(x_1,x_2)=\mathrm P^{(\lambda)}(x_2,x_1)$.
(b) is easy to prove by induction.
\end{proof}
\end{lemma}

\begin{remark}\plabel{rem:cgener}
These are also implied by the formula
$\mathrm P^{(\lambda)}(\Theta^{(\lambda)}(x_1),\ldots,\Theta^{(\lambda)}(x_n))=\Theta^{(\lambda)}(x_1+\ldots+x_n)$,
cf. Theorem \ref{thm:rresdecomp} later.
\qedremark
\end{remark}

Consider $\mu^{(\lambda)}_{k+n}(X_1,\ldots,X_k,Y_1,\ldots,Y_n)$.
Tentatively, it has shape
\begin{equation}\sum \lambda^{\#XY\text{-steps}} (\lambda-1)^{\#YX\text{-steps}}
\ldots\mu_{k_s}(X_{i_s},\ldots)\mu_{n_{s+1}}(Y_{j_{s+1}},\ldots)
\mu_{k_{s+2}}(X_{i_{s+1}},\ldots) \ldots
\plabel{eq:resdecomp}
\end{equation}
where $X$-blocks and $Y$-blocks follow each other, the variables are with increasing indexing in each
$X$-block and $Y$-block, respectively; and, multiplicatively, every variable has multiplicity $1$.
So, this decomposition, in its patterns, corresponds to $\mathrm M^{(\lambda)}(X,Y)$.
Similar picture applies in the case when the variables are partitioned to
$p$ many segments. Then it is related to $\mathrm M^{(\lambda)}$ with $p$ many variables.
\begin{theorem} \plabel{th:resdecomp}
Suppose that $\phi_1,\phi_2$ are $\mathfrak A$-valued ordered measures of finite variation.
Then
\[\underbrace{\sum_{n=1}^\infty \left|\mu^{(\lambda)}_{n,\mathrm R}(\phi_1\boldsymbol.\phi_2)\right|
}_{\equiv\Theta^{(\lambda)}(\phi_1\boldsymbol.\phi_2)}\leq
\mathrm P^{(\lambda)}_\real\Biggl(
\underbrace{\sum_{n=1}^\infty \left|\mu^{(\lambda)}_{n,\mathrm R}(\phi_1)\right|}_{\equiv\Theta^{(\lambda)}(\phi_1)}
,\underbrace{\sum_{n=1}^\infty \left|\mu^{(\lambda)}_{n,\mathrm R}(\phi_2)\right|}_{\equiv\Theta^{(\lambda)}(\phi_2)}\Biggr).\]
If the RHS is finite, then
\[\mu_{\mathrm R}^{(\lambda)}(\phi_1\boldsymbol.\phi_2)
=\mathrm M^{(\lambda)}\left(\mu_{\mathrm R}^{(\lambda)}(\phi_1),\mu_{\mathrm R}^{(\lambda)}(\phi_2)\right) ,\]
meaning so that the sum obtained from the formal power series of  $\mathrm M^{(\lambda)}$ is absolute convergent.
Analogous statements hold decompositioned to several parts.
\begin{proof}
This corresponds to the decomposition \eqref{eq:resdecomp} in the integrand.
The multidecomposition case can be proven analogously.
\end{proof}
\end{theorem}
\begin{theorem}\plabel{thm:rresdecomp} As formal commutative power series,
\[\Theta^{(\lambda)}(x+y)=\mathrm P^{(\lambda)}\left(\Theta^{(\lambda)}(x),\Theta^{(\lambda)}(y)  \right).\]
Consequently,
\[\Theta_\real^{(\lambda)}(x+y)=\mathrm P_\real^{(\lambda)}\left(\Theta_\real^{(\lambda)}(x),\Theta_\real^{(\lambda)}(y)  \right).\]
Analogous statements hold with several variables.
\begin{proof}
Let $\Theta^{(\lambda)}(x,y)=\sum_{k,n=0}^\infty\Theta^{(\lambda)}_{k,n}x^ky^n$, where $\Theta^{(\lambda)}_{k,n}$ is $\frac1{k!n!}$
times the sum of the absolute value of the coefficients in $\mu^{(\lambda)}(X_1,\ldots,X_k,Y_1,\ldots,Y_n)$.
Based on the placement of $X_1$ in the corresponding monomials, one has the formal differential equation / IVP
\[\frac{\partial}{\partial x}\Theta^{(\lambda)}(x,y)
=(1+|\lambda| \Theta^{(\lambda)}(x,y) )(1+|1-\lambda| \Theta^{(\lambda)}(x,y) ), \]
\[\Theta^{(\lambda)}(0,y)=\Theta^{(\lambda)}(y).\]
Comparison to \eqref{eq:Thetaeq} yields $\Theta^{(\lambda)}(x,y)=\Theta^{(\lambda)}(x+y)$.
On the other hand, the nature of decomposition \eqref{eq:resdecomp} implies
$\Theta^{(\lambda)}(x,y)=\mathrm P^{(\lambda)}\left(\Theta^{(\lambda)}(x),\Theta^{(\lambda)}(y) \right) $
(the multinomial identity is applies with respect to the distribution of the $X_i$ and $Y_j$, respectively).
The multivariable case follows by induction.
\end{proof}
\begin{proof}[Alternative proof]
Let $x,y\geq 0$ but $x+y<2$.
Let us apply Theorem \ref{thm:rresdecomp} to $\phi_1=\mathrm Z^1_{[0,x)}$ and  $\phi_2=\mathrm Z^1_{[x,x+y)}$.
Taking norm at both sides yields the equality
$\Theta_\real^{(\lambda)}(x+y)=\mathrm P_\real^{(\lambda)}\left(\Theta_\real^{(\lambda)}(x),\Theta_\real^{(\lambda)}(y)\right) $
for $x+y<2$.
This, however, already implies the equality of the formal series.
The multivariable case can be proven analogously.
\end{proof}
\begin{proof}[Third proof] We know the $\Theta^{(\lambda)}$ and $\mathrm P^{(\lambda)}$ in question explicitly, see
\eqref{eq:ThetaLdef} and Theorem \ref{th:Euler}, thus we can check the statement arithmetically.
The multivariable case follows by induction.
\end{proof}
\end{theorem}
\begin{cor}\plabel{cor:rresdecomp}
\[\mathrm P^{(\lambda)}\left(\Theta^{(\lambda)}(x), y \right)=
\Theta^{(\lambda)}\left(x+\left(\Theta^{(\lambda)}\right)^{-1}(y)\right);\]
and
\[\mathrm P_\real^{(\lambda)}\left(\Theta_\real^{(\lambda)}(x), y \right)=
\Theta_\real^{(\lambda)}\left(x+\left(\Theta_\real^{(\lambda)}\right)^{-1}(y)\right).\]
\begin{proof}In Theorem \ref{thm:rresdecomp},
substitute $\left(\Theta^{(\lambda)}\right)^{-1}(y)$ or $\left(\Theta_\real^{(\lambda)}\right)^{-1}(y)$
to the place of $y$, respectively.
\end{proof}
\end{cor}
\begin{remark}\plabel{rem:invtheta1}
For $\lambda\in[0,1]$,
\[
\left(\Theta^{(\lambda)}\right)^{-1}(y)=
\begin{cases}
\dfrac{\log\dfrac{1+\lambda y}{1+(1-\lambda)y}}{2\lambda-1}
=\dfrac{2\artanh\dfrac{(2\lambda-1)y}{2+y}}{2\lambda-1}&\text{if}\quad \lambda\in[0,1]\setminus \{\frac12\},\\\\
\dfrac y{1+\frac 12y }&\text{if}\quad \lambda=\frac12.
\end{cases}
\]
More generally, for $\lambda\in\mathbb C$,
\[\left(\Theta^{(\lambda)}\right)^{-1}(y)=\frac{\left(\Theta^{\left(\frac{|\lambda|}{
|\lambda|+|1-\lambda|}\right)}\right)^{-1}((|\lambda|+|1-\lambda|)y)}{|\lambda|+|1-\lambda|}.\eqedremark\]
\begin{commentx}
In particular,
\[\left(\Theta^{(\lambda)}\right)^{-1}(y)=
\dfrac{\log\dfrac{1+|\lambda| y}{1+|1-\lambda|y}}{\sgn(2\lambda-1)}
=2\artanh\dfrac{y}{2+|2\lambda-1|y}\qquad\text{if}\quad \lambda\in\mathbb R\setminus[0,1].
\eqedremark
\]
\end{commentx}
\end{remark}
\snewpage

Using the ``decomposition method'', we can see that any ``delay''
in the norm growth of the resolvent leads to improvement convergence:
\begin{theorem} \plabel{the:resdelay}
Suppose that $\phi_1,\phi_2,\phi_3$ are $\mathfrak A$-valued ordered measures of finite variation; $N\in[0,+\infty)$.
If
\[\Theta^{(\lambda)}(\phi_2)\equiv
\sum_{n=1}^\infty \left|\mu^{(\lambda)}_{n,\mathrm R}(\phi_2)\right|\leq N,\] then
\[\Theta^{(\lambda)}(\phi_1\boldsymbol.\phi_2\boldsymbol.\phi_3)\equiv
\sum_{n=1}^\infty \left|\mu^{(\lambda)}_{n,\mathrm R}(\phi_1\boldsymbol.\phi_2\boldsymbol.\phi_3)\right|\leq
\Theta^{(\lambda)}_\real\left(\int |\phi_1|+\int |\phi_3|+\left(\Theta^{(\lambda)}\right)^{-1}(N)\right).\]
\begin{proof} According Theorem \ref{th:resdecomp} (multidecomposition version),
\begin{align}
\sum_{n=1}^\infty \left|\mu^{(\lambda)}_{n,\mathrm R}(\phi_1\boldsymbol.\phi_2\boldsymbol.\phi_3)\right|&\leq
 \mathrm P^{(\lambda)}_\real\left(\Theta^{(\lambda)}_\real\left(\int |\phi_1|\right),N,\Theta^{(\lambda)}_\real
 \left(\int |\phi_3|\right) \right)\notag\\
 &=  \mathrm P^{(\lambda)}_\real\left( \mathrm P^{(\lambda)}_\real
 \left(\Theta^{(\lambda)}_\real\left(\int |\phi_1|\right),\Theta^{(\lambda)}_\real\left(\int |\phi_3|\right) \right),N \right)
\notag\\
&= \mathrm P^{(\lambda)}_\real\left( \Theta^{(\lambda)}_\real
\left(\int |\phi_1|+\int |\phi_3|\right) ,N \right).\plabel{eq:ensai}
\end{align}
Taking Corollary \ref{cor:rresdecomp} into account, we obtain the statement.
\end{proof}
\end{theorem}

\begin{theorem}
 \plabel{cor:resdelay}
Suppose that $\phi_1,\phi_2,\phi_3$ are $\mathfrak A$-valued ordered measures of finite variation,
$\phi=\phi_1\boldsymbol.\phi_2\boldsymbol.\phi_3$; $N\in[0,+\infty)$.
If
\[\Theta^{([0,1]\mathrm m)}(\phi_2)\equiv\sup_{\lambda\in[0,1]}\sum_{n=1}^\infty \left|\mu^{(\lambda)}_{n,\mathrm R}(\phi_2)\right|\leq N,\]
and
\[\int |\phi_1|+\int |\phi_3|<\frac{4}{2+N},\]
then
\begin{multline*}
\Theta^{([0,1]\mathrm m)}(\phi_1\boldsymbol.\phi_2\boldsymbol.\phi_3 )\equiv
\sup_{\lambda\in[0,1]}
\sum_{n=1}^\infty \left|\mu^{(\lambda)}_{n,\mathrm R}(\phi_1\boldsymbol.\phi_2\boldsymbol.\phi_3)\right|
\leq
\\
\leq\Theta^{(1/2)}_\real\left(\int |\phi_1|+\int |\phi_3|+\left(\Theta^{(1/2)}\right)^{-1}(N)\right)<+\infty;
\end{multline*}
consequently, the Magnus expansion is convergent and the logarithmic Magnus formula holds for $\phi$.
\begin{proof}
We continue the previous proof from \eqref{eq:ensai}.
For $\lambda\in[0,1]$,
$\Theta^{(\lambda)}(x)\stackrel{\forall x}\leq \Theta^{(1/2)}(x)$ (see \eqref{eq:ThetaMest}) and
$\mathrm P^{(\lambda)}(x,y)\stackrel{\forall x,y}\leq \mathrm P^{(1/2)}(x,y)$
(commutative variables, it is easy to see) imply
\[\mathrm P^{(\lambda)}_\real\left( \Theta^{(\lambda)}_\real
\left(\int |\phi_1|+\int |\phi_3|\right) ,N \right)
\leq \mathrm P^{(1/2)}_\real\left( \Theta^{(1/2)}_\real
\left(\int |\phi_1|+\int |\phi_3|\right) ,N \right). \]
Corollary \ref{cor:rresdecomp} and the concrete shape of $\Theta^{(1/2)}(x)$ imply the statement.
\end{proof}
\end{theorem}

For the sake of an alternative formulation; we set
\[N^{(\lambda)}(\phi)=\left(\Theta_{\real}^{(\lambda)}\right)^{-1}\left(\Theta^{(\lambda)}(\phi)\right),\]
where $\left(\Theta_{\real}^{(\lambda)}\right)^{-1}\left(+\infty\right)=\mathrm C^{(\lambda)}_\infty$.
One can see that the map
$(\lambda,x)\in\mathbb C\times[0,+\infty]\times [0,+\infty]\rightarrow\left(\Theta_{\real}^{(\lambda)}\right)^{-1}(x)\in [0,+\infty]$
is continuous (with the appropriate topology for $[0,+\infty]$).
This is a back-normalized ``resolvent mass''.
An immediate consequence of our earlier discussions is that
\[N^{(\lambda)}(\phi)\leq\min\left(\int|\phi|, \mathrm C^{(\lambda)}_\infty\right).\]
\begin{theorem} \plabel{th:ressubadd}
Suppose that $\phi_1,\ldots,\phi_n$ are $\mathfrak A$-valued ordered measures of finite variation.
Then
\[N^{(\lambda)}(\phi_1\boldsymbol.\ldots\boldsymbol.\phi_n)\leq
N^{(\lambda)}(\phi_1 )
+\ldots+
N^{(\lambda)}( \phi_n).\]
In particular, if $N^{(\lambda)}(\phi_1 )+\ldots+N^{(\lambda)}( \phi_n)<\mathrm C^{(\lambda)}_\infty$,
then $\Theta^{(\lambda)}(\phi_1\boldsymbol.\ldots\boldsymbol.\phi_n)<+\infty$.
\begin{proof}
This is an immediate consequence of Theorem \ref{th:resdecomp} and Theorem \ref{thm:rresdecomp}.
\end{proof}
\end{theorem}
This generalizes and reformulates Theorem \ref{the:resdelay}.
A complementary viewpoint is to consider the ``delay''
\[\Delta^{(\lambda)}(\phi)=\left(\int|\phi|\right)-N^{(\lambda)}(\phi).\]
Now the wisdom is that if we can exhibit nonzero delays we can improve the convergence
of the Magnus expansion.

There are possible variations toward similar but cumulative statements
(cf.~Theorem \ref{cor:resdelay}), but the $\lambda$-wise estimates (as above) are the most important.

\begin{theorem}\plabel{lem:mag2criter}
 Suppose that $\phi$ is a  $\mathfrak A$-valued measure of finite variation,
which arises as a Lebesgue-Bochner integrable function times the Lebesgue measure on an interval.
Assume $\lambda\in\mathbb C\setminus\left( (-\infty,0]\cup[1,+\infty) \right)$ and $\int|\phi|>0$.
Then the following hold:

(a) If $k\geq2$, then
\[\left|\mu^{(\lambda)}_{k,\mathrm R}(\phi)\right|<\Theta_k^{(\lambda)}\cdot\left(\int|\phi|\right)^k.\]

(b)
If $\int|\phi|\leq\mathrm C^{(\lambda)}_\infty$, then
\begin{equation}
\Theta^{(\lambda)}(\phi)\equiv
\sum_{n=1}^\infty \left|\mu^{(\lambda)}_{n,\mathrm R}(\phi)\right|<\Theta^{(\lambda)}_{\real}\left(\int|\phi|\right).
\plabel{eq:rana}
\end{equation}
In particular, the LHS is finite.
\begin{proof}
(a)
We prove only the case $k=2$, the higher cases follow the same principles.
Suppose that $\phi(t)=h(t)\mathbf 1_I$.
Let $\varepsilon>0$.
Then $h(t)$ can be approximated by step-functions (in $L^1$ norm) arbitrarily well.
This implies that there is a nontrivial subinterval $I'\subset I$ of nonzero length and a nonzero element $a\in\mathfrak A$ such that
$\int_{I'} |h(t)-a|\leq \varepsilon|a|\cdot|I'|$. Let $k(t)=h(t)-a$.
Note that
\[\mu^{(\lambda)}_2(X_1,X_2)=\lambda X_1X_2+(\lambda-1)X_2X_1.\]
Thus
\begin{align}
\Bigl|\int_{t_1\leq t_2\in I'}&\mu^{(\lambda)}_2(h(t_1),h(t_2))\,\mathrm dt_1\,\mathrm dt_2\Bigr|\leq
\notag\\
\leq& \left|\left(\lambda-\frac12\right)\int_{t_1,t_2\in I'} a^2 \,\mathrm dt_1\,\mathrm dt_2\right|
\notag\\
&+\max\left(|\lambda|,|\lambda-1|\right)\int_{t_1,t_2\in I'} |a|\,|k(t_2)|+|k(t_1)|\,|a|+|k(t_1)|\,|k(t_2)|  \,\mathrm dt_1\,\mathrm dt_2
\notag\\
\leq& |a|^2|I'|^2\left( \left|\lambda-\frac12\right|+\max\left(|\lambda|,|\lambda-1|\right)\cdot\left(2\varepsilon+\varepsilon^2\right) \right).
\notag
\end{align}

On the other hand, in the formal estimate we count at least
\begin{align}
&  \frac{|\lambda|+|\lambda-1|}2\int_{t_1,t_2\in I'} |a^2| \,\mathrm dt_1\,\mathrm dt_2
\notag\\
&-\max\left(|\lambda|,|\lambda-1|\right)\int_{t_1,t_2\in I'} |a|\,|k(t_2)|+|k(t_1)|\,|a|+|k(t_1)|\,|k(t_2)|  \,\mathrm dt_1\,\mathrm dt_2
\notag\\
\geq& |a|^2|I'|^2\left(\frac{|\lambda|+|\lambda-1|}2
 - \max\left(|\lambda|,|\lambda-1|\right)\cdot\left(2\varepsilon+\varepsilon^2\right) \right).
\notag
\end{align}

Between the two estimates there is a difference
\[|a|^2|I'|^2\left(\frac{|\lambda|+|\lambda-1|}2-\left|\lambda-\frac12\right|
 - 2\max\left(|\lambda|,|\lambda-1|\right)\cdot\left(2\varepsilon+\varepsilon^2\right) \right).\]

Assume that $\varepsilon>0$ was chosen originally so that
\[\frac{|\lambda|+|\lambda-1|}2-\left|\lambda-\frac12\right|
 - 2\max\left(|\lambda|,|\lambda-1|\right)\cdot\left(2\varepsilon+\varepsilon^2\right) =\delta>0.\]
Then,
\[\left|\mu^{(\lambda)}_{2,\mathrm R}(\phi|_{I'})\right|
+ |a|^2|I'|^2\delta\leq\Theta_2^{(\lambda)}\cdot\left(\int|\phi_{I'}|\right)^2.\]
Extending to estimate to $I$ trivially, we find
\[\left|\mu^{(\lambda)}_{2,\mathrm R}(\phi)\right|+ |a|^2|I'|^2\delta\leq\Theta_2^{(\lambda)}\cdot\left(\int|\phi|\right)^2.\]

(b) Adding further terms, we obtain
\begin{equation}
\sum_{n=1}^\infty \left|\mu^{(\lambda)}_{n,\mathrm R}(\phi)\right|+ |a|^2|I'|^2\delta\leq\Theta^{(\lambda)}_{\real}\left(\int|\phi|\right).
\plabel{eq:prerana}
\end{equation}

If  $\int|\phi|<\mathrm C^{(\lambda)}_\infty$, then $\Theta^{(\lambda)}_{\real}\left(\int|\phi|\right)<+\infty$,
 and \eqref{eq:prerana} implies \eqref{eq:rana}.
Note, in this case $N^{(\lambda)}(\phi)<\int|\phi|$.

Finally, the critical case $\int|\phi|=\mathrm C^{(\lambda)}_\infty$ remains to be treated.
Then, we can take a decomposition $\phi=\phi_1\boldsymbol.\phi_2$ such that
 $N^{(\lambda)}(\phi_1)<\int|\phi_1|$ and $N^{(\lambda)}(\phi_2)<\int|\phi_2|$ hold.
Consequently, $N^{(\lambda)}(\phi)=N^{(\lambda)}(\phi_1\boldsymbol.\phi_2)
=N^{(\lambda)}(\phi_1)+N^{(\lambda)}(\phi_2)< \int|\phi_1|+\int|\phi_2|=\mathrm C^{(\lambda)}_\infty$.
Then $N^{(\lambda)}(\phi)<\mathrm C^{(\lambda)}_\infty$ implies
$\Theta^{(\lambda)}(\phi)<+\infty=\Theta_{\real}^{(\lambda)}( \mathrm C^{(\lambda)}_\infty)$.
\end{proof}
\end{theorem}

\begin{theorem}\plabel{thm:ressub} Suppose that $\phi$ is a  $\mathfrak A$-valued measure,
 which arises as a Lebesgue-Bochner integrable function times the Lebesgue measure on an interval.
Suppose that $\smallint |\phi|=r$. (As we know, $\spec \Rexp(\phi)\subset \spec \Rexp(r\cdot \mathrm Z^1_{[0,1)})$.)
Then
\[\left(\spec \Rexp(\phi)\right) \,\cap\ \partial\left(\spec \Rexp(r\cdot \mathrm Z^1_{[0,1)}) \right)
\subset \{\mathrm e^r,\mathrm e^{-r}\}.\]
\begin{proof}
By Theorem \ref{lem:mag2criter}, there are ``delays'' for
$\lambda\in\mathbb C\setminus\left( (-\infty,0]\cup[1,+\infty) \right)$.
In terms of the ordinary spectrum that translates to delays for $s\in\mathbb C\setminus[0,+\infty)$.
\end{proof}
\end{theorem}

\begin{theorem}\plabel{thm:mag2crit} Suppose that $\phi$ is a  $\mathfrak A$-valued measure,
 which arises as a Lebesgue-Bochner integrable function times the Lebesgue measure on an interval.
Suppose that $\smallint |\phi|\leq2$.
Then $\phi$ is $M$-controlled.
In particular, the Magnus series is absolute convergent and the logarithmic Magnus formula holds.
\begin{proof}
By Theorem \ref{lem:mag2criter}(b), $\Theta^{(\lambda)}(\phi)<+\infty$ for any $\lambda\in[0,1]$.
By Theorem \ref{th:control}((viii)) we obtain that $\phi$ is $M$-controlled.
In particular, Theorem \ref{th:logMagnus} applies.
\end{proof}
\end{theorem}

\begin{remark}\plabel{lem:magleb1}
Technicalities aside, (the proof of) Theorem  \ref{lem:mag2criter} says that any internal commutativity in the Magnus integrand
 leads to some decrease in the ``resolvent mass''.
 In Theorem \ref{lem:mag2criter} and Theorem \ref{thm:mag2crit}, it is sufficient that
$\phi$ is Lebesgue-Bochner integrable times the Lebesgue measure only on a subinterval but with nonzero variation there.
This implies that these statements extend automatically to the case of noncontinuous measures:
There the blowup has a subinterval with the required properties.
With similar assumptions in Corollary \ref{cor:resdelay}, the boundedness of
$\Theta^{(\lambda)}(\phi)$ ($\lambda\in[0,1]$), and thus the convergence of the Magnus expansion
and the validity of the logarithmic Magnus formula,
extends to the case $\int |\phi_1|+\int |\phi_3|=\frac{4}{2+N}$.
\qedremark
\end{remark}

\begin{remark}\plabel{lem:magleb3}
However, if we ask only that $\phi$ is a Lebesgue-Bochner integral not in uniform but in strong sense,
 then one cannot make a similar statement.
Indeed, let us consider the measures
 $\Xi_{n}= X_{1}\underbrace{\mathbf 1_{[0,1/2^n)}\boldsymbol.\ldots\boldsymbol. X_{2^n}\mathbf 1_{[0,1/2^n)}}_{2^n\text{ terms}},$
 where the $X_j$ are formal noncommutative variables with $\ell_1$ norm.
If $c\geq 0$, then $\smallint |c\Xi_{n}|=c$, and it is easy to show that $n\nearrow\infty$ yields
 $|\mu_{\mathrm R}{}(c\Xi_{n})|\nearrow\Theta(c)$.
Taking an infinite direct sum of the $2\cdot\Xi_n$'s, we obtain a  strong Lebesgue-Bochner construction
 of variation $2$,
 whose Magnus expansion, however, cannot converge due to the  Banach-Steinhaus theorem.
\qedremark
\end{remark}
At this point we could argue that $2$ can be improved as the cumulative convergence radius of the Baker--Campbell--Hausdorff expansion.
Indeed, there is much ``internal commutativity'' there which can be quantified on the universal level.
However, this will be done more precisely in Section \ref{sec:DiscRes} by direct estimates.
\snewpage

Next, we address an issue for later reference.
In general, computing with $\Theta^{(\lambda)}$ ($\lambda\in[0,1]$) can be complicated.
It may be useful to reduce certain computations to  $\Theta^{(1/2)}$.

\begin{lemma}\plabel{lem:depp1}
%(Increasing delay principle.)
Assume that $\lambda\in[0,1]$, $x\in[0,+\infty)$, and $0<y<x$.
Then
\[
\Theta^{(\lambda)}(x)-\Theta^{(\lambda)}(y)
\leq
\Theta^{(1/2)}(x)-\Theta^{(1/2)}(y).
\]
Equality holds only for $\lambda=1/2$.
\begin{proof}
For $\lambda\neq\frac12$,
considering the ODEs describing $\Theta^{(\lambda)}$ and $\Theta^{(1/2)}$, it is clear that $\Theta^{(1/2)}$ can grow more from $y$ to $x$,
being even the initial value at $y$ higher.
\end{proof}
\end{lemma}

\begin{lemma}\plabel{lem:depp2}
%(Decreasing delay principle.)
Assume that $\lambda\in[0,1]$, $x\in[0,+\infty)$, and $0<r<\Theta^{(\lambda)}(x)<+\infty$.
Then
\[
x-\left(\Theta^{(\lambda)}\right)^{-1}\left(\Theta^{(\lambda)}(x)-r \right)
\geq
x-\left(\Theta^{(1/2)}\right)^{-1}\left(\Theta^{(1/2)}(x)-r \right)
\]
Equality holds only for $\lambda=1/2$.
\begin{proof}
Let $y$ be so that $\Theta^{(1/2)}(x)-\Theta^{(1/2)}(y)=r$.
By the previous lemma, for $\lambda\neq\frac12$, we know
$\Theta^{(\lambda)}(x)-\Theta^{(\lambda)}(y)<r$.
Therefore $\Theta^{(\lambda)}(x)-r<\Theta^{(\lambda)}(y)$, and thus
\[\left(\Theta^{(\lambda)}\right)^{-1}\left(\Theta^{(\lambda)}(x)-r  \right)<y
\equiv \left(\Theta^{(1/2)}\right)^{-1}\left(\Theta^{(1/2)}(x)-r  \right).\]
And this is equivalent to the statement to prove.
\end{proof}
\end{lemma}
\begin{theorem}(``Delay estimate reduction principle''.)
Assume that  $\lambda\in[0,1]$ and
 $0<r<\Theta^{(\lambda)}\left(\int|\phi|\right)<+\infty$, and
\[\Theta^{(\lambda)}\left(\phi\right)\leq \Theta^{(\lambda)}\left(\int|\phi|\right)-r.\]
Then
\[N^{(\lambda)}(\phi)\equiv
\left(\Theta^{(\lambda)} \right)^{-1} \left(\Theta^{(\lambda)}\left(\int|\phi|\right)\right)
\leq \left(\Theta^{(1/2)}\right)^{-1}\left(\Theta^{(1/2)}\left(\int|\phi|\right)-r  \right),\]
or equivalently,
\[\Delta^{(\lambda)}(\phi)\equiv
\left(\int|\phi|\right)-N^{(\lambda)}(\phi)\geq x-\left(\Theta^{(1/2)}\right)^{-1}\left(\Theta^{(1/2)}\left(\int|\phi|\right)-r  \right).\]
\begin{proof}
This is an immediate consequence of the previous lemma.
\end{proof}
\end{theorem}
This is, however, not a perfect tool.
In general delays are not very significant for $\lambda\sim0$ and $\lambda\sim1$ while convergence is otherwise good.

\snewpage
\begin{commenty}

\subsection{Holomorphic calculus through the resolvent}\plabel{ss:resCalc}
~\\

Here we use the negative spectral cut to discuss powers.
\begin{commentx}
\begin{example}\plabel{rem:MPsq}
The spectral square root function
\[\sqrt A=\frac1\pi\int_{\lambda=0}^1\frac{A}{\lambda+(1-\lambda)A}\frac{\mathrm d\lambda}{\sqrt{\lambda(1-\lambda)}}
=1+\frac1\pi\int_{\lambda=0}^1\frac{A-1}{\lambda+(1-\lambda)A}\sqrt{\frac{\lambda}{1-\lambda}}\,\mathrm d\lambda
\]
\[\left(=\int_{s=-\infty}^0\frac{1}{\pi\sqrt{-s}}\frac{A}{(A-s)}\,\mathrm ds
=1+\int_{s=-\infty}^0\frac{\sqrt{-s}}{\pi}\frac{(A-1)}{(1-s)(A-s)}\,\mathrm ds\right)\]
is quite similar to $\log$.
One can proceed with the calculus of $\sqrt{\exp_{\mathrm R}(\phi)}$ similarly.
\qedexer
\end{example}
\end{commentx}
One can check that for $\Rea\alpha\in(-1,1)$,
\begin{multline*}
A^\alpha
=1+ \frac{ \sin(\alpha\pi)}{\pi}\int_{\lambda=0}^1\frac{A-1}{\lambda+(1-\lambda)A}\left(\frac{\lambda}{1-\lambda}\right)^\alpha\,\mathrm d\lambda
\qquad \rechoicecomm{\\}
\left(
=1+\int_{s=-\infty}^0\frac{(-s)^\alpha\sin(\alpha\pi)}{\pi}\frac{(A-1)}{(1-s)(A-s)}\,\mathrm ds\right).
\end{multline*}
(The details are spelled out in Mart\'inez Carracedo, Sanz Alix \cite{MS} or Haase \cite{Haa}.
Together with the spectral conditions, this can be used as the definition for  $\Rea\alpha\in(-1,1)$,
then it can be extended by the prescription $A^{\alpha+n}=A^{\alpha}A^n$ with $n\in\mathbb Z$.
It is holomorphic in $\alpha$.)
As a formal consequence, this allows to obtain
\[\mu_{k}^{\langle\langle\alpha\rangle\rangle}(X_1,\ldots,X_k)
=\int_{\lambda=0}^1\frac{\sin(\alpha\pi)}\pi\left(\frac{\lambda}{1-\lambda}\right)^\alpha\mu_{k}^{(\lambda)}(X_1,\ldots,X_k)\,\mathrm d\lambda
\]
for $k\geq1$, $\Rea\alpha\in(-1,1)$.
In practice the integral is again just the beta function identity again; yielding an alternative proof
for Theorem \ref{th:powex} (for $\Rea\alpha\in(-1,1)$, but then the statement extends naturally).

Similarly, for $k\geq1$, standard spectral calculus (a  simple application of the Sokhotski–Plemelj formula) yields
\begin{multline}
(\log A)^k=\int_{\lambda=0}^1
 \frac{
\left(\log\left(\dfrac{\lambda}{1-\lambda}\right)+\mathrm i\pi\right)^k
-\left(\log\left(\dfrac{\lambda}{1-\lambda}\right)-\mathrm i\pi\right)^k
 }{2\mathrm i\pi}
\cdot \frac{A-1}{\lambda +(1-\lambda)A }\,d\lambda
\\\qquad\left(=\int^{0}_{s=-\infty}
\frac{
\left(\log(-s)+\mathrm i\pi\right)^k -\left(\log(-s)-\mathrm i\pi\right)^k
 }{2\mathrm i\pi}
\cdot\frac{A-1}{(1-s)(A-s)}\, \mathrm ds\right).
\plabel{eq:logpow}
\end{multline}

As a more compact notation, for $n\geq1$  and $\lambda\in(0,1)$, let
\begin{equation}
E_n(\lambda)= \frac{
\left(\log\left(\dfrac{\lambda}{1-\lambda}\right)+\mathrm i\pi\right)^n
-\left(\log\left(\dfrac{\lambda}{1-\lambda}\right)-\mathrm i\pi\right)^n
 }{2\mathrm i\pi}.
\plabel{eq:logpowco}
\end{equation}
Moreover, let
\[G_n(\lambda)=G_n(\lambda,\lambda-1).\]

$G_n(\lambda)$ is a polynomial of order $n-1$ of $\lambda$; and
$E_n(\lambda)$ is a polynomial of order $n-1$ of $\log\left(\frac{\lambda}{1-\lambda}\right)$.
The polynomials $G_n(\lambda)$ form the ``Goldberg basis''.
The functions $E_n(\lambda)$ form ``the dual Goldberg basis'':
\begin{lemma}\plabel{th:lowcan}
For $k\geq1$, $n\geq1$,
\[\int_{\lambda=0}^1 E_k(\lambda)G_n(\lambda)\,\mathrm d\lambda=\delta_{k,n}.\]
\begin{proof}
Applying \eqref{eq:logpow}, \eqref{eq:logpowco} to $A=\exp(X)$ formally, we obtain
\[X^k=\int_{\lambda=0}^1 E_k(\lambda) \mathcal R^{\lambda}(\exp X)\,\mathrm d\lambda. \]
Taking $\mathcal R^{\lambda}(\exp X)=\sum_{n=1}^\infty G_{n}(\lambda)X^n$ into consideration, this yields the statement.
\end{proof}
\end{lemma}
\begin{commentx}
\begin{remark}
\plabel{eq:dualgoldgen}
For $n\geq1$,
\[\frac{\mathrm d}{\mathrm d\lambda}E_n(\lambda)=\frac{n}{\lambda(1-\lambda)}E_{n-1}(\lambda)\]
holds (with $E_0(\lambda)=0$). This is a dual statement to Lemma \ref{lem:goldgen}.
\qedremark
\end{remark}
\end{commentx}
Applying \eqref{eq:logpow} formally, we find
\begin{equation}
\mu^{[k]}_n(X_1,\ldots,X_n)=\frac1{k!}\int_{\lambda=0}^1
E_k(\lambda)
\cdot \mu^{(\lambda)}_n(X_1,\ldots,X_n)\,d\lambda.
\plabel{eq:highcan}
\end{equation}
(But this is also obtainable from Lemma \ref{th:lowcan} and Theorem \ref{th:gemcan}.)
It is, however, not very practical algebraically, as the integration of the terms
$\lambda^n\left(\log\dfrac{\lambda}{1-\lambda}\right)^m$ is cumbersome.
\end{commenty}

\snewpage
\section{The discrete resolvent approach}
\plabel{sec:DiscRes}

%The discrete resolvent approach can be considered  the resolvent method applied to mBCH measures.
%However, the emphases are different.
Here we take the ``resolvent Goldbergian'' approach.
We try to make this discussion independent from the previous section, as much it is possible,  in order to avoid the measure terminology;
although we do not redevelop the functions $\mathrm M^{(\lambda)}$ and $\mathrm P^{(\lambda)}$ here.

\subsection{Exact and formal decompositions of the resolvent}
\begin{lemma}\plabel{lem:dres1}
In a real or complex algebra, if $\lambda$ is a scalar, then
\[R=\frac{A-1}{\lambda+(1-\lambda)A}\qquad\Leftrightarrow\qquad A=\frac{1+\lambda R}{1-(1-\lambda) R};\]
and in those cases
\[\frac{1}{\lambda+(1-\lambda)A}=1-(1-\lambda)R\qquad\text{and}\qquad \frac{1}{1-(1-\lambda)R}=\lambda+(1-\lambda)A.\]
\begin{proof}Elementary computation.\end{proof}
\end{lemma}
\begin{lemma}\plabel{lem:dres2}
(Resolvent decomposition, exact.)
Assume that $\mathcal R^{(\lambda)}(A)$ and $\mathcal R^{(\lambda)}(B)$ exist. Then

(i) $\mathcal R^{(\lambda)}(AB)$ exists

(ii) $\mathcal R^{(\lambda)}(BA)$ exists

(iii) $\left(1+\lambda(1-\lambda)\mathcal R^{(\lambda)}(A)\mathcal R^{(\lambda)}(B)   \right)^{-1}$ exists

(iv) $\left(1+\lambda(1-\lambda)\mathcal R^{(\lambda)}(B)\mathcal R^{(\lambda)}(A)   \right)^{-1}$ exists
\\ are equivalent to each other.

In those cases,
\begin{multline}
\left(1+\lambda(1-\lambda)\mathcal R^{(\lambda)}(A)\mathcal R^{(\lambda)}(B)\right)^{-1}=
(\lambda+(1-\lambda)B)\cdot\\\cdot\left(1-(1-\lambda)\mathcal R^{(\lambda)}(AB) \right)(\lambda+(1-\lambda)A),
\plabel{eq:enso}
\end{multline}
and
\begin{align}
\mathcal R^{(\lambda)}(AB)
=&\mathcal R^{(\lambda)}(A)(1-\lambda (\lambda-1)\mathcal R^{(\lambda)}(B)\mathcal R^{(\lambda)}(A))^{-1}\plabel{eq:bires}\\
&+\mathcal R^{(\lambda)}(B)(1-\lambda(\lambda-1) \mathcal R^{(\lambda)}(A)\mathcal R^{(\lambda)}(B))^{-1}\notag\\
&+\lambda \mathcal R^{(\lambda)}(A)\mathcal R^{(\lambda)}(B)(1-\lambda(\lambda-1) \mathcal R^{(\lambda)}(A)
\mathcal R^{(\lambda)}(B))^{-1}\notag\\
&+(\lambda-1) \mathcal R^{(\lambda)}(B)\mathcal R^{(\lambda)}(A)(1-\lambda(\lambda-1) \mathcal R^{(\lambda)}(B)
\mathcal R^{(\lambda)}(A))^{-1}.
\notag
\end{align}
\renewcommand{\qedsymbol}{}
\begin{proof}[Note]
Using the identities
$ \mathcal R^{(\lambda)}(AB) A=A\mathcal R^{(\lambda)}(BA) $
and
$ \mathcal R^{(\lambda)}(BA) B=B\mathcal R^{(\lambda)}(AB) $
one can write \eqref{eq:enso} in various forms.
Using the identities
\[\mathcal R^{(\lambda)}(A)(1-\lambda(\lambda-1) \mathcal R^{(\lambda)}(B)\mathcal R^{(\lambda)}(A))^{-1}
=(1-\lambda(\lambda-1) \mathcal R^{(\lambda)}(A)\mathcal R^{(\lambda)}(B))^{-1}\mathcal R^{(\lambda)}(A)\]
and
\[\mathcal R^{(\lambda)}(B)(1-\lambda(\lambda-1) \mathcal R^{(\lambda)}(A)\mathcal R^{(\lambda)}(B))^{-1}
=(1-\lambda(\lambda-1) \mathcal R^{(\lambda)}(B)\mathcal R^{(\lambda)}(A))^{-1}\mathcal R^{(\lambda)}(B),\]
one can write \eqref{eq:bires} somewhat flexibly.
\end{proof}
\renewcommand{\qedsymbol}{$\Box$}
\begin{proof} (i)$\Leftrightarrow$(iii) follows from the identity
\begin{multline}
1+\lambda(1-\lambda)\mathcal R^{(\lambda)}(A)\mathcal R^{(\lambda)}(B)=
(\underbrace{1-(1-\lambda)\mathcal R^{(\lambda)}(A)}_{(\lambda+(1-\lambda)A)^{-1}} ) \\
\Biggl( \lambda+ (1-\lambda)
\underbrace{\frac{1+\lambda\mathcal R^{(\lambda)}(A) }{1-(1-\lambda)\mathcal R^{(\lambda)}(A) }}_A
\underbrace{\frac{1+\lambda\mathcal R^{(\lambda)}(B) }{1-(1-\lambda)\mathcal R^{(\lambda)}(B) }}_B
\Biggr)
(\underbrace{1-(1-\lambda)\mathcal R^{(\lambda)}(B)}_{(\lambda+(1-\lambda)B)^{-1}} ).
\plabel{eq:downbires}
\end{multline}
(ii)$\Leftrightarrow$(iv) is similar.
(iii)$\Leftrightarrow$(iv) follows from the general identity
\[(1-HG)^{-1}=1-H(1-GH)^{-1}G.\]
Similarly, the identities in the remark are consequences of the general identity
\[(1-GH)^{-1}G=G(1-HG)^{-1}.\]

Equality \eqref{eq:enso} follows from \eqref{eq:downbires}.
It remains to prove \eqref{eq:bires}. Now,
\begin{multline}
\mathcal R^{(\lambda)}(A)+\mathcal R^{(\lambda)}(B)+(2\lambda-1)\mathcal R^{(\lambda)}(A)\mathcal R^{(\lambda)}(B)=
(\underbrace{1-(1-\lambda)\mathcal R^{(\lambda)}(A)}_{(\lambda+(1-\lambda)A)^{-1}} ) \\
\Biggl( -1+
\underbrace{\frac{1+\lambda\mathcal R^{(\lambda)}(A) }{1-(1-\lambda)\mathcal R^{(\lambda)}(A) }}_A
\underbrace{\frac{1+\lambda\mathcal R^{(\lambda)}(B) }{1-(1-\lambda)\mathcal R^{(\lambda)}(B) }}_B
\Biggr)
(\underbrace{1-(1-\lambda)\mathcal R^{(\lambda)}(B)}_{(\lambda+(1-\lambda)B)^{-1}} ).
\plabel{eq:upbires}
\end{multline}

Comparing \eqref{eq:upbires} and \eqref{eq:downbires}, we find
\begin{align}
(1-(1-\lambda)\mathcal R^{(\lambda)}(A) )\mathcal R^{(\lambda)}(AB)=\,&
(\mathcal R^{(\lambda)}(A)+\mathcal R^{(\lambda)}(B)+(2\lambda-1)\mathcal R^{(\lambda)}(A)\mathcal R^{(\lambda)}(B))\notag\\
&(1+\lambda(1-\lambda)\mathcal R^{(\lambda)}(A)\mathcal R^{(\lambda)}(B))^{-1}\notag\\
&(1-(1-\lambda)\mathcal R^{(\lambda)}(A) ).\notag
\end{align}

Expanding the to $3\times1\times2=6$ terms, it yields
\begin{align}
=\,&+\mathcal R^{(\lambda)}(A)(1-\lambda(\lambda-1) \mathcal R^{(\lambda)}(A)\mathcal R^{(\lambda)}(B))^{-1}\tag{gi}\\
&+\mathcal R^{(\lambda)}(B)(1-\lambda(\lambda-1) \mathcal R^{(\lambda)}(A)\mathcal R^{(\lambda)}(B))^{-1}\tag{b}\\
&+(2\lambda-1)\mathcal R^{(\lambda)}(A)\mathcal R^{(\lambda)}(B)(1-\lambda(\lambda-1) \mathcal R^{(\lambda)}(A)
\mathcal R^{(\lambda)}(B))^{-1}\tag{cf}\\
&-(1-\lambda)\mathcal R^{(\lambda)}(A)\mathcal R^{(\lambda)}(A)(1-\lambda(\lambda-1) \mathcal R^{(\lambda)}(B)
\mathcal R^{(\lambda)}(A))^{-1}\tag{e}\\
&-(1-\lambda)\mathcal R^{(\lambda)}(B)\mathcal R^{(\lambda)}(A)(1-\lambda(\lambda-1) \mathcal R^{(\lambda)}(B)
\mathcal R^{(\lambda)}(A))^{-1}\tag{d}\\
&-(1-\lambda)(2\lambda-1)\mathcal R^{(\lambda)}(A)\mathcal R^{(\lambda)}(B)\mathcal R^{(\lambda)}(A)(1-\lambda(\lambda-1)
\mathcal R^{(\lambda)}(B)\mathcal R^{(\lambda)}(A))^{-1}\tag{ahj}.
\end{align}

This can be decomposed further as
\begin{align}
=\,&+\mathcal R^{(\lambda)}(A)(1-\lambda(\lambda-1) \mathcal R^{(\lambda)}(B)\mathcal R^{(\lambda)}(A))^{-1}\tag{a}\\
&+\mathcal R^{(\lambda)}(B)(1-\lambda(\lambda-1) \mathcal R^{(\lambda)}(A)\mathcal R^{(\lambda)}(B))^{-1}\tag{b}\\
&+\lambda \mathcal R^{(\lambda)}(A)\mathcal R^{(\lambda)}(B)(1-\lambda(\lambda-1) \mathcal R^{(\lambda)}(A)
\mathcal R^{(\lambda)}(B))^{-1}\tag{c}\\
&+(\lambda-1) \mathcal R^{(\lambda)}(B)\mathcal R^{(\lambda)}(A)(1-\lambda(\lambda-1) \mathcal R^{(\lambda)}(B)
\mathcal R^{(\lambda)}(A))^{-1}\tag{d}\\
&-(1-\lambda)\mathcal R^{(\lambda)}(A)\mathcal R^{(\lambda)}(A)(1-\lambda(\lambda-1) \mathcal R^{(\lambda)}(B)
\mathcal R^{(\lambda)}(A))^{-1}\tag{e}\\
&-(1-\lambda)\mathcal R^{(\lambda)}(A)\mathcal R^{(\lambda)}(B)(1-\lambda(\lambda-1) \mathcal R^{(\lambda)}(A)
\mathcal R^{(\lambda)}(B))^{-1}\tag{f}\\
&-(1-\lambda)\lambda\mathcal R^{(\lambda)}(A) \mathcal R^{(\lambda)}(A)\mathcal R^{(\lambda)}(B)(1-\lambda(\lambda-1)
\mathcal R^{(\lambda)}(A)\mathcal R^{(\lambda)}(B))^{-1}\tag{g}\\
&-(1-\lambda)(\lambda-1) \mathcal R^{(\lambda)}(A)\mathcal R^{(\lambda)}(B)\mathcal R^{(\lambda)}(A)(1-\lambda(\lambda-1)
\mathcal R^{(\lambda)}(B)\mathcal R^{(\lambda)}(A))^{-1}\tag{h}\notag\\
&+A\tag{i}\notag\\
&-A\tag{j}\notag.
\end{align}
(For example, line (ahj) decomposes to lines (a),(h),(j).)
Let temporarily $\mathcal R^{(\lambda)}(A,B)$ denote the  RHS of \eqref{eq:bires}.
Then lines (i) and (j) cancel each other, while lines (a)--(h) yield the $2\times4=8$ terms of
$(1-(1-\lambda)\mathcal R^{(\lambda)}(A) )\mathcal R^{(\lambda)}(A,B)$.
Thus,
\[(1-(1-\lambda)\mathcal R^{(\lambda)}(A) )\mathcal R^{(\lambda)}(AB)=(1-(1-\lambda)\mathcal R^{(\lambda)}(A) )\mathcal R^{(\lambda)}(A,B),\]
which implies $\mathcal R^{(\lambda)}(AB)=\mathcal R^{(\lambda)}(A,B)$.
\end{proof}
\end{lemma}

\begin{lemma}\plabel{lem:powres}
(Resolvent decomposition, formal.)
If $X,Y$ are formal variables, then
\begin{align}
\mathcal R^{(\lambda)}((\exp X)(\exp Y))=\sum_{k=0}^\infty\biggl(
&\lambda^k(\lambda-1)^k\mathcal R^{(\lambda)}(\exp Y)(\mathcal R^{(\lambda)}(\exp X)\mathcal R^{(\lambda)}(\exp Y))^{k}\notag\\
+&\lambda^k(\lambda-1)^k\mathcal R^{(\lambda)}(\exp X)(\mathcal R^{(\lambda)}(\exp Y)\mathcal R^{(\lambda)}(\exp X))^{k}\notag\\
+&\lambda^{k+1}(\lambda-1)^{k}(\mathcal R^{(\lambda)}(\exp X)\mathcal R^{(\lambda)}(\exp Y))^{k+1} \notag\\
+&\lambda^{k}(\lambda-1)^{k+1}(\mathcal R^{(\lambda)}(\exp Y)\mathcal R^{(\lambda)}(\exp X))^{k+1} \qquad\,\,\biggr).\notag
 \end{align}
One can replace $\exp X$ and $\exp Y$ by other formal perturbations of $1$.
\begin{proof}
$\mathcal R^{(\lambda)}(\exp X)=X+\ldots$ is purely formal, so the previous Lemma can be applied.
$(1-\lambda(\lambda-1) \mathcal R^{(\lambda)}(\exp X)\mathcal R^{(\lambda)}(\exp Y))^{-1}$
expand as
\[\sum_{k=0}^\infty\lambda^k(\lambda-1)^k(\mathcal R^{(\lambda)}(\exp X)\mathcal R^{(\lambda)}(\exp Y))^k,\]
etc. The general case is similar, or follows by a change of variables.
(Remark: A direct formal proof is considerably simpler than the proof of Lemma \ref{lem:dres2}.)
\end{proof}
\end{lemma}
\begin{theorem}\plabel{th:powres}(Formal multivariable resolvent decomposition.)
Suppose that $X_1,\ldots,X_n$ are formal variables. Then
\begin{multline}
\mathcal R^{(\lambda)}((\exp X_1)\cdot\ldots\cdot(\exp X_n))=\\
\\=\sum_{\substack{i_1,\ldots,i_s\in\{1,\ldots,n\}\\s\geq1,\, i_j\neq i_{j+1}}}
\lambda^{\asc(i_1,\ldots,i_s)}(\lambda-1)^{\des(i_1,\ldots,i_s)}
\mathcal R^{(\lambda)}(\exp X_{i_1})\cdot\ldots\mathcal R^{(\lambda)}(\exp X_{i_s}).
\plabel{eq:thpowres}
\end{multline}
One can replace $\exp X$ and $\exp Y$ by  other formal perturbations of $1$.
\begin{proof}
Splitting up $(\exp X_1)\cdot\ldots\cdot(\exp X_n)=\left((\exp X_1)\cdot\ldots\cdot(\exp X_{n-1})\right)\cdot (\exp X_n)$,
let us apply Lemma \ref{lem:powres} inductively.
Again, the general case is similar, or follows by a change of variables.
\end{proof}
\end{theorem}

\begin{remark}\plabel{rem:powres}
We can easily prove Theorem \ref{th:powres}  based on the ideas of Section \ref{sec:MagnusResolvent}:
In a formal setting, we can apply Theorem \ref{th:resdecomp} to $X_1\mathbf 1_{[0,1)}\bo.\ldots\bo.X_n\mathbf 1_{[0,1)}$.
(We can do this in the $\mathrm F^1_{\mathbb R}$ topology,
 or, alternatively, we can give $X_i$ arbitrarily small norm; or we can use a truly formal, degree-wise, setup.)
Conversely, Theorem \ref{th:powres} implies  Theorem \ref{th:powexr}.
\qedremark
\end{remark}

\snewpage

\subsection{Some algebraic consequences for the BCH expansion}\plabel{ss:algco}
~\\

As an illustration for the use of the formal resolvent expansion,
we give a somewhat peculiar proof for
\begin{theorem}[F. Schur \cite{Sch1} (1890), \cite{Sch2},  Poincar\'e \cite{PH} (1899) ] \plabel{thm:schurres}
Let us consider the (formal) power series
\[\beta(x)=\frac x{\mathrm e^x-1}=\sum_{n=0}^\infty \beta_nx^n.\]

Then, using commutator notation,
\begin{equation}
\BCH(X,Y)_{\deg_{X,Y}=(n,1)}=(-1)^n\beta_n \underbrace{[X,[X,\ldots [X,}_{n \text{ times}}Y]\ldots]];
%\tag{\theequation R}
\plabel{0R}
\end{equation}
and
\begin{equation}
\BCH(X,Y)_{\deg_{X,Y}=(1,n)}=(-1)^n\beta_n [[\ldots[X\underbrace{,Y]\ldots ,Y],Y]}_{n \text{ times}}.
%\tag{\theequation L}
\plabel{0L}
\end{equation}
 (Here $(-1)^n\beta_n$ makes difference to $\beta_n$ only as $(-1)^1\beta_1=\frac12$.)

\begin{proof}
I.~e., we have to prove that the coefficient of $X^{n_1}YX^{n-n_1}$ in $\BCH(X,Y)$ is
\[(-1)^n\beta_n \cdot (-1)^{n-n_1}\binom{n}{n_1};\]
and the coefficient of $Y^{n_1}XY^{n-n_1}$ in $\BCH(X,Y)$ is
\[(-1)^n\beta_n \cdot (-1)^{n_1}\binom{n}{n_1}.\]

We will prove only the first statement.
Let us consider $\BCH(X,Y)$; replace the terms $X^{n_1}YX^{n-n_1}$ by $x_1^{n_1}x_2^{n-n_1}$, and ignore the other terms.
Then, according Lemma \ref{lem:powres} and the logarithm formula \eqref{eq:logdef}, the generating function to obtain is
\[\int_{\lambda=0}^1\left( 1
+\lambda(\lambda-1)\mathcal R^{(\lambda)}(\exp x_1)
 \mathcal R^{(\lambda)}(\exp x_2)
 +\lambda\mathcal R^{(\lambda)}(\exp x_1)
 +(\lambda-1)\mathcal R^{(\lambda)}(\exp x_2)   \right)   \,\mathrm d\lambda\]
\[=\int_{\lambda=0}^1  \frac{1}{\lambda \mathrm e^{-x_1}+(1-\lambda)} \frac{1}{\lambda+(1-\lambda)\mathrm e^{x_2}} \,\mathrm d\lambda
=\left[\frac{1}{\mathrm e^{x_2-x_1} -1} \log\left( \frac{\lambda+(1-\lambda)\mathrm e^{x_1}}{ \lambda+(1-\lambda)\mathrm e^{x_2}}
\right)\right]_{\lambda=0}^1\]
\[=\frac{x_2-x_1}{\mathrm e^{x_2-x_1}-1}=\beta(x_2-x_1)=\sum_{n=0}^\infty(-1)^n\beta_n\cdot (x_1-x_2)^n.\]
The last expression for the generating function expresses the statement to prove.

(In the formula above, the transition from term 2 to term 4 may look suspicious, as it is not necessarily clear that term 3 is well-defined.
It is, but if one does not want to dissect that formal power series, then one can argue as follows:
If the sequence of terms 2, 3, 4 gets multiplied by $\mathrm e^{x_2-x_1} -1$, then the whole transition becomes uncontroversial,
and the equality of terms 2 and 4 follows the from zero-divisor-freeness of the formal power series.)
\end{proof}
\renewcommand{\qedsymbol}{$\triangle$}

\begin{proof}[Note]
Schur \cite{Sch1} develops the equivalents of \eqref{0R} and \eqref{0L} in the setting of Lie group theory.
Later, Campbell, Pascal,  Baker and Hausdorff also have it and use it in their investigations
(cf.~Achilles, Bonfiglioli \cite{AB}) but
Poincar\'e \cite{PH} is notable because he is first one the who discusses these statements in the formal / analytic settings
 much considered nowadays.
For the sake of simplicity, we call all these results as `Schur's formulae'.
\end{proof}
\begin{proof}[Note]
Goldberg \cite{G} contains more on generating functions.
Mielnik, Pleba\'nski \cite{MP} have (several) combinatorial arguments of greater generality.
These complement Reutenauer \cite{R} in the direction of analytical combinatorics.
\end{proof}
\renewcommand{\qedsymbol}{$\Box$}
\end{theorem}
%\snewpage

\snewpage

\subsection{Analytic consequences for the BCH expansion}\plabel{ss:anaco}
~\\

In this subsection we apply the formal resolvent decompositions
in order to obtain convergence estimates with respect to (m)BCH expanions.

Let us recall the formal power series
$\mathrm M^{(\lambda)}(X_1,\ldots,X_n)$ and $\mathrm P^{(\lambda)}(X_1,\ldots,X_n)$ from
point \ref{po:deco}.
With this terminology, \eqref{eq:thpowres} reads as
\[\mathcal R^{(\lambda)}((\exp X_1)\cdot\ldots\cdot(\exp X_n))
=\mathrm M^{(\lambda)}(\mathcal R^{(\lambda)}(\exp X_1),\ldots,\mathcal R^{(\lambda)}(\exp X_n)).\]

Let $\Gamma^{(\lambda)}(X_1,\ldots,X_n)$ denote $\mathcal R^{(\lambda)}((\exp X_1)\cdot\ldots\cdot(\exp X_n))$
but where every coefficient of the monomials $X_{j_1}\cdot\ldots\cdot X_{j_s}$ is replaced by its absolute value.
Then, we find
\begin{equation}\Gamma^{(\lambda)}(X_1,\ldots,X_n)=\mathrm P^{(\lambda)}(\Gamma^{(\lambda)}(X_1),\ldots,\Gamma^{(\lambda)}(X_n)).
\plabel{eq:dresest}\end{equation}
(Here we have used the fact that there are no consecutive multiplicities in $\mathrm M^{(\lambda)}$.)

The connection to the Baker--Campbell--Hausdorff expansion is as follows.
Formally (just taking the logarithm),
\begin{equation}
\BCH(X_1,\ldots,X_n)=\int_{\lambda=0}^1 \mathcal R^{(\lambda)}((\exp X_1)\cdot\ldots\cdot(\exp X_n))\,\mathrm d\lambda
.
\plabel{eq:logBCH}
\end{equation}
Thus, taking absolute values in the power series,
\[\Gamma(X_1,\ldots,X_n)\stackrel{\forall X_1,\ldots,X_n}\leq \int_{\lambda=0}^1 \Gamma^{(\lambda)}(X_1,\ldots,X_n) \,\mathrm d\lambda.\]
In particular, applying this to nonnegative variables,
\begin{equation}
\Gamma_\real(x_1,\ldots,x_n)\leq \int_{\lambda=0}^1 \Gamma^{(\lambda)}_\real(x_1,\ldots,x_n) \,\mathrm d\lambda.
\plabel{eq:gameq}
\end{equation}

In fact, if $\Gamma^{(\lambda)}_\real(x_1,\ldots,x_n)$ is bounded on $\lambda\in[0,1]$,
then the BCH expansion not only converges, but the resolvents exists, and for $|X_i|\leq x_i$
the formula \eqref{eq:logBCH} actually yields
\[\BCH(X_1,\ldots,X_n)=\log((\exp X_1)\ldots(\exp X_n))\]
in strict analytical sense.
\snewpage
\begin{remark}
\plabel{rem:subconv}
Even at this point, we can improve on the convergence range of the BCH formula easily:

For $\lambda\in[0,1]$, let $\rho^{(\lambda)}(x)$ be the power series expansion of $\mathcal R^{(\lambda)}(\exp x)$ around $x=0$.
 $\mathcal R^{(\lambda)}(\exp x)$ satisfies the recursion / formal IVP
\[\rho^{(\lambda)\prime}(x)=1-(1-2\lambda)\rho^{(\lambda)}(x)-\lambda(1-\lambda)\rho^{(\lambda)}(x)^2,\]
\[\rho^{(\lambda)}(0)=0.\]
Thus, thinking combinatorially, we can estimate the absolute value of the coefficients in $\rho^{(\lambda)}(x)$
 by the coefficients of the solution of the  recursion / formal IVP
\[g^{(\lambda)\prime}(x)=1+|1-2\lambda|g^{(\lambda)}(x)+\lambda(1-\lambda)g^{(\lambda)}(x)^2,\]
\[g^{(\lambda)}(0)=0.\]
The actual solution (until it blows up) is given by
\begin{commentx}
\newcommand{\tanexpr}{\dfrac{\tan\left(\dfrac{x}2\sqrt{8\lambda(1-\lambda)-1}\right)}{\sqrt{8\lambda(1-\lambda)-1}}}
\newcommand{\tanhexpr}{\dfrac{\tanh\left(\dfrac{x}2\sqrt{1-8\lambda(1-\lambda)}\right)}{\sqrt{1-8\lambda(1-\lambda)}}}
\[g^{(\lambda)}(x)=\begin{cases}
\dfrac{2\tanexpr}{1-|1-2\lambda|\tanexpr}
&
\begin{matrix}
\text{if }1-8\lambda(1-\lambda)<0\\
\qquad(\text{i.~e.~}|\lambda-\frac12|<2^{-3/2}),
\end{matrix}
\\\\
\dfrac{x}{1-\dfrac{\sqrt2}4x}&
\begin{matrix}
\text{if }1-8\lambda(1-\lambda)=0\\
\qquad(\text{i.~e.~}|\lambda-\frac12|=2^{-3/2}),
\end{matrix}
\\\\
\dfrac{2\tanhexpr}{1-|1-2\lambda|\tanhexpr}
&
\begin{matrix}
\text{if }1-8\lambda(1-\lambda)>0\\
\qquad(\text{i.~e.~}|\lambda-\frac12|>2^{-3/2}).
\end{matrix}
\end{cases}\]
Note: Using a more compact notation,
\end{commentx}
\[g^{(\lambda)}(x)=\frac{x\Tan\left(\dfrac{x^2}4({8\lambda(1-\lambda)-1})\right)}{1-|1-2\lambda|
\frac12x\Tan\left(\dfrac{x^2}4({8\lambda(1-\lambda)-1})\right) },\]
where
\begin{equation}
\Tan(z)=
\begin{cases}
\dfrac{\tan{\sqrt{z}}}{\sqrt{z}}&\text{if }z>0,\\
1&\text{if }z=0,\\
\dfrac{\tanh{\sqrt{-z}}}{\sqrt{-z}}&\text{if }z>0.
\end{cases}
\plabel{eq:defTan}
\end{equation}
\begin{commentx}
(The domain conditions correspond to $|\lambda-\frac12|<2^{-3/2}$ -- elliptic,
  $|\lambda-\frac12|=2^{-3/2}$ -- parabolic, $|\lambda-\frac12|>2^{-3/2}$
-- hyperbolic, respectively.)
\end{commentx}

Thus, according to the previous discussion,
\[\Gamma^{(\lambda)}(x) \stackrel{\forall x}\leq g^{(\lambda)}(x).\]
(This is exact for $\lambda=0,\frac12,1$, cf. Lemma \ref{lem:Gcomputer}.)

Consequently,
\[\Gamma^{(\lambda)}(x,y) \stackrel{\forall x,y}\leq\mathrm P^{(\lambda)}(g^{(\lambda)}(x),g^{(\lambda)}(y))
=\frac{g^{(\lambda)}(x)+g^{(\lambda)}(y) +g^{(\lambda)}(x)g^{(\lambda)}(y) }{1-\lambda(1-\lambda)g^{(\lambda)}(x)g^{(\lambda)}(y) }.\]

Let
\begin{commentx}
\[h(\lambda)=\begin{cases}
\dfrac{2\arctan\left(\dfrac{\sqrt{8\lambda(1-\lambda)-1}}{|1-2\lambda|}\right)}{\sqrt{8\lambda(1-\lambda)-1}}
&
\begin{matrix}
\text{if }1-8\lambda(1-\lambda)<0\\
\qquad(\text{i.~e.~}|\lambda-\frac12|<2^{-3/2}),
\end{matrix}
\\\\
2\sqrt2&
\begin{matrix}
\text{if }1-8\lambda(1-\lambda)=0\\
\qquad(\text{i.~e.~}|\lambda-\frac12|=2^{-3/2}),
\end{matrix}
\\\\
\dfrac{2\artanh\left(\dfrac{\sqrt{1-8\lambda(1-\lambda)}}{|1-2\lambda|}\right)}{\sqrt{1-8\lambda(1-\lambda)}}
&
\begin{matrix}
\text{if }1-8\lambda(1-\lambda)>0\\
\qquad(\text{i.~e.~}|\lambda-\frac12|>2^{-3/2}).
\end{matrix}
\end{cases}\]
Note: Using a more compact notation,
\end{commentx}
\[h(\lambda)=\begin{cases}
\dfrac2{|1-2\lambda|}\ATT\left(\dfrac{{8\lambda(1-\lambda)-1}}{|1-2\lambda|^2}\right)&\text{if }\lambda\neq0,\\
\pi&\text{if }\lambda=1/2,
\end{cases}\]
where
\begin{equation}
\ATT(x)=\begin{cases}
\dfrac{\arctan\sqrt z}{\sqrt z}&\text{if }z>0,\\
1&\text{if }z=0,\\
\dfrac{\artanh\sqrt{- z}}{\sqrt{ -z}}&\text{if }-1<z<0
\end{cases}
\qquad
\begin{pmatrix} \text{remark:}\\
\text{analitically} \\ \text{extendible to} \\ z\in\mathbb C\setminus(-\infty,1]
\end{pmatrix}
.
\plabel{eq:ATdef}
\end{equation}
%\begin{commentx}
Or, even more compactly,
\[h(\lambda)=\frac1{\sqrt{\lambda(1-\lambda)}}\AC\left(\frac{|2\lambda-1|}{2\sqrt{\lambda(1-\lambda)}}\right),\]
where
\begin{equation}
\AC(z)=\begin{cases}
\dfrac{\arccos z}{\sqrt{1-z^2}}&\text{if }-1< z<1,\\[3mm]
1&\text{if }z=1,\\[1mm]
\dfrac{\arcosh z}{\sqrt{z^2-1}}\qquad&\text{if } 1<z\\
\end{cases}
\qquad
\begin{pmatrix} \text{remark:}\\
\text{analitically} \\ \text{extendible to} \\ z\in\mathbb C\setminus(-\infty,1]
\end{pmatrix}
.
\plabel{eq:ACdef}
\end{equation}
%\end{commentx}

One can check that if $x,y\geq0$ and $x+y=h(\lambda)$, then
\[ \lambda(1-\lambda)g^{(\lambda)}\left(x\right)g^{(\lambda)}\left(y\right)=1.\]

Let
\[\widetilde{\mathrm C}_2=\min_{x\in[0,1]} h(\lambda)=\min_{x\in[0,1/2]} h(\lambda).\]
Now,
\[\widetilde{\mathrm C}_2=2.701428513\ldots,\]
and
\[\widetilde{\mathrm C}_2=h(\lambda_0),\]
where $\lambda_0=0.258744139\ldots$. Here $\lambda_0$ is the unique solution
\[h(\lambda)\cdot 4\lambda(1-\lambda)|1-2\lambda|=1\]
for $\lambda\in[0,1/2]$.
(As such, it can be computed with Newton iterations very efficiently.)

Consequently, if $x,y\geq0$, $x+y<\widetilde{\mathrm C}_2$, then
\[\sup_{\lambda\in[0,1]} \lambda(1-\lambda)g^{(\lambda)}\left(x\right)g^{(\lambda)}\left(y\right)<1.\]

Thus,
\[\Gamma^{(\lambda)}_\real(x,y)\leq
 \mathrm P_\real^{(\lambda)}(g^{(\lambda)}(x),g^{(\lambda)}(y))
 =\frac{g^{(\lambda)}(x)+g^{(\lambda)}(y) +g^{(\lambda)}(x)g^{(\lambda)}(y)
 }{1-\lambda(1-\lambda)g^{(\lambda)}(x)g^{(\lambda)}(y) }.\]

Hence, if $X,Y\in \mathfrak A$, and
\[|X|+|Y|<\widetilde{\mathrm C}_2,\]
then the (logarithmic) BCH formula holds.

This $\widetilde{\mathrm C}_2=2.701428513\ldots$ is already a  significant improvement compared to $2$.
\qedremark
\end{remark}

\snewpage
\begin{lemma}
\plabel{lem:Rover}
For $\lambda\in[0,1]$, the convergence radius of $\mathcal R^{(\lambda)}(\exp X)$ around $X=0$ is given by
\[{\mathrm C}_\infty^{(\lambda),\boldsymbol\varepsilon}=\text{``$\left|\log \frac\lambda{\lambda-1}\right|$''}=
\begin{cases}
 \sqrt{\pi^2+\left(\log\frac\lambda{1-\lambda}\right)^2}&\text{if}\quad\lambda\in(0,1),\\
 \\
+\infty&\text{if}\quad\lambda\in\{0,1\}.
\end{cases}\]
For $\lambda\in(0,1)$, taken in complex sense, $X\mapsto\mathcal R^{(\lambda)}(\exp X)$ has just two simple poles on the boundary
of the convergence disk.

${\mathrm C}_\infty^{(\lambda),\boldsymbol\varepsilon}$ is a strictly convex, nonnegative function in $\lambda\in(0,1)$,
symmetric for $\lambda\mapsto1-\lambda$;  its minimum is ${\mathrm C}_\infty^{(1/2)}=\pi$.
In particular, in $\lambda\in[0,1]$, it yields a $[\pi,+\infty]$-valued strictly convex continuous function.
For $\lambda\in[0,1]$, let
\[w^{(\lambda),\boldsymbol\varepsilon}=1/{\mathrm C}_\infty^{(\lambda),\boldsymbol\varepsilon}.\]
In $\lambda\in[0,1]$, $w^{(\lambda),\boldsymbol\varepsilon}$ is a $[0,1/\pi]$-valued strictly convex continuous function,
symmetric for $\lambda\mapsto1-\lambda$;
its maximum is $w^{(1/2),\boldsymbol\varepsilon}=1/\pi$.
\begin{proof}
The case $\lambda\in\{0,1\}$ is trivial. For $0<\lambda<1$, formally,
\begin{align}
\mathcal R^{(\lambda)}(\exp X)&=\frac{(\exp X)-1}{\lambda+(1-\lambda)(\exp X)}\notag\\
&=\frac{1}{\sqrt{\lambda(1-\lambda)}}\frac{\sinh\frac X2}{\cosh \frac12\left( X -\log\frac\lambda{1-\lambda} \right)}\notag\\
&=\frac{2\lambda-1}{2\lambda(1-\lambda)}+\frac{1}{2\lambda(1-\lambda)}\tanh \frac12\left( X -\log\frac\lambda{1-\lambda} \right).\notag
\end{align}
Analytically, this function has simple poles at $X=\log\frac\lambda{1-\lambda} +(2k+1)\pi\mathrm i$, $k\in\mathbb Z $
(with residue $\frac1{\lambda(1-\lambda)}$).
The rest is elementary calculus.
\end{proof}
\end{lemma}
\begin{lemma}\plabel{lem:Gcomputer}
For $\lambda=0,\frac12,1$, we have
\begin{align}
\mathcal R^{(0)}(\exp X)&=1-(\exp -X);&\Gamma^{(0)}(X)&=\exp X-1;\notag\\
\mathcal R^{(1/2)}\left(\exp X\right)&=2\tanh\frac X2;&\Gamma^{(1/2)}(X)&=2\tan\frac X2;\notag\\
\mathcal R^{(1)}(\exp X)&=(\exp X)-1;&\Gamma^{(1)}(X)&=\exp X-1.\notag
\end{align}
\begin{proof}
It follows from the power series expansions.
\end{proof}
\end{lemma}

Firstly, we draw conclusions regarding the Cayley transform case, i.~e.~$\lambda=\frac12$, where
$-\frac12\mathcal R^{(1/2)}(A)$ is the real involutive Cayley transform $\frac{1-A}{1+A}$.

\snewpage
\begin{theorem}\plabel{thm:cayleyest}
If $x_1,\ldots,x_n\in[0,\infty)$ ($n\geq2$) and
\[x_1+\ldots+x_n\leq 2n\arctan \frac1{n-1},\]
then
\[ \Gamma^{(1/2)}_\real(x_1,\ldots,x_n)\leq \frac {2n\tan\frac{x_1+\ldots+x_n}{2n}}{1-(n-1)\tan\frac{x_1+\ldots+x_n}{2n}}<+\infty.\]

However
\[  \Gamma^{(1/2)}_\real\Bigl( \underbrace{2\arctan \frac1{n-1} ,\ldots,2\arctan \frac1{n-1}}_{\text{$n$ terms}}\Bigr)=+\infty. \]

(For $n=1$,   $2\arctan \frac1{n-1}=\pi$ could be taken.)
\begin{proof}
Lemma \ref{lem:Mcomputer} and Lemma \ref{lem:Gcomputer} allow us to compute $\Gamma^{(1/2)}(X_1,\ldots,X_n)$ formally.
In particular, we can compute $\Gamma^{(1/2)}_\real(X_1,\ldots,X_n)$ explicitly
if $n=2$ or $x_1=\ldots=x_n$.

Simple analysis shows that the if $c\in[0,\pi]$, then the function
$x\mapsto \Gamma^{(1/2)}(x,c-x)$ is concave on the interval $x\in[0,c]$.
Using the associativity of $\mathrm P^{(\lambda)}$, this allows us to make the estimate
\[\Gamma^{(1/2)}_\real(x_1,\ldots,x_n)\leq \Gamma^{(1/2)}_\real\left(\frac{x_1+\ldots+x_n}n,\ldots,\frac{x_1+\ldots+x_n}n\right);\]
and the RHS  is known according to the previous discussion.
\end{proof}
\end{theorem}
\begin{theorem}\plabel{th:cayleyexp}
If $|X_1|+\ldots+|X_n|<2n\arctan\frac1{n-1}$, then \eqref{eq:thpowres} specified to $\lambda=\frac12$, i. e.
\begin{multline}
 \mathcal R^{(1/2)}\left((\exp X_1)\cdot\ldots\cdot(\exp X_n)\right)=\\
\\=\sum_{\substack{i_1,\ldots,i_s\in\{1,\ldots,n\}\\s\geq1,\, i_j\neq i_{j+1}}}
(-1)^{\des(i_1,\ldots,i_s)}\frac1{2^{s-1}}
\mathcal R^{(1/2)}\left(\exp X_{i_1}\right)\cdot\ldots\mathcal R^{(1/2)}\left(\exp X_{i_s}\right),
\notag
\end{multline}
is well-defined and absolute convergent.
This is, however, not necessarily true if
$ |X_1|=\ldots=|X_n|=2\arctan\frac1{n-1} $, then the LHS might not even exist.
\begin{proof}
The convergence statement follows from the previous theorem imediately.
For the non-convergence statement, let us consider the
Banach algebra of noncommutative power series of $Y_1,\ldots,Y_n$ with $\ell^1$ norm of the coefficients,
and apply $X_n=2\arctan\frac1{n-1}Y_n$.
The power series expansion of the LHS is formally determined, yet we know that it would have norm $+\infty$.
\end{proof}
\end{theorem}
\begin{cor}\plabel{cor:biexp}
If $|X|+|Y|<\pi$, then $(\exp X)+(\exp Y)$ is invertible.
\begin{proof}
The previous statement yields
the existence of
$\frac12-\frac14\mathcal R^{(1/2)}((\exp X)(\exp Y))=
\frac12-\frac12 \frac{(\exp X)(\exp Y)-1}{(\exp X)(\exp Y)+1}=
\frac1{(\exp X)(\exp Y)+1}$,
that is $((\exp X)(\exp Y)+1)^{-1}$.
Substituting $Y=-Y$, and multiplying by $(\exp Y)^{-1}$ on the left yields $((\exp X)+(\exp Y))^{-1}$.
\end{proof}
\end{cor}
The convergence of the BCH formula requires more work:

\begin{lemma}\plabel{lem:DRESgamma}
(a) The function
\[\widetilde\Gamma_\real^{[1]}:[0,1]\times[0,+\infty]\rightarrow[0,+\infty]\]
\[(\lambda,x)\mapsto \Gamma^{(\lambda)}_\real(x)\]
is continuous.
$\widetilde\Gamma_\real^{[1]}(\lambda,x)$ is monotone increasing in the variable $x$, and it is strictly
increasing in the variable $x$ on $\left( \Gamma_\real^{[1]}\right)^{-1}([0,+\infty))$.

(b) $\widetilde\Gamma_\real^{[1]}(\lambda,x)<+\infty$ if and only
$x<\sqrt{\pi^2+\left(\log\frac\lambda{1-\lambda}\right)^2}=$``$\left|\log \frac\lambda{\lambda-1}\right|$''
(the latter one is defined as $+\infty$ for $\lambda=0,1$).
Some special values are given by
\[\widetilde\Gamma_\real^{[1]}(\lambda,0)=0;\]
and
\begin{align}
\widetilde\Gamma_\real^{[1]}(0,x)&=(\exp x)-1&&\text{if\quad} x<+\infty;\notag\\
\widetilde\Gamma_\real^{[1]}\left(\frac12,x\right)&=2\tan\frac x2&&\text{if\quad} x<\pi;\notag\\
\widetilde\Gamma_\real^{[1]}(1,x)&=(\exp x)-1&&\text{if\quad} x<+\infty.\notag
\end{align}
\begin{proof}
(b)
Then the finiteness statement follows from Lemma \ref{lem:Rover} and simple principles of complex analysis.
The rest is straightforward, cf. Lemma \ref{lem:Gcomputer}.

(a) The Cauchy formula can be applied to give locally uniform (in $\lambda$) estimates
for the power series coefficients. This leads to continuity on the open set $\left( \Gamma_\real^{[1]}\right)^{-1}([0,+\infty))$.
The monotonicity statements are immediate from the nonnegative power series expansion and
that $\Gamma^{(\lambda)}(x)=x+\ldots$.
Separate continuity (i.~e.~continuity for fixed $\lambda$) is straightforward from the pole structure above.
Separate continuity, continuity on the finite part and monotonicity implies global continuity.
\end{proof}
\end{lemma}

\begin{lemma}\plabel{lem:DRESM}
(a) The function
\[\mathrm P^{[2]}_\real:[0,1]\times[0,+\infty]\times[0,+\infty]\rightarrow[0,+\infty] \]
\[(\lambda,x,y)\mapsto\mathrm P^{(\lambda)}_\real(x,y)=
\begin{cases}
\dfrac{x+y+xy}{1-\lambda(1-\lambda)xy}&\text{if $x,y<+\infty,\,\lambda(1-\lambda)xy<1$},\\
+\infty&\text{otherwise}
\end{cases}\]
is continuous.
$\mathrm P^{[2]}_\real(\lambda,x,y)$ is monotone increasing in the variables $x,y$;
and $\mathrm P^{[2]}_\real(\lambda,x,y)$ is strictly increasing in the variables $x,y$ on $(\mathrm P^{[2]}_\real)^{-1}([0,+\infty)$.

(b) Special values are given by
\[\mathrm P^{[2]}_\real(\lambda,x,0)=x,\qquad\mathrm P^{[2]}_\real(\lambda,0,y)=y,\]
\[\mathrm P^{[2]}_\real(\lambda,+\infty,y)=\mathrm P^{[2]}_\real(\lambda,x,+\infty)=+\infty.\]
\begin{proof}
This is just elementary analysis. However, we remark the following.
Consider the function $\iota:[0,1]\rightarrow[0,+\infty]$ such that
\[\iota( x)=\begin{cases} \dfrac{x}{1-x}&\text{if }0\leq x<1, \\ +\infty&\text{if } x=1. \end{cases}\]
Then reparametrizing the second and third variables and the range in $\mathrm P^{[2]}_\real$, it yields
\[\iota^*\mathrm P^{[2]}_\real(\lambda,x,y)=\min\left(1,\frac{x+y-xy}{1-\lambda(1-\lambda)xy}\right), \]
which is particularly transparent regarding continuity.
\end{proof}
\end{lemma}

\begin{lemma}\plabel{the:DRESgamma}
The function
\[\widetilde\Gamma^{[n]}_\real:[0,1]\times\underbrace{[0,+\infty]\times\ldots\times[0,+\infty]}_{\text{$n$ terms}}  \rightarrow[0,+\infty]\]
\[(\lambda,x_1,\ldots,x_n)\mapsto \Gamma^{(\lambda)}_\real(x_1,\ldots,x_n)\]
is continuous.
$\widetilde\Gamma^{[n]}_\real(\lambda,x_1,\ldots,x_n)$ is monotone increasing in the variables $x_1,\ldots,x_n$;
and $\widetilde\Gamma^{[n]}_\real(\lambda,x_1,\ldots,x_n)$ is strictly increasing in the variables $x_1,\ldots,x_n$
on $(\widetilde\Gamma^{[n]}_\real)^{-1}([0,+\infty))$.
\begin{proof}
This follows from combining Lemma \ref{lem:DRESgamma}.a and Lemma \ref{lem:DRESM}.a.
\end{proof}
\end{lemma}
For $x_1,\ldots,x_n\in[0,+\infty]$, let
\[\overline\Gamma_\real(x_1,\ldots,x_n) =\sup_{\lambda \in[0,1]}\,\Gamma^{(\lambda)}_\real(x_1,\ldots,x_n).\]

(In contrast to previous definitions, this is not derived from a power series but of real function theoretic nature.)

\begin{cor}\plabel{cor:DRESsup}
The function
\[\overline\Gamma^{[n]}_\real:\underbrace{[0,+\infty]\times\ldots\times[0,+\infty]}_{\text{$n$ terms}}  \rightarrow[0,+\infty]\]
\[(x_1,\ldots,x_n)\mapsto \overline\Gamma_\real(x_1,\ldots,x_n) =\sup_{\lambda \in[0,1]}\,\Gamma^{(\lambda)}_\real(x_1,\ldots,x_n).\]
is continuous.
$\overline\Gamma^{[n]}_\real(x_1,\ldots,x_n)$ is monotone increasing in the variables $x_1,\ldots,x_n$;
and $\overline\Gamma^{[n]}_\real(x_1,\ldots,x_n)$ is strictly increasing in the variables $x_1,\ldots,x_n$
on  $(\overline\Gamma^{[n]}_\real)^{-1}([0,+\infty))$.
\begin{proof}
This follows from general properties of continuous functions on compact sets.
\end{proof}
\end{cor}
\snewpage

\begin{lemma}\plabel{lem:DRESest}
Suppose that $x_1,\ldots,x_n\in[0,+\infty)$, $x_1+\ldots+ x_{k-1}>0$,  $x_k+\ldots +x_{n}>0$.
Then the following conditions are equivalent:

(i) $\overline\Gamma_\real(x_1,\ldots,x_n)=+\infty$, but for $0\leq y_j\leq x_j$,
\[y_1+\ldots+y_n<x_1+\ldots+x_n\quad \Rightarrow\quad \overline\Gamma_\real(y_1,\ldots,y_n)<+\infty;\]

(ii) $\overline\Gamma_\real(x_1,\ldots,x_n)=+\infty$, but
\[0\leq t<1\quad \Rightarrow\quad \overline\Gamma_\real(tx_1,\ldots,tx_n)<+\infty;\]

(iii)
\[\sup_{\lambda\in[0,1]} \lambda(1-\lambda)\Gamma^{(\lambda)}_\real(x_1,\ldots,x_{k-1})\Gamma^{(\lambda)}_\real(x_k,\ldots,x_{n})=1. \]
\begin{proof}
(i)$\Rightarrow$(ii) is obvious, and (ii)$\Rightarrow$(iii)$\Rightarrow$(i)  follows from the formal identity
\begin{multline}
\Gamma^{(\lambda)}(x_1,\ldots,x_{n})=\\=
\dfrac{\Gamma^{(\lambda)}(x_1,\ldots,x_{k-1})+\Gamma^{(\lambda)}(x_k,\ldots,x_{n})+\Gamma^{(\lambda)}(x_1,\ldots,x_{k-1})
\Gamma^{(\lambda)}(x_k,\ldots,x_{n})}{1-\lambda(1-\lambda)\Gamma^{(\lambda)}(x_1,\ldots,x_{k-1})\Gamma^{(\lambda)}(x_k,\ldots,x_{n})}
\notag
\end{multline}
and the monotonicity properties.
\end{proof}
\end{lemma}

An immediate consequence of \eqref{eq:gameq} is
\[\Gamma_\real(x_1,\ldots,x_n)\leq \overline\Gamma_\real(x_1,\ldots,x_n).\]

Consider a sequence $\Xi_1,\Xi_2,\ldots$ with values in $\Dbar(0,1)\subset\mathbb C$.
We define an operation $\Xi_\leq$ on formal power series as follows.
Consider the monomial $X_{i_1}^{m_1}\cdot\ldots\cdot X_{i_s}^{m_s}$, where $i_j\neq i_{j+1}$, and $m_j\geq 1$.
On this monomial the effect of  $\Xi_\leq$ is given by
\[\Xi_\leq(X_{i_1}^{m_1}\cdot\ldots\cdot X_{i_s}^{m_s})=(-1)^{\des(i_1,\ldots, i_s)} \Xi_{m_1}
\cdot\ldots\cdot \Xi_{m_s} X_{i_1}^{m_1}\cdot\ldots\cdot X_{i_s}^{m_s}.\]
Then extend $\Xi_\leq$ linearly.
Notice that every monomial gets multiplied by an element of $\Dbar(0,1)$.
In particular, this operation is ``contractive''.

Let us apply $\Xi_\leq$ to  $\mathcal R^{(\lambda)}((\exp X_1)\cdot\ldots\cdot(\exp X_n))$ in the light of \eqref{eq:thpowres}.
We find that the negative signs from the descents vanish, and $\mathcal R^{(\lambda)}(\exp X_j)$
gets replaced by $\Xi_\leq\mathcal R^{(\lambda)}(\exp X_j)$. Thus, formally,
\[\Xi_\leq\mathcal R^{(\lambda)}((\exp X_1)\cdot\ldots\cdot(\exp X_n))=
\mathrm P^{(\lambda)}(\Xi_\leq\mathcal R^{(\lambda)}(\exp X_1),\ldots,\Xi_\leq\mathcal R^{(\lambda)}(\exp X_n)).\]
Let ${}^\Xi\Gamma^{(\lambda)}(X_1,\ldots,X_n)$ denote the power series
$\Xi_\leq\mathcal R^{(\lambda)}((\exp X_1)\cdot\ldots\cdot(\exp X_n))$.
Then, due to the contractivity property, we find that $x_1,\ldots,x_n\in\mathbb R$ the value
${}^\Xi\Gamma^{(\lambda)}(x_1,\ldots,n_n)$ is well-defined (absolute convergent),
as long as $\Gamma_{\real}^{(\lambda)}(|x_1|,\ldots,|x_n|)<+\infty$.

\begin{theorem}\plabel{lem:DRESover}
\[{}^\Xi\Gamma^{[n]}:  (\Gamma^{[n]})^{-1}([0,+\infty))\rightarrow\mathbb C\]
\[(\lambda,x_1,\ldots,x_n)\mapsto {}^\Xi\Gamma^{(\lambda)}(x_1,\ldots,x_n)\]
is a smooth function. It satisfies
\[\left| {}^\Xi\Gamma^{[n]}(\lambda,x_1,\ldots,x_n)\right|\leq \Gamma^{[n]}(\lambda,x_1,\ldots,x_n).\]
\begin{proof}
First, we see that the statement holds for $n=1$, as we have the locally uniform (in $\lambda$) estimates
for the coefficients of the power series expansion of not only for $\mathcal R^{(\lambda)}(\exp X)$
but for $\frac{\mathrm d \mathcal R^{(\lambda)}(\exp X) }{\mathrm d \lambda}=(\mathcal R^{(\lambda)}(\exp X) )^2$,
and other higher derivatives.
After that one can use induction to extend the result.
\end{proof}
\end{theorem}

\begin{theorem}\plabel{thm:conver}
 Assume that $x_1,\ldots, x_n\in[0,\infty)$, and at least two them is positive. Then
\[\overline\Gamma_{\real}(x_1,\ldots,x_n)<+\infty\]
if and only if
\[\Gamma_{\real}(x_1,\ldots,x_n)<+\infty.\]
\begin{proof}
$(\Rightarrow)$ is obvious; we have prove  $(\Leftarrow)$.
We assume that $x_i\neq0$ ($n\geq2$).
Assume, indirectly, that $\Gamma(x_1,\ldots,x_n)<+\infty$, but $\overline\Gamma(x_1,\ldots,x_n)=+\infty$.
We can also assume that for any $t\in[0,1)$, $\overline\Gamma^{(\sup)}(tx_1,\ldots,tx_n)<+\infty$.
Equivalently, this means that for $t\in[0,1)$
\[\sup_{\lambda\in[0,1]}\lambda(1-\lambda)\Gamma^{(\lambda)}_{\real}(tx_1)\Gamma^{(\lambda)}_{\real}(tx_2,\ldots,tx_n) <1,\]
however, the set
\[\Lambda_0=\{\lambda \,:\, \lambda(1-\lambda)\Gamma^{(\lambda)}_{\real}(x_1)\Gamma^{(\lambda)}_{\real}(x_2,\ldots,x_n) =1 \}\]
is non-empty.
Consider, say, $\lambda_0=\inf \Lambda_0\in\Lambda_0$.
Let
\[\Xi_k=\sgn \left( \text{the coefficient of $X^k$ in $\mathcal R^{(\lambda)}(\exp X,\lambda_0)$}\right).\]
Then
\[f(\lambda,t)=\lambda(1-\lambda){}^\Xi\Gamma^{(\lambda)}(tx_1){}^\Xi\Gamma^{(\lambda)}(tx_2,\ldots,tx_n) \]
is smooth on $[0,1]\times[0,1]$, and $|f(\lambda,t)|\leq1$.
Due to its definition, $f(\lambda_0,1)=1$.

Notice that $f(\lambda,t)=1$ implies $t=1$, $\lambda\in\Lambda$, and
${}^\Xi\Gamma^{(\lambda)}(x_1),{}^\Xi\Gamma^{(\lambda)}(x_2,\ldots,x_n)>0$.
Indeed, $t=1$ is implied by the supremum property; then
\begin{equation}
\lambda(1-\lambda){}^\Xi\Gamma^{(\lambda)}(tx_1){}^\Xi\Gamma^{(\lambda)}(tx_2,\ldots,tx_n)=1.
\plabel{eq:crit}
\end{equation}
We have to rule out $\Gamma^{(\lambda)}(x_1),\Gamma^{(\lambda)}(x_2,\ldots,x_n)<0$.
As the starting terms in the power series of $\Gamma^{(\lambda)}(x_1)$
and $\Gamma^{(\lambda)}(x_2,\ldots,x_n)$ are $x_1$ and $x_2+\ldots+x_n$, respectively; negative signs imply
$|{}^\Xi\Gamma^{(\lambda)}(x_1)|<\Gamma^{(\lambda)}_{\real}(x_1)$
and $|{}^\Xi\Gamma^{(\lambda)}(x_2,\ldots,x_n)|<\Gamma^{(\lambda)}_{\real}(x_2,\ldots,x_n)$
in contradiction to \eqref{eq:crit}.

Consider
\[F(\lambda,t)=\dfrac{\Gamma^{(\lambda)}(tx_1)+{}^\Xi\Gamma^{(\lambda)}(tx_2,\ldots,tx_n)+\Gamma^{(\lambda)}(tx_1)
{}^\Xi\Gamma^{(\lambda)}(tx_2,\ldots,tx_n)}{1-\lambda(1-\lambda)
{}^\Xi\Gamma^{(\lambda)}(tx_1){}^\Xi\Gamma^{(\lambda)}(tx_2,\ldots,tx_n)}.\]
It is well-defined and equal to ${}^\Xi\Gamma^{(\lambda)}(tx_1,tx_2,\ldots,tx_n)$ for $0\leq t<1$.
However, according to the previous discussion it extends continuously to $t=1$
by setting $f(\lambda,1)=+\infty$ for $f(\lambda,1)=1$, and the
standard arithmetical value otherwise; the negative part is bounded. Then
\[\lim_{t\nearrow1} {}^\Xi\Gamma(tx_1,\ldots,tx_n)
\equiv\lim_{t\nearrow1} \int_{\lambda=0}^1 F(\lambda,t)\,\mathrm d\lambda=\int_{\lambda=0}^1 F(\lambda,1)\,\mathrm d\lambda.\]
On the other hand, the RHS is $+\infty$ because the smoothness of
$1-f(\lambda,1)$ at $\lambda=\lambda_0$ implies a nonintegrable singularity in $F(\lambda,1)$.
Thus
\[\lim_{t\nearrow1} {}^\Xi\Gamma(tx_1,\ldots,tx_n)=+\infty.\]
But ${}^\Xi\Gamma(tx_1,\ldots,tx_n)\leq
\Gamma_\real(tx_1,\ldots,tx_n)\leq\Gamma_\real(x_1,\ldots,x_n)$ implies $\Gamma_\real(x_1,\ldots,x_n)=+\infty$.
\end{proof}
\end{theorem}

\begin{theorem}\plabel{th:BCHdiv}
 There exits $\varepsilon>0$, such that
\[  \Gamma_\real\Bigl( \underbrace{\frac1{1+\varepsilon}2\arctan \frac1{n-1} ,\ldots,
\frac1{1+\varepsilon}2\arctan \frac1{n-1}}_{\text{$n$ terms}}\Bigr)=+\infty. \]
\begin{proof}
It is easy to show that
\[  \Gamma_\real\Bigl( \underbrace{2\arctan \frac1{n-1} ,\ldots,2\arctan \frac1{n-1}}_{\text{$n$ terms}}\Bigr)=+\infty. \]
Indeed,
$\Gamma^{(1/2)}_{\real}(\ldots)=+\infty\,\Rightarrow\,\overline\Gamma_{\real}(\ldots)=
+\infty\,\Rightarrow\,\Gamma_{\real}(\ldots)=+\infty$.

Next, it is sufficient to show that the supremum of
\begin{equation}
\lambda(1-\lambda)\Gamma^{(\lambda)}_\real
\Bigl( \underbrace{2\arctan \frac1{n-1} ,\ldots,2\arctan \frac1{n-1}}_{\text{$n-1$ terms}}\Bigr)
\Gamma^{(\lambda)}_\real\Bigl( 2\arctan \frac1{n-1} \Bigr)
\plabel{eq:suparg}
\end{equation}
$(\lambda\in[0,1])$ is not taken at $\lambda=\frac12$. (The actual value is $1$ for $\lambda=1/2$.)
Now,
\begin{equation}
\frac{\mathrm d}{\mathrm d\lambda}\mathcal R^{(\lambda)}(\exp X)\Bigl|_{\lambda=\frac12}=
\left(2\tanh\frac X2\right)^2=4\left(1-\frac{\mathrm d}{\mathrm dX}\left(2\tanh\frac X2\right)\right).
\plabel{eq:tanhdiff}
\end{equation}
Note that \eqref{eq:tanhdiff} is supported in even degrees, in contrast to
$\mathcal R^{(\lambda)}(\exp X)\Bigl|_{\lambda=\frac12}=2\tanh\dfrac X2$,
which is supported in odd degrees.
Knowing the power series expansion of $\tanh X$, it is a matter of elementary analysis to show that for $0\leq x<\pi$,
\[\lim_{\lambda \rightarrow \frac12\pm}\frac{\Gamma^{(\lambda)}_\real(x)- \Gamma^{(1/2)}_\real(x)}{\lambda-1/2}
=\pm\left(2\tan\frac x2\right)^2.\]

As \eqref{eq:suparg} is built up from $\Gamma^{(\lambda)}_\real(x)$ (with $x=2\arctan \frac1{n-1}$) and $\lambda(1-\lambda)$
as a positive infinite series, and $\dfrac{\mathrm d}{\mathrm d\lambda}\lambda(1-\lambda)\Bigl|_{\lambda=\frac12}=0$; we see that
\eqref{eq:suparg} has a strict local minimum at $\lambda=1/2$.

Once we know that the supremum of \eqref{eq:suparg} is beyond $1$, then, using continuity, we can slightly contract the variables to
the same effect.
\end{proof}
\end{theorem}
Thus, the (cumulative) convergence radius of the $n$-term BCH-formula is at most $2n\arctan \frac1{n-1}$.
As $n\nearrow\infty$, we have $2n\arctan \frac1{n-1}\searrow2$.
In particular, this demonstrates that the convergence radius of the Magnus expansion cannot be improved beyond
$2$ in the general setting of Banach algebras.
\snewpage

\begin{defin} Let $\mathrm C_1=+\infty$. For $n\geq2$, let
\[\mathrm C_n=\min \{x_1+\ldots+x_n\,:\, x_1,\ldots,x_n\in[0,+\infty), \Gamma_{\real}(x_1,\ldots, x_n)=+\infty\}  ,\]
(equivalent with  $\overline\Gamma_{\real}(x_1,\ldots, x_n)=+\infty$ and $\overline\Gamma_{\real}^{[n]}$ is continuous),
that is the cumulative convergence radius of the Baker--Campbell--Hausdorff formula with $n$ terms.
\end{defin}
\begin{theorem}\plabel{th:BCHconv}
Let $\phi=X_1\mathbf 1_{[0,1)}\boldsymbol.\ldots\boldsymbol. X_n\mathbf 1_{[0,1)}$, $n\geq2$, and $\lambda\in[0,1]$.
Then
\[
\sum_{n=1}^\infty\left|\mu^{(\lambda)}_{n,\mathrm R}(\phi)\right|
\leq \Gamma^{(\lambda)}_\real(|X_1|,\ldots,|X_n|)\leq\overline \Gamma_\real(|X_1|,\ldots,|X_n|),
\]
where every term is well-defined, and they are finite for $|X_1|+\ldots+|X_n|<\mathrm C_n$.

Thus, if $|X_1|+\ldots+|X_n|<\mathrm C_n$ holds, then the conditions of Theorem \ref{th:logMagnus} are satisfied,
and the (stronger) logarithmic version of the Baker-Campbell-Hausdorff formula (with $n$ variables) holds.

On the other hand, there is a Banach algebra $\mathfrak A$ with $X_1,\ldots,X_n\in\mathfrak A$ such that
$|X_1|+\ldots+|X_n|=\mathrm C_n$, but
there is no element $X$ such that $\exp X=(\exp X_1)\ldots(\exp X_n)$ holds.
\begin{proof}
This a consequence of the previous discussion.
\end{proof}
\end{theorem}

\begin{theorem}\plabel{thm:BCHest}
If $n\geq2$, then $2<\mathrm C_n< 2n\arctan \frac1{n-1}$.
As $n$ increases, $\mathrm C_n$ strictly decreases.
\begin{proof}

Regarding the lower bound, we always have convergence for $x_1+\ldots+ x_n\leq2$ according to Theorem \ref{thm:mag2crit}.
Thus, thus, on that domain, $\overline\Gamma(x_1,\ldots,x_n)$ is finite.
Then, by continuity, also in a neighborhood.
Regarding the upper bound, this follows from Theorem \ref{th:BCHdiv}.

Regarding the monotonicity of $\mathrm C_n$, note that if $0<\Gamma^{(\lambda)}_\real( x_k)<+\infty$ and
\[ \lambda(1-\lambda)\Gamma^{(\lambda)}_\real( x_1,\ldots,x_{k-1})\Gamma^{(\lambda)}_\real( x_k)=1, \qquad x_k>0\]
implies
\[ \lambda(1-\lambda)\Gamma^{(\lambda)}_\real( x_1,\ldots,x_{k-1})\Gamma^{(\lambda)}_\real\left( \frac{x_k}2,\frac{x_k}2\right)>1; \]
because in this case $\lambda\neq0,1$, thus
$\Gamma^{(\lambda)}_\real\left( \frac{x_k}2,\frac{x_k}2\right)>\Gamma^{(\lambda)}_\real( x_k)$.
(See $\mathcal R^{(\lambda)}((\exp X)(\exp Y))=\lambda XY+(\lambda-1)YX+\ldots$ for cancellations.)
\end{proof}
\end{theorem}

For $n=2$, this gives  $2<\mathrm C_2<\pi$ for the (cumulative) convergence radius of the BCH expression.
One may wonder what is the actual value.
Due to Theorem \ref{thm:conver}, this can be considered as a numerical problem.
\snewpage

\begin{theorem}\plabel{thm:c2val}
\[\mathrm C_2=2.89847930\ldots\quad.\]

Furthermore, if $\overline\Gamma_\real^{(\lambda)}(x,y)=+\infty$,
$x+y=\mathrm C_2$, then $x=y=\frac12\mathrm C_2$, and
\[0.35865<\min(\lambda,1-\lambda)<0.35866\quad.\]
\end{theorem}
\begin{proof}

First, we would like estimate the value of $\Gamma^{(\lambda)}(x)$.
\begin{lemma}\plabel{lem:DREStech}
 Let $0<\lambda<1$.
We will use the notation
$\left|\log \frac\lambda{\lambda-1}\right|=\sqrt{\pi^2+\left(\log\frac\lambda{1-\lambda}\right)^2}$.
Let
\[\gamma^{(\lambda)}(X)=\left(2+|2\lambda-1|X+4\frac{    \dfrac{X^2}{
\left|\log \frac\lambda{\lambda-1}\right|^2\left(\left|\log
\frac\lambda{\lambda-1}\right|-X\right)}+ \dfrac{X^2}{6\pi^2(3\pi-X)}     }{\sqrt{\lambda(1-\lambda)}}\right)\sinh\frac X2.\]
Then
\[\Gamma^{(\lambda)}(X)\stackrel{\forall X}\leq \gamma^{(\lambda)}(X).\]
\begin{proof} Due to
\[\mathcal R^{(\lambda)}(\exp X)
=\frac{1}{\sqrt{\lambda(1-\lambda)}}\frac{\sinh\frac X2}{\cosh \frac12\left( X -\log\frac\lambda{1-\lambda} \right)},\]
it is sufficient to estimate the $\sech$ term.
Using the complex analytic formula
\[\frac1{\cosh\frac12x}
=\sum_{N=1}^\infty \left( \frac{(-1)^N2\mathrm i}{x-(2N-1)\pi\mathrm i} - \frac{(-1)^N2\mathrm i}{x+(2N-1)\pi\mathrm i} \right),\]
we find that
\begin{align}
\frac{1}{\cosh \frac12\left( X -\log\frac\lambda{1-\lambda} \right)}
\stackrel{\forall X}\leq\,&2\sqrt{\lambda(1-\lambda)}+|2\lambda-1|\sqrt{\lambda(1-\lambda)}\,X\notag\\
&+4\sum_{r=2}^\infty \frac{X^r}{\sqrt{\pi^2+\left(\log\frac\lambda{1-\lambda}\right)^2}^{r+1} }
+4\sum_{r=2}^\infty \sum_{N=2}^\infty\frac{X^r}{((2N-1)\pi)^{r+1}}.\notag
\end{align}
Using the estimate
\[\sum_{N=2}^\infty\frac{1}{((2N-1)\pi)^{r+1}}=\frac{(1-2^{r+1})\zeta(r+1)-1}{\pi^{r+1}}\leq\frac{3/2}{(3\pi)^{r+1}},\]
and combining the expressions, we arrive to the statement.
\end{proof}
\end{lemma}
For $x\in[0,\pi)$, this allows to make the upper and lower estimates
\begin{equation}
\Gamma^{(\lambda)}_{\text{$k$-th $X$-Taylor}}(x)\leq
\Gamma^{(\lambda)}_\real(x)\leq
\Gamma^{(\lambda)}_{\text{$k$-th $X$-Taylor}}(x)-\gamma^{(\lambda)}_{\text{$k$-th $X$-Taylor}}(x) +\gamma^{(\lambda)}_\real(x).
\plabel{eq:estim}
\end{equation}
One can, of course, use other series instead of $\gamma^{(\lambda)}(X)$.
Also note that in this manner we can estimate $\frac{\mathrm d}{\mathrm dx}\Gamma^{(\lambda)}_\real(x)$,  etc.
\snewpage

Having good estimates allows to zero in on $\mathrm C_2$:
\begin{lemma}\plabel{lem:DRESc2}
\[2.8984<\mathrm C_2<2.8985\quad.\]
Moreover, if $\overline\Gamma_\real^{(\lambda)}(x,y)=+\infty$,
i. e.
\[\sup_{\lambda\in[0,1]} \lambda(1-\lambda)\Gamma_\real^{(\lambda)}(x)\Gamma_\real^{(\lambda)}(y)\geq1;\]
and $x+y=\mathrm C_2$,
(that is $x+y$ is minimal form the property above, and in particular,  the supremum is actually $1$),
then
\[0.358<\min(\lambda,1-\lambda)<0.359,\]
and
\[1.445<x,y<1.455\quad.\]
\begin{proof}
This is a series of finer and finer estimates locating critical behaviour.
\end{proof}
\end{lemma}

\begin{lemma}\plabel{lem:DRESmax}
 If
\[2.8984<c<2.8985\]
and
\[0.358<\lambda<0.359\]
and
\[1.445<x<1.455,\]
then the function
\[x\rightarrow \Gamma_\real^{(\lambda)}(x)\Gamma_\real^{(\lambda)}(c-x)\]
is concave in $x$. In particular, due to symmetry, its maximum is taken at $x=c/2$.
\begin{proof}
For a fixed $\lambda$, the function $\Gamma_\real^{(\lambda)}(x)$ is smooth in $x$; and the first and second derivatives
can be estimated
as in \eqref{eq:estim}.
The second derivative of  $\Gamma_\real^{(\lambda)}(x)\Gamma_\real^{(\lambda)}(c-x)$ is
\[ \Gamma_\real^{(\lambda)\prime\prime}(x)\Gamma_\real^{(\lambda)}(c-x)-2\Gamma_\real^{(\lambda)\prime}(x)
\Gamma_\real^{(\lambda)\prime}(c-x)+\Gamma_\real^{(\lambda)}(x)\Gamma_\real^{(\lambda)\prime\prime}(c-x). \]
Thus, we should use upper estimators in the first and third terms, and we should use lower estimators in the middle term.
When we do this, the second derivative turns out to be smaller than $-1$ on the given range.
\end{proof}
\end{lemma}

Theorem \ref{thm:c2val} but with weaker constants follows from a combination of the previous two lemmas.
However, once we know that the maximum is taken for $x=y=\frac12\mathrm C_2$, we can search for $\mathrm C_2$ more effectively.
[End of proof for Theorem \ref{thm:c2val}.]
\end{proof}

\begin{remark}\plabel{rem:DRESinc}
Note that $\mathrm C_2$ is not taken at $\lambda=1/2$.
In fact, combining arguments from Theorem \ref{thm:cayleyest} and Theorem \ref{thm:BCHest},
one can show that this is neither the case for higher $\mathrm C_n$.
On the other hand, for a fixed $x\geq0$,
\[\lim_{\substack{x_1+\ldots+x_n=x, x_i\geq 0,\\
\max\{x_1,\ldots,x_n\}\rightarrow0}}\Gamma_\real^{(\lambda)}(x_1,\ldots,x_n)=\Theta_\real^{(\lambda)}(x).\]
From monotonicity with respect to refinements, and the actual shape of $\Theta_\real^{(\lambda)}(x)$, one can conclude that
the critical $\lambda$ must converge to $1/2$.
\qedremark
\end{remark}

\begin{commenty}
\subsection{Other algebraic consequences of the discrete resolvent decomposition}\plabel{ss:othalg}
~\\

Here we give a formula for the coproduct of the infinitesimal resolvent expression $\mu^{(\lambda)}_n(X_1,\ldots,X_n)$.

Recall that $\mu_n^{(\lambda)}(X_1,\ldots,X_n)$ is just the part of $\mathcal R^{(\lambda)}(\exp(X_1)\cdot\ldots\cdot\exp(X_n) )$
 which is universally $1$-homogeneous in its variables.
Due to this $1$-homogeneity, the coproduct can be computed in the following manner:
First separate the set of indices by $ I=\{i_1,\ldots,i_k\}$ ($i_1<\ldots<i_k$) and $J=\{j_1,\ldots,j_{n-k}\}$ ($j_1<\ldots<j_{n-k}$)
 such that $I\,\dot\cup\,J=\{1,\ldots,n\}$.
Next, we replace $X_s$ by $\tilde X_s$ where $\tilde X_s=X_s\otimes1$ if $s\in I$, and $\tilde X_s=1\otimes X_s$ if $s\in J$.
Then we have to take the part of
 $\mathcal R^{(\lambda)}\Bigl(\exp(\tilde X_1)\cdot\ldots\cdot \exp(\tilde X_n)\Bigr)$ which is uniformly $1$-homogeneous in the variables $X_1,\ldots,X_n$.
Finally, we have to sum all of these $1$-homogeneous parts for all separating pairs $I,J$.

The separation of variables makes commuting sets of variables, thus
\begin{equation}
\exp(\tilde X_1)\cdot\ldots\cdot \exp(\tilde X_n)
=\left( \exp(\tilde X_{i_1})\cdot\ldots\cdot \exp(\tilde X_{i_k}) \right)
\left( \exp(\tilde X_{j_1})\cdot\ldots\cdot \exp(\tilde X_{j_{n-k}}) \right)
\plabel{eq:decexp}
\end{equation}
holds.
Then the discrete resolvent decomposition can be applied.
Note that, in terms of the computation,  \eqref{eq:decexp} provides extra flexibility, as
 the original ordering of the variables is conveniently changed for the application of discrete resolvent decomposition.

Furthermore, we set
\[{}^{[r]}\mu^{(\lambda)}_n(X_1 ,\ldots, X_n )=\sum_{
\substack{
I_1\,\dot\cup\,\ldots\,\dot\cup\,I_r=\{1,\ldots,r\}
\\
I_s=\{i_{s,1},\ldots, i_{s,|I_s|}\}\neq\emptyset
\\
i_{s,1}<\ldots< i_{s,|I_s|}
}
}
\mu^{(\lambda)}_{|I_1|}(X_{i_{1,1}},\ldots,X_{i_{1,|I_1|}})
\cdot\ldots\cdot
 \mu^{(\lambda)}_{|I_r|}(X_{i_{r,1}},\ldots,X_{i_{r,|I_r|}}).
\]
(Here $r\geq1$; the result vanishes for $r>n$.)
Note that
\begin{equation}
\frac1{(r-1)!}\left(\frac{\mathrm d}{\mathrm d\lambda}\right)^{r-1} \mu^{(\lambda)}_n(X_1,\ldots, X_n)
={}^{[r]}\mu^{(\lambda)}_n(X_1,\ldots, X_n).
\plabel{eq:itder}
\end{equation}
This is easy to see from \eqref{eq:muresdef} noting that the derivative of the connecting term
\[\lambda^{\asc(\sigma(u),\sigma(u+1) )}\lambda^{\des(\sigma(u),\sigma(u+1) )}=
\begin{cases}
\lambda&\text{ if }\sigma(u)<\sigma(u+1)\\
\lambda-1&\text{ if }\sigma(u)>\sigma(u+1)
\end{cases}
\]
 is $1$ anyway;
thus taking derivatives just eliminates the connecting terms randomly.
Another way to argue is by noting that
\[\frac1{(r-1)!}\left(\frac{\mathrm d}{\mathrm d\lambda}\right)^{r-1}\mathcal R^{(\lambda)}(\exp(X),\lambda)
=\mathcal R^{(\lambda)}(\exp(X),\lambda)^r;\]
applying to $X=\log(\exp(X_1)\cdot\ldots\cdot\exp(X_n))$; and taking the usual multiplicity argument in $X_1,\ldots,X_n$.
\snewpage
%\begin{commentx}
%~\\~\\~\\~\\
%\end{commentx}
\begin{theorem}
\plabel{th:murescoproduct}
It yields
 \begin{align}
\mu^{( \lambda)}_n&(X_1\otimes1+1\otimes X_1,\ldots, X_n\otimes1+1\otimes X_n)=
\notag
\\\notag=
 &\mu^{( \lambda)}_n(X_1 ,\ldots, X_n )\otimes1+1\otimes\mu^{( \lambda)}_n(X_1 ,\ldots, X_n )
\\
&+\sum_{\substack{  I\,\dot\cup\,J=\{1,\ldots,n\}\\
I=\{i_1,\ldots,i_{|I|}\}\neq\emptyset,\,\,i_1<\ldots<i_{|I|}\\
J=\{j_1,\ldots,j_{|J|}\}\neq\emptyset,\,\,j_1<\ldots<j_{|J|}
 }} \Biggl(
 \notag%\plabel{eq:murescoproduct}
 \\\notag
&\sum_{k=1}^{\min(|I|,|J|)}
(2\lambda-1)\lambda^{k-1}(\lambda-1)^{k-1}\cdot{}^{[k]}\mu^{( \lambda)}_n(X_{i_1} ,\ldots, X_{i_{|I|}} )
\otimes
 {}^{[k]}\mu^{( \lambda)}_n(X_{j_1} ,\ldots, X_{j_{|J|}} )
\\\notag
&\quad+\sum_{k=1}^{\min(|I|-1,|J|)}
\lambda^k(\lambda-1)^k\cdot{}^{[k+1]}\mu^{( \lambda)}_n(X_{i_1} ,\ldots, X_{i_{|I|}} )
\otimes
 {}^{[k]}\mu^{( \lambda)}_n(X_{j_1} ,\ldots, X_{j_{|J|}} )
\\\notag
&\quad+\sum_{k=1}^{\min(|I|,|J|-1)}
  \lambda^k(\lambda-1)^k\cdot{}^{[k]}\mu^{( \lambda)}_n(X_{i_1} ,\ldots, X_{i_{|I|}} )
\otimes
 {}^{[k+1]}\mu^{( \lambda)}_n(X_{j_1} ,\ldots, X_{j_{|J|}} )
\Biggr)
\\\notag=
 &\mu^{( \lambda)}_n(X_1 ,\ldots, X_n )\otimes1+1\otimes\mu^{( \lambda)}_n(X_1 ,\ldots, X_n )
\\
&+\sum_{\substack{  I\,\dot\cup\,J=\{1,\ldots,n\}\\
I=\{i_1,\ldots,i_{|I|}\}\neq\emptyset,\,\,i_1<\ldots<i_{|I|}\\
J=\{j_1,\ldots,j_{|J|}\}\neq\emptyset,\,\,j_1<\ldots<j_{|J|}
 }}\frac{\mathrm d}{\mathrm d\lambda} \Biggl(
 \sum_{k=1}^{\min(|I|,|J|)}
 \frac{ \lambda^k(\lambda-1)^k}k\cdot
 \notag%\plabel{eq:murescoproductver}
\\\notag
&
\frac1{(k-1)!}\left(\frac{\mathrm d}{\mathrm d\lambda}\right)^{k-1}\mu^{( \lambda)}_n(X_{i_1} ,\ldots, X_{i_{|I|}} )
\otimes
\frac1{(k-1)!}\left(\frac{\mathrm d}{\mathrm d\lambda}\right)^{k-1}\mu^{( \lambda)}_n(X_{j_1} ,\ldots, X_{j_{|J|}} )
\Biggr)
 .
\end{align}
(In the $\sum$'s, instead of the $\min(\ldots,\ldots)$'s we could have taken $\infty$ for the upper limit, as the omitted terms vanish.)
\begin{proof}
The first equality is the direct transcription of Theorem \ref{th:powres} (cf. Lemma \ref{lem:powres}) with respect to the order of variables
as in the right side of \eqref{eq:decexp} and then summed up.
The result can be rewritten using \eqref{eq:itder}, yielding the second equality.
%For the second equality, using \eqref{eq:itder}, we can write \eqref{eq:murescoproduct} into a more suggestive form.
\end{proof}
\end{theorem}

\begin{cor}
\plabel{cor:murescoface}
 For $n=1$, $\mu^{(\lambda)}_n(1)=1$. For $n\geq2$,
\begin{align*}
\mu^{(\lambda)}_n&(X_1,\ldots,X_s,1,X_{s+1},\ldots,X_{n-1})=
\\=&
(2\lambda-1)\mu^{(\lambda)}_{n-1}(X_1,\ldots,X_{n-1})+\lambda(\lambda-1)\,{}^{[2]}\mu^{(\lambda)}_{n-1}(X_1,\ldots,X_{n-1})
\\=&\frac{\mathrm d}{\mathrm d\lambda}
\left(\lambda(\lambda-1)\mu^{(\lambda)}_{n-1}(X_1,\ldots,X_{n-1})\right).
\end{align*}
\begin{proof}
This follows taking the $(1,n-1)$-homogeneous terms in %\eqref{eq:murescoproduct} or \eqref{eq:murescoproductver}.
 the previous theorem.
\end{proof}
\end{cor}
\snewpage
\begin{remark}
\plabel{rem:murescoface}
A formal consequence of the Corollary \ref{cor:murescoface} is, for $n\geq2$,
\[\frac{\mathrm d}{\mathrm dX}\mu^{(\lambda)}_n(X ,\ldots,X )=
n\frac{\mathrm d}{\mathrm d\lambda}\left(\lambda(\lambda-1)\mu^{(\lambda)}_{n-1}(X,\ldots,X)\right).\]
As $\mu^{(\lambda)}_n(X ,\ldots,X )=n!G(\lambda,\lambda-1)X^n$, this shows that Corollary \ref{cor:murescoface}  is
a generalization of Lemma \ref{lem:goldgen}.
\qedremark
\end{remark}

\begin{theorem}
\plabel{th:Fmu}
For the co-shuffle $F$,
 \begin{align}
F(\mu^{( \lambda)}_n&(X_1,\ldots,  X_n))=
\notag
\\\notag=
 &2\mu^{( \lambda)}_n(X_1 ,\ldots, X_n )
 \\\notag
&+\sum_{k=1}^{\infty}
(2\lambda-1)\lambda^{k-1}(\lambda-1)^{k-1}\cdot{}^{[2k]}\mu^{( \lambda)}_n(X_1 ,\ldots, X_n )
\\\notag
&+\sum_{k=1}^{\infty}
2\lambda^k(\lambda-1)^k\cdot{}^{[2k+1]}\mu^{( \lambda)}_n(X_1 ,\ldots, X_n )
\\\notag=
 &2\mu^{( \lambda)}_n(X_1 ,\ldots, X_n )
 \\\notag
&+\sum_{k=1}^{\infty}
(2\lambda-1)\lambda^{k-1}(\lambda-1)^{k-1}\cdot\frac1{(2k-1)!}\left(\frac{\mathrm d}{\mathrm d\lambda}\right)^{2k-1}
\mu^{( \lambda)}_n(X_1 ,\ldots, X_n )
\\\notag
&+\sum_{k=1}^{\infty}
2\lambda^k(\lambda-1)^k\frac1{(2k)!}\left(\frac{\mathrm d}{\mathrm d\lambda}\right)^{2k}
\mu^{( \lambda)}_n(X_1 ,\ldots, X_n )
\\\notag=
 &2\mu^{( \lambda)}_n(X_1 ,\ldots, X_n )
\\\notag
&+\frac{\mathrm d}{\mathrm d\lambda}\sum_{k=1}^{\infty}
2\lambda^k(\lambda-1)^k\frac1{(2k)!}\left(\frac{\mathrm d}{\mathrm d\lambda}\right)^{2k-1}
\mu^{( \lambda)}_n(X_1 ,\ldots, X_n )
 .
\end{align}
%(In the $\sum$'s, we could have taken $\infty$ for the upper limit, as the omitted terms vanish.)
\begin{proof}
This follows from the previous theorem, taking \eqref{eq:itder} into account.
\end{proof}
\end{theorem}
\begin{cor}
\plabel{cor:FrieSpec}
For polynomials of $\lambda$, we define the co-shuffle differential operator as
 \begin{align}
\notag\Delta^F=
 &2\Id
 \\\notag
&+\sum_{k=1}^{\infty}
(2\lambda-1)\lambda^{k-1}(\lambda-1)^{k-1}\cdot\frac1{(2k-1)!}\left(\frac{\mathrm d}{\mathrm d\lambda}\right)^{2k-1}
\\\notag
&+\sum_{k=1}^{\infty}
2\lambda^k(\lambda-1)^k\frac1{(2k)!}\left(\frac{\mathrm d}{\mathrm d\lambda}\right)^{2k}
 .
\end{align}
Then $G_n(\lambda)\equiv G_n(\lambda,\lambda-1)$ is an eigenvector of eigenvalue $2^n$
(leading to an eigendecomposition of the space of the polynomials).
\begin{proof}
This follows by comparing Theorems \ref{th:gemcan}, \ref{th:Frie}, and  \ref{th:Fmu}.
\end{proof}
\end{cor}

\end{commenty}
\snewpage

\section{The Magnus expansion in the adjoint representation}
\plabel{sec:adj}
In the previous section our experience was that the resolvent method is more capable than the
 original combinatorial method.
Here we will see a particular setting where the resolvent method does not do that well.
This will be the case of the adjoint representation.
See Appendix \ref{sec:sc} for the general spectral considerations of adjoint actions.

\subsection{The adjoint representation of the Magnus expansion}\plabel{ss:adjaction}
~\\

As an analytic application of the spectral calculus associated to the resolvent expression,
a case of some practical importance is as follows.
Let us consider the function
\[\beta(x)=\frac{x}{\mathrm e^x-1}\]
understood as a meromorphic function on $\mathbb C$; analytic, in particular for $|\Ima x|<2\pi$.
Recall that (our spectral) $\log A$, whenever defined, has spectrum lying in $\{x\in\mathbb C\,:\,|\Ima x|<\pi\}$.
Therefore $\beta(\log A)$ are $\beta(-\log A)$ are automatically defined.
In fact,
\[
\beta(\log A)=\int_{\lambda=0}^1 \frac{1 }{\lambda +(1-\lambda)A }\,\mathrm d\lambda
\qquad
\text{and}
\qquad
\beta(-\log A)=\int_{\lambda=0}^1 \frac{A }{\lambda +(1-\lambda)A }\,\mathrm d\lambda.
\]
(Here the constant term `$1$' of $\beta$ goes to $1_{\mathfrak A}$ in the spectral calculus.)
\begin{lemma}\plabel{lem:buu}
For $\int |\phi|<2$, or $\int |\phi|\leq2$ and $\phi$ is the Lebesgue measure times a Lebesgue--Bochner
integrable function, or even if $\phi$ is just $M$-controlled,
\[\beta(\log\Rexp(\phi))=1_{\mathfrak A}+\sum_{k=1}^\infty \int_{\lambda=0}^1 (\lambda-1)\,\mu_{k,\mathrm R}^{(\lambda)}(\phi)\,\mathrm d\lambda \]
and
\[\beta(-\log\Rexp(\phi))=1_{\mathfrak A}+\sum_{k=1}^\infty \int_{\lambda=0}^1 \lambda\,\mu_{k,\mathrm R}^{(\lambda)}(\phi)\,\mathrm d\lambda .\]
Moreover, in these circumstances, `$\beta$' can be replaced by $\beta^{\mathrm{pow}}$.
\begin{proof}
This immediately follows from the resolvent integrals
\[\beta(X)=1+\int_{\lambda=0}^1 (\lambda-1)\,\mathcal R^{(\lambda)}(\exp X)\,\mathrm d\lambda
\qquad\text{and}\qquad
\beta(-X)=1+\int_{\lambda=0}^1 \lambda\,\mathcal R^{(\lambda)}(\exp X)\,\mathrm d\lambda.
\]
The comment about the replacement of `$\beta$' follows from  $\mathrm r\left(\Rexp(\phi)\right)<\pi$.
\end{proof}
\end{lemma}

We can consider the adjoint action $\ad:\mathfrak A\rightarrow \mathcal B(\mathfrak A)$.
Here $(\ad X)Y=[X,Y]\equiv XY-YX$.
We can take the operator norm $\|\cdot\|_{\mathfrak A}$ on $\mathcal B(\mathfrak A)$.
Then $\|\ad X\|_{\mathfrak A}\leq2|X|$.
It is also easy  to see that $\Rexp( \ad \phi) Y=\Rexp(\phi) Y \Rexp(\phi)^{-1}$.
For the sake of the next statement we use the fact  that $\mu_n(X_1,\ldots, X_n)$
 is formally a commutator polynomial of the $X_i$ (cf. Theorem \ref{th:magcomm}); therefore
\begin{equation}
\ad \mu_n(X_1,\ldots, X_n)=\mu_n(\ad X_1,\ldots,\ad  X_n).
\plabel{eq:admu}
\end{equation}

We also note that whenever $\log A$ is defined,
 the spectrum of $\ad \log A$ lies in $\{x\in\mathbb C\,:\,|\Ima x|<2\pi\}$ (cf. Lemma \ref{lem:commspec}).
Therefore $\beta(\ad \log A)$ and $\beta(-\ad \log A)$ are automatically defined.
(Here the constant term `$1$' of $\beta$ goes to $\Id_{\mathfrak A}$ in the spectral calculus.)
\snewpage

\begin{theorem}\plabel{th:buu}
If $\int |\phi|<2$,  or $\int |\phi|\leq2$ and $\phi$ is the Lebesgue measure times a Lebesgue--Bochner
integrable function,  or even if $\phi$ is just $M$-controlled, then
\begin{equation}
\ad \log\Rexp(\phi)=\sum_{k=1}\mu_{\mathrm R,k}(\ad \phi)
\equiv\sum_{k=1}^\infty \int_{\lambda=0}^1 \mu_{k,\mathrm R}^{(\lambda)}(\ad \phi)\,\mathrm d\lambda,
\plabel{eq:magnusadj}
\end{equation}
and
\begin{equation}
\beta(\ad \log\Rexp(\phi))
=\Id_{\mathfrak A}+\sum_{k=1}^\infty \int_{\lambda=0}^1 (\lambda-1)\,\mu_{k,\mathrm R}^{(\lambda)}(\ad \phi)\,\mathrm d\lambda ,
\plabel{eq:betaconv1}
\end{equation}
and
\begin{equation}\beta(-\ad \log\Rexp(\phi))
=\Id_{\mathfrak A}+\sum_{k=1}^\infty \int_{\lambda=0}^1 \lambda\,\mu_{k,\mathrm R}^{(\lambda)}(\ad \phi)\,\mathrm d\lambda;
\plabel{eq:betaconv2}
\end{equation}
the sums are decreasing geometrically. (The point is validity for the domain $\int |\phi|<2$.)

Moreover, in these circumstances, `$\beta$' can be replaced by $\beta^{\mathrm{pow}}$.
\begin{proof}
Equation \eqref{eq:admu} and the absolute convergence of the Magnus expansion implies \eqref{eq:magnusadj}.
Next, we prove  \eqref{eq:betaconv1}.
If  $\int |\phi|<1$, then $\int \|\ad \phi\|_{\mathrm{op}}<2$,
 $\beta(\ad \log\Rexp(\phi))=\beta(\log\Rexp(\ad \phi))$, and the latter term expands as indicated.
Otherwise we should replace $\phi$ by $t\cdot\phi$, and consider it for $t\in\intD(0,1+\varepsilon)$.
As the LHS defines an analytic function on this domain, the power series representation extends.
 \eqref{eq:betaconv2} is similar.
The comment about the replacement of `$\beta$' follows from  $\mathrm r\left(\ad\log\Rexp(\phi)\right)<2\pi$.
\end{proof}
\end{theorem}
\begin{remark}\plabel{rem:schurformalpro}
In the setting of the previous lemma, continuing the lines \eqref{eq:betaconv1} and \eqref{eq:betaconv2} further,
\begin{multline*}
\beta(\ad \log\Rexp(\phi))Y
=\lim_{t\searrow 0}\frac{\log\Rexp(Y\mathbf 1_{[0,t)}\boldsymbol.\phi) -\log\Rexp(\phi)} {t}\\
=\int_{\lambda=0}^1
\frac1{\lambda+(1-\lambda)\Rexp(\phi)}Y\Rexp(\phi)\frac1{\lambda+(1-\lambda)\Rexp(\phi)} \,\mathrm d\lambda,
\end{multline*}
and
\begin{multline*}
\beta(-\ad \log\Rexp(\phi))Y
=\lim_{t\searrow 0}\frac{\log\Rexp(\phi\boldsymbol.Y\mathbf 1_{[0,t)}) -\log\Rexp(\phi)} {t}\\
=\int_{\lambda=0}^1
\frac1{\lambda+(1-\lambda) \Rexp(\phi)}\Rexp(\phi)Y\frac1{\lambda+(1-\lambda)\Rexp(\phi)} \,\mathrm d\lambda.
\end{multline*}
These are the left and right ``chronological derivatives'' of $\log\Rexp(\phi)$.
\qedremark
\end{remark}
\snewpage
This leads, more generally, to the ODE viewpoint.
First, we formulate it in the formal case where it is uncontroversial.
Assume that $\phi$ is supported the interval $I$, and $T$ is formal commutative variable.
Then
\begin{multline*}
\mu_{\mathrm R}(T\cdot\phi)= \int_{t\in I}\beta(\ad\log\Rexp(T\cdot \phi|_{I\cap(t,+\infty)}) )\,(T\cdot\phi(t))
\\\equiv \int_{t\in I}\beta(\ad\mu_{\mathrm R}(T\cdot \phi|_{I\cap(t,+\infty)}) )\,(T\cdot\phi(t))
\end{multline*}
and
\begin{multline*}
\mu_{\mathrm R}(T\cdot\phi)= \int_{t\in I}\beta(-\ad\log\Rexp(T\cdot \phi|_{I\cap(-\infty,t)}) )\,(T\cdot\phi(t))
\\\equiv
\int_{t\in I}\beta(-\ad\mu_{\mathrm R}(T\cdot \phi|_{I\cap(-\infty,t)}) )\,(T\cdot\phi(t)).
\end{multline*}

Indeed, this is just the rewriting of the integration simplex-wise, using $t_1$ as $t$ or $t_n$ as $t$ respectively,
and meanwhile applying \eqref{eq:betaconv1} and \eqref{eq:betaconv2} formally.

More generally, using the notation
\[\Omega^{\phi}_{\mathrm R+}(t)=\mu_{\mathrm R}(\phi|_{I\cap(t,+\infty)})\]
and
\[\Omega^{\phi}_{\mathrm R-}(t)=\mu_{\mathrm R}(\phi|_{I\cap(-\infty,t)}),\]
we can write
\begin{equation}
\Omega^{T\cdot\phi}_{\mathrm R+}(s)
=\int_{t\in I\cap(s,+\infty)}\beta\left(\ad \Omega^{T\cdot\phi}_{\mathrm R+}(t)\right)(T\cdot\phi(t))
\plabel{eq:arr1}\end{equation}
and
\begin{equation}\Omega^{T\cdot\phi}_{\mathrm R-}(s)
=\int_{t\in I\cap(-\infty,s)}\beta\left(-\ad \Omega^{T\cdot\phi}_{\mathrm R-}(t)\right)(T\cdot\phi(t)).
\plabel{eq:arr2}\end{equation}
If $\phi(t)$ is sufficiently nice, say $\phi(t)=f(t)\,\mathrm dt$, where $f$ is continuous, and $I=[a,b]$ then
these equations can be written in differential form
\[\dfrac{\mathrm d\Omega^{T\cdot\phi}_{\mathrm R+}(t)}{\mathrm dt}=
\beta\left(\ad \Omega^{T\cdot\phi}_{\mathrm R+}(t)\right)(T\cdot f(t))\qquad\qquad \left( \Omega^{T\cdot\phi}_{\mathrm R+}(b)=0 \right)\]
and
\[\dfrac{\mathrm d\Omega^{T\cdot\phi}_{\mathrm R-}(t)}{\mathrm dt}=
\beta\left(-\ad \Omega^{T\cdot\phi}_{\mathrm R-}(t)\right)(T\cdot f(t))\qquad\qquad \left( \Omega^{T\cdot\phi}_{\mathrm R-}(a)=0 \right).\]
(What makes the ``ODE viewpoint'' is not necessarily the differential formulation but the way of accounting the integration on the
simplices.)
\snewpage

For completeness, we also give the formulae in the ``L'' formalism:
For
\[\Omega^{\phi}_{\mathrm L+}(t)=\mu_{\mathrm L}(\phi|_{I\cap(t,+\infty)})\]
and
\[\Omega^{\phi}_{\mathrm L-}(t)=\mu_{\mathrm L}(\phi|_{I\cap(-\infty,t)}),\]
we have
\begin{equation}
\Omega^{T\cdot\phi}_{\mathrm L+}(s)
=\int_{t\in I\cap(s,+\infty)}\beta\left(-\ad \Omega^{T\cdot\phi}_{\mathrm L+}(t)\right)(T\cdot\phi(t))
\plabel{eq:arr3}\end{equation}
and
\begin{equation}
\Omega^{T\cdot\phi}_{\mathrm L-}(s)
=\int_{t\in I\cap(-\infty,s)}\beta\left(\ad \Omega^{T\cdot\phi}_{\mathrm L-}(t)\right)(T\cdot\phi(t)).
\plabel{eq:arr4}\end{equation}

Using \eqref{eq:arr1}, \eqref{eq:arr2}, \eqref{eq:arr3}, \eqref{eq:arr4},
 $\Omega^{T\cdot\phi}_{**}(t)$ (and thus the Magnus expansion) can be built up in degrees of $T$ inductively.
Apart from the presence of this formal $T$,
the ODE formalism is a traditional viewpoint of Magnus expansion, which indeed
continues the Schur--Poincaré and Campbell--Baker--Hausdorff lines of arguments (cf. Discussion \ref{disc:tradBCH}).
For example, in the paper of Magnus \cite{M}, $\Omega(t)\equiv\Omega^{ \text{`$x\mapsto$'} A(x)\,\mathrm dx}_{\mathrm L-}(t)$
and $Y(t)\equiv\Lexp(\text{`$x\mapsto$'} A(x)\,\mathrm dx|_{I\cap(-\infty,t)})$ is used.
From our viewpoint, the Magnus expansion is a series in the ``homogeneity of $\phi$''
(here: in the variable $T$), the ODE approach also considers this homogeneity accounting but introduces
 an additional, chronological, accounting,  an ``attack of direction'', which is not quite
 necessary, but it might be very useful.
(Without the homogeneity accounting , we have just a ``continuation of logarithm''.)
Of course, we have already encountered a kind of ODE approach with the chronological decompositions and the ``delay'' method.

Again, the analytical question is whether we can put $T=1$ into the formal expressions.
Some convergence estimates can be obtained using the formal case,
 but here, in the general Banach algebraic context, we are already beyond that.

\begin{theorem}\plabel{th:conex}
If $\int |\phi|<2$,  or $\int |\phi|\leq2$ and $\phi$ is the Lebesgue measure times a Lebesgue--Bochner
integrable function,  or even if $\phi$ is just $M$-controlled, then
\[\Omega^{\phi}_{\mathrm R+}(s)
=\int_{t\in I\cap(s,+\infty)}\beta\left(\ad \Omega^{\phi}_{\mathrm R+}(t)\right)(\phi(t)),\]
\[\Omega^{\phi}_{\mathrm R-}(s)
=\int_{t\in I\cap(-\infty,s)}\beta\left(-\ad \Omega^{\phi}_{\mathrm R-}(t)\right)(\phi(t)),\]
etc.
\begin{proof}
For $\int|\phi|<1$, these statements hold by naive arguments using absolute convergence.
In the general  $M$-controlled case, one should replace $\phi$ by $\tau\cdot\phi$, and use
analytical continuation from $\tau=0$ to $\tau=1$ (we only use that the expressions are well-defined).
\end{proof}
\end{theorem}

It is natural to ask whether \eqref{eq:magnusadj} holds in logarithmic sense,
i. e. $\ad \log\Rexp(\phi)= \log\Rexp(\ad\phi)$ holds or not.
Or, rephrased, whether $\int$ and $\sum$ can be interchanged in \eqref{eq:magnusadj} or not.
The general answer turns out to be negative:
For example, we cannot guarantee the existence of
$\mathcal R^{(1/2)}(\Rexp(\ad\phi)) $ (and, in particular,
the absolute convergence of $\sum_{k=1}^\infty   \mu_{k,\mathrm R}^{(1/2)}(\ad \phi)$) for $\int |\phi|\geq\sqrt2$.
In order to deal with this further, we need some technical tools.
\snewpage

\subsection{The formal adjoint representation and the convergence problem}\plabel{ss:formaladj}
~\\

In another approach to the adjoint actions, we can extend $\mathfrak A$ by a universal element $W$ of norm $1$,
and then, we can consider the adjoint actions on the universal element $W$.
This we call the ``extended setting''.
Technically, that means that we consider the projective tensor algebra $\mathfrak A^{1}[W]$,
and consider $\ad^*:\mathfrak A\rightarrow \mathcal B(\mathfrak A^{1}[W])$
and $\Ad^*:\mathfrak A\rightarrow \mathcal B(\mathfrak A^{1}[W])$.
Here we consider the operator norm $\|\cdot\|_{\mathfrak A}^*$ on $\mathfrak A^{1}[W]$,
but the relevant norms can be tested on the single element $W$.
(I. e. $\ad X$ is $\ad_{\mathfrak A} X$ and  $\ad^* X$ is $\ad_{\mathfrak A^{1}[W]} X$ properly.
We  also use the notation $\|\cdot\|_{\mathfrak A}^*$ instead of  $\|\cdot\|_{\mathfrak A^{1}[W]}$.
We may still use the notation $\ad X$ if the ambient algebra is specified explicitly.)
Then, for  $X\in\mathfrak A$, we have $\|\ad X\|_{\mathfrak A} \leq \|\ad^* X\|_{\mathfrak A}^*\leq 2|X|_{\mathfrak A}$, etc.
(Here we also have the extended $X_{\mathrm L}^*$ and $X_{\mathrm R}^*$.)
In fact, due to the universality of the choice of $W$, we see that $ \|\ad^* X\|_{\mathfrak A}^*$
is the most pessimistic choice for the operator norm of $\ad X$ in any ambient algebra of $\mathfrak A$.
Here we have similar results for the spectrum $\ad^*$ and $\Ad^*$ as before,
 which is reasonable, as we have simply change
from the Banach algebra $\mathfrak A$ to  $\mathfrak A^{1}[W]$.
(For example, throughout the previous subsection, $\ad$ can be replaced by $\ad^*$.
Also, the extension will be isospectral with respect to elements of $\mathfrak A$:
Any resolvental inverse would remain such after substituting $W=0$.)
We note that if $\mathfrak A=\mathrm F^1_{\mathbb K}(\Omega,\mathfrak F,\omega)$, then
the structure of $\mathfrak A^{1}[W]$ is really simple: all we have to do is to extend the measure space $(\Omega,\mathfrak F,\omega)$
by a Dirac delta measure.

We will use the fact that $\ad$ and $\Ad$ are two-sided multiplication operators,
and in the extended setting this can be transferred to dealing with a special algebra.

Ultimately, the interest of the following discussion will be the case of formal power series with continuous variables,
 but, for the sake of simplify we will start with the case of ordinary formal power series with $\ell^1$ norm.
If $P\in\mathrm F^{1}[Y_1,Y_2,\ldots]$, then the adjoint action $\ad P$ is given by $(\ad P)W=PW-WP$.
This $W$ can be in $\mathrm  F^{1}[Y_1,Y_2,\ldots]$ arbitrary, but it does not matter if $W$ is a newly adjoined
 variable to $\mathrm  F^{1}[Y_1,Y_2,\ldots]$.
(By the universal properties of these algebras, the norm of $(\ad P)Z$ can be used for testing the operator norm for $ \ad P $.)
In order not spend time with searching such a $W$, the adjoint action can be imagined to be two-step process:
First we assign to $P$ the element
\[\widetilde\delta(P)=P\otimes 1^{\mathrm{op}}-1\otimes P^{\mathrm{op}}\in
\mathrm  F^{1}[Y_1,Y_2,\ldots] \otimes \mathrm  F^{1}[Y_1,Y_2,\ldots]^{\mathrm{op}}  \]
 (the opposite algebra and the commutative tensor product is used).
This gets further represented on $W$ by the recipe  $(P_1\otimes P_2^{\mathrm{op}})W=P_1WP_2$.
Similarly, $\Ad P$ can be represented by
 \[\widetilde\Delta(P)=P\otimes (P^{-1})^{\mathrm{op}}\in
\mathrm  F^{1}[Y_1,Y_2,\ldots] \otimes \mathrm  F^{1}[Y_1,Y_2,\ldots]^{\mathrm{op}}.  \]
The advantage is that on $\mathrm  F^{1}[Y_1,Y_2,\ldots] \otimes \mathrm  F^{1}[Y_1,Y_2,\ldots]^{\mathrm{op}} $
 the $\ell^1$ norm can be used, and we do not have to adjoin a new variable $W$ in order to test the operator norm
 of the adjoint action.
Here it is actually reasonable to pass to the complete(d) algebra
 $\mathrm  F^{1}[Y_1,Y_2,\ldots] \otimes_{\ell^1} \mathrm  F^{1}[Y_1,Y_2,\ldots]^{\mathrm{op}} $.
This passage will automatically be understood in the future.
\snewpage

More generally, we will have an action $U$ on $\mathfrak A^1[W]$ which sends $W$ to an $U(W)$ which is $1$-homogeneous in $W$,
and the general $U(S)$ is $U(W)$ but $S$ is substituted to the place of $W$.
This $U$ can be represented linearly by an element of $\mathfrak A\otimes\mathfrak A$ (projective tensor product,
corresponding to $U(W)$), but as an algebra element as of $\mathfrak A\otimes\mathfrak A^{\mathrm{op}}=:\widetilde{\mathbb I}\, \mathfrak A$.

Let $\LR$ denote the action discussed above.
Then $\LR(\mathfrak A\otimes\mathfrak A^{\mathrm{op}})$ is the set of elements in $\mathcal B(\mathfrak A^1[W])$ obtained
in this manner.
\begin{lemma}\plabel{lem:resfone}
If $U\in \LR(\mathfrak A\otimes\mathfrak A^{\mathrm{op}})$ and $U^{-1}$ exist in $\mathcal B(\mathfrak A^1[W])$,
 then $U^{-1}\in \LR(\mathfrak A\otimes\mathfrak A^{\mathrm{op}})$.
Similar observation applies to other resolvent expressions.
\begin{proof}
Assume that $U=\LR(\theta)$.
Let us consider $U^{-1}(W)$.
Its component of $\deg_W=1$ induces an element $\eta\in \mathfrak A\otimes\mathfrak A^{\mathrm{op}}$.
Then $W=U(U^{-1}(W))$ implies that $\theta\eta=1_{\mathfrak A\otimes\mathfrak A^{\mathrm{op}}}$.
Taking the representation $\LR$, this implies $U \LR(\eta)=\Id_{\mathfrak A^1[W]}$.
Multiplying by $U^{-1}$ on the left, we obtain $\LR(\eta)=U^{-1}$.
As the identity is also represented in this manner, as $\LR(1_{\mathfrak A}\otimes1_{\mathfrak A}^{\mathrm{op}})$,
 the extension to resolvent expressions is natural.
\end{proof}
\end{lemma}

\begin{lemma}\plabel{lem:startfone}
(a) If $X\in\mathfrak A$, then
\[\spec(\ad^* X)=\spec(\widetilde\delta X)=\spec(X^*_{\mathrm L})+\spec(-X^*_{\mathrm R})= \spec(X)-\spec(X).\]

(b) If $X\in\mathfrak A$ is invertible, then
\[\spec(\Ad^* X)=\spec(\widetilde\Delta X)=\spec(X^*_{\mathrm L})\cdot\spec((X^*_{\mathrm R})^{-1})= \spec(X)/\spec(X).\]
\begin{proof}
Assume that $\lambda_1\in\spec(X)$.
Then the closure of $\mathfrak A(X-\lambda_1)$ cannot be the whole $\mathfrak A$;
indeed it cannot intersect even the open unit ball around $1_{\mathfrak A}$.
Using the Banach--Hahn theorem, there is a weakly contractive linear functional
$\rho_1$ which vanishes on $\overline{\mathfrak A(X-\lambda_1)}$ and it takes the value $1$ on $1_{\mathfrak A}$.
Similarly, if $\lambda_2\in\spec(X)$, then there is a weakly contractive linear functional
$\rho_2$ which vanishes on $\overline{(X-\lambda_2)\mathfrak A}$ and it takes the value $1$ on $1_{\mathfrak A}$.
Then  $\overline{\mathfrak A(X-\lambda_1) \otimes  \mathfrak A + \mathfrak A \otimes(X-\lambda_2)\mathfrak A}$
is a proper subspace of $\mathfrak A\otimes \mathfrak A$.
Indeed, $\rho_1\otimes\rho_2$ vanishes on the subspace, while it takes $1$ on $1_{\mathfrak A}\otimes 1_{\mathfrak A}$.
$X_1\otimes X_2^{\mathrm{op}}\in\mathfrak A\otimes \mathfrak A^{\mathrm{op}}$ acts on
$\mathfrak A\otimes\mathfrak A$ by $(X_1\otimes X_2^{\mathrm{op}})(A_1\otimes A_2)=X_1A_1\otimes A_2X_2$;
and so on the factor space $(\mathfrak A\otimes\mathfrak A) / \overline{\mathfrak A(X-\lambda_1) \otimes  \mathfrak A + \mathfrak A \otimes(X-\lambda_2)\mathfrak A} $.
Let $\Omega$ be the image of $1_{\mathfrak A}\otimes 1_{\mathfrak A}$ in the latter factor space;
we know that $\Omega\neq0$. In case (a)
\[(X\otimes 1-1\otimes X^{\mathrm{op}}-(\lambda_1-\lambda_2)\Id )\Omega=0;\]
and in case (b)
\[(X\otimes (X^{-1})^{\mathrm{op}}-(\lambda_1/\lambda_2)\Id )\Omega=0;\]
contradicting to the invertibility of those elements in $\mathfrak A\otimes\mathfrak A^{\mathrm{op}}$.
\end{proof}
\end{lemma}

\snewpage

Returning to the free $\ell^1$ case, we can avoid the use of the opposite algebra with an appropriate
(norm-compatible) formal star operation $^\ddag$.
In that setting, we assign to $P$ the element
\[\widehat\delta(P)=P\otimes 1-1\otimes P^\ddag\in
\mathrm  F^{1}[Y_1,Y_2,\ldots] \otimes \mathrm  F^{1}[Y_1,Y_2,\ldots] ; \]
which can be represented on $Z$   by the recipe  $(P_1\otimes P_2 )W=P_1WP_2^\ddag$; etc.
We can use the formal star operation given by the prescription $(Y_i)^\ddag=-Y_i$.
This has the advantage that $\widehat\delta(Y_i)=Y_i\otimes 1+1\otimes Y_i$.
In what follows, beside
$\widetilde{\mathbb I}\mathrm  F^{1}[Y_1,Y_2,\ldots] \equiv\mathrm  F^{1}[Y_1,Y_2,\ldots] \otimes \mathrm  F^{1}[Y_1,Y_2,\ldots]^{\mathrm{op}}$
we also use the notation
$\widehat{\mathbb I}\mathrm  F^{1}[Y_1,Y_2,\ldots] \equiv\mathrm  F^{1}[Y_1,Y_2,\ldots] \otimes \mathrm  F^{1}[Y_1,Y_2,\ldots] $.
The locally convex algebras $\widetilde{\mathbb I}\mathrm  F^{1,\mathrm{loc}}[Y_1,Y_2,\ldots] $ and
 $\widehat{\mathbb I}\mathrm  F^{1,\mathrm{loc}}[Y_1,Y_2,\ldots] $
 can also be considered.
These are induced by (finite) seminorms through the $\ell^1$  norms of the algebras truncated up to given homogeneity degrees.
Then there is an ``overall $\ell^1$ norm'' which may take infinite values.
The finite cores of the locally convex algebras are the ``plain'' algebras.
Mutatis mutandis, the discussion also extends to algebras of formal power series with continuous variables.

Let $T$ be a formal commutative variable.
Let us consider the formal power series
\begin{multline*} {}^\delta\Theta^{(\lambda)}(T)=\left|\mathcal R^{(\lambda)}\left(\Rexp\left(
 T\cdot\widehat\delta\mathrm Z^1_{[0,1)}\right)\right)\right|^{\forall T}_{\ell^1}
=\sum_{k=1}^\infty T^k\cdot\left|\mu^{(\lambda)}_{\mathrm R, k}\left(
 \widehat\delta\mathrm Z^1_{[0,1)} \right)\right|_{\ell^1}
\\
=\sum_{k=1}^\infty T^k\cdot\frac1{k!}\left|\mu^{(\lambda)}_{ k}
\left(\widehat\delta(Y_1),\ldots,\widehat\delta(Y_k)  \right)\right|_{\ell^1}
\equiv \sum_{k=1}^\infty {}^\delta\Theta^{(\lambda)}_k T^k.
\end{multline*}
This has the universal property that
\[\sum_{k=1}^\infty \left\|\mu^{(\lambda)}_{\mathrm R, k}\left(
 \ad\phi \right)\right\|_{\mathrm{op}} \leq {}^\delta\Theta^{(\lambda)}_{\mathrm{real}}(\textstyle\int|\phi|).\]
In particular, $\mathcal R^{(\lambda)}(\Rexp(\phi))$ exists if ${}^\delta\Theta^{(\lambda)}_{\mathrm{real}}(\textstyle\int|\phi|)<+\infty$.
On the other hand, if ${}^\delta\Theta^{(\lambda)}_{\mathrm{real}}(r)=+\infty$, then
$\phi=\ad\left( r\cdot\mathrm Z^1_{[0,1)} \right)$ over $\mathrm F^1([0,2))$ will yield a counterexample to the existence of
$\mathcal R^{(\lambda)}(\Rexp(\phi))W$ with $W=Z_{[1,2)}$.
But $\phi=\ad\left( r\cdot\mathrm Z^1_{[0,1)} \right)$ over $\mathrm F^1([0,1))$ would have trouble with boundedness.

Ultimately, the resolvent problem is governed by ${}^\delta\Theta^{(\lambda)}(T)$.
The coefficients of this series can be computed from the coproduct expression of
 $  \mu^{(\lambda)}_{ k}\left( Y_1,\ldots, Y_k   \right)$.
This expression is, however, complicated combinatorially.
We will recover some information about the convergence radius of ${}^\delta\Theta^{(\lambda)}(T)$ less directly, by spectral methods.
\snewpage

\subsection{The resolvent method in the adjoint representation}\plabel{ss:magnusadj}
~

\begin{theorem}\plabel{thm:resadjmagnus}
If $\int |\phi|<\sqrt2$, or even if
$\int |\phi|\leq\sqrt2$ and $\phi$ is the Lebesgue measure times a Bochner--Lebesgue integrable function, then
$\ad\phi$ is $M$-controlled.
Moreover, for $\lambda\in[0,1]$,
\begin{equation}
\mathcal R^{(\lambda)}(\Rexp(\ad\phi))=\sum_{k=1}^\infty\mu_{\mathrm R,k}^{(\lambda)}(\ad \phi).
\plabel{eq:agora}
\end{equation}
Also, $\log  \Rexp(\ad\phi)$ exists and
\[\log  \Rexp(\ad\phi)=\ad\log  \Rexp(\phi).\]
\begin{proof}
In general, $\Rexp(\ad \phi)=\Bigl(\Rexp( \phi) \Bigr)_{\mathrm L} \Bigl(\Rexp( \phi)^{-1} \Bigr)_{\mathrm R} $
(i. e. multiplication by $\Rexp( \phi)$ on the left and multiplication by $\Rexp( \phi)^{-1}$ on the right).
Applying Lemma \ref{lem:mulspec} with $X=\Rexp( \phi)$, we obtain the existence
of $\mathcal R^{(\lambda)}(\Rexp(\ad \phi))$.
In particular as this applies in the universal case, we obtain ${}^\delta\Theta^{(\lambda)}_{\mathrm{real}}(r)<+\infty $
 if $r<\sqrt2$ and  $\lambda\in [0,1]$ hold.
By the discussion in the previous subsection, this leads to \eqref{eq:agora}.
The same formula in the critical case  $\int |\phi|=\sqrt2$ follows by analytic continuation.
The conclusion of about the logarithm follows from analytic controllability and the unicity of analytic continuation.
\end{proof}
\end{theorem}
\snewpage

However,
\begin{theorem}\plabel{thm:adjturnback}
$\mathcal R^{(1/2)}(\Rexp(\ad\phi))$ does not exists in
$\mathcal B( \mathrm F^{1}([0,1)) )$ if we consider $\phi=\sqrt2\cdot \mathrm Z^1_{[0,1)}$ in $ \mathrm F^1([0,1))$.
This is a counterexample with $\int|\phi|=\sqrt2$.
There are counterexamples of multiple Baker--Campbell--Hausdorff type with $\int|\phi|=\sqrt2+\varepsilon$.
\begin{proof}
We want to show that the inverse of
 $ \Rexp(\ad\phi)+\Id =  \Ad\Rexp(\phi)+\Id$ does not exists in $\mathcal B( \mathrm F^{1}([0,1)) )$.
Firstly, let us consider the extended algebra $\mathrm F^{1}([0,1)[W] $, and we try to solve the equation
\begin{equation}
\left(\Ad\Rexp(\phi)+\Id\right)Q=W
\plabel{eq:romme}
\end{equation}
 for $Q$.
Here $\Ad\Rexp(\phi)+\Id=2\Id +$higher terms of $\ad\phi$.
Therefore \eqref{eq:romme} can be solved uniquely in up to all degrees in the truncation algebras,
 i. e. for some element $Q\in \mathrm F^{1,\mathrm{loc}}([0,1)[W] $.
We claim, however, that $Q\notin \mathrm F^{1}([0,1))[W]$, i. e. the norm of $Q$ is $+\infty$.
Indeed, otherwise $Q$ would yield a representative element in
 $\mathrm F^{1}([0,1))\otimes \mathrm F^{1}([0,1))^{\mathrm{op}}$
 inverting $\Rexp(\phi)\otimes (\Rexp(\phi)^{-1})^{\mathrm{op}}+1\otimes1^{\mathrm{op}}$.
This is, however, impossible: We know that the spectrum of
 $\Rexp(\phi)=\Rexp( \sqrt2\cdot \mathrm Z^1_{[0,1)})$ contains $\pm\mathrm i$;
 and by Lemma \ref{lem:startfone}(b), the spectrum of
 $\Rexp(\phi)\otimes (\Rexp(\phi)^{-1})^{\mathrm{op}}$ contains $\mathrm i/(-\mathrm i)=-1$.
Therefore $|Q|=+\infty$ is established.
In other terms,
\[\sum_{k=1}|Q_{[k]}|=+\infty,\]
 where $_{[k]}$ denotes the $k$-homogeneous part.
We cannot use, however $W$ in $\mathrm F^{1}([0,1))$.
However, we can use $W_s=\frac1sZ_{[0,s)}$ for $s\in(0,1)$.
Let $Q_s=\left(\Ad\Rexp(\phi)+\Id\right)^{-1} W_s $.
Then it is easy to see that
\[\lim_{s\searrow0}|(Q_s)_{[k]}|= |Q_{[k]}|.\]
This, however, implies that
\[\lim_{s\searrow0}| Q_s |= +\infty.\]
Therefore some $Q_s$ either cannot exist in $\mathrm F^{1}([0,1))$, or they  just
 contradict to the boundedness of $\left(\Ad\Rexp(\phi)+\Id\right)^{-1}$ anyway.

The counterexamples of mBCH type can be obtained by piecewise constant replacements.
\end{proof}
\end{theorem}
This means that, in general, the existence of the logarithm and the resolvent method fail for $\ad\phi$
 if the cumulative norm of $\phi$ reaches $\sqrt 2$,  meanwhile the Magnus expansion converges absolutely
 up to cumulative norm $2$ for $\phi$.
~\\

In fact, there is another familiar situation when the resolvent approach is not that effective in the study
 of the Magnus expansion.
This is when $\mathfrak A$ is a commutative Banach algebra.
Then the Magnus expansion is convergent in any case.
However we can guarantee the existence of $\mathcal R^{(\lambda)}(\Rexp \phi)$ for $\lambda\in[0,1]$
 in terms of the cumulative norm only for
 $\int|\phi|<$``$\left|\log\frac\lambda{\lambda-1}\right|$''$=\sqrt{\pi^2
 +\left(\log\frac\lambda{1-\lambda}\right)^2}$;
 and the existence of $\log(\Rexp \phi)$ for $\int|\phi|<\pi$.
(This we know by commutative spectral calculus, and counterexamples are already by $\mathfrak A=\mathbb C$.)
While this behaviour is different from the one above, it is not that different,
 even if much more special.

\snewpage
\appendix

\section{Ordered measures}
\plabel{sec:mea}
\textbf{General measures.}
Suppose that $\mathfrak a$ is a Banach space.
An $\mathfrak a$-valued measure on a bounded interval $I$ is a $\sigma$-additive $\mathfrak a$-valued
function on the set of (possible degenerate, even empty) subintervals of $I$.
If $\phi$ is  $\mathfrak a$-valued measure on an interval $I$,
and $J\subset I$ is a subinterval, then we can define the variation by
\[|\phi|(J):=\sup\left\{\sum_{s\in S} |\phi(J_s)|  \,:\, \text{the $J_s$ are finitely many disjoint
intervals in $J$}  \right\}.\]
This yields a $[0,\infty]$-valued measure on $I$.
We will be interested only in measures $\phi$ of finite variation, i.~e.~for which
$\int|\phi|\equiv |\phi|(I)<+\infty$.

An $\mathfrak a$-valued measure $\phi$ on the interval $I$ is absolutely continuous if there
exists a nonnegative Lebesgue-integrable function $H$ on $I$ such that
$|\phi|(J)\leq \int_{t\in J} H(t)\,\mathrm dt$ for any subinterval $J\subset I$.
For us, absolute continuity has little importance in itself.
However, the most notable examples are absolutely continuous:

\begin{example}\plabel{ex:mea1}
Suppose that  $I\subset\mathbb R $ is an interval, $h$ is an $\mathfrak a$-valued
Lebesgue-Bochner integrable function on $I$.
Let $\phi$ be the measure which is $h$ times the Lebesgue measure restricted to $I$;
i.~e.~for any subinterval $J\subset I $
\[\phi(J)=\int_{t\in J}h(t)\,\mathrm dt.\]
Then $\phi$ is an $\mathfrak a$-valued measure of finite variation;
moreover $|h|$ is Lebesgue-integrable on $I$, and for any subinterval $J\subset I $
\[ |\phi|(J)=\int_{t\in J} |h(t)|\,\mathrm dt. \]
(Or, in other terms, $|\phi|=|h|\mathbf 1_I$.)
\qedexer
\end{example}
\snewpage
\begin{example}\plabel{ex:mea2}
Let $\mathfrak A=\mathfrak B(\mathfrak X)$ be the set of bounded operators on a Banach space $\mathfrak X$.
Suppose $h$ is an $\mathfrak A$-valued
strongly Lebesgue-Bochner integrable function on $I$; i.~e.~for any $y\in\mathfrak X $
the function
$hy:I\rightarrow \mathfrak X\,:\,  t\in I \mapsto h(t)y\in\mathfrak X $ is Lebesgue-Bochner integrable.
We define $\phi(J)$ as the operator  given by
\[\phi(J)y=\int_{t\in J}h(t)y\,\mathrm dt.\]
(Notice that that $\phi(J)$ is not necessarily in $\mathfrak A=\mathfrak B(\mathfrak X)$ .)
\textit{Assume} that
\begin{equation}\sup\left\{\sum_{s\in S} |\phi(J_s)|  \,:\, \text{the $J_s$ are finitely many disjoint
intervals in $J$ } \right\}<+\infty.
\plabel{eq:finvar}\end{equation}
Then $\phi$ is an $\mathfrak A$-valued measure of finite variation.
Moreover, there is a nonnegative Lebesgue-integrable function $\lceil h\rceil$ on $I$ such that for any subinterval $J\subset I$,
\[ |\phi|(J)=\int_{t\in J} \lceil h\rceil(t)\,\mathrm dt. \]
In fact, $\lceil h\rceil(t)$ is the supremum of the set of
functions $\{ |h(t)y|\,:\,y\in \mathfrak X,|y|=1\}$   in
$\leq_{\textrm{almost everywhere}}$ sense (which is an object whose existence is not a priori guaranteed).

Conversely, if $h$ is an $\mathfrak A$-valued
strongly Lebesgue-Bochner integrable function on $I$,
and  $H$ is a nonnegative Lebesgue integrable function such that for any $y\in\mathfrak X$, $|y|=1$
the inequality
$|h(t)y| \leq_{\textrm{ almost everywhere in $t$}} H(t)$ holds, then \eqref{eq:finvar} holds,
and for any subinterval $J\subset I$, the inequality $|\phi|(J)\leq \int_{t\in J}H(t)\, \mathrm dt$ holds.
\qedexer
\end{example}
An $\mathfrak a$-valued measure $\phi$ on the interval $I$ is continuous if for any $\varepsilon>0$, and any subinterval $J\subset I$,
$J$ is a countable union of subintervals $J_{\xi}$, such that $|\phi|(J_{\xi})<\varepsilon$.
(Thus, continuity is an absolute notion and absolute continuity is a relative notion.)

Whenever $\phi$ is an $\mathfrak a$-valued measure on  the interval $I$, and $f$ is, say,
a real valued function on $I$,  and $A\subset I$, then by $\int_{A}f\, \phi\equiv\int_{t\in A}f(t)\, \phi(t) $ we mean
$\int_{I}(\chi_Af)\, \phi\equiv\int_{t\in I}(\chi_A(t)f(t))\, \phi(t) $, where
$\chi_A$ is the characteristic function of the set $A$.
I. e., we integrate over the original measure and not over a restriction, or an extension, or both.
(Although this will not  matter in our cases.)

In the previous discussion, $I$ can replaced by an abstract set $\Omega$,
and the system of subintervals can be replaced by a semiring $\mathfrak F$.
For the sake of simplicity, we require that $\Omega$ should be a countable union
of elements from $\mathfrak F$. (I. e. $\mathfrak F$ is $\sigma$-finite.)
Then $\phi:\mathfrak F\rightarrow\mathfrak a$ is should be $\sigma$-additive.
The Lebesgue-measure can be replaced by any measure (i.~e.~$\sigma$-additive function)
$\omega:\mathfrak F\rightarrow[0,\infty)$.
Lebesgue-integrable and Lebesgue--Bochner-integrable should be replaced
by  Lebesgue-integrable and Lebesgue--Bochner-integrable with respect to $\omega$.
Otherwise, the same applies.

The only case where we need this greater generality is when
it rises from the product construction from the interval case.
A generalized interval $I\subset\mathbb R^n$ is a direct product $I_1\times\ldots\times I_n$ of intervals  $I_j\subset\mathbb R$.
(Bounded) generalized intervals of a generalized interval form a semiring, etc.
Suppose that $\phi_j$ is an $\mathfrak a_j$-valued measure of finite variation on the (generalized) interval $I_j$ for
$j=1,\ldots,n$.
Also assume that $F:\mathfrak a_1\times\ldots\times\mathfrak a_n\rightarrow \mathfrak a$ is an $n$-linear continuous function
with norm $|F|$.
Then we can define $\phi=F(\phi_1,\ldots,\phi_n)$ by
\[\phi(J_1\times\ldots\times J_n)=F( \phi_1(J_1),\ldots,\phi_n(J_n)).\]
This yields an $\mathfrak a$-valued measure $\phi$ of finite variation
$\int |\phi|  \leq|F|\cdot(\int |\phi_1|)\cdot\ldots\cdot(\int |\phi_n|)  $.
(The product construction itself can also be carried out in greater generality. )
\snewpage

\textbf{Ordered measures.}
An $\mathfrak a$-valued measure $\phi$ of finite variation on the interval $I$ is continuous
if for any $x\in I$, the equality $\phi(\{x\})=0$ holds.
In general, if $\phi$ is of finite variation, then  $\phi(\{x\})=0$ holds for all except countably many $x\in I$.
These points, quite canonically, can be blown up to
closed intervals with continuous uniform distribution, leading to a measure $\phi^*$ on an interval $I^*$.
In many arguments, it is advantageous to pass from $\phi$ to $\phi^*$.

For an ordered $n$-tuple $(t_1,\ldots,t_n)$,
we define its multiplicity sequence $\mul(t_1,\ldots,t_n)$ as the sequence $(m_1,\ldots, m_s)$ of positive integers,
where
\[m_1+\ldots+m_s=n\]
and
\[t_{i}\neq t_{t+1}\text{ if and only if } i\in\{m_1,m_1+m_2,\ldots,m_1+m_2+\ldots+m_{s-1}\}.\]
We also use the multiindex factorial $(m_1,\ldots, m_s)!=m_1!\cdot m_2!\cdot\ldots\cdot m_s!$.

If $\phi$ is an $\mathfrak a$-valued measure of finite variation on the interval $I\subset \mathbb R$,
and $F:\mathfrak a^n\rightarrow \mathfrak a$
is a continuous $n$-linear function, then we can consider
the right-ordered evaluation
\[F_{\mathrm R}(\phi):=\int_{\substack{(t_1,\ldots,t_n)\in I^n\\ t_1\leq\ldots
\leq t_n}}\frac1{\mul(t_1,\ldots,t_n)!}F(\phi(t_1),\ldots,\phi(t_n)); \]
and the left-ordered evaluation
\[F_{\mathrm L}(\phi):=\int_{\substack{(t_n,\ldots,t_1)\in I^n\\ t_1\leq\ldots
\leq t_n}}\frac1{\mul(t_1,\ldots,t_n)!}F(\phi(t_n),\ldots,\phi(t_1)). \]

The reason behind this  definition is as follows.
The definition is uncontroversial if $\phi$ is continuous.
Indeed, in this case  $\mul(t_1,\ldots,t_n)!=1$ except on negligible subset.
On the other hand, if $\phi$ is not necessarily continuous,
but $\phi^*$ is a continuous blowup of $\phi$, then
 $F_{\mathrm R}(\phi)=F_{\mathrm R}(\phi^*)$, and $F_{\mathrm L}(\phi)=F_{\mathrm L}(\phi^*)$
holds, i.~e.~the  definition behaves naturally with respect to blow-ups.
This allows to pass between $\phi$ and $\phi^*$ in many arguments.

In general, we consider the  case of continuous measures as the basic one, where expressions
are slightly simpler;
however, in some  situations, it is  instructive to see how the general case formulates.

Even more generally, $I$ can be replaced by a fully ordered set $\Omega$ which
contains a countable dense set in the order topology;
and $\mathfrak F$ can be replaced by a semiring of some
possibly degenerate subintervals of $\Omega$.
($\mathfrak F$ is still required to be $\sigma$-finite.)
This generality is required only if we want to concatenate intervals abstractly.
Otherwise, one should restrict to half-open intervals, or to measures which are continuous at the endpoints.
In those cases the concatenation of measures on intervals is unproblematic.
\snewpage

\section{$L^1$ spaces of noncommutative power series}
\plabel{sec:ellone}
%\textbf{$L^1$ spaces of formal power series.}
Let $\mathbb K$ be $\mathbb R$ or $\mathbb C$.
In the case of the BCH expansion, the most important Banach algebras with respect to our investigations
are the Banach algebras of formal noncommutative power series with the $\ell^1$ norm:

If $A$ is a formal power  noncommutative series with variables from $\{Y_\lambda:\lambda\in\Lambda\}$ with coefficients from $\mathbb K$,
then let $|A|_{\ell^1}$ denote the sum of the absolute value of the coefficients.
This may be infinite, but the elements $A$ with $|A|_{\ell^1}<+\infty$ form the Banach algebra
$\mathrm F^1_{\mathbb K}[Y_\lambda\,:\,\lambda\in\Lambda]$.
A slightly bigger space but still an algebra is
$\mathrm F^{1,\mathrm{loc}}_{\mathbb K}[Y_\lambda\,:\,\lambda\in\Lambda]$;
this contains those formal power series $A$ such that $|\pi_k A|_{\ell^1}<+\infty$ for any $k$,
where $\pi_k$ denotes the projection to the component of formal degree $k$.
$\mathrm F^{1,\mathrm{loc}}_{\mathbb K}[Y_\lambda\,:\,\lambda\in\Lambda]$ is a locally convex algebra
with seminorms $|\pi_{\leq k}\cdot|_{\ell^1}$.
This is a useful space, because some  expressions might converge to an
$A$ in $\mathrm F^{1,\mathrm{loc}}_{\mathbb K}[Y_\lambda\,:\,\lambda\in\Lambda]$
but with $|A|_{\ell^1}=+\infty$, which precludes convergence in $\mathrm F^1_{\mathbb K}[Y_\lambda\,:\,\lambda\in\Lambda]$.
If $\widehat{\phantom{k}}:\mathrm F^1_{\mathbb K}[Y_\lambda\,:
\,\lambda\in\Lambda]\rightarrow \mathrm F^1_{\mathbb K}[Y_\lambda\,:\,\lambda\in\Lambda]$
is the map which replaces power series coefficients with absolute values, then it is norm-preserving and
weakly contractive ($1$-Lipschitz), i.~e.~$|\widehat A-\widehat B|_{\ell^1} \leq | A- B|_{\ell^1} $.
If $\sim$ is an equivalence relation on $\Lambda$, and
$\widetilde{\phantom{k}}:\mathrm F^1_{\mathbb K}[Y_\lambda\,:\,\lambda\in\Lambda]\rightarrow
\mathrm F^1_{\mathbb K}[Y_{\tilde\lambda}\,:\,\tilde\lambda\in\Lambda/\sim]$
is the map which identifies variables according to $\sim$, then it is also
weakly contractive ($1$-Lipschitz).
One can also factorize by simply sending some variables to $0$.
Similar observations apply to  $\mathrm F^{1,\mathrm{loc}}_{\mathbb K}[Y_\lambda\,:\,\lambda\in\Lambda]$.

In the case of the Magnus expansion, a continuous version of the previous example plays a similar role:
Consider the noncommutative polynomial algebra
\[\mathrm F^*_{\mathbb K}(\mathbb R)=\mathrm F_{\mathbb K}[Z_{[a,b)}\,:\,a<b\in\mathbb R ],\]
i.~e.~the free noncommutative algebra generated by the symbols $Z_{[a,b)}$ where $a,b\in\mathbb R$, $a<b$.
For a noncommutative polynomial $P(Z_{[a_i,b_i)}\,:\,1\leq i\leq s)$
we define its norm as follows.
Take any finite set $C=\{c_1\ldots,c_k\}\subset \mathbb R$ such that $c_1<\ldots< c_k$ and
\[\{a_1,\ldots,a_s,b_1,\ldots,b_s\}\subset C.\]
We can assume that $a_i=c_{A_i}$ and $b_i=c_{B_i}$.
Then we define
\begin{equation}
|P(Z_{[a_i,b_i)}\,:\,1\leq i\leq s)|_{\ell^1*}:=
\left|P\left(\sum_{r=A_i}^{B_i-1}(c_{r+1}-c_r)Y_{[c_r,c_{r+1})}\,:\,1\leq i\leq s\right)\right|_{\ell^1},
\plabel{eq:LF}
\end{equation}
where the latter $\ell^1$ is understood with respect to $\mathrm F^1_{\mathbb K}[  Y_{[c_i,c_{i+1})}\,:\, 1\leq i< k]$.
It is easy to see that further refinements of $C$ do not change $|\cdot|_{\ell^1*}$.
Thus it is well-defined and $|Z_{[a,b)}|_{\ell^1*}=b-a$.
It is also easy to see that $|\cdot|_{\ell^1*}$ is a seminorm on $\mathrm F^*_{\mathbb K}(\mathbb R)$
with the submultiplicative property.
However, it is not a norm, because $|\cdot|_{\ell^1*}$ cannot distinguish between
$Z_{[a,b)}$ and $Z_{[a,c)}+Z_{[c,b)}$ for $a<c<b$; not even if they are part of larger polynomials.
The seminorm $|\cdot|_{\ell^1*}$ actually descends to the factor algebra
\[\mathrm F_{\mathbb K}(\mathbb R)=\mathrm F_{\mathbb K}[Z_{[a,b)}\,:
\,a<b\in\mathbb R ]/\{Z_{[a,b)}=Z_{[a,c)}+Z_{[c,b)}\,:\, a<c<b\in\mathbb R\}.\]
There it induces a norm $|\cdot|_{\ell^1}$, because if \eqref{eq:LF} yields $0$, then it means that, using an appropriate $C$, the
refined expansion of $P(Y_{[a_i,b_i)}\,:\,1\leq i\leq s)$ yields $0$, thus it is $0$.
This normed space can be completed (still $|\cdot|_{\ell^1}$), yielding the Banach space $\mathrm F^1_{\mathbb K}(\mathbb R)$.
More generally, it can be considered as $\mathrm F^*_{\mathbb K}(\mathbb R)$ factorized by elements of seminorm $0$,
and completed.
This construction is effectively based on the semiring of the intervals $[a,b)$, and common refinements of them, and
the Lebesgue measure (interval length).

This construction can be formulated in greater generality.
Suppose that $\Omega$ is a set, $\mathfrak F$ is semiring on it, and $\omega: \mathfrak F\rightarrow[0,+\infty)$ is a measure
(i.~e.~$\sigma$-additive).
Then consider the noncommutative polynomial algebra
\[\mathrm F_{\mathbb K}^*(\Omega,\mathfrak F,\omega)=\mathrm F_{\mathbb K}[Z_A\,:\,A\in\mathfrak F].\]
For a noncommutative polynomial $P(Z_{A_i}\,:\,1\leq i\leq s)$,
take a finite set $C\subset \mathfrak F$, such that any $C$ contains pairwise disjoint, nonempty sets,
and any element $A_i$, $1\leq i\leq s$, is the disjoints union of some elements from $C$.
Then we define
\begin{equation}
|P(Z_{A_i}\,:\,1\leq i\leq s)|_{\ell^1*}:=\left|P\left(\sum_{A\in C,\,A\subset A_i}\omega(A)Y_A\,:\,1\leq i\leq s\right)\right|_{\ell^1},
\plabel{eq:LF2}
\end{equation}
where the latter $\ell^1$ is understood with respect to $\mathrm F^1_{\mathbb K}[  Y_A\,:\, A\in C]$.
Further refinements of $C$ do not change $|\cdot|_{\ell^1*}$,
thus it is well-defined. In particular, $|Z_{A_i} |_{\ell^1*}=\omega(A)$.
Again, $|\cdot|_{\ell^1*}$ is a seminorm on
$\mathrm F_{\mathbb K}^*(\Omega,\mathfrak F,\omega)$
with the submultiplicative property.
Then we factorize $\mathrm F_{\mathbb K}^*(\Omega,\mathfrak F,\omega)$ with those
elements whose seminorm is $0$, and then we complete it to a Banach algebra $\mathrm F^1_{\mathbb K}(\Omega,\mathfrak F,\omega)$.
Actually, we can simply describe the algebra after the factorization.
It is given by
\begin{multline}
\mathrm F_{\mathbb K}(\Omega,\mathfrak F,\omega)=\mathrm F_{\mathbb K}[Z_A\,:\,A\in\mathfrak F]/\{
Z_A=0 \text{ if }\omega(A)=0, \\ Z_A=Z_{A_1}+\ldots+Z_{A_s}\text{ if }
A=A_1\dot\cup\ldots\dot\cup A_s\}.
\plabel{eq:genring}
\end{multline}
Indeed, these factorizations are valid, yet if  \eqref{eq:LF2} yields $0$, then it means
that the corresponding polynomial can be expanded to $0$ using them.
Thus $\mathrm F_{\mathbb K}^1(\Omega,\mathfrak F,\omega)$ is the completion of $\mathrm F_{\mathbb K}(\Omega,\mathfrak F,\omega)$.
(We still use the notation $|\cdot|_{\ell^1}$ for the norm.)
We will use the notation $Z^1_A$ for the corresponding generators of the ring \eqref{eq:genring}.
The resulted spaces are naturally graded, and for any $A\in \mathrm F^1_{\mathbb K}(\Omega,\mathfrak F,\omega)$,
\[\left|A\right|_{\ell^1}=\sum_{k=0}^\infty\left|\pi_kA\right|_{\ell^1},\]
where $\pi_kA$ is the component of grade $k$ of $A$.

If we use the construction where $\Omega$ is $\mathbb R$, $\mathfrak F=\mathfrak I_{[)}$ is the set
of bounded intervals $[a,b)$, and $\omega$ is the restriction of the Lebesgue measure (i.~e.~interval length), then
this yields the construction of the previous example.
If we use the semiring $\mathfrak F'=\mathfrak I$ of the possibly degenerate (not necessarily half-open) intervals,
then it yields the same construction.
In fact, they are the same even on the pre-completed algebraic level   $\mathrm F_{\mathbb K}(\Omega,\mathfrak F,\omega)$.
Thus we will use the notation $\mathrm F^1_{\mathbb K}(\mathbb R)$ without particularly mentioning the semiring.
The same applies for $\mathrm F^1_{\mathbb K}([0,1])$, $\mathrm F^1_{\mathbb K}([0,1))$, etc.
In fact, we can even use the semiring $\mathfrak F''$ of the subsets of $\mathbb R$ of finite Lebesgue measure.
This will be different on the algebraic level, but essentially the same on the completed level;
see the following  measure theoretic interpretation:
\snewpage

Let
$(\Omega^{\in\mathbb N},\mathfrak F^{\in\mathbb N},\omega^{\in\mathbb N})
=\bigcup^{\circ}_{n\in\mathbb N}(\Omega^{n},\mathfrak F^{n},\omega^{n})$.
Then $\mathrm F^1_{\mathbb K}(\Omega,\mathfrak F,\omega)$ is naturally isomorphic to
$L^1(\Omega^{\in\mathbb N},\mathfrak F^{\in\mathbb N},\omega^{\in\mathbb N})$ as a base space with the appropriate norm;
the grading corresponds to the power decomposition,
the product is induced from taking exterior direct of function in homogeneous components,
and $Z_A$ $(A\in\mathfrak F)$ corresponds to the characteristic function of $A$ on $\Omega$.
Indeed, if we map $Z_{A_1}\ldots Z_{A_s}$ into the class of
$(t_1,\dots,t_n)\mapsto \chi_{A_1}(t_1)\ldots\chi_{A_s}(t_s)$,
then it  $\mathrm F^*_{\mathbb K}(\Omega,\mathfrak F,\omega)$ (and $\mathrm F_{\mathbb K}(\Omega,\mathfrak F,\omega)$)
gets mapped into a dense subset of  $L^1(\Omega^{\in\mathbb N},\mathfrak F^{\in\mathbb N},\omega^{\in\mathbb N})$
isometrically.
We now that the latter one is a Banach algebra, thus completion yields an isomorphism between
$\mathrm F^1_{\mathbb K}(\Omega,\mathfrak F,\omega)$ and  $L^1(\Omega^{\in\mathbb N},\mathfrak F^{\in\mathbb N},\omega^{\in\mathbb N})$.
Thus the ``completion picture'' and the ``exterior $L^1$ algebra picture'' are equivalent to each other.
At this point the second one looks simpler because it is constructed using pre-built elements.
However, if we are to take more sophisticated norms on formal power series, the first one
is simpler to handle. In general, we use the terminology of the completion picture.
(Such spaces often occur in functional analysis and probability theory, but with `$\otimes$' as the product.
This is to indicate the we deal with a ``free noncommutative'' construction.
But, this latter notation for the product is quite unnecessary.)

We define the tautological measure $\mathrm Z^1_{(\Omega,\mathfrak F,\omega)}:\mathfrak F\rightarrow
\mathrm F^1_{\mathbb K}(\Omega,\mathfrak F,\omega)$ as the
$\mathrm F^1_{\mathbb K}(\Omega,\mathfrak F,\omega)$-valued measure given by
\[\mathrm Z^1_{(\Omega,\mathfrak F,\omega)}(A)=Z^1_A\in \mathrm F^1_{\mathbb K}(\Omega,\mathfrak F,\omega)\]
for $A\in\mathfrak F$.
The variation measure of $\mathrm Z_{(\Omega,\mathfrak F,\omega)}$ is the measure $\omega$.
(In particular, it is absolutely continuous with respect to $\omega$.)
If $f:\Omega^n\rightarrow \mathbb K$ is Lebesgue integrable with respect to  $\omega^n$, then
\[\int_{(t_1,\ldots,t_n)\in\Omega^n}f(t_1,\ldots,t_n)\,\mathrm Z^1_{(\Omega,\mathfrak F,\omega)}(t_1)\,\ldots\,
\mathrm Z^1_{(\Omega,\mathfrak F,\omega)}(t_n)  \in \mathrm F^1_{\mathbb K}(\Omega,\mathfrak F,\omega)\]
exists and corresponds to ``$f$'' in the exterior $L^1$ algebra picture.
In particular, it has norm $|f|_{L^1(\Omega^n,\mathfrak F^n,\omega^n)}$.
We use the notation $\mathrm Z^1_{\mathbb R}$, $\mathrm Z^1_{[0,1]}$ etc in the case of $\mathrm F^1_{\mathbb K}(\mathbb R)$,
$\mathrm F^1_{\mathbb K}([0,1])$, etc., respectively.
These latter measure has approximately the same level of individuality as the Lebesgue measure.
Thus it is fair to say that $c\cdot\mathrm Z^1_{[0,1]}$, with $c\in(0,+\infty)$,
is \textit{the} totally noncommutative continuous mass of norm $c$.

The space  $\mathrm F^{1,\mathrm{loc}}_{\mathbb K}(\Omega,\mathfrak F,\omega)$ can be defined similarly, as in the discrete case.
There is a map
 $\widehat{\phantom{k}}:\mathrm F^1_{\mathbb K}(\Omega,\mathfrak F,\omega)\rightarrow \mathrm F^1_{\mathbb K}(\Omega,\mathfrak F,\omega)$
induced by replacing power series coefficients with absolute values.
(In the exterior $L^1$ algebra picture it is taking the absolute value of the representing function.)
This is also norm-preserving and weakly contractive.
Regarding factorizations, we are interested only in the simplest cases.
Assume that $A_1\dot\cup\ldots\dot\cup A_s\subset\Omega$, $A_i\in\mathfrak F$.
Then there is a map
$\widetilde{\phantom{k}}:\mathrm F^1_{\mathbb K}(\Omega,\mathfrak F,\omega)\rightarrow \mathrm F^1_{\mathbb K}[Y_{A_i}\,:\,1\leq i\leq s]$
induced by
$Z^1_{A}\mapsto\sum_{i=1}^s\omega(A\cap A_i)Y_{A_i}$.
This corresponds to restricting to $A_1\dot\cup\ldots\dot\cup A_s$, and coarsening into finitely many points.
This also yields a weak contraction.
Similar observations apply to  $\mathrm F^{1,\mathrm{loc}}_{\mathbb K}(\Omega,\mathfrak F,\omega) $.

\snewpage
\section{The spectral calculus of $\ad$ and $\Ad$}
\plabel{sec:sc}

Here the most general principles of holomorphic spectral calculus will be recalled but specialized to very concrete settings.

In the setting of a Banach algebra $\mathfrak A$,
we can consider the adjoint action $\ad:\mathfrak A\rightarrow \mathcal B(\mathfrak A)$,
where $(\ad X)Y=[X,Y]\equiv XY-YX$.
For $\mathcal B(\mathfrak A)$ we use the operator norm $\|\cdot\|_{\mathfrak A}$.
It is a trivial estimate that $\|\ad X\|_{\mathfrak A}\leq 2|X|_{\mathfrak A}$.
We note that the left multiplication $X_{\mathrm L}$ by $X$ on $\mathfrak A$
and right multiplication $X_{\mathrm R}$ by $X$ on $\mathfrak A$ are isospectral to $X$.
In this terminology, $\ad X=X_{\mathrm L}-X_{\mathrm R}$.
Also, $X_{\mathrm R}$ and $X_{\mathrm L}$ commute.

\begin{lemma}\plabel{lem:commspec}
 For the spectrum of $\ad X$,
\[\spec(\ad X)\subset \spec(X_{\mathrm L})+\spec(-X_{\mathrm R}) = \spec(X)-\spec(X) .\]
In particular, if $\spec X\subset \{z\in\mathbb  C\,:\,|\Ima z|<r \}$, then
$\spec(\ad X)\subset \{z\in\mathbb  C\,:\,|\Ima z|<2r \} $.

\begin{proof}
Suppose that $\gamma=\bigcup_i\gamma_i$ is a finite disjoint union of simple cycles
encircling, altogether, every point of $\spec(X_{\mathrm L})$ exactly $+1$ times;
and $\delta=\bigcup_j\delta_j$ is a finite disjoint union of simple cycles
encircling, altogether, every point of $\spec(-X_{\mathrm R})$ exactly $+1$ times.
Assume that $\gamma$ and $\delta$ bound $R_1$ and $R_2$, respectively; i. e.
$\gamma=\partial_{\mathfrak{or}}R_1$ and $\delta=\partial_{\mathfrak{or}}R_2$.
Then, for $\lambda\notin R_1+ R_2$,
\begin{equation}
\frac1{X_{\mathrm L}+(-X)_{\mathrm R}-\lambda\Id}=\frac1{(2\pi\mathrm i)^2}\int_{w:\delta}\int_{z:\gamma}
\frac1{z+w-\lambda}\,\frac{\mathrm dz}{z\Id-X_{\mathrm L}}\,\frac{\mathrm dw}{w\Id-(-X_{\mathrm R})}
\plabel{eq:commspec}
\end{equation}
holds.
Indeed, let us use the shorthand (i. e. non-symbolic) notation
$\mathcal F(\mathrm{expr})\equiv \frac1{(2\pi\mathrm i)^2}\int_{w:\delta}\int_{z:\gamma}
\mathrm{expr} \cdot\frac1{z+w-\lambda}\,\frac{\mathrm dz}{z\Id-X_{\mathrm L}}\,\frac{\mathrm dw}{w\Id-(-X_{\mathrm R})}$.
Then the calculation
\begin{multline*}
(X_{\mathrm L}+(-X)_{\mathrm R}-\lambda\Id) \mathcal F(1)=
\mathcal F(X_{\mathrm L}+(-X)_{\mathrm R}-\lambda\Id)=\mathcal F(X_{\mathrm L}+w-\lambda )
=\\=\mathcal F(z+w-\lambda )=
\frac1{(2\pi\mathrm i)^2}\int_{w:\delta}\int_{z:\gamma}
\frac{\mathrm dz}{z\Id-X_{\mathrm L}}\,\frac{\mathrm dw}{w\Id-(-X_{\mathrm R})}=\Id
\end{multline*}
is valid as we have integrands which are holomorphic in $z$ in a neighborhood of $\spec(X_{\mathrm L})$
and holomorphic in $w$ in a neighborhood of $\spec(-X_{\mathrm R})$ separately.
(For example, in order to check the third equality above, we have to show
$\mathcal F(z\Id-X_{\mathrm L} )=0$. And, indeed,
\[\mathcal F(z\Id-X_{\mathrm L} )=\frac1{(2\pi\mathrm i)^2}\int_{w:\delta}\int_{z:\gamma}
(z\Id-X_{\mathrm L} ) \cdot\frac1{z+w-\lambda}\,\frac{\mathrm dz}{z\Id-X_{\mathrm L}}\,\frac{\mathrm dw}{w\Id-(-X_{\mathrm R})}\]
\[=\frac1{(2\pi\mathrm i)^2}\int_{w:\delta}\left(\int_{z:\gamma}
  \frac1{z+w-\lambda}\, \mathrm dz \right)\,\frac{\mathrm dw}{w\Id-(-X_{\mathrm R})}
  =\frac1{(2\pi\mathrm i)^2}\int_{w:\delta}0\,\frac{\mathrm dw}{w\Id-(-X_{\mathrm R})}=0;\]
the main point is $z\mapsto\frac1{z+w-\lambda}$ being holomorphic near $R_1$.)

Therefore the arithmetic resolvent property checks out, characterising the resolvent.

Due to the compactness of the spectra and the continuity of multiplication,
 $\gamma$ and $\delta$ can be chosen in a $\varepsilon_1$-neighborhood of $\spec (X_{\mathrm L})$
 and in a $\varepsilon_2$-neighborhood of $\spec( -X_{\mathrm R}  )$, respectively, such that $R_1+R_2$
 will be in an $\varepsilon$-neighborhood of $\spec( X_{\mathrm L}) +\spec (-X_{\mathrm R}  )$.
Thus, the complement of $ R_1+ R_2 $ will exhaust the complement of $ \spec( X_{\mathrm L}) +\spec (-X_{\mathrm R}  ) $.
\end{proof}
\end{lemma}
\snewpage

\begin{remark}\plabel{rem:adform}
(a) The fact that $\spec(X)\subset\Dbar(0,r)$ implies $\spec(X)\subset\Dbar(0,2r)$ follows easily from the spectral radius limit formula.

(b) Assume that $\spec(X)\subset\left\{z\in\mathbb C\,:\, -\frac{|\Ima \lambda|}2<\Ima z<\frac{|\Ima \lambda|}2\right\}$.
Then
\begin{equation}
\frac1{X_{\mathrm L}+(-X)_{\mathrm R}-\lambda\Id}
=\frac{\sgn\Ima\lambda}{2\pi\mathrm i}\int_{u\in\mathbb R}
\frac{1}{(\frac\lambda2+u)\Id-X_{\mathrm L}}
\cdot
\frac1{(\frac\lambda2-u)\Id+X_{\mathrm R}}
\,\mathrm du .
\plabel{eq:spectconv}
\end{equation}
(This formula, via linear transformations $A\rightsquigarrow cA+d$, can easily be transcribed to bands not around the real axis.)
Indeed, this follows from  \eqref{eq:commspec},
considering certain (infinite) rectangles $R$ and $(1-\varepsilon)R$, resolving the  integral along the inner contour,
and pushing the loose boundary piece to the infinity.

In fact, \eqref{eq:spectconv} can also be established by a more arithmetical argument
\begin{align*}
&(X_{\mathrm L}+(-X)_{\mathrm R}-\lambda\Id)\cdot \frac{1}{2\pi\mathrm i}\int_{u\in\mathbb R}
\frac{1}{(\frac\lambda2+u)\Id-X_{\mathrm L}}\cdot\frac1{(\frac\lambda2-u)\Id+X_{\mathrm R}}\,\mathrm du\\
&=\frac{1}{2\pi\mathrm i}\int_{u\in\mathbb R}
\left(-\frac{1}{(\frac\lambda2+u)\Id-X_{\mathrm L}}-\frac1{(\frac\lambda2-u)\Id+X_{\mathrm R}}\right)\,\mathrm du\\
&=\lim_{N\rightarrow+\infty}\left(\frac{1}{2\pi\mathrm i}\int_{u=-N}^N
-\frac{1}{(\frac\lambda2+u)\Id-X_{\mathrm L}}\,\mathrm du\right)
+\lim_{N\rightarrow+\infty}\left(\frac{1}{2\pi\mathrm i}\int_{u=-N}^N
-\frac1{(\frac\lambda2-u)\Id+X_{\mathrm R}}\,\mathrm du\right)
\\
&=(\sgn\Ima\lambda)\frac12\Id_{\mathrm L}+(\sgn\Ima\lambda)\frac12\Id_{\mathrm R}=(\sgn\Ima\lambda)\Id.
\end{align*}
This is equivalent to \eqref{eq:spectconv}.

We remark that
considering the linear transforms, this special case is already sufficient to establish
$\spec(\ad X)\subset\conv(\spec(X_{\mathrm L})+\spec(-X_{\mathrm R}))=\conv(\spec(X)-\spec(X)).$
\qedremark
\end{remark}

If $X\in\mathfrak A$ is invertible, then let $\Ad X$ denote the operator
$\Ad X\,:\,\mathfrak A\rightarrow\mathfrak A,$ given by $Y\mapsto XYX^{-1}$.
Then $\Ad X=(X_{\mathrm L})(X_{\mathrm R})^{-1}$.

\begin{lemma}\plabel{lem:mulspec}
 If $X$ is invertible, then for the spectrum of $\Ad X$,
\[\spec(\Ad X)\subset \spec(X_{\mathrm L})\cdot\spec((X_{\mathrm R})^{-1})=\spec(X)/\spec(X).\]
In particular, if $\spec X\subset \{z\in\mathbb  C\,:\,\Rea z>0 \}$, then
$\spec(\Ad X)\subset \mathbb C\setminus [0,+\infty)$.
\begin{proof}
This is entirely similar to the Lemma \ref{lem:commspec} but with respect to the formula
\begin{equation}
\frac1{X_{\mathrm L}(X_{\mathrm R})^{-1}-\lambda\Id}=\frac1{(2\pi\mathrm i)^2}\int_{w:\delta}\int_{z:\gamma}
\frac1{z w-\lambda}\,\frac{\mathrm dz}{z\Id-X_{\mathrm L}}\,\frac{\mathrm dw}{w\Id- (X _{\mathrm R})^{-1}}.
\plabel{eq:mulspec}
\qedhere
\end{equation}
\end{proof}
\end{lemma}
\begin{remark}\plabel{rem:mulform}
Assume that $\spec X\subset \{z\in\mathbb  C\,:\,\Rea z>0 \}$.
Then the inverse of $X_{\mathrm L}(X_{\mathrm R})^{-1}$ is $(X_{\mathrm L})^{-1}X_{\mathrm R}$, while for $\lambda\in(-\infty,0)$,
\[\frac1{X_{\mathrm L}(X_{\mathrm R})^{-1}-\lambda\Id}=\int_{u\in\mathbb R}\frac1{2\pi\sqrt{-\lambda}}\cdot
\frac{1}{u\mathrm i\sqrt{-\lambda}\Id-X_{\mathrm L}}\cdot
\frac{X_{\mathrm R}}{-u\mathrm i\frac1{\sqrt{-\lambda}}\Id-X_{\mathrm R}}\,\mathrm du.
\]
(This formula can also be used to deal with a slightly different setup in sectoriality.)
\qedremark
\end{remark}

(For finer aspects of commutative holomorphic multivariable calculus see Waelbroeck \cite{Wae},
and Arens, Calderón \cite{AC} and subsequent developments.)

\snewpage
\section{The Banach algebraic form of Schur's formulae}
\plabel{sec:schuranal}
In many approaches to the Magnus / Baker--Campbell--Hausdorff expansions,
Schur's formulae play the role of a stepping stone.
The formal versions have already been encountered in Theorem \ref{thm:schurres}.
Here we give the form valid in general Banach algebras.
A standard proof is given here, using functional calculus;
it is independent from the  combinatorial resolvent discussion in Section \ref{sec:DiscRes}
and from the analytic resolvent discussion in Section \ref{sec:adj}.
As Schur's formulae can be proven quite simply anyway; the point here is more like reviewing the spectral aspects.
Here $\mathfrak A$ will be a Banach algebra, $X\in\mathfrak A$.

Consider the meromorphic function
\begin{equation}\beta(x)=\frac x{\mathrm e^x-1}=\sum_{j=0}^\infty \beta_j x^j.\plabel{eq:schur}\notag\end{equation}
Note that this function has poles at $2\pi\mathrm i(\mathbb Z\setminus\{0\})$.
Now, $\beta(\ad X)$ can be constructed as long as $\spec(\ad X)$ is disjoint
$2\pi\mathrm i(\mathbb Z\setminus\{0\})$.
On the other hand, for the entire function
\begin{equation}\alpha(x)=\frac {\mathrm e^x-1}x=\sum_{j=0}^\infty \frac1{(j+1)!} x^j,\plabel{eq:antischur}\notag\end{equation}
the map $\alpha(\ad X)$ is well-constructed anyway.
The significance of $\beta(\ad X)$ is that it provides an inverse to $\alpha(\ad X)$ (for a good $X$).
The significance of $\alpha(\ad X)$, in turn, comes from the derivative the exponential function.
Using the power series of $\exp$, it easy to see that
\[ \left.\frac{\mathrm d}{\mathrm dt}\exp(X+tZ)\right|_{t=0}\exp(-X)=\int_{s=0}^1\exp(sX)Z\exp(-sX)\,\mathrm ds=\alpha(\ad X)Z\]
and
\[ \exp(-X)\left.\frac{\mathrm d}{\mathrm dt}\exp(X+tZ)\right|_{t=0}=\int_{s=0}^1\exp(-sX)Z\exp(sX)\,\mathrm ds=\alpha(-\ad X)Z.\]
\snewpage
\begin{theorem}[F. Schur \cite{Sch1} (1890), \cite{Sch2}, Poincar\'e \cite{PH} (1899) ] \plabel{thm:schuranal}
If $|X|<\pi$, or $\spec(X)\subset\{z\in\mathbb C\,:\,|z|<\pi\}$, or $\spec(X)\subset\{z\in\mathbb C\,:\,|\Rea z|<\pi\}$ , then
$|\ad X|_{\mathfrak A}<2\pi$, or  $\spec(\ad X)\subset\{z\in\mathbb C\,:\,| z|<2\pi\}$,
or  $\spec(\ad X)\subset\{z\in\mathbb C\,:\,|\Rea z|<2\pi\}$, respectively.
In particular,
$\beta(\ad X):\mathfrak A\rightarrow \mathfrak A$ makes sense as an absolute convergent power series (first two cases)
 or as a homomorphic function of  $\ad X$.

In these cases, for $Y\in\mathfrak A$,
%\stepcounter{equation}
\begin{equation}\frac{\mathrm d}{\mathrm dt}\log(\exp(tY)\exp(X))\Bigl|_{t=0}=\beta(\ad X)Y
%\tag{\theequation L}
\plabel{eq:schurL}\end{equation}
and
\begin{equation}\frac{\mathrm d}{\mathrm dt}\log(\exp(X)\exp(tY))\Bigl|_{t=0}=\beta(-\ad X)Y
%\tag{\theequation R}
\plabel{eq:schurR}\end{equation}
hold;
with the usual $\log$ branch cut along the negative real axis.
\end{theorem}

\begin{proof}
The arguments regarding the norms are trivial, and the more general spectral behaviour
follows from Lemma \ref{lem:commspec}.
Regarding \eqref{eq:schurL} and \eqref{eq:schurR}:
In those spectral ranges,  $\exp A$ (defined for $\spec(A)\subset \{z\in\mathbb C\,:\,|\Rea z|<\pi\}$)
and   $\log A$ (defined for $\spec(A)\subset \mathbb C\setminus(-\infty,0]$) are smooth and inverses
of each other.

More specifically, for fixed $X$ with $\spec(X)\subset \{z\in\mathbb C\,:\,|\Rea z|<\pi\}$),
the maps
\[\hat Z\mapsto \hat Y=\log(\exp(\hat Z)\exp(-X))\]
and
\[\hat Y\mapsto \hat Z=\log(\exp(\hat Y)\exp(X))\]
are inverses of each other (and smooth) for $\hat Z\sim X$ and $\hat Y\sim 0$.
For the first map, as the differential of $\log$ at $1$ is the identity map,
the linear derivative at $\hat Z=X$ is given by $Z\mapsto \alpha(\ad X)Z$.
For the second map, this forces the linear derivative  at $\hat Y=0$ to be the inverse map $Y\mapsto \beta(\ad X)Y$.
This establishes \eqref{eq:schurL}.
Analogously, the maps
\[\hat Z\mapsto \hat Y=\log(\exp(-X)\exp(\hat Z))\]
and
\[\hat Y\mapsto \hat Z=\log(\exp(X)\exp(\hat Y))\]
are inverses of each other (and smooth) for $\hat Z\sim X$ and $\hat Y\sim 0$.
For the first map, the linear derivative at $\hat Z=X$ is given by $Z\mapsto \alpha(-\ad X)Z$.
For the second map, this forces the linear derivative at $\hat Y=0$ to be the inverse map $Y\mapsto \beta(-\ad X)Y$.
This establishes \eqref{eq:schurR}.
\end{proof}
\begin{commentx}
\begin{remark}\plabel{rem:schurformal}
The inverse relationship between the tangent map of the exponential and the log can be expressed by the fact that
for $A$ $\log$-able,
\[\mathrm T\exp_{\log A}:\mathfrak A\rightarrow\mathfrak A\qquad Y\mapsto\int_{t=0}^1 A^{t}YA^{1-t}\,\mathrm dt\]
and
\[\mathrm T\log_A:\mathfrak A\rightarrow\mathfrak A\qquad Y\mapsto\int_{\lambda=0}^1
\frac1{\lambda+(1-\lambda)A}Y\frac1{\lambda+(1-\lambda)A} \,\mathrm d\lambda\]
are inverses of each other.
\qedremark
\end{remark}
\end{commentx}
\begin{cor}
[F. Schur \cite{Sch1} (1890), \cite{Sch2},  Poincar\'e \cite{PH} (1899) ]
\plabel{thm:schurformal}
If $X$ and $Y$ are formal noncommutative variables, then
%\stepcounter{equation}
\begin{equation}
\log(\exp(Y)\exp(X))_{\text{the multiplicity of $Y$ is $1$}}=\beta(\ad X)Y
%\tag{\theequation L}
\plabel{eq:schur2L}
\end{equation}
and
\begin{equation}
\log(\exp(X)\exp(Y))_{\text{the multiplicity of $Y$ is $1$}}=\beta(-\ad X)Y ;
%\tag{\theequation R}
\plabel{eq:schur2R}
\end{equation}
where $\beta(\ad X)$ is understood in the sense of formal power series.
\begin{proof}
This follows even from the general matrix case of Theorem \ref{thm:schuranal}.
\end{proof}
\end{cor}
By this, we have rederived Theorem \ref{thm:schurres}.
(In turn, we could have obtained Theorem \ref{thm:schuranal} from
Theorem \ref{thm:schurres} / Corollary \ref{thm:schurformal} by analytic extension;
but the discussion of the spectrally allowed domain is unavoidable in any case.)

\begin{discussion}
\plabel{disc:tradBCH}
From \eqref{eq:schurL} and $\eqref{eq:schurR}$, one can immediately deduce that
using the abbreviation $\log(\exp(U)\exp(V))\equiv\BCH(U,V)$,
\begin{equation}
\BCH(Y,X)=X+\int_{t=0}^1\beta(\ad \BCH(tY,X))Y\,\mathrm dt
\plabel{eq:CBHL}
\end{equation}
and
\begin{equation}
\BCH(X,Y)=X+\int_{t=0}^1\beta(-\ad \BCH(X,tY))Y \,\mathrm dt.
\plabel{eq:CBHR}
\end{equation}
Thinking about $\beta$ as $\beta^{\mathrm{pow}}$, and iterating \eqref{eq:CBHL} and \eqref{eq:CBHR} (separately or in a mixed manner),
after sufficiently  many invocation of the formulae above, the first $n$ orders of $\BCH(X,Y)$ clear to  commutator polynomials.
This is essentially the original Campbell--Baker--Hausdorff approach, cf. Achilles, Bonfiglioli \cite{AB},
except they deal with the original, more complicated, Lie group setting.
\qedremark
\end{discussion}

\begin{remark}\plabel{rem:tradBCH}
The traditional approach above is also applicable to the analytic study of BCH expansion.
It seems to be less effective in the Banach algebraic case, therefore it is not discussed here.
See, however, Part III for the Banach--Lie setting. \qedremark
\end{remark}
\snewpage
\begin{remark}\plabel{rem:Bernoulli}
The Bernoulli numbers $B_j$ are defined by the expansion
\[\beta(x)=\frac x{\mathrm e^x-1}=\sum_{j=0}^\infty\beta_j x^j=\sum_{j=0}^\infty\frac {B_j}{j!} x^j\qquad(|x|<2\pi).\]
Then
\begin{equation}B_j=\sum_{k=0}^j\begin{pmatrix} j\\ k\end{pmatrix}B_k\qquad(j\geq2);\plabel{a}\end{equation}
\begin{equation}B_j=-\begin{pmatrix} j\\ k\end{pmatrix}\sum_{k=0}^{j-1}\frac{B_k}{j+1-k}\qquad(j\geq2);\plabel{b}\end{equation}
\begin{equation}B_{2j+1}=0\qquad(j\geq1);\plabel{c}\end{equation}
\begin{equation}\tan x=\sum_{j=1}^\infty  2^{2j}(2^{2j}-1)\frac{(-1)^{j+1}B_{2j}}{(2j)!}x^{2j-1} \quad(|x|<\pi/2); \plabel{d}\end{equation}
\begin{equation}\sgn B_{2j}=(-1)^{j+1}\qquad(j\geq1);\plabel{e}\end{equation}
\begin{equation}\tanh x=\sum_{j=1}^\infty  2^{2j}(2^{2j}-1)\frac{B_{2j}}{(2j)!}x^{2j-1} \quad(|x|<\pi/2); \plabel{dd}\end{equation}
\begin{equation}B_0=1,\quad\! B_1=-\frac12,\quad\!  B_2=\frac16,\quad\!   B_4=-\frac1{30},
\quad\!   B_6=\frac1{42},\quad\!   B_8=-\frac1{30},\quad\!   B_{10}=\frac5{66};\plabel{f} \end{equation}
\begin{equation}B_{2j}=(-1)^{j+1}\frac{(2j)!}{(2\pi)^{2j}}\,2\zeta(2j)
\equiv(-1)^{j+1}\frac{(2j)!}{(2\pi)^{2j}}\,\,2\sum_{N=1}^\infty\frac1{N^{2j}}
\quad (j\geq 1).\plabel{g}\end{equation}
(The power series above are to be interpreted analytically.)
Properties (\ref{a}--\ref{f}) are relatively straightforward to prove using elementary methods, in particular, \eqref{c} from
\[\beta(x)-\beta(-x)=-x;\]
while
(\ref{g}) follows from Euler's formula
\[\pi\cot \pi x=\frac 1x +\sum_{\substack{N=-\infty\\N\neq 0}}^\infty \left( \frac1{x-N}+\frac1{N} \right).\]

Thus, the Bernoulli numbers turn up in some relevant power series.
\qedremark
\end{remark}

\begin{commenty}

\snewpage
\section{Perturbations of the resolvent}\plabel{sec:resolvid}

For the reader's convenience,
here we collect the basic facts concerning perturbations of the resolvent but rewritten to terms of the ``multiplicative format'' \eqref{eq:resdef}.

\begin{lemma}
\plabel{lem:perres1}
[Resolvent perturbations of first type.]

(a) If $\lambda,\nu$ are scalars, and  $\mathcal R^{(\lambda)}(A)$ exists, moreover $|\nu-\lambda|\,|\mathcal R^{(\lambda)}(A)|<1$,
then $\mathcal R^{(\nu)}(A)$ exists and
\[\mathcal R^{(\nu)}(A)=\frac{ \mathcal R^{(\lambda)}(A)}{1-(\nu-\lambda)\mathcal R^{(\lambda)}(A)}=
\sum_{k=0}^\infty (\nu-\lambda)^k\left(\mathcal R^{(\lambda)}(A)\right)^{k+1}.\]

(b) If $\lambda,\nu$ are scalars, and  $\mathcal R^{(\lambda)}(A)$ does not exist, then
\[\mathcal R^{(\nu)}(A)\text{ does not exists \qquad or \qquad}  |\mathcal R^{(\nu)}(A)|\geq \frac1{|\nu-\lambda|}.\]

\begin{proof}[Note] In the statement, the occurrences of
`$|\mathcal R^{(\lambda)}(A)|$' can be replaced by   $ \mathrm r\left(\mathcal R^{(\lambda)}(A)\right) $.
\qedno
\end{proof}
\begin{proof}
(a) This is a standard argument using Neumann series.
(b) is a consequence  of (a) with the role of $\lambda$ and $\nu$ interchanged.
\end{proof}
\end{lemma}

\begin{lemma}
\plabel{lem:perres2} [Resolvent perturbations of second type]
If $\lambda$ is a scalar,  $\mathcal R^{(\lambda)}(A)$ exists, $ B\in \mathfrak A$, and either
$\left|  (\lambda-1)(B-A)(1+(\lambda-1) \mathcal R^{(\lambda)}(A)) \right|<1$ or
\rechoicecomm{\\ \phantom{a}\qquad}
  $\left|  (\lambda-1)(1+(\lambda-1) \mathcal R^{(\lambda)}(A))\boldsymbol(B-A) \right|<1$ holds, then
\begin{equation*}
\mathcal R^{(\lambda)}(B)=\mathcal R^{(\lambda)}(A)+\sum_{n=1}^\infty
(\lambda-1)^{n-1}\cdot
\bigl(1+(\lambda-1) \mathcal R^{(\lambda)}(A)\bigr)\left((B-A)\bigl(1+(\lambda-1) \mathcal R^{(\lambda)}(A)\bigr)\right)^n.
%\plabel{eq:formresser}
\end{equation*}
\begin{proof}[Note]
 $\bigl(1+(\lambda-1) \mathcal R^{(\lambda)}(A)\bigr)\equiv(\lambda+(1-\lambda)A)^{-1} $.
 \qedno
\end{proof}
\begin{proof}[Note] In the statement, `$\left|  (\lambda-1)(B-A)(1+(\lambda-1) \mathcal R^{(\lambda)}(A)) \right|$'
and
\rechoicecomm{\\ \phantom{a}\qquad}
`$\left|  (\lambda-1)(1+(\lambda-1) \mathcal R^{(\lambda)}(A))\boldsymbol(B-A) \right|$'
can be replaced by
\rechoicecomm{\\ \phantom{a}\qquad}
$\mathrm r\left(  (\lambda-1)(B-A)(1+(\lambda-1) \mathcal R^{(\lambda)}(A)) \right)
 \equiv
\mathrm r\left(  (\lambda-1)(1+(\lambda-1) \mathcal R^{(\lambda)}(A))(B-A)\right) $.
\qedno
 \end{proof}
\begin{proof}
Again, this is an argument using Neumann series.
\end{proof}
\end{lemma}
\renewcommand{\qedsymbol}{$\triangle$}

\end{commenty}
\snewpage

\end{document}